\documentclass[11pt, oneside]{amsart}

\usepackage{tabularx, ifthen}
\usepackage{amssymb} \usepackage{amsfonts} \usepackage{amsmath}
\usepackage{amsthm} \usepackage{epsfig, subfig}
\usepackage{ amscd, amsxtra, latexsym,pb-diagram,pb-xy}
\usepackage[all]{xy}
\usepackage{caption}
\usepackage{enumerate}
\usepackage{color}

\usepackage[verbose,hyperpageref]{backref}
\backrefsetup{verbose=false}

\renewcommand*{\backrefalt}[4]{%
\ifcase #1 %
No citations.%
\or
(Cited p. #2.)%
\else
(Cited pp. #2.)%
\fi
}

\usepackage[pagebackref=true]{hyperref}

\usepackage{paralist}

\usepackage{epigraph}

\newtheorem{lemma}{Lemma}[section]
\newtheorem{thm}[lemma]{Theorem}
\newtheorem{prop}[lemma]{Proposition}
\newtheorem{cor}[lemma]{Corollary}

\newtheorem{thmintro}{Theorem}
\newtheorem{conjintro}[thmintro]{Conjecture}
\newtheorem{questintro}[thmintro]{Question}
\newtheorem{remintro}[thmintro]{Remark}
\newtheorem{corintro}[thmintro]{Corollary}

\newtheorem{claim}[lemma]{Claim}

\theoremstyle{definition}
\newtheorem{defn}[lemma]{Definition}

\newtheorem{rem}[lemma]{Remark}
\newtheorem*{rem*}{Remark}
\newtheorem{conv}[lemma]{Convention}
\theoremstyle{definition}
\newtheorem{notation}[lemma]{Notation}

\definecolor{darkgreen}{cmyk}{1,0,1,.2}

\newcommand{\showcomments}{yes}
\renewcommand{\showcomments}{no}

\newsavebox{\commentbox}
{\ifthenelse{\equal{\showcomments}{yes}}%
{\footnotemark
        \begin{lrbox}{\commentbox}
        \begin{minipage}[t]{1.25in}\raggedright\sffamily\tiny
        \footnotemark[\arabic{footnote}]}
{\begin{lrbox}{\commentbox}}}%
{\ifthenelse{\equal{\showcomments}{yes}}%
{\end{minipage}\end{lrbox}\marginpar{\usebox{\commentbox}}}
{\end{lrbox}}}

\newcommand{\stabilizer}{\mathrm{Stab}}

\newcommand{\Act}{{\rm Act}}

\newcommand{\nest}{\sqsubseteq}
 \usepackage{mathabx}
\newcommand{\propnest}{\sqsubsetneq}
\newcommand{\orth}{\bot}
\newcommand{\transverse}{\pitchfork}

\newcommand{\duaug}[2]{{#1}^{+{#2}}}
\newcommand{\link}{\mathrm{Lk}}
\newcommand{\sat}{\mathrm{Sat}}
\newcommand{\naturals}{\mathbb N}
\newcommand{\dist}{\mathrm{d}}
\newcommand{\diam}{\mathrm{diam}}
\newcommand{\neb}{\mathcal N}
\newcommand{\image}{\mathrm{im}}

\newcommand{\dimension}{\mathrm{dim}}
\newcommand{\fix}{\mathrm{Fix}}
\newcommand{\pstab}{\mathrm{PStab}}
\newcommand{\act}{\mathrm{Act}}
\newcommand{\actminus}[1]{\act({#1})\setminus \{{#1}\}}

\newcommand{\cuco}[1]{{\mathcal #1}}
\newcommand{\tup}[1]{\vec{#1}}
\newcommand{\tsh}[1]{\left\{\kern-.7ex\left\{#1\right\}\kern-.7ex\right\}}
\newcommand{\Tsh}[2]{\tsh{#2}_{#1}}
\newcommand{\ignore}[2]{\Tsh{#2}{#1}}

\newcommand{\lift}[1]{\widehat{#1}}

\newcommand{\frakS}{\mathfrak S}

\makeatletter
\@tfor\next:=abcdefghijklmnopqrstuvwxyzABCDEFGHIJKLMNOPQRSTUVWXYZ\do{%
  \def\command@factory#1{%
    \expandafter\def\csname cal#1\endcsname{\mathcal{#1}}
    \expandafter\def\csname frak#1\endcsname{\mathfrak{#1}}
    \expandafter\def\csname scr#1\endcsname{\mathscr{#1}}
    \expandafter\def\csname bb#1\endcsname{\mathbb{#1}}
    \expandafter\def\csname rm#1\endcsname{\mathrm{#1}}
      \expandafter\def\csname bf#1\endcsname{\mathbf{#1}}
  }
 \expandafter\command@factory\next
}
\makeatother

\newcommand{\fontact}{{\mathcal C}}

\def\co{\colon}
\newcommand{\genus}{\mathrm{Gen}}

\begin{document}

\title[Combinatorial HHS \& quotients of MCG]{A combinatorial 
take on hierarchical hyperbolicity\\ 
and applications to quotients of mapping class groups}
\author[J Behrstock]{Jason Behrstock}
\address{Lehman College and The Graduate Center, CUNY, New York, New York, USA}
\email{jason.behrstock@lehman.cuny.edu} 
\author[M Hagen]{Mark Hagen}
\address{School of Mathematics, University of Bristol, Bristol, UK}
\email{markfhagen@posteo.net}
\author[A Martin]{Alexandre Martin}
	\address{Maxwell Institute and Department of Mathematics, Heriot-Watt University, Edinburgh, UK}
\email{alexandre.martin@hw.ac.uk}
\author[A Sisto]{Alessandro Sisto}
	\address{Maxwell Institute and Department of Mathematics, Heriot-Watt University, Edinburgh, UK}
	\email{a.sisto@hw.ac.uk}

	\thanks{To appear in the Journal of Topology.}
	
\maketitle

\begin{abstract}
We give a simple combinatorial criterion, in terms of an action on a
hyperbolic simplicial complex, for a group to be hierarchically
hyperbolic.  We apply this to show that quotients of mapping class
groups by large powers of Dehn twists are hierarchically hyperbolic
(and even relatively hyperbolic in the genus 2 case). In genus at 
least three, there are no known infinite hyperbolic quotients of
mapping class groups. However, using the hierarchically hyperbolic quotients
we construct, we show, under a residual finiteness assumption, that
mapping class groups have many non-elementary hyperbolic quotients.
Using these quotients, we relate questions of Reid and
Bridson--Reid--Wilton about finite quotients of mapping class groups
to residual finiteness of specific hyperbolic groups.
\end{abstract}

 \tableofcontents

\section*{Introduction}

Hierarchically hyperbolic spaces/groups (HHS/HHGs) were introduced in
\cite{HHS_I} as a common framework for studying the coarse geometry of
mapping class groups and CAT(0) cube complexes.  Since then, this
notion has found applications in a variety of flavors, including 
results of a 
coarse geometric \cite{HHS_III,HHS:quasiflats,RST} and of an algebraic
nature
\cite{Haettel,HHS:boundary,ABD,AbbottBehrstock:conjugator,AbbottNgSpriano:uniform},
with many of these new even for well-studied examples such as mapping
class groups.  Also, the class of groups known to be hierarchically hyperbolic groups has expanded 
considerably beyond the motivating examples of mapping class groups, compact special groups, and hyperbolic 
groups, and there are many ways to produce new hierarchically hyperbolic groups from old: 
suitable graphs of groups~\cite{HHS_II,BerlaiRobbio,RobbioSpriano}, graph 
products~\cite{BerlyneRussell,BerlaiRobbio}, and, somewhat in the spirit of the 
present paper, certain ``small-cancellation'' quotients, including, say, quotients of mapping class groups 
by normal subgroups generated by high powers of a pseudo-Anosov~\cite{HHS_III}.  The machinery built in the 
present paper gives a method for producing even more new examples.  Specifically, we show in 
Theorem~\ref{thmi:MCG}, that quotients of mapping class groups by subgroups generated by 
suitable powers of (all) Dehn twists are hierarchically hyperbolic.

The main drawback of the theory has been that,
despite the simplification in the definition provided by \cite{HHS_II},
verifying that any particular space is an HHS requires a lot of work, and
understanding of the HHS machinery.  We remedy this by providing a
combinatorial sufficient condition for a space/group to be hierarchically hyperbolic.  This criterion is 
simpler to state than the definition of hierarchical hyperbolicity, and it is (in principle and in examples) 
easier to verify for a given space/group.

We first need a definition.  Let a group $H$ act on
a (possibly disconnected) simplicial complex $Y$.  We say that $Y$ is a \emph{hyperbolic
$H$--space} if it becomes hyperbolic upon addition of finitely many $H$--orbits of edges (see Definition \ref{defn:hyp_H_space}).

A
simplified, but still powerful, version of the criterion is the 
following theorem which, informally, states that if $G$ acts cocompactly
with finite stabilizers of maximal simplices on a hyperbolic
simplicial complex with hyperbolic links, then $G$ is
hierarchically hyperbolic under a geometric condition (quasi-isometric
embedding in condition \ref{item:stab_hyp_intro}) and a combinatorial
condition on intersections of links (condition
\ref{item:propaganda_join_intro}).

\begin{thmintro}\label{thm:main_intro}
  Let the group $G$ act cocompactly on the flag simplicial complex $X$, and suppose that maximal simplices have finite stabilizers.  Suppose that:
  
\begin{enumerate}[(A)]
\item\label{item:stab_hyp_intro} For every simplex $\Delta$ of $X$, 
its link, $\link(\Delta)$, is a hyperbolic 
$Stab_G(\link(\Delta))$-space quasi-isometrically embedded in $X-\bigcup_{\link(\Sigma)=\link(\Delta)}\Sigma$.

\item\label{item:propaganda_join_intro} For all simplices $\Delta,\Sigma$ of $X$ there exist simplices $\Pi,\Pi'$ of $\link(\Delta)$ such that
$\link(\Delta)\cap\link(\Sigma)=\link(\Delta\star\Pi)\star\Pi'.$

\item \label{item:link_conn_intro} If the simplex $\Delta$ is not a co-dimension 1 face of a maximal simplex, then $\link(\Delta)$ is connected.

\end{enumerate}

Then $G$ is a hierarchically hyperbolic group.
\end{thmintro}

In the statement, $\star$ denotes the simplicial join, and we are allowing the empty simplex in the various conditions, so that in particular $X=\link(\emptyset)$ is hyperbolic. For more explanation of the statement and comments on the various conditions, see Section \ref{sec:propaganda}.

A more general version is given in Theorem \ref{thm:hhs_links} below, and both Theorem \ref{thm:hhs_links} and 
Theorem \ref{thm:propaganda} give further details about the hierarchically hyperbolic structure.

The setup of Theorem~\ref{thm:main_intro}, and in particular the idea 
of considering a hyperbolic complex with hyperbolic links, is inspired by
the curve complex $\mathcal C(S)$ of a surface $S$. Links of
simplices in $\mathcal C(S)$ are related to curve complexes of
subsurfaces of $S$ by the following observations, which hold with some
low-complexity exceptions.  The link of a vertex is the curve graph of
the complement of the corresponding curve.  More generally, for each
subsurface $U$ of $S$, the curve complex $\mathcal C(U)$ arises as the
link of a simplex in $\mathcal C(S)$.  Conversely, the link of each
simplex is either a nontrivial join (hence bounded), or is the curve
graph of a subsurface.

However, curve graphs are the hyperbolic complexes that witness the hierarchical hyperbolicity of 
pants graphs, not mapping class groups. This is because annular curve graphs are not links of simplices 
in the curve graph; in fact, $\mathcal C(S)$ does not even contain the vertex set of $\mathcal C(U)$ when $U$ 
is an annulus. For a discussion of the hyperbolic complex that witnesses hierarchical hyperbolicity of the 
mapping class group, see Section~\ref{subsec:blow_up}.  

For the pants graph, checking the conditions of Theorem \ref{thm:hhs_links} (which are similar to those 
of Theorem \ref{thm:main_intro}) only uses hyperbolicity of curve graphs \cite{MM_I}, the fact 
\cite[Lemma 2.3]{MM_II} that subsurface projections are coarsely Lipschitz where defined (this is to check the 
quasi-isometric embedding condition), and arguments involving subsurfaces filled by multicurves to check 
the combinatorial conditions.  We leave the details to the reader since pants graphs are already known to be 
hierarchically hyperbolic~\cite[Theorem G]{HHS_I}.

Some applications of Theorem~\ref{thm:main_intro} are discussed below, but we expect it (and the more 
general version, Theorem~\ref{thm:hhs_links}) to have many other applications.  In fact, 
Theorem~\ref{thm:hhs_links} is used in \cite{Veech:hhs} to verify hierarchical hyperbolicity of
surface extensions of naturally-occurring subgroups of mapping class groups, namely lattice Veech groups. Moreover, in \cite{HHS:artin}, Theorem~\ref{thm:hhs_links} is used  to show that Artin groups of extra-large type are hierarchically 
hyperbolic groups, and in \cite{HHS:graph-man}, the combinatorial HHS viewpoint 
is used to prove that fundamental groups of graph manifolds are hierarchically hyperbolic \emph{groups}, not 
merely hierarchically hyperbolic \emph{spaces} (as was previously shown in~\cite{HHS_II}).

We emphasize that the combinatorial HHS viewpoint encapsulated in Theorem~\ref{thm:hhs_links} operates in 
tandem with the definition of hierarchical hyperbolicity from~\cite{HHS_II}, rather than replacing it.  By this 
we mean the following:
\begin{itemize}
 \item Theorem~\ref{thm:hhs_links} is a good way to prove that a space/group is hierarchically hyperbolic, but 
many of the consequences of hierarchical hyperbolicity --- for example, finite asymptotic 
dimension~\cite{HHS_III}, control of quasiflats~\cite{HHS:quasiflats}, ``coarse 
rank-rigidity'' and the omnibus subgroup theorem~\cite{HHS:boundary,PetytSpriano}, the Tits 
alternative~\cite{HHS:boundary,HHS:corrigendum}, uniform exponential growth~\cite{AbbottNgSpriano:uniform}, 
etc. rely on the coarse geometric machinery built on the original definition.  The combinatorial viewpoint 
does not currently allow us to prove geometric results of this type directly.  Even establishing much more 
basic geometric properties of a combinatorial HHS --- as in Theorem~\ref{thm:main_intro} and its proof --- 
takes 
considerable work.
\item More fundamentally, the main tools of hierarchical hyperbolicity --- notably, the distance formula --- 
aren't readily extractable from the combinatorial viewpoint on its own.  
\item The property of being a combinatorial HHS seems to be strictly stronger than the property of being an 
HHS.  This can be illustrated using fairly artificial CAT(0) square complex examples to show that the natural 
candidate combinatorial HHS structure one might try to build from an HHS structure does not work without 
additional hypotheses on the HHS.  A partial converse of the form ``hierarchically hyperbolic groups satisfying natural additional conditions are combinatorially hierarchically hyperbolic" is the subject of ongoing work, so we will not deal in detail with the above mentioned example here.  Suffice it to say that it involves an infinite CAT(0) square complex with trivial automorphism group, which is a hierarchically hyperbolic space in view of \cite{HHS_I} but for which the natural candidate combinatorial hierarchically hyperbolic structure involves an underlying simplicial complex in which the set of links of simplices, partially ordered by inclusion, contains arbitrarily long chains.
\end{itemize}

Nonetheless, someone wishing to establish a property of some group $G$ that is already known for 
hierarchically hyperbolic groups can (with luck) now do so with no engagement whatsoever with the HHS 
machinery: they could instead try to build a combinatorial HHS on which $G$ acts.  A guide on how to use 
Theorem~\ref{thm:hhs_links} for further application appears in Section~\ref{subsec:user} and 
Section~\ref{subsec:blow_up}.

\subsection*{(Hierarchically) hyperbolic quotients of mapping class groups}
The main application of Theorem~\ref{thm:main_intro} presented in this paper is the study of certain quotients 
of mapping class groups.  Below, we first discuss the natural quotients obtained by modding out high powers of 
Dehn twists.  We then discuss a construction of non-elementary hyperbolic quotients, which works under the 
assumption that specific hyperbolic groups encountered in the construction are residually finite (and without 
this assumption in genus $2$).

\subsubsection*{Quotients by powers of Dehn
twists}\label{subsubsec:intro_kill_dehn} We first apply
Theorem~\ref{thm:main_intro} to quotients of mapping class groups by
large powers of Dehn twists, further advancing the technology developed by Dahmani in
\cite{D_PlateSpinning} to resolve Ivanov's deep relations question~\cite[Section 12]{Ivanov_list}, 
and the extensions of that technology from~\cite{DFDT}.

We now briefly survey the history of the study of quotients of mapping class groups by powers of Dehn twists. This dates to at least 1974 where it appears in 
Birman's classic monograph  
\cite{Birman:Braids}. In that text, Birman notes that for 
the closed genus two surface the normal closure of the squares of 
Dehn twists is of index $6!$ in the mapping class group; she then 
asked whether the index is finite or infinite for arbitrary genus. 
Humphries later resolved this for the normal subgroups generated by squares 
or by cubes of Dehn twists (finite for 
the closed or once punctured surfaces of genus two or three; 
otherwise infinite), see \cite{Humphries:powersdehn}. Humphries also 
showed that for the genus two surface with any number of punctures 
and any power of Dehn twists greater than 3, the corresponding 
quotient group was infinite. Eventually, Funar proved that as long as 
the genus is at least 3, then the quotient by the normal subgroups 
along powers of Dehn twists are infinite, as long as they are at 
least 13th powers \cite{Funar}. Funar was interested in these 
subgroups because of their connection with TQFT representations. See 
also \cite{AramayonaFunar}.

Another, more recent, motivation for studying these quotients comes
from the algebraic counterpart of Thurston's Dehn filling theorem in
the context of relatively hyperbolic groups
\cite{Os-perfill,GrMa-perfill}. This algebraic version has numerous important
applications, including a role in the proof of the Virtual Haken conjecture
\cite{vhak} and in the solution of the isomorphism problem for certain
relatively hyperbolic groups \cite{DG:recognize_DF, DT-DecidingIsom}.
Mapping class groups are not non-trivially relatively hyperbolic
except in very low complexity \cite{AAS:MCG_not_rel_hyp, BDM:thick},
but, in a number of ways, the subgroups generated by Dehn
twists around curves in a pants decomposition play a role analogous 
to that of peripheral subgroups.

In \cite{DFDT}, it is proven that quotients of mapping class groups by
large powers of Dehn twists are acylindrically hyperbolic, providing
an analogue of the Dehn filling theorem.  However, while acylindrical
hyperbolicity has many interesting consequences, it only captures 
part of the geometry of mapping class groups.  One would hope that
performing Dehn fillings preserves much more of the hierarchically 
hyperbolic structure of the mapping class group --- indeed, in Theorem 
\ref{thmi:MCG} below we will establish that the quotients are also 
hierarchically hyperbolic.

Given a surface $S$, we denote by $DT_K$ the normal subgroup generated
by all $K$-th powers of Dehn twists.  Using Theorem
\ref{thm:main_intro}, we prove:

\begin{thmintro}\label{thmi:MCG}
 Let $S$ be a finite-type surface. Then there exists $K_0>0$ so that for any non-zero multiple $K$ of $K_0$, 
$MCG(S)/DT_K$ is an infinite hierarchically hyperbolic group.
\end{thmintro}

In fact, we provide an explicit HHS structure, which is described in
Theorem \ref{thm:MCG_mod_Dehn}.  Corollaries of hierarchical 
hyperbolicity for these groups include: 
finiteness of the asymptotic dimension of $MCG(S)/DT_K$
\cite{HHS_III}, uniform exponential growth
\cite{AbbottNgSpriano:uniform}, and (using the description of the HHS
structure) that the maximal dimension of quasiflats is $\lfloor
(3g-3)/2\rfloor$ and each quasi-flat of that dimension is a union of
finitely many standard orthants \cite{HHS_I,HHS:quasiflats}.  The
latter result might be useful to prove quasi-isometric rigidity, which
we formulate as a question:

\begin{questintro}\label{quest:qi_rigid}
 For $K$ as in Theorem \ref{thmi:MCG}, is $MCG(S)/DT_K$ quasi-isometrically rigid?
\end{questintro}

A strategy to give a positive answer to Question \ref{quest:qi_rigid} could involve adapting arguments 
from \cite[Section 5]{HHS:quasiflats} and from \cite{Bowditch:WP_rigid}, 
and proving a combinatorial rigidity result whereby one shows that the automorphism group of a suitable simplicial complex coincides with the desired group.\footnote{After this paper was originally circulated, this approach has been successfully implemented for punctured spheres in \cite{MangioniSisto}.} Using the lifting techniques from 
\cite{DFDT} that we develop further in Subsection \ref{subsec:lifting_gons}, it might be possible to prove 
such a result by reducing it to combinatorial rigidity of a suitable complex on which $MCG(S)$ acts. We pointed 
to \cite{Bowditch:WP_rigid} because we expect the structure of $MCG(S)/DT_K$ to more closely resemble that of 
the Weil-Petersson metric than that of the mapping class group itself, since, roughly speaking, we made the 
annular curve graphs bounded and hence irrelevant.

\subsubsection*{Hyperbolic quotients}\label{subsubsec:intro_hyperbolic_quotient}
The case of the genus 2 surface is notable in that quotients by powers of Dehn twists are not only hierarchically hyperbolic, but they are in fact relatively hyperbolic, as we will see by applying a result of Russell \cite{Russell:rel_HHS} to the HHS structure on the quotients:

\begin{thmintro}\label{thmi:genus_2}
  There exists $K_0\ge1$ so that for all non-zero multiples $K$ of $K_0$, the following holds. The 
quotient $MCG(\Sigma_2)/DT_K$ is hyperbolic relative to an infinite index subgroup commensurable to the product of two 
$C'(1/6)$-groups.
\end{thmintro}

We note that the theorem can be used in conjunction with relatively hyperbolic Dehn filling \cite{Os-perfill,GrMa-perfill} to produce many hyperbolic quotients of $MCG(\Sigma_2)$, as stated in the following corollary, which we deduce from Theorem \ref{thmi:genus_2} in Remark
\ref{rem:Francesco_trick}.  In said remark, we also point out a
different construction relying on results of a very different nature
and on a trick suggested to us by Francesco Fournier Facio. 

\begin{corintro}\label{MCG_2FRH}
 $MCG(\Sigma_2)$ is fully residually non-elementary hyperbolic.
\end{corintro}

Recall that a group is \emph{fully residually $\mathcal P$} if for every finite collection of elements, there 
is a  quotient satisfying $\mathcal P$ into which the collection injects. 

For higher genus, we cannot apply the relatively hyperbolic Dehn filling. However, we now outline how one could construct
hyperbolic quotients using HHS machinery instead.  After modding out powers of Dehn
twists, we are left with a hierarchically hyperbolic group that has
strictly lower complexity (in the HHS sense, the exact meaning of
which is immaterial for this discussion).  At the bottom of the
hierarchy, we find a collection of hyperbolic groups.  The idea to
further reduce complexity is to repeat the previous procedure, namely
modding out deep finite-index subgroups of these hyperbolic groups.
We proceed inductively, reducing complexity, until we are left with a
hierarchically hyperbolic group of complexity $1$.  A general fact
about hierarchically hyperbolic groups is that when the complexity is
$1$, the group is hyperbolic.  We show that this construction works
provided that all hyperbolic groups encountered at the various stages
are residually finite. In Theorem \ref{thmi:MCG_quotient_HHG_intro}
below we formulate this under the assumption that all
hyperbolic groups are residually finite, which is much stronger than what we strictly need. We do not know
whether there are enough residually finite hyperbolic groups to run our 
construction. Theorem \ref{thmi:MCG_quotient_HHG_intro}, as
well as the other theorems below, can be seen as a route towards establishing the existence of non-residually finite hyperbolic
groups or as an invitation to find a suitable class of hyperbolic
groups that are residually finite.  We discuss this further below.

\begin{thmintro}\label{thmi:MCG_quotient_HHG_intro}
   Let $S$ be a connected orientable surface of finite type of complexity at least 2.  If all 
hyperbolic groups are residually finite, then $MCG(S)$ is fully residually non-elementary hyperbolic.
\end{thmintro}

As seen in Corollary~\ref{MCG_2FRH}, our techniques prove that 
the mapping class group of 
the closed genus--$2$ surface is fully
residually non-elementary hyperbolic without a residual finiteness
assumption.

\begin{remintro}\label{rem:conditional}
In the above theorem, and the next two theorems, we condition the conclusion on residual 
finiteness of all hyperbolic groups.  This is for simplicity in formulating the theorems.  In reality, during 
the proofs of these theorems, one encounters specific hyperbolic groups, and it is only for these particular groups whose 
residual finiteness is necessary.  Therefore, one can interpret these theorems as potential ways of proving 
the existence of a non-residually finite hyperbolic group (by exhibiting a mapping class group with no 
non-elementary hyperbolic quotients, say), but we prefer to see these theorems as an invitation to study the 
residual finiteness question for the specific hyperbolic groups encountered in the proofs.
\end{remintro}

 We use the flexibility of the construction 
described above to obtain the necessary multitude of hyperbolic quotients, as 
described in Theorem~\ref{thm:MCG_quotient}.

\begin{remintro}\label{rem:strong_T}
The following was pointed out by Dawid Kielak.  If $MCG(S)$ admits a
non-elementary hyperbolic quotient, then $MCG(S)$ admits an affine
isometric action on an $L^p$ space with unbounded
orbits~\cite{Nica,Yu} for some $1<p<+\infty$, and in particular it
does not have property $F_{L^p}$ for said $p$.  (In contrast, 
$SL_n(\mathbb Z)$ does have $F_{L^p}$ for all
$1<p<+\infty$ \cite{BFGM}.)  In view of
Theorem~\ref{thmi:MCG_quotient_HHG_intro}, if some mapping class group
has $F_{L^p}$, say for all $1<p<+\infty$, then there is a
non-residually finite hyperbolic group.
\end{remintro}

The flexibility of our construction of quotients can also be exploited to prove further results, as we now discuss.  Recall 
from~\cite{FarbMosher} that $H\leq MCG(S)$ is 
\emph{convex-cocompact} if some $H$--orbit in the Teichm\"uller 
space of $S$ is quasiconvex. There are several characterizations of 
convex-cocompactness; see~\cite{KentLeininger:convex_cocompact, Hamenstadt, 
DurhamTaylor:stability}. One reason this notion is interesting is its connection with hyperbolicity of 
fundamental groups of surface bundles over surfaces \cite{FarbMosher,Hamenstadt}.

Reid posed the question of whether convex-cocompact subgroups of 
$MCG(S)$ are separable \cite[Question 3.5]{Reid:separable}. Recall that a subgroup $H<G$ is separable if for every $x\in G-H$ there 
exists a finite group $F$ a surjective homomorphism $\phi:G\to F$ with $\phi(x)\notin \phi(H)$.

Note that in general $MCG(S)$ contains non-separable subgroups. In fact this is already the 
case for $MCG(S_{0,5})$~\cite{LeiningerMcReynolds}. (Nonetheless, various 
geometrically natural subgroups, e.g. curve-stabilizers, are known to be 
separable in $MCG(S)$~\cite{LeiningerMcReynolds}.) 

Our techniques reduce Reid's question to residual finiteness of certain hyperbolic groups, which we formulate 
as:

\begin{thmintro}\label{thmi:separable}
     Let $S$ be a connected orientable surface of finite type of complexity at least 2.  If all hyperbolic 
groups are residually finite, then 
every convex-cocompact subgroup of $MCG(S)$ is separable.
\end{thmintro}

The proof of Theorem~\ref{thmi:separable} relies on the hyperbolic 
quotients of $MCG(S)$ arising in the proof of 
Theorem~\ref{thmi:MCG_quotient_HHG_intro}. The assumption about residual finiteness is 
invoked again to apply a result of~\cite{AGM} (namely, if all hyperbolic groups are residually finite then all 
hyperbolic groups are QCERF).  Again, we do not really need residual finiteness of all hyperbolic groups, 
just the ones encountered in our construction and in the (iterated Dehn filling) construction from~\cite{AGM}.

The next application relates to a question of Bridson--Reid--Wilton.  In fact, we reduce~\cite[Question 
5.1]{BridsonReidWilton} to the questions of residual finiteness of certain hyperbolic groups and of the congruence subgroup property for mapping class groups.  

\begin{thmintro}\label{thmi:basically_omnipotence}
Let $S$ be a connected orientable surface of finite type of complexity at least 2.  If all hyperbolic groups 
are residually finite, then the 
following holds. Let $g,h\in MCG(S)$ be pseudo-Anosovs with no common proper power, and let $q\in \mathbb 
Q_{>0}$. Then there 
exists a finite group $G$ and a homomorphism $\psi:MCG(S)\to G$ so that $ord(\psi(g))/ord(\psi(h))=q$, where 
$ord$ denotes 
the order.
\end{thmintro}

The property established in the theorem is called \emph{omnipotence for pseudo-Anosovs}. In \cite{BridsonReidWilton}, the authors study profinite rigidity of $3$--manifold groups using the notion of \emph{$\pi_1(\Sigma)$-congruence} omnipotence, where there is the additional requirement that the finite quotients are congruence quotients. It is not known whether mapping class groups have the congruence subgroup property, see \cite[Problem 2.10]{Kirby_list}\cite[Conjecture on page 75]{Ivanov_list}, but if so, then the two notions of omnipotence for pseudo-Anosovs are equivalent.

In \cite{BridsonReidWilton} it is shown that a positive answer to their Question 5.1 implies the following.  Let $M$ be a closed hyperbolic 3-manifold with 
first Betti number $1$.  Let $N$ be a compact 3-manifold so that $\pi_1M$ and $\pi_1N$ have isomorphic 
profinite completions. Then $M$ and $N$ have a common finite cyclic cover.

A heuristic discussion of the proof of Theorem~\ref{thmi:MCG_quotient_HHG_intro} is given in 
Section~\ref{subsec:MCG_proof_outline}. The proof of Theorem~\ref{thmi:MCG_quotient_HHG_intro} is essentially self-contained, only using the statement of Theorem \ref{thm:main_intro}.

\subsection*{Speculations}

As mentioned above, Theorems \ref{thmi:MCG_quotient_HHG_intro}, \ref{thmi:separable}, and 
\ref{thmi:basically_omnipotence} hold provided that ``sufficiently many'' hyperbolic groups are residually 
finite, and therefore there are two natural research directions to explore. The first is to use those results 
to show that there exist hyperbolic groups that are not residually finite. Consider, for example, the 
following question:

\begin{questintro}\label{quest:quotient}
 Do all mapping class groups of closed oriented surface of genus at least 1 have an infinite hyperbolic quotient?
\end{questintro}

In view of Theorem \ref{thmi:MCG_quotient_HHG_intro}, a negative 
answer to said question shows that there exists some 
non-residually-finite hyperbolic group. Similarly, the same is true if either the question by Reid ~\cite[Question 3.5]{Reid:separable} or that by Bridson--Reid--Wilton \cite[Question 5.1]{BridsonReidWilton} have a negative answer.

(A priori, having an infinite hyperbolic quotient is weaker than being fully residually non-elementary hyperbolic. However, in our context these properties are equivalent in view of the argument in Remark \ref{rem:Francesco_trick}, since mapping class groups do not have two-ended quotients.)

In the other direction, one might try to find a suitable class of hyperbolic groups that can be shown to be 
residually finite, and are sufficient to prove the theorems above without additional assumptions. This would 
mirror the developments that led to the proof of the virtual Haken
conjecture.  Indeed, if all hyperbolic groups were residually finite,
then the virtual Haken conjecture would follow from all hyperbolic
groups being in fact QCERF \cite{AGM}, the existence of quasiconvex
surface subgroups \cite{KahnMarkovic}, and the connection between
separability and embedding into a finite cover
\cite{scott1978subgroups}.  Also, although the proof does not follow
exactly this template, the conjecture was eventually proven by showing
that cubulated hyperbolic groups are virtually special \cite{vhak},
which implies that they are (QCERF, whence) residually
finite \cite{HaglundWise:special}.

Further, provided that residual finiteness issues are resolved, it would be interesting to determine whether 
the hyperbolic quotients obtained via our construction can be CAT(0). We believe that this is a natural 
question given that actions of mapping class groups on CAT(0) spaces are very constrained 
\cite{Bridson:MCG_not_CAT(0),KapovichLeeb:actions}. One can even ask:

\begin{questintro}
 Can an infinite hyperbolic quotient of the mapping class group of a closed oriented surface of sufficiently high genus be CAT(0)?
\end{questintro}

\subsection*{Largest hyperbolic quotients}
For the purposes of the following discussion, assume that all hyperbolic groups are residually finite.  It is 
natural to wonder whether the hyperbolic groups we construct in the proof of Theorem~\ref{thm:MCG_quotient} are 
the ``largest possible'', meaning that any hyperbolic quotient of the mapping class group is a quotient of one 
of them. 

We do not believe this to be true, because, in our quotients, too many stabilizers of subsurfaces have finite 
image.  However, we formulate a conjecture related to this below. 

Fix a closed surface $S$ of genus at least $2$, and let $\mathcal Y$ be the collection of all (isotopy classes of) essential subsurfaces $Y$ of $S$ so that
there exists a mapping class $g$ with $Y$ and $gY$ disjoint and not
isotopic.  For example, an annulus around a non-separating
curve is in $\mathcal Y$, while the annulus around the ``middle curve'' of
the genus-2 closed surface is not.  

\begin{remintro}\label{rem:terminology}
The set $\mathcal Y$ coincides with the set of non-\emph{$MCG(S)$--overlapping} subsurfaces in the sense of~\cite{CMM:normal}, as well as with the set of \emph{nondisplaceable} subsurfaces in the sense 
of~\cite[Theorem 1]{HQR:big}.
\end{remintro}

For a subsurface $Y\subseteq S$,
denote by $MCG(Y|S)$ the subgroup of $MCG(S)$ consisting of all
mapping classes supported on $Y$.  The conjecture is:

\begin{conjintro}\label{conj:duplicable}
 Let $S$ be a closed surface of genus at least 3.  Then there are epimorphisms $\phi:MCG(S)\to G$, with $G$ 
hyperbolic, such 
that the following holds for any essential subsurface $Y\subseteq S$: the group $\phi(MCG(Y|S))$ is finite if 
and only if 
$Y\in\mathcal Y$.
\end{conjintro}

The statement of the conjecture might require some adjustments.  We believe the conjecture to be at least 
morally correct 
provided that no residual finiteness issues arise.

The conjecture is inspired by the observations below, which show that if $Y\in \mathcal Y$,  then 
$MCG(Y|S)$ 
becomes virtually cyclic in every hyperbolic quotient of $MCG(S)$. 

Since every infinite virtually cyclic group surjects onto 
$\mathbb Z/2\mathbb Z$, this implies that $MCG(Y|S)$ actually becomes finite at least if $Y$ has genus at least 3, but 
we believe that 
with additional arguments this can also be shown in lower genus.

Let $Y\in\mathcal Y$, with corresponding mapping class
$g$.  Suppose that the epimorphism $\phi:MCG(S)\to G$, where $G$ is
hyperbolic, is so that $H=\phi(MCG(Y|S))$ is infinite (otherwise we
are done).  Then $H$ contains an infinite order element, say $h$.
Moreover, $\phi(g) H\phi(g)^{-1}$ is contained in the centralizer of
$h$ (notice that $gMCG(Y|S)g^{-1}=MCG(gY|S)$, and that $MCG(Y|S)$
commutes with $MCG(gY|S)$).  The centralizer of $h$ is virtually
cyclic, and hence so is $H$.

\medskip

\textbf{Outline}
In Section \ref{sec:setup_statement}, we introduce the notions needed to state Theorem \ref{thm:hhs_links}, in 
particular the notion of a combinatorial HHS and state the theorem. 
We conclude the section with some remarks 
that might be of use to the reader wishing to apply the theorem, and we 
also include a simple example.

In Section \ref{sec:background} we recall the definition of HHS, and the (very few) results needed for this paper.

\textbf{Sections \ref{sec:y_hyp}--\ref{sec:propaganda} contain the proofs of Theorems \ref{thm:hhs_links} and \ref{thm:main_intro}, but only the statement of Theorem \ref{thm:main_intro} is used in the subsequent sections.} 

In Section \ref{sec:y_hyp}, we prove that the spaces used to define HHS projections are hyperbolic. This is 
crucial to prove that the candidate projections to the various hyperbolic spaces behave as expected.  In Section \ref{sec:hieromorphisms}, we study induced combinatorial HHS structures on links, which will enable inductive arguments. In Section \ref{sec:check_axioms}, we complete the proof of Theorem \ref{thm:hhs_links} by checking the HHS axioms. At this stage, the hardest part of the proof will be the Uniqueness axiom.

In Section \ref{sec:propaganda}, we show that the conditions in Theorem \ref{thm:main_intro} imply those in 
Theorem \ref{thm:hhs_links}.  At this point, we have proved Theorem~\ref{thm:main_intro}, and we move on to 
the mapping class group applications.

\textbf{Sections \ref{sec:MCG} and \ref{sec:proof_quotients} focus on quotients of mapping class groups.}

In Section \ref{sec:MCG}, we state Theorem~\ref{thm:MCG_quotient}, about the hierarchically hyperbolic 
structures on our quotients of mapping class groups and deduce Theorem~\ref{thmi:MCG} 
(Theorem~\ref{thm:MCG_mod_Dehn}), Theorem~\ref{thmi:genus_2} (Corollary~\ref{cor:genus_2_non_intro}),  
Theorem~\ref{thmi:MCG_quotient_HHG_intro} (Corollary~\ref{cor:fully_residually_hyperbolic}), 
Theorem~\ref{thmi:separable} (Corollary~\ref{cor:separable}), and Theorem~\ref{thmi:basically_omnipotence} 
(Corollary~\ref{cor:basically_omnipotence}).  

In Section~\ref{sec:proof_quotients}, we prove Theorem~\ref{thm:MCG_quotient}.  Here, we improve the lifting 
technology of \cite{DFDT} and combine it with Theorem \ref{thm:main_intro}. In fact, we expect many of the 
new technical lemmas in this section to be useful for other Dehn filling-type theorems, possibly even outside 
the hierarchically hyperbolic context.  Section~\ref{sec:proof_quotients} begins with a fairly detailed 
heuristic outline of the proof.

\medskip

\textbf{Remark.}
This paper arose from two separate projects, which are naturally linked and we therefore merged. The 
results in Sections~\ref{sec:y_hyp}--\ref{sec:propaganda}, in particular Theorems \ref{thm:hhs_links} and 
\ref{thm:main_intro}, are due to MH, AM, and AS. The remaining results are due to all four authors.

\subsection*{Acknowledgements}
We are very grateful to Fran\c{c}ois Dahmani for numerous essential discussions.  A similar approach to the 
construction of further quotients of mapping class groups was discussed with him and was instrumental to this 
work, as have been his explanations of technical points in his paper \cite{D_PlateSpinning}.  We thank Ian Agol, Matt Clay, Francesco Fournier Facio, Dawid Kielak, Dan 
Margalit, Michah Sageev, and Alex Wright for useful comments and suggestions. We also thank Carolyn Abbott, 
Daniel Berlyne, Thomas Ng, Alex Rasmussen, and Jacob Russell who together 
carefully read through an early version of the paper with a very 
sharp eye and provided us with extremely helpful feedback.  We are also extremely grateful to the referees for a very large number of comments, and 
 several corrections, that have enormously improved the paper.

The authors would like to thank the Isaac Newton Institute for Mathematical Sciences, Cambridge, for support 
and hospitality during the programme \emph{Non-positive curvature, group actions, and cohomology} where work on this paper 
was undertaken. This work was supported by EPSRC grant no EP/K032208/1. We thank the International Centre for 
Mathematical Sciences (ICMS) for their hospitality during a 
Research-in-Groups in which much of this work was done.  Behrstock was supported by NSF grant DMS-1710890. Hagen was 
partially supported by EPSRC New Investigator Award  EP/R042187/1.  Martin was partially supported by EPSRC 
New Investigator Award EP/S010963/1. Sisto  was partially supported by the Swiss National 
Science Foundation (grant \#182186).

\section{Combinatorial HHS}\label{sec:setup_statement}
We first state Theorem~\ref{thm:hhs_links} and supply tools for its 
proof. In Subsection~\ref{subsec:user}, we illustrate how these tools 
work, via an explicit example (Section~\ref{subsub:AFP}); it might be instructive for the reader 
to refer to that for motivation.  

\subsection{Basic definitions}\label{subsec:setup} Let $X$ be a flag simplicial 
complex.   

\begin{defn}[Join, link, star]\label{defn:join_link_star}
Given disjoint simplices $\Delta,\Delta'$ of $X$, we let $\Delta\star\Delta'$ denote the simplex spanned by 
$\Delta^{(0)}\cup\Delta'^{(0)}$, if it exists.  More generally, if $K,L$ are disjoint induced subcomplexes of 
$X$ such that every vertex of $K$ is adjacent to every vertex of $L$, then 
$K\star L$ is the induced subcomplex with vertex set $K^{(0)}\cup L^{(0)}$.  We refer to $K\star L$ as the \emph{join} of $K$ and $L$.  

For each simplex $\Delta$, the \emph{link} $\link(\Delta)$ is the union of 
all simplices $\Sigma$ of $X$ such that $\Sigma\cap\Delta=\emptyset$ and $\Sigma\star\Delta$ is a simplex of $X$.  Observe 
that $\link(\Delta)=\emptyset$ if $\Delta$ is a maximal simplex.  Conversely, if $\link(\Delta)=\emptyset$, then $\Delta$ is 
not properly contained in a simplex, i.e., $\Delta$ is maximal.  More generally, if $K$ is an induced subcomplex of $X$, then 
$\link(K)$ is the union of all simplices $\Sigma$ of $X$ such that $\Sigma\cap K=\emptyset$ and $\Sigma\star K$ is a 
subcomplex of $X$.

The \emph{star} of $\Delta$ is $Star(\Delta)=\link(\Delta)\star\Delta$, i.e., the union of all simplices of $X$ 
that contain 
$\Delta$.

We often refer to $0$--simplices as \emph{vertices} and $1$--simplices as \emph{edges}, and make no distinction between $1$--dimensional simplicial complexes and simplicial graphs.  In a (not necessarily simplicial) graph or a simplicial complex $Y$, we use the term \emph{open star} of a vertex $v$ to refer to the union of $\{v\}$ with all open simplices (or open edges) of $Y$ whose closures contain $v$.  (By \emph{open simplex} or \emph{open edge}, we mean the image of the restriction to the interior of a cell of the appropriate characteristic map, which need not be open in $X$.)  We sometimes refer to \emph{removing the open star} of $v$, the result of which is an induced subcomplex of $Y$ consisting of exactly those simplices (or edges) that do not contain $v$.

We emphasize that $\emptyset$ is a simplex of $X$, whose link is all of $X$ and whose star is all of $X$.  
\end{defn}

\begin{defn}[$X$--graph, $W$--augmented complex]\label{defn:X_graph}
An \emph{$X$--graph} is a graph $W$ whose vertex set is the set of all maximal simplices of 
$X$.  

For a flag complex $X$ and an $X$--graph $W$, the \emph{$W$--augmented graph} 
$\duaug{X}{W}$ is the graph defined as follows:
\begin{itemize}
     \item the $0$--skeleton of $\duaug{X}{W}$ is $X^{(0)}$;
     \item if $v,w\in X^{(0)}$ are adjacent in $X$, then they are adjacent in $\duaug{X}{W}$; 
	
	 \item if two vertices in $W$ are adjacent, then we consider 
	 $\sigma,\rho$, the associated maximal simplices of $X$, and 
	 in $\duaug{X}{W}$ we connect each vertex of $\sigma$ to each vertex 
	 of $\rho$.
\end{itemize}
We equip $W$ with the usual path-metric, in which each edge has unit length, and do the same for $\duaug{X}{W}$.
\end{defn}

We are aiming to construct a hierarchically hyperbolic structure $(W,\mathfrak S)$.  The actually hierarchically hyperbolic space will be the graph $W$, equipped with the usual path-metric.  The ``curve graph'' (i.e., the hyperbolic space associated to the unique $\nest$--maximal element of 
$\mathfrak S$) will be $X^{+W}$.

\begin{rem*}[Connectedness of $W$]
A priori, there is no assumption that $W$ is connected.  In practice, connectedness of $W$ will be deduced using the other parts of the definition of a combinatorial HHS (Definition~\ref{defn:combinatorial_HHS} below).  Specifically, during the proof of Theorem~\ref{thm:hhs_links}, we verify that the links of $X$ provide a hierarchically hyperbolic structure for $W$, and, during the part of that proof where the ``uniqueness axiom'' for hierarchically hyperbolic spaces (Definition~\ref{defn:HHS}.\eqref{item:dfs_uniqueness}) is checked, we verify that $W$ is connected.  This relies on the fact --- coming from Definition~\ref{defn:combinatorial_HHS} --- that various auxiliary graphs related to $W$ and $X$ are hyperbolic when given the usual path-metric, and in particular connected.  It also uses induction on the ``complexity'' $n$ of $X$ (see Definition~\ref{defn:combinatorial_HHS}.\eqref{item:chhs_flag} and Definition~\ref{defn:CHHS-finite-complexity}), to say that various subgraphs $W^\Delta$ of $W$ associated to links of simplices $\Delta$ in $X$ are hierarchically hyperbolic and, in particular, connected).

Probably the strongest reason not to simply hypothesize connectedness of $W$ is that, in support of the above induction, we will need to verify that for each nonempty non-maximal simplex $\Delta$ of $X$, the subgraph $W^\Delta$ of $W$ spanned by maximal simplices of the form $\sigma\star\Delta$ is a $\link(\Delta)$--graph that combines with the complex $\link(\Delta)$ to form a combinatorial HHS of strictly lower complexity.  This is Proposition~\ref{prop:compatibility_1}.  By not including connectedness of $W$ in the definition of an $X$--graph, we avoid having to verify connectedness of $W^\Delta$ when proving Proposition~\ref{prop:compatibility_1}, and instead are able to assume it as an inductive hypothesis when proving connectedness of $W$ later, in Theorem~\ref{thm:hhs_links}.
\end{rem*}

\begin{defn}[Equivalent simplices, saturation]\label{defn:simplex_equivalence}
For $\Delta,\Delta'$ simplices of $X$, we write $\Delta\sim\Delta'$ to mean
$\link(\Delta)=\link(\Delta')$. We denote by $[\Delta]$ the $\sim$--equivalence class of $\Delta$.  Let 
$\sat(\Delta)$ denote 
the set of vertices $v\in X$ for which there exists a simplex $\Delta'$ of $X$ such that $v\in\Delta'$ and 
$\Delta'\sim\Delta$, i.e., $$\sat(\Delta)=\left(\bigcup_{\Delta'\in[\Delta]}\Delta'\right)^{(0)}.$$
We refer to $\sat(\Delta)$ as the \emph{saturation} of $\Delta$.
We denote by $\mathfrak S$ the set of $\sim$--classes of \textbf{non-maximal} simplices in 
$X$.
\end{defn}

\begin{rem}\label{rmk:lk_sat}
 Notice that $\sat(\Delta)\subseteq \link(\link(\Delta))^{(0)}$ (any vertex in $\sat(\Delta)$ is contained in a simplex all 
of whose vertices are connected to any vertex in $\link(\Delta)$).
 Also, we have $\link(\sat(\Delta))=\link(\Delta)$ (indeed, since $\Delta^{(0)}\subseteq \sat(\Delta)$, we have the inclusion 
``$\subseteq$'', while on the other hand any vertex connected to all vertices of $\Delta$ is connected to all elements of 
$\sat(\Delta)$, giving the other inclusion).
\end{rem}

\begin{defn}[Complement, link subgraph]\label{defn:complement}
Let $W$ be an $X$--graph.  For each simplex $\Delta$ of $X$, the \emph{complement subgraph} $Y_\Delta$ is the subgraph of $\duaug{X}{W}$ induced by the set
$(\duaug{X}{W})^{(0)}-\sat(\Delta)$ of vertices.

The \emph{augmented link} $\mathcal C(\Delta)$ of $\Delta$ is the induced subgraph of $Y_\Delta$ spanned by $\link(\Delta)^{(0)}$.  Note that 
$\mathcal 
C(\Delta)=\mathcal C(\Delta')$ whenever $\Delta\sim\Delta'$.  We emphasize 
that we are taking links in $X$, not in $X^{+W}$, and then considering the 
subgraphs of $Y_\Delta$ induced by those links.  
\end{defn}

(The notation $\mathcal C(\Delta)$ is chosen since these spaces will be the underlying hyperbolic spaces in a hierarchically hyperbolic structure for $W$; compare Definition~\ref{defn:HHS}.  The use of ``$\mathcal C$'' was originally motivated by usage in concrete examples from~\cite{HHS_I}: \textbf{c}urve graphs in the setting of mapping class groups, and \textbf{c}ontact graphs in the setting of CAT(0) cube complexes.)

\begin{defn}\label{defn:CHHS-finite-complexity}
The simplicial 
complex $X$ has \emph{finite complexity} if there exists $n\in\naturals$ so that any chain $\link(\Delta_1)\subsetneq\dots\subsetneq\link(\Delta_i)$, where each $\Delta_j$ is a simplex of $X$, has length at 
most $n$; the minimal such $n$ is the \emph{complexity} of $X$.
\end{defn}

\begin{rem}[Complexity versus dimension]\label{rem:dimension}
Suppose that $X$ has complexity $n$.  Let $\Delta$ be a $t$--simplex of $X$.  Then $\Delta$ contains a chain 
$\emptyset\subsetneq\Delta_0\subsetneq\Delta_1\subsetneq\ldots\subsetneq\Delta_t$ with each $\Delta_i$ an 
$i$--simplex.  
Observe that $\link(\Delta_{i+1})\subsetneq\link(\Delta_i)$ for all $i$, so $t+2\leq n$, i.e., $\dimension X\leq n-2$.  

However, for a general simplicial complex $X$, the complexity cannot be bounded in terms of the dimension.  
Indeed, let $X$ be the following $1$--dimensional 
simplicial complex.  Let $X^{(0)}=\{v_n\}_{n\ge0}\sqcup\{h_n\}_{n\ge0}$, and, for each $i\ge0$, join $v_i$ by an 
edge to $h_0,\ldots,h_i$.  Then $\link(v_i)\subsetneq \link(v_{i+1})$ for all $i$, i.e., the complexity is infinite.
\end{rem}

Our basic object is a \emph{combinatorial hierarchically hyperbolic space}:

\begin{defn}[Combinatorial HHS]\label{defn:combinatorial_HHS}
 A \emph{combinatorial HHS}, abbreviated \emph{CHHS}, $(X,W)$ consists of a flag simplicial 
 complex $X$ and an $X$--graph $W$ satisfying all of the following conditions.
 \begin{enumerate}
  \item \label{item:chhs_flag}$X$ has complexity $n<+\infty$.
  
  \item \label{item:chhs_delta}There is a constant $\delta$ so that for each non-maximal simplex $\Delta$, the subgraph 
$\mathcal C(\Delta)$ is 
$\delta$--hyperbolic and $(\delta,\delta)$--quasi-isometrically embedded in the complement subgraph 
$Y_\Delta$, which was defined in Definition~\ref{defn:complement}.
\item \label{item:chhs_join}Let $\Delta$ and $\Sigma$ be non-maximal simplices such that there exists 
a non-maximal simplex $\Gamma$ with the following properties:
\begin{itemize}
    \item $\link(\Gamma)\subseteq\link(\Delta)$,
    \item $\link(\Gamma)\subseteq \link(\Sigma)$, and
    \item $diam(\mathcal C 
(\Gamma))\geq \delta$.
\end{itemize}

Then there exists a 
simplex $\Pi$ in the link of $\Sigma$ such that $\link(\Sigma\star\Pi)\subseteq \link(\Delta)$ and all non-maximal simplices $\Gamma$ satisfying the above three itemized properties also satisfy $\link(\Gamma)\subseteq\link(\Sigma\star\Pi)$.
\item \label{item:C_0=C} If $v,w$ are distinct non-adjacent vertices of $\link(\Delta)$, for some simplex $\Delta$ of $X$, and $v,w$ are contained in 
$W$-adjacent maximal simplices, then they are contained in $W$-adjacent simplices of the form $\Delta\star\Sigma$.
\end{enumerate}
Sometimes we use the notation $(X,W,\delta,n)$ when we have to keep track of the constants.
\end{defn}

\begin{rem}\label{rem:join_nonmax}
The simplex $\Sigma\star\Pi$ in Definition~\ref{defn:combinatorial_HHS}.\eqref{item:chhs_join} is necessarily 
non-maximal.  Indeed, its link is nonempty since it contains $\link(\Gamma)$ for some non-maximal $\Gamma$.
\end{rem}

\subsection{On the various parts of the definition}

We regard the first two conditions of Definition \ref{defn:combinatorial_HHS} as the most important ones, 
and the ones with 
solid theoretical reasons to be there. 

As a side note, the quasi-isometric embedding part of Condition \ref{item:chhs_delta} 
can be viewed as an analogue  of Bowditch's fineness condition for relative hyperbolicity \cite[Proposition 2.1.(F5)]{Bow:rel_hyp}.  We discuss this in the context of examples in Section~\ref{subsec:user}.

While the first two conditions do not seem to be sufficient to yield an HHS, we expect that 
the last two conditions can be replaced with ``better'' ones. We do not have compelling reasons for those 
properties to be required; our best explanations are as follows. Our heuristic justification for Condition 
\ref{item:chhs_join} is that it seems to be what is needed to perform arguments that in the context of the 
curve graph would require the use of tight geodesics, see in particular the proof of the ``Uniqueness'' axiom 
(i.e., the verification that Definition~\ref{defn:HHS}.\eqref{item:dfs_uniqueness} is satisfied). The heuristic 
justification for Condition \ref{item:C_0=C} is that, without it, it might be possible to ``move'' between 
places in the link of some simplex without this being doable within the link itself. 

In the interest of the reader who might need alternative conditions, or who might be interested in finding the ``right'' ones, we list where the last two conditions get used:

\begin{rem}[Potentially replacing conditions~\eqref{item:chhs_join} and \eqref{item:C_0=C}]\label{rem:join_heuristic}
Here are the only places where Condition ~\eqref{item:chhs_join} is used:
\begin{itemize}
\item the proof of Proposition~\ref{prop:Y_delta_hyperbolic}, via Lemma~\ref{lem:Y_delta_containment},
 \item the proof of Lemma \ref{lem:curve_graph_iota_v2}.\ref{item:containment_Y},
 \item the proof of Proposition \ref{prop:compatibility_1} (only to check the analogous condition for a link),
 \item the proof of Theorem \ref{thm:hhs_links}, in Section~\ref{sec:check_axioms}, in the ``Consistency for 
nesting'', ``Large links'', and ``Uniqueness'' (Case 1) parts.
\end{itemize}
The only uses of Condition \eqref{item:C_0=C} are in Section~\ref{sec:y_hyp}, in the proofs of Lemma~\ref{lem:Y_delta_hyperbolic} and Lemma~\ref{lem:electric-intersection}, and in Section \ref{sec:hieromorphisms}, via Lemma 
\ref{lem:C_0=C}.
\end{rem}

Here are some further comments on Definition~\ref{defn:combinatorial_HHS}.\eqref{item:chhs_join}. To understand the point of Definition~\ref{defn:combinatorial_HHS}.\eqref{item:chhs_join}, we consider a more intuitive, and strictly stronger, version of the condition: one could insist that intersections of links are always links, i.e., if $\link(\Delta)\cap\link(\Sigma)$ is nonempty, then there exists a (necessarily non-maximal) simplex $\Pi$ such that $\link(\Sigma)\cap\link(\Delta)=\link(\Pi)$.  

This captures the right intuition. Indeed, our goal is to show that the set of equivalence classes of non-maximal simplices will give a hierarchically hyperbolic structure, where nesting is containment of links.  For each $\Delta$, the associated hyperbolic space will be $\mathcal C(\Delta)$, and if $\link(\Pi)\subsetneq\link(\Delta)$ for some non-maximal simplex $\Pi$, we need $\Pi$ to correspond to a bounded subset of $\mathcal C(\Delta)$, because this is required by the definition of a hierarchically hyperbolic space (Definition~\ref{defn:HHS}.\eqref{item:dfs_nesting}).  This much always works: $\sat(\Pi)$ has to intersect $Y_\Delta$, so $\link(\Pi)$ is bounded in $Y_\Delta$ and hence in $\mathcal C(\Delta)$ in view of Definition~\ref{defn:combinatorial_HHS}.\eqref{item:chhs_delta}. 

Now, if $\Sigma$ and $\Delta$ have intersecting links, but there is no containment between their links, then they should be transverse elements in the HHS structure (see Definition~\ref{defn:HHS}.\eqref{item:dfs_transversal}).  This again requires that $\Sigma$ correspond to a bounded set in $\mathcal C(\Delta)$.  The naive hypothesis $\link(\Sigma)\cap\link(\Delta)=\link(\Pi)$ would achieve this, since $\Pi$ would be nested in $\Delta$.

This hypothesis is too strong to accommodate desirable examples, including the case where $X$ is the curve graph of a surface (so that simplices are multicurves).  See Remark~\ref{rem:condition-B}.

So, we define our bounded sets in $\mathcal C(\Delta)$ as follows.  When $\Sigma$ is nested in, or transverse to, $\Delta$, we observe as above that $\sat(\Sigma)\cap Y_\Delta$ is always nonempty and bounded in $Y_\Delta$, i.e., $\mathcal C(\Sigma)\cap Y_\Delta$ is ``coned off" in $Y_\Delta$.  

To define the required bounded subset of $\mathcal C(\Delta)$ associated to $\Sigma$, which is done in Definition~\ref{defn:projections}, we apply the coarse closest-point projection $Y_\Delta\to \mathcal C(\Delta)$ to $\sat(\Sigma)\cap Y_\Delta$.  This projection obviously gives a uniformly bounded set if $\mathcal C(\Delta)$ has uniformly bounded diameter.  If not, we verify in Section~\ref{sec:y_hyp} that $Y_\Delta$ is hyperbolic, from which Definition~\ref{defn:chhs_action}.\eqref{item:chhs_delta} implies that the coarse projection has the necessary properties. 

The need for a condition like Definition~\ref{defn:combinatorial_HHS}.\eqref{item:chhs_join} then makes itself felt in proving hyperbolicity of $Y_\Delta$.  Recall that we only need to do this when $\mathcal C(\Delta)$ has diameter at least the threshold $\delta$.  The point in the proof of hyperbolicity (Lemma~\ref{lem:Y_delta_hyperbolic}) where this is needed is: via Lemma~\ref{lem:SC_nesting}, and its consequence, Lemma~\ref{lem:Y_delta_containment}, if $\Delta$ is nested in some $\Pi$, then $Y_\Pi\subset Y_\Delta$.  The hypothesis we have chosen seems to be a suitable way to arrange this while being weak enough to cover natural examples.  As noted above, we also use it in similar ways in a few other places.

\subsection{Projections to links}\label{subsec:projections_to_links}

In this section we relate combinatorial objects to HHS objects, the connection being justified by Theorem \ref{thm:hhs_links}. The reader who is not interested in the details of the hierarchical structure obtained but only wishes to use it as a simple criterion to prove the hierarchical hyperbolicity of a space/group can skip this section and go directly to Section~\ref{subsec:main_theorem_statement}.

\begin{defn}[Nesting, orthogonality, transversality, complexity]\label{defn:nest_orth}
Let $X$ be a simplicial complex.  Let $\Delta,\Sigma$ be non-maximal simplices of $X$.  Then:
\begin{itemize}
     \item $[\Delta]\nest[\Sigma]$ if $\link(\Delta)\subseteq\link(\Sigma)$, and we say $[\Delta]$ is \emph{nested} in $[\Sigma]$;
     \item $[\Delta]\orth[\Sigma]$ if $\link(\Sigma)\subseteq \link(\link(\Delta))$, and we say $[\Delta]$ and $[\Sigma]$ are \emph{orthogonal}.
\end{itemize}
If $[\Delta]$ and $[\Sigma]$ are neither $\orth$--related nor $\nest$--related, we write 
$[\Delta]\transverse[\Sigma]$, and say $[\Delta]$ and $[\Sigma]$ are \emph{transverse}.

Note that $[\emptyset]$ is the unique $\nest$--maximal $\sim$--class of simplices in $X$ and that $\nest$ is a partial 
ordering on the set of $\sim$--classes of simplices in $X$.  Notice that the simplicial 
complex $X$ has finite complexity if there exists $n\in\naturals$ so that any $\nest$--chain has length at 
most $n$; the minimal such $n$ is the complexity of $X$.  
\end{defn}

\begin{rem}
 The definition of $\Sigma\orth\Delta$ is equivalent to saying that any vertex in the link of $\Sigma$ is joined by an edge to any vertex in the link of $\Delta$.
\end{rem}

One might be tempted to think of nesting as being equivalent to inclusion of simplices, but this only works in one direction, namely:

\begin{rem}\label{rem:set_theory}
Let $\Delta,\Delta'$ be simplices of $X$. If $\Delta\subseteq \Delta'$, then $[\Delta']\nest [\Delta]$.  

Similarly, if $\sat(\Delta)\subset\sat(\Delta')$, then any vertex in $\link(\Delta')$ is adjacent to every vertex in $\sat(\Delta)$, so $\link(\Delta')\subset\link(\Delta)$, i.e., $[\Delta']\nest[\Delta]$.  Again, the converse does not hold, although Lemma~\ref{lem:Y_delta_containment} gives a partial converse.
\end{rem}

Notice that Definition~\ref{defn:combinatorial_HHS}.\eqref{item:chhs_join} can be rephrased as follows:

\begin{itemize}
 \item Whenever $\Delta$ and $\Sigma$ are non-maximal simplices for which there exists 
a non-maximal simplex $\Gamma$ such that $[\Gamma]\nest[\Delta]$, $[\Gamma]\nest[\Sigma]$, and $diam(\mathcal C 
(\Gamma))\geq \delta$, then there exists a 
simplex $\Pi$ in the link of $\Sigma$ such that $[\Sigma\star\Pi]\nest [\Delta]$ and all $[\Gamma]$ as above satisfy $[\Gamma]\nest[\Sigma\star\Pi]$.
\end{itemize}

Also, note that if $\Pi$ is the simplex associated to $\Sigma$ and $\Delta$, provided by 
Definition~\ref{defn:combinatorial_HHS}.\eqref{item:chhs_join}, then since $\Pi\subset\link(\Sigma)$, the 
simplex $\Pi\star\Sigma$ exists automatically and, moreover, $[\Pi\star\Sigma]\nest[\Sigma]$. 

We note the following special case of Condition~\eqref{item:chhs_join} for later use:

\begin{lemma}
\label{lem:SC_nesting}
Suppose that 
$[\Sigma]\nest[\Delta]$ and $diam(\mathcal C(\Sigma))\geq\delta$.  Then $[\Sigma]=[\Delta\star\Pi]$ for some 
simplex $\Pi$ of $\link(\Delta)$.
\end{lemma}

\begin{proof}
Definition~\ref{defn:combinatorial_HHS}.\eqref{item:chhs_join} provides a simplex $\Pi$ of $\link(\Delta)$ such that 
$[\Delta\star\Pi]\nest[\Sigma]$.  Since $diam(\mathcal C(\Sigma))\geq\delta$, we also have (setting 
$\Gamma=\Sigma$)
$[\Sigma]\nest[\Delta\star\Pi]$, so $[\Sigma]=[\Delta\star\Pi]$.
\end{proof}

Our next goal is to define projections from $W$ to $\mathcal C([\Delta])$ for $[\Delta]\in\mathfrak S$.  This will 
use the following lemma:

\begin{lemma}\label{lem:proj_defined}
Let $X$ be a flag simplicial complex, let $\Delta$ be a non-maximal simplex, and let $\Sigma$ be a maximal simplex.  Then 
$\Sigma\cap Y_\Delta$ is nonempty and has diameter at most $1$.
\end{lemma}

\begin{proof}
Let $Z$ be the subcomplex of $X$ spanned by $\sat(\Delta)$.  Then for each maximal simplex $\Pi$ of $Z$, we have 
$\link(\Pi)\supseteq\link(\Delta)$.  If $\Sigma\cap Y_\Delta=\emptyset$, then $\Sigma^{(0)}\subseteq\sat(\Delta)$, 
so $\Sigma\subset Z$.  Moreover, by maximality, we have $\link(\Sigma)\supseteq\link(\Delta)$.  But 
$\link(\Sigma)=\emptyset$, by maximality of $\Sigma$, while $\link(\Delta)\ne\emptyset$, by non-maximality of $\Delta$.  
Hence $\Sigma\cap Y_\Delta\ne\emptyset$.

Since the vertices of $\Sigma$ are pairwise-adjacent in $X$, they are pairwise-adjacent in $\duaug{X}{W}$, so since 
$Y_\Delta$ is an induced subgraph, the vertices of $\Sigma\cap Y_\Delta$ are pairwise-adjacent in $Y_\Delta$, as required.
\end{proof}

\begin{defn}[Projections]\label{defn:projections}
Let $(X,W,\delta,n)$ be a combinatorial HHS.  

Fix $[\Delta]\in\mathfrak S$ and define a map $\pi_{[\Delta]}:W\to 2^{\mathcal C([\Delta])}$ as 
follows.  Let $p:Y_\Delta\to2^{\mathcal C([\Delta])}$ be the coarse closest point projection, i.e.,
$$p(x)=\{y\in\mathcal C([\Delta]):\dist_{Y_\Delta}(x,y)\le\dist_{Y_\Delta}(x,\mathcal C([\Delta]))+1\}.$$

Suppose that $w$ is a vertex of $W$, so $w$ corresponds to a unique simplex $\Sigma_w$ of $X$.  Since 
$\Sigma_w$ 
is maximal (by Definition~\ref{defn:X_graph}), and $\Delta$ is non-maximal (by the definition of 
$\mathfrak S$), the graph $\Sigma_w\cap Y_\Delta$ is nonempty and has diameter at most $1$, by Lemma~\ref{lem:proj_defined}. 
 Define $$\pi_{[\Delta]}(w)=p(\Sigma_w\cap Y_\Delta).$$

We have thus defined  $\pi_{[\Delta]}:W^{(0)}\to 2^{\mathcal C([\Delta])}$.   If $v,w\in W$ are joined by an 
edge $e$ of $W$, 
then $\Sigma_v,\Sigma_w$ are joined by edges in $\duaug{X}{W}$, and we let 
$\pi_{[\Delta]}(e)=\pi_{[\Delta]}(v)\cup\pi_{[\Delta]}(w)$.

Now let $[\Delta],[\Delta']\in\mathfrak S$ satisfy $[\Delta]\transverse[\Delta']$ or $[\Delta']\propnest [\Delta]$. 
Let $$\rho^{[\Delta']}_{[\Delta]}=p( 
\sat(\Delta')\cap Y_\Delta),$$ where $p:Y_\Delta\to 2^{\mathcal C([\Delta])}$ is 
coarse closest-point projection. 

Let $[\Delta]\propnest 
[\Delta']$. Let $\rho^{[\Delta']}_{[\Delta]}:\mathcal C([\Delta'])\to 
2^{\mathcal C([\Delta])}$ be defined as follows.  On $\mathcal C([\Delta'])\cap 
Y_\Delta$, it is the restriction of $p$ to $\mathcal C([\Delta'])\cap 
Y_\Delta$. Otherwise, it takes the value $\emptyset$.
\end{defn}

\begin{rem}[See the future]\label{rem:see_the_future}
In Lemma~\ref{lem:Y_delta_hyperbolic}, we will show that $Y_\Delta$ is $\delta_0$--hyperbolic, for some uniform $\delta_0$, provided $\diam(\mathcal C(\Delta))\geq\delta$. Since $\mathcal C([\Delta])$ is 
$(\delta,\delta)$--quasi-isometrically embedded, we will then have that $\mathrm{diam}(p(\Sigma_w\cap Y_\Delta))$ is bounded 
in terms of $\delta,\delta_0$, i.e., $\pi_{[\Delta]}(w)$ is a nonempty, uniformly bounded set.  When $\diam(\mathcal C(\Delta))\leq\delta$, then the same conclusion is immediate, with no need for hyperbolicity of $Y_\Delta$.  Once we have established that either $Y_\Delta$ is hyperbolic or $\mathcal C(\Delta)$ is uniformly bounded, then the coarse closest-point projection will send points to uniformly bounded sets, and we will, when convenient and when we are only concerned about distances up to uniformly bounded error, think of $p$ as a map. 
\end{rem}

\subsection{Statement of Theorem~\ref{thm:hhs_links}}\label{subsec:main_theorem_statement}
Our main theorem about combinatorial HHS is Theorem~\ref{thm:hhs_links}.  See Section~\ref{sec:background} for the 
definition of a hierarchically hyperbolic space (HHS) and a hierarchically hyperbolic group (HHG).

\begin{thm}[HHS structures from $X$--graphs]\label{thm:hhs_links}
Let $(X,W)$ be a combinatorial HHS.

Let $\mathfrak S$ be as in Definition~\ref{defn:simplex_equivalence}, define nesting and orthogonality 
relations on $\mathfrak S$ as in Definition~\ref{defn:nest_orth}, let the associated hyperbolic spaces be as in 
Definition~\ref{defn:combinatorial_HHS}, and define projections as in Definition~\ref{defn:projections}. 

Then $(W,\mathfrak S)$ is a hierarchically hyperbolic space, and the HHS constants only depend on $\delta,n$ as in Definition 
\ref{defn:combinatorial_HHS}.  

Moreover, suppose that $G$ is a group acting cocompactly on $X$.  Suppose that the $G$--action on the set of 
maximal simplices of $X$ extends to an action on $W$ which is 
metrically proper and cobounded. Then $(G,\mathfrak S)$ is a hierarchically hyperbolic group.
\end{thm}

\begin{rem}
In the ``moreover'' part of the statement, one actually only needs something weaker than cocompactness of the $G$--action on 
$X$.  Specifically, the exact property we need is that $G$ acts on $X$ with finitely many 
orbits of subcomplexes of the form $\link(\Delta)$, where $\Delta$ is a non-maximal simplex of $X$.
\end{rem}

As in Definition~\ref{defn:chhs_action} below, we say that the group $G$ \emph{acts} on the combinatorial HHS 
$(X,W)$ if $G$ acts by simplicial automorphisms on $X$ and the action on $W^{(0)}$ induced by the action on 
$X$ extends to an action on the whole graph $W$ (i.e., it preserves $W$--adjacency).  In this language, we can 
rephrase the latter part of Theorem~\ref{thm:hhs_links}:

\begin{cor}\label{cor:exportable}
 Let $G$ act on the combinatorial HHS $(X,W)$.  Suppose that the action of $G$ on $X$ is cocompact and the 
action on $W$ is proper and cocompact.  Then $G$ is a hierarchically hyperbolic group.
\end{cor}

\begin{rem}
 Notice that, under the assumptions of Corollary \ref{cor:exportable}, we have that the action of $G$ on $X$ is acylindrical in view of \cite[Theorem K]{HHS_I}.
\end{rem}

We fix the notation of Theorem \ref{thm:hhs_links} from now on.  The proof of Theorem~\ref{thm:hhs_links} is 
in Section~\ref{subsec:verify_axioms_assuming_hyperbolic}, after some necessary preparation.

\subsection{User's guide and simple examples}\label{subsec:user}

We make some remarks that could be useful for the reader interested in applying Theorem \ref{thm:hhs_links} to establish hierarchical hyperbolicity in their example of interest. First, Theorem \ref{thm:propaganda} provides a simpler set of conditions that do not involve the $X$-graph $W$, and the reader is advised to first check whether that theorem applies in their situation. Typical obstructions to using the simplified version arise from bounded links, which are treated more flexibly in Theorem \ref{thm:hhs_links}.  

If not, it has to be noted that, in situations where there is a natural $X$ to consider, the $X$ might actually have to be changed within its quasi-isometry class to satisfy the fine geometry constraints.  For example, in the amalgamated free product example just below, we see that the natural candidate, the Bass-Serre tree, may not work as our $X$, and we need to ``blow up'' the vertices of the tree to stars before proceeding.

One strategy is to build the correct ``model'' bottom-up, meaning starting from the hyperbolic complexes for the expected sub-HHS (e.g., in a tree of HHS, one might want to suitably combine the hyperbolic complexes for the various vertex spaces, see below).

It might also be useful to note that taking direct products at the level of complexes corresponds to taking joins at the level of links.  Also, relative hyperbolicity corresponds to disjoint unions of cones over the hyperbolic complexes for the peripherals, as we discuss a bit more at the end of this subsection.

\subsubsection{Amalgamated free product example}\label{subsub:AFP}  We now give an example of a combinatorial HHS. Let us consider an amalgamated product $G=A*_C B$ of hyperbolic 
groups over a common quasiconvex almost-malnormal subgroup $C$, so that $G$ is hyperbolic, and hence hierarchically hyperbolic, by the Bestvina--Feighn combination theorem \cite{BF:combination}. 

We will define the simplicial complex $X$, 
which will be quasi-isometric to the Bass-Serre tree for $G$. Notice that the Bass-Serre tree itself is not the 
right complex to consider since there is no link ``encoding $C$'', meaning that $C$ does not act on any link in such a way that, say, $C$ has unbounded orbits if it is infinite, as one would expect from the right complex in view of the distance formula. 

Let us now construct 
$X$, as follows. The vertex set of $X$ is $G\sqcup \{v_{gA}:gA\in G/A\}\sqcup \{v_{gB}:gB\in G/B\}$ (where 
$G/H$ denotes the set of left cosets of $H$ in $G$). Edges of $X$ correspond to either containment of an 
element of $G$ in a coset of $A$ or $B$, or to pairs of cosets of $A$ and $B$ intersecting non-trivially; edges of the latter type correspond to edges of the Bass-Serre tree.

\begin{figure}
	\begin{center}
		\scalebox{1}{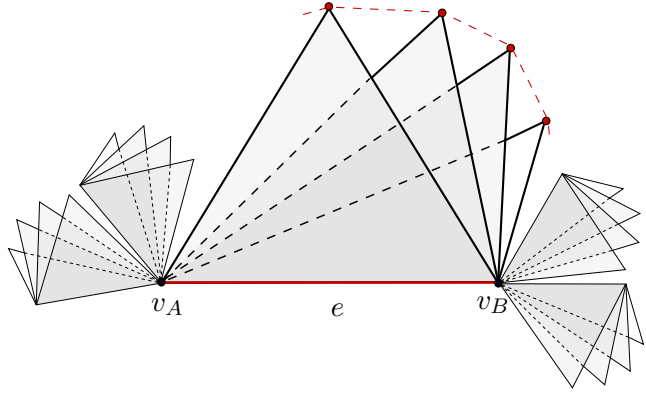}
		\caption{A portion of the complex $X$. The vertices in red correspond to elements of $C$, and all have the same link (red edge $e$  in the picture). The link of that edge in $X$ is a discrete set in bijection with $C$. However, due to the choice of $W$, the augmented link $\mathcal C(e)$ is the Cayley graph of $C$ with respect to $S_C$ (red dotted lines in the picture).}
		\label{fig_example}
	\end{center}
\end{figure}

Finally, we let $X$ be the flag complex with the $1$-skeleton we just described. Roughly speaking, $X$ is the 
Bass-Serre tree corresponding to the splitting, but where every edge is now contained in triangles indexed by 
$C$ (see Figure \ref{fig_example}). Notice that maximal simplices have vertex set of the form 
$g,v_{gA},v_{gB}$, and are in bijection with $G$. The link of $g\in G$ is a single edge, while the vertex set 
of the link of $v_A$ is in bijection with $A\sqcup A/C$.

Examples of saturations are that the saturation  of 
$g\in G$ is $gC$, while the saturation of the edge with endpoints $1,A$ is $C\cup\{v_A\}$.

We now define $W$ as any Cayley graph of $G$ corresponding to a generating set $S_A\cup S_B\cup S_C$ with $S_H$ a generating set of $H$ for $H\in \{A,B,C\}$, and $S_A\cap C=S_B\cap C=\emptyset$. Then it can be checked that $X^{+W}$ is quasi-isometric to $X$ and to the Bass-Serre tree of $G$, that $\mathcal C (v_A)$ is the Cayley graph of $A$ with respect to $S_A$ with the cosets of $C$ coned-off (and similarly for $B$), and that $\mathcal C(e)$, where $e$ is the edge with endpoints $v_A,v_B$, is the Cayley graph of $C$ with respect to $S_C$.  The link of the edge $e'$ joining, say, $1\in G$ to $v_A$ is the single vertex $v_B$.

\subsubsection{Relative hyperbolicity}
We now discuss relative hyperbolicity, and the analogy between the notion of a combinatorial HHS and Bowditch's \emph{fine graphs} \cite[Proposition 2.1]{Bow:rel_hyp}.

First, consider infinite hyperbolic groups $A,B$ and let $G=A* B$.  Let $X$ be the Bass-Serre tree.  The vertices are of the form $v_{gA}$ or $v_{gB}$, for $g\in G$, i.e., they are indexed by left cosets of $A$ and $B$.  Links of vertices are discrete, and any two intersect in at most one vertex, so there is no proper containment of links beyond the fact that all links are contained in the link of $\emptyset$.

There is a natural $G$--equivariant bijection $G\to W^{(0)}$, where $W^{(0)}$ is the set of edges (maximal simplices) of $X$: each $g\in G$ appears in precisely one coset $gA$ and one coset $gB$, and the intersection of these cosets corresponds to an edge of $X$.  Fixing a finite generating set for $G$ consisting of the disjoint union of generating sets of $A$ and $B$, we join $x,y\in W^{(0)}$ (edges of $X$) by an edge of $W$ if and only if the corresponding elements of $G$ are adjacent in the Cayley graph of $G$.  So, $W$ is a copy of the Cayley graph of $G$.  The graph $X^{+W}$ is naturally quasi-isometric to the coned-off Cayley graph of $G$, with cones over each $gA,gB$, so it is hyperbolic.  The link of each vertex $v$ becomes, upon addition of $W$--edges, a copy of a Cayley graph of $A$ or $B$, and hence hyperbolic.  One verifies that this augmented link is quasi-isometrically embedded in $X^{+W}-\{v\}$ by constructing a Lipschitz retract from $X^{+W}-\{v\}$ to the link (see, for instance, \cite{HHS:graph-man} for a more general version of such a construction in the context of Bass-Serre trees).

This is a simple example of a relatively hyperbolic group --- thinking of $A$ and $B$ as the peripheral subgroups --- where an associated hyperbolic fine graph, namely $X$, also functions as the underlying complex of a CHHS structure.  As we saw above, obtaining CHHS structures from splittings normally requires ``blowing up'' the Bass-Serre tree; that it doesn't in this example is an artifact of the example being a free product.

However, coned-off Cayley graphs for relatively hyperbolic groups do relate closely to combinatorial HHS.  Indeed, let $G$ be hyperbolic relative to a subgroup $P$.  Fix a Cayley graph $K_0$ of $G$ associated to a finite generating set containing generators for $P$, and let $K$ be obtained from $K_0$ by adding a vertex $v_{gP}$, joined by an edge to each element of $gP$, for each left coset $gP$.  So, $K$ is the standard hyperbolic fine graph witnessing relative hyperbolicity.  

On the other hand, suppose that $(X_P,W_P)$ is a combinatorial HHG structure for $P$, i.e., a combinatorial HHS, with a simplicial $P$--action on $X_P$ such that, for simplicity, $W_P$ with the induced $P$--action is equivariantly isomorphic to the Cayley graph of $P$ with the finite generating set mentioned above.  For each $gP$, let $X'_{gP}$ be obtained from a copy $X_{gP}$ of $X_P$ by adding a vertex $v_{gP}$, joining it to each vertex of $X_{gP}$, and taking the induced flag complex.  We let $W_{gP}$ be a copy of $W_P$; note that the vertices of $W_{gP}$ correspond naturally to the elements of $gP$.  Let $X$ be the disjoint union of the $X_{gP}'$.  Let $W$ be formed from the disjoint union of the $W_{gP}$ as follows: first, note that the vertex set of $W$ is naturally in bijection with $G$, since the cosets $gP$ partition $G$.  Join vertices of $W$ by an edge if the corresponding elements of $G$ are adjacent in our Cayley graph.  The pair $(X,W)$ is a combinatorial HHS with $G$ acting geometrically on $W$.  For example, $X^{+W}$ is quasi-isometric to $K$, so it is hyperbolic.  This is a sketch of a proof that a group hyperbolic relative to combinatorially hierarchically hyperbolic groups is again combinatorially hierarchically hyperbolic, and justifies the analogy with fine hyperbolic graphs.

\subsection{Mapping class groups and blow-ups}\label{subsec:blow_up}
This subsection contains further discussion on applications of
Theorem~\ref{thm:hhs_links} and is not used elsewhere in the paper.

Although $MCG(S)$ is known to be a hierarchically hyperbolic group ---
the index set $\mathfrak S$ consists of isotopy classes of (possibly
disconnected) subsurfaces, and the associated hyperbolic spaces are
curve graphs~\cite[Section 11]{HHS_II} --- one cannot apply
Theorem~\ref{thm:hhs_links} to the curve graph $\fontact S$ in order
to realize $MCG(S)$ as a combinatorial HHS. The reason is that annular
curve graphs --- projections to which need to appear one way or
another in any HHS structure --- do not arise as links of simplices in
the curve complex of $S$.  When we take $X=\fontact S$ and $W$
to be the pants graph of $S$, then applying Theorem~\ref{thm:hhs_links}
yields the HHS structure on the pants graph (i.e., on
Teichm\"uller space with the Weil-Petersson metric) as a combinatorial
HHS. Note that this standard HHS structure on the Weil-Petersson metric, described in \cite[Theorem G]{HHS_I} or \cite[Theorem 1.1, Example 2.3]{Vokes}, 
has index set consisting of the \emph{non-annular}
subsurfaces, which do correspond to links in $\fontact S$.

However, one can most likely modify $\fontact S$ by something like the following ``blow-up'' construction to 
get a combinatorial HHS 
structure on the 
marking graph (which is quasi-isometric to the mapping class group).  For each vertex $\gamma$ of $\fontact S$, let $B(\gamma)$ be 
the graph 
obtained from the vertex set of the annular curve graph $\fontact\gamma$ by adding a vertex $\gamma$ joined by 
an edge to 
each vertex of $\fontact\gamma$.  Let $X$ be the flag complex on the graph obtained from 
$\bigsqcup_{\gamma\in\fontact 
S^{(0)}}B(\gamma)$ by adding an edge joining every vertex of $B(\gamma)$ to every vertex of $B(\alpha)$ 
whenever 
$\gamma,\alpha$ are $\fontact S$--adjacent.  (So, in particular, $X$ is quasi-isometric to $\fontact S$: just 
collapse each 
$B(\gamma)$ to the vertex $\gamma$ and each subgraph $B(\gamma)\star B(\alpha)$ to the edge $[\alpha,\gamma]$.) 
 The idea would then be to associate maximal simplices in $X$ to markings, take $W$ to be a suitably chosen 
version of the marking graph, and verify that $(X,W)$ is a combinatorial HHS.  

It's instructive to see how annular curve graphs arise as links in this setup.  Let $\gamma$ be a curve, and 
let 
$\gamma_1,\ldots,\gamma_n$ be a maximal collection of disjoint curves with $\gamma_1=\gamma$.  For $i\geq 2$, choose a point 
$x_i\in\fontact\gamma_i$ (regarded as a subgraph of $X$).  Let $\Sigma$ be the simplex with vertex set 
$\gamma_1,\ldots,\gamma_n,x_2,\ldots,x_n$.  Then $\link_X(\Sigma)$ is exactly $\fontact\gamma$.

Also, let $\Sigma$ be a simplex of the form $[\gamma_1,x_1],\ldots,[\gamma_n,x_n]$, where each $[\gamma_i,x_i]$ 
joins a 
curve $\gamma_i$ to $x_i\in\fontact\gamma_i$, and the $\gamma_i$ are all disjoint.  Then $\link_X(\Sigma)$ is 
the union of 
the subgraphs $B(\alpha)$, as $\alpha$ varies over all the curves in the complement of 
$\gamma_1,\ldots,\gamma_n$.  In 
particular, if $Y$ is a subsurface, then taking $\gamma_1,\ldots,\gamma_n$ to consist of the boundary curves of 
$Y$ that 
are essential in $S$, together with pants decompositions of the components of $S-Y$, we find that 
$\link_X(\Sigma)$ 
corresponds to $\fontact Y$.  So, except for some bounded links, the links in $X$ correspond to the curve 
graphs of the 
various subsurfaces in the HHS structure on the marking graph from \cite[Section 11]{HHS_II}.

We expect that there is a wide range of contexts in which similar
blow-up constructions can be used to construct combinatorial HHS
structures.  For example, we expect that for (many) right-angled Artin
groups, one can ``blow up'' the Kim-Koberda \emph{extension
graph}~\cite{KK1,KK2} to obtain a combinatorial HHS structure, and 
whence an alternate proof of hierarchical hyperbolicity for these groups via
Theorem~\ref{thm:hhs_links}.

In \cite{HHS:artin}, hierarchical hyperbolicity of extra-large-type Artin groups is established by 
constructing an 
appropriate 
version of the extension graph, blowing it up, and verifying that this blown-up graph gives the $X$ such that 
$(X,W)$ is a 
combinatorial HHS, where $W$ is an appropriately-chosen Cayley graph of the Artin group.  In fact, we expect that such a blow-up construction can be used to prove a partial converse to 
Theorem~\ref{thm:hhs_links} 
under some extra conditions on the hierarchically hyperbolic space.

\section{Background on hierarchical hyperbolicity}\label{sec:background}

\subsection{Axioms}
We recall from~\cite{HHS_II} the definition of a hierarchically hyperbolic space.

\begin{defn}[Hierarchically hyperbolic space]\label{defn:HHS}
The $q$--quasigeodesic space  $(\cuco X,\dist_{\cuco X})$ is a \emph{hierarchically hyperbolic space} if there 
exists $\delta\geq0$, an index set $\mathfrak S$, and a set $\{\mathcal C U:U\in\mathfrak S\}$ of 
$\delta$--hyperbolic spaces $(\mathcal C U,\dist_U)$,  such that the following conditions are 
satisfied:\begin{enumerate}
\item\textbf{(Projections.)}\label{item:dfs_curve_complexes} There is
a set $\{\pi_U: \cuco X\rightarrow2^{\mathcal C U}\mid U\in\mathfrak S\}$
of \emph{projections} sending points in $\cuco X$ to sets of diameter
bounded by some $\xi\geq0$ in the various $\mathcal C U\in\mathfrak S$.
Moreover, there exists $K$ so that for all $U\in\mathfrak S$, the coarse map $\pi_U$ is $(K,K)$--coarsely
Lipschitz and $\pi_U(\cuco X)$ is $K$--quasiconvex in $\mathcal C U$.

 \item \textbf{(Nesting.)} \label{item:dfs_nesting} $\mathfrak S$ is
 equipped with a partial order $\nest$, and either $\mathfrak
 S=\emptyset$ or $\mathfrak S$ contains a unique $\nest$--maximal
 element; when $V\nest U$, we say $V$ is \emph{nested} in $U$.  (We
 emphasize that $U\nest U$ for all $U\in\mathfrak S$.)  For each
 $U\in\mathfrak S$, we denote by $\mathfrak S_U$ the set of
 $V\in\mathfrak S$ such that $V\nest U$.  Moreover, for all $V,U\in\mathfrak S$
 with $V\propnest U$ there is a specified subset
 $\rho^V_U\subset\mathcal C U$ with $\diam_{\mathcal C U}(\rho^V_U)\leq\xi$.
 There is also a \emph{projection} $\rho^U_V: \mathcal C
 U\rightarrow 2^{\mathcal C V}$.  (The similarity in 
notation is justified by viewing $\rho^V_U$ as a coarsely constant map $\mathcal C
 V\rightarrow 2^{\mathcal C U}$.)
 
 \item \textbf{(Orthogonality.)} 
 \label{item:dfs_orthogonal} $\mathfrak S$ has a symmetric and
 anti-reflexive relation called \emph{orthogonality}: we write $V\orth
 U$ when $V,U$ are orthogonal.  Also, whenever $V\nest U$ and $U\orth
 A$, we require that $V\orth A$.  We require that for each
 $T\in\mathfrak S$ and each $U\in\mathfrak S_T$ for which
 $\{V\in\mathfrak S_T\mid V\orth U\}\neq\emptyset$, there exists $B\in
 \mathfrak S_T-\{T\}$, so that whenever $V\orth U$ and $V\nest T$, we
 have $V\nest B$.  Finally, if $V\orth U$, then $V,U$ are not
 $\nest$--comparable.
 
 \item \textbf{(Transversality and consistency.)}
 \label{item:dfs_transversal} If $V,U\in\mathfrak S$ are not
 orthogonal and neither is nested in the other, then we say $V,U$ are
 \emph{transverse}, denoted $V\transverse U$.  There exists
 $\kappa_0\geq 0$ such that if $V\transverse U$, then there are
  sets $\rho^V_U\subseteq\mathcal C W$ and
 $\rho^U_V\subseteq\mathcal C V$ each of diameter at most $\xi$ and 
 satisfying: $$\min\left\{\dist_{
 U}(\pi_U(x),\rho^V_U),\dist_{
 V}(\pi_V(x),\rho^W_V)\right\}\leq\kappa_0$$ for all $x\in\cuco X$.
 
 For $V,U\in\mathfrak S$ satisfying $V\nest U$ and for all
 $x\in\cuco X$, we have: $$\min\left\{\dist_{
 U}(\pi_U(x),\rho^V_U),\diam_{\mathcal C
 V}(\pi_V(x)\cup\rho^U_V(\pi_U(x)))\right\}\leq\kappa_0.$$ 
 
 The preceding two inequalities are the \emph{consistency inequalities} for points in $\cuco X$.
 
 Finally, if $U\nest V$, then $\dist_U(\rho^U_T,\rho^V_T)\leq\kappa_0$ whenever $T\in\mathfrak S$ satisfies either $V\propnest T$ or $V\transverse T$ and $T\not\perp U$.
 
 \item \textbf{(Finite complexity.)} \label{item:dfs_complexity} There exists $n\geq0$, the \emph{complexity} of $\cuco X$ (with respect to $\mathfrak S$), so that any set of pairwise--$\nest$--comparable elements has cardinality at most $n$.
  
 \item \textbf{(Large links.)} \label{item:dfs_large_link_lemma} There
exist $\lambda\geq1$ and $E\geq\max\{\xi,\kappa_0\}$ such that the following holds.
Let $U\in\mathfrak S$ and let $x,x'\in\cuco X$.  Let
$N=\lambda\dist_{U}(\pi_U(x),\pi_U(x'))+\lambda$.  Then there exists $\{T_i\}_{i=1,\dots,\lfloor
N\rfloor}\subseteq\mathfrak S_U-\{U\}$ such that for all $T\in\mathfrak
S_U-\{U\}$, either $T\in\mathfrak S_{T_i}$ for some $i$, or $\dist_{
T}(\pi_T(x),\pi_T(x'))<E$.  Also, $\dist_{
U}(\pi_U(x),\rho^{T_i}_U)\leq N$ for each $i$.

 \item \textbf{(Bounded geodesic image.)}
 \label{item:dfs:bounded_geodesic_image} There exists $E>0$ such that 
 for all $U\in\mathfrak S$,
 all $V\in\mathfrak S_U-\{U\}$, and all geodesics $\gamma$ of
 $\mathcal C U$, either $\diam_{\mathcal C V}(\rho^U_V(\gamma))\leq E$ or
 $\gamma\cap\neb_E(\rho^V_U)\neq\emptyset$.
 
 \item \textbf{(Partial Realization.)} \label{item:dfs_partial_realization} There exists a constant $\alpha$ with the following property. Let $\{V_j\}$ be a family of pairwise orthogonal elements of $\mathfrak S$, and let $p_j\in \pi_{V_j}(\cuco X)\subseteq \mathcal C V_j$. Then there exists $x\in \cuco X$ so that:
 \begin{itemize}
 \item $\dist_{V_j}(\pi_{V_j}(x),p_j)\leq \alpha$ for all $j$,
 \item for each $j$ and 
 each $V\in\mathfrak S$ with $V_j\nest V$, we have 
 $\dist_{V}(\pi_V(x),\rho^{V_j}_V)\leq\alpha$, and
 \item for each $j$ and 
 each $V\in\mathfrak S$ with $V_j\transverse V$, we have $\dist_V(\pi_V(x),\rho^{V_j}_V)\leq\alpha$.
 \end{itemize}

\item\textbf{(Uniqueness.)} For each $\kappa\geq 0$, there exists
$\theta_u=\theta_u(\kappa)$ such that if $x,y\in\cuco X$ and
$\dist_{\cuco X}(x,y)\geq\theta_u$, then there exists $V\in\mathfrak S$ such
that $\dist_V(\pi_V(x),\pi_V(y))\geq \kappa$.\label{item:dfs_uniqueness}
\end{enumerate}
 
We often refer to $\mathfrak S$, together with the nesting
and orthogonality relations, and the projections as a \emph{hierarchically hyperbolic structure} for the space $\cuco
X$.  Observe that $\cuco X$ is hierarchically hyperbolic with respect
to $\mathfrak S=\emptyset$, i.e., hierarchically hyperbolic of
complexity $0$, if and only if $\cuco X$ is bounded.  Similarly,
$\cuco X$ is hierarchically hyperbolic of complexity $1$ with respect
to $\mathfrak S=\{\cuco X\}$, if and only if $\cuco X$ is hyperbolic.
\end{defn}

\begin{rem}\label{rem:redundancy}
Jacob Russell has pointed out that the ``$\dist_{
W}(\pi_W(x),\rho^{T_i}_W)\leq N$'' requirement in Definition~\ref{defn:HHS}.\eqref{item:dfs_large_link_lemma} 
follows from the consistency and bounded geodesic image axioms, and is therefore redundant; see~\cite[Remark 
2.10]{Russell:rel_HHS}.
\end{rem}

\begin{notation}\label{notation:suppress_pi}
Where it will not cause confusion, given $U\in\mathfrak S$, we will often suppress the projection
map $\pi_U$ when writing distances in $\mathcal C U$, i.e., given $x,y\in\cuco X$ and
$p\in\mathcal C U$  we write
$\dist_U(x,y)$ for $\dist_U(\pi_U(x),\pi_U(y))$ and $\dist_U(x,p)$ for
$\dist_U(\pi_U(x),p)$. Note that when we measure distance between a 
pair of sets (typically both of bounded diameter) we are taking the minimum distance 
between the two sets. 
Given $A\subset \cuco X$ and $U\in\mathfrak S$ 
we let $\pi_{U}(A)$ denote $\cup_{a\in A}\pi_{U}(a)$.
\end{notation}

\begin{defn}[Hierarchically hyperbolic group]\label{defn:HHG}
The group $G$ is a \emph{hierarchically hyperbolic group (HHG)} if there exists a hierarchically hyperbolic space 
$(Z,\mathfrak S)$ such that the following hold:
\begin{itemize}
     \item $G$ acts metrically properly and coboundedly by isometries on the quasigeodesic space $Z$.
     \item $G$ acts on $\mathfrak S$ with finitely many orbits, and the $G$ action preserves the relations 
$\nest,\orth,\transverse$.
\item For all $U\in\mathfrak S$ and $g,h\in G$, there is an isometry $g:\fontact U\to\fontact gU$ such that 
the isometry $(gh):\fontact U\to\fontact ghU$ is the composition of the isometries $g:\fontact hU\to\fontact ghU$ and 
$h:\fontact U\to\fontact hU$.
\item For all $U\in\mathfrak S,g\in G,z\in Z$, we have $\pi_{gU}(gz)=g(\pi_U(z))$.
\item For all $U,V\in\mathfrak S$ such that $U\transverse V$ or $U\propnest V$, and all $g\in G$, we have 
$\rho^{gU}_{gV}=g(\rho^U_V)$.
\end{itemize}
From the first bullet and Milnor--Schwarz, $G$ is finitely generated and, when $G$ is equipped with any 
word-metric, composing the projections $\pi_U$ with any orbit map $G\to Z$ shows that we can take $Z=G$ in the above 
definition.  In particular, $(G,\mathfrak S)$ is an HHS.  When we 
wish to emphasize the particular HHS structure on which 
$G$ is acting, we say that $(G,\mathfrak S)$ is an HHG.
\end{defn}

\begin{rem}
In the definition, we have asked that $G$ act \emph{metrically properly and coboundedly}, rather than properly and 
cocompactly, since it is sometimes convenient to check that $G$ is an HHG by constructing an action on an HHS $(Z,\mathfrak 
S)$ where $Z$ is not proper.
\end{rem}

\subsection{Useful facts}

We now recall results from \cite{HHS_II} that will be useful later on.  To avoid some technicalities, we will 
assume that, given an HHS $(\cuco X,\mathfrak S)$, the maps $\pi_U,U\in\mathfrak S$ are uniformly coarsely 
surjective, which can always be arranged (see~\cite[Remark 1.3]{HHS_II}).

\begin{defn}[Consistent tuple]\label{defn:consistent_tuple}
Let $\kappa\geq0$ and let $\tup b\in\prod_{U\in\mathfrak S}2^{\mathcal C U}$ be a tuple such that for each $U\in\mathfrak S$, 
the $U$--coordinate  $b_U$ has diameter $\leq\kappa$.  Then $\tup b$ is \emph{$\kappa$--consistent} if for all 
$V,W\in\mathfrak S$, we have $$\min\{\dist_V(b_V,\rho^W_V),\dist_W(b_W,\rho^V_W)\}\leq\kappa$$ whenever $V\transverse W$ and 
$$\min\{\dist_W(b_W,\rho^V_W),\diam_V(b_V\cup\rho^W_V(b_W))\}\leq\kappa$$ whenever $V\propnest W$.
\end{defn}

The following is \cite[Theorem~3.1]{HHS_II}:

\begin{thm}[Realization]\label{thm:realization}
Let $(\cuco X,\mathfrak S)$ be a hierarchically hyperbolic space. Then for each $\kappa\geq1$, there exists 
$\theta=\theta(\kappa)$ so that,  for any $\kappa$--consistent tuple $\tup b\in\prod_{U\in\mathfrak S}2^{\mathcal C U}$, 
there exists $x\in\cuco X$ such that $\dist_V(x,b_V)\leq\theta$ for all $V\in\mathfrak S$.
\end{thm}

\noindent Observe that uniqueness (Definition~\eqref{item:dfs_uniqueness}) implies that the \emph{realization 
point} $x$ for 
$\tup b$ provided by Theorem~\ref{thm:realization} is coarsely unique.  

The following is 
\cite[Theorem~4.5]{HHS_II}:

\begin{thm}[Distance formula]\label{thm:distance_formula}
Let $(\cuco X,\mathfrak S)$ be a hierarchically hyperbolic space.  Then
there exists $s_0$ such that for all $s\geq s_0$, there exist $C,K$ so
that for all $x,y\in\cuco X$,
$$\dist(x,y)\asymp_{_{K,C}}\sum_{U\in\mathfrak
S}\ignore{\dist_U(x,y)}{s}.$$
\end{thm}

\noindent (The notation $\ignore{A}{B}$ denotes the quantity which is $A$ if $A\geq B$ and $0$ otherwise. The notation $A\asymp_{\lambda,\lambda} B$ means $A\leq \lambda B+\lambda$ and $B\leq \lambda A+\lambda$.)

We will  use the following variation, which is well-known to experts:

\begin{thm}[Distance Formula$+\epsilon$]\label{thm:DF_plus_epsilon}
 Let $(\cuco X,\mathfrak S)$ be an HHS. For every $\lambda\geq 1$, there exist $T\geq 2\lambda$ and $\kappa\geq 
1$, depending only on the HHS constants and $\lambda$,  with the following property. Let $x,y\in \cuco X$, 
and consider for every $Y\in\mathfrak S$ some $h_Y$ with $h_Y\asymp_{_{\lambda,\lambda}} \dist_Y(x,y)$. Then
 $$\dist_{X}(x,y)\asymp_{_{\kappa,\kappa}} \sum_{Y\in\mathfrak S} \ignore{h_Y}{T}.$$
\end{thm}

\begin{proof}
The proof follows from manipulating the thresholds of the usual distance formula.

 Let $DF_T(x,y)=\sum_{Y\in\mathfrak S} \ignore{\dist_Y(x,y)}{T}$, and 
 let $H_T(x,y)$ be the sum in the 
statement of the theorem we are proving. Then, for $T\ge 2\lambda$, we claim that we 
have:
 $$ \frac{1}{2\lambda}DF_{\lambda(\lambda+T)}(x,y)\leq H_T(x,y)\leq 2\lambda \ DF_{T/\lambda-1}(x,y).$$
 These inequalities, combined with the distance formula (Theorem~\ref{thm:distance_formula}), prove the 
required statement.
 
  Let us show the second inequality. If $h_Y\geq T$ (so that $h_Y$ contributes to $H_T(x,y)$), then we have 
$\dist_Y(x,y)\geq h_Y/\lambda -1\geq h_Y/(2\lambda)$ since $T\geq 2\lambda$.
 
Hence:
 
 $$H_T(x,y)=\sum_{Y: h_Y\geq T} h_Y\leq 2\lambda \sum_{Y: h_Y\geq T} \dist_Y(x,y)\leq 2\lambda \ DF_{T/\lambda 
-1}(x,y),$$
 as required.
 
 Let us show the first inequality, with the same method. If $\dist_Y(x,y)\geq \lambda(\lambda+T)$ then $h_Y\geq 
\dist_Y(x,y)/\lambda - 1 \geq T$, and also $h_Y\geq \dist_Y(x,y)/(2\lambda)$, since $T$ is sufficiently 
large. 
Hence:
 $$DF_{\lambda(\lambda+T)}(x,y)=\sum_{ Y: \dist_Y(x,y) \geq \lambda(\lambda+T) } \dist_Y(x,y) \leq 2\lambda \sum_{ Y: \dist_Y(x,y) \geq \lambda(\lambda+T) } h_Y\leq 2\lambda \ H_T(x,y),$$
 as required.
\end{proof}

The following will be used for one of the corollaries on quotients of mapping class groups in Section \ref{sec:proof_quotients}, namely Corollary \ref{cor:separable}, and in the proof of Theorem~\ref{thm:MCG_quotient}.\eqref{item:poc}. This result appears as the (3) implies (2) implication of {\cite[Theorem B]{ABD}}.

\begin{lemma}[{\cite{ABD}}]\label{lem:conv_cocpt=cobounded}
 Let $(G,\mathfrak S)$ be an HHG, and let $Q<G$ be such that orbit maps to $\mathcal CS$ are quasi-isometric 
embeddings, where $S$ is the $\nest$--maximal element of $\mathfrak S$. Then there exists $\kappa$ such that 
for all $Y\in\mathfrak S-\{S\}$ we have $\dist_Y(g,h)\leq \kappa$ for all $g,h\in Q$.
\end{lemma}

\begin{proof}
 A detailed proof appears in \cite{ABD}, here is a sketch.
 
 In view of the bounded geodesic image axiom, if $\dist_Y(g,h)$ is large, then $\rho^Y_S$ lies uniformly close 
to a geodesic from $\pi_S(g)$ to $\pi_S(h)$. Again in view of the bounded geodesic image axiom, along such a 
geodesic the projection to $\mathcal C(Y)$ is coarsely constant far from $\rho^Y_S$, and similarly for a 
neighborhood of the geodesic. Moreover, such a geodesic stays close to the quasiconvex subset $\pi_S(Q)$ of 
$\mathcal C(S)$, and in particular we see that $Q$ contains elements $g',h'$ so that:
 \begin{itemize}
  \item $\pi_S(g'),\pi_S(h')$ lie within uniformly bounded distance in $\mathcal C(S)$,
  \item $\dist_Y(g,h)$ differs from $\dist_Y(g',h')$ by a uniformly bounded amount.
 \end{itemize}
 
In view of the first item there are finitely many possible pairs $(g',h')$ up to the $Q$--action. In 
particular, there is a bound on $\dist_Y(g',h')$, whence a bound on $\dist_Y(g,h)$ as required.
\end{proof}

\section{Hyperbolicity of $Y_\Delta$}\label{sec:y_hyp}
Let $(X,W,\delta,n)$ be a combinatorial HHS.  The goal of this section is to show that $Y_\Delta$ is uniformly hyperbolic 
whenever $\Delta$ is a non-maximal simplex of $X$ for which $\diam(\mathcal C(\Delta))\geq \delta$.  We need this for several reasons in the proof of Theorem~\ref{thm:hhs_links}.  The most fundamental reason is: to construct an HHS structure on $W$, we will need the projections $\pi_{\Delta}:W\to\mathcal C(\Delta)$ from Definition~\ref{defn:projections} to be well-defined, coarsely Lipschitz coarse maps.  This follows from the fact that $\mathcal C(\Delta)\to Y_\Delta$ is a uniform quasi-isometric embedding, once we show that $Y_\Delta$ is hyperbolic.

We say that a constant $K$ is 
\emph{uniform} if $K$ depends only on $\delta$ and $n$.

\subsection{Preliminary lemma about hyperbolic graphs}\label{subsec:prelim}
Given a graph $Z$, we say that a set $\mathcal V$ of vertices of $Z$ is \emph{discrete} if no two elements of 
$\mathcal V$ are joined by an edge of $Z$.

The following lemma will be used inductively to prove hyperbolicity of $Y_\Delta$. The ``moreover'' part will not be used to prove hyperbolicity of $Y_\Delta$, but rather an additional statement, namely Lemma~\ref{lem:moreover}, which will be required in Section~\ref{sec:hieromorphisms}.

There are various ways to prove Lemma~\ref{lem:discrete_hyperbolic}; we chose the way that serves as a warm-up for the proof of Theorem \ref{thm:hhs_links}.

\begin{lemma}\label{lem:discrete_hyperbolic}
 For every $\delta$ there exists $\delta'$ with the following property. Let $Z$ be a $\delta$--hyperbolic graph 
and let $\mathcal V$ be a discrete collection of vertices. For each $v\in\mathcal V$, let $Z_v$ be the induced 
subgraph of $Z$ with vertex set $Z^{(0)}-\{v\}$.

Suppose that the following hold for all $v\in\mathcal V$:
\begin{itemize}
 \item $\link(v)$ is $\delta$--hyperbolic;
 \item $\link(v)$ is $(\delta,\delta)$--quasi-isometrically embedded in $Z_v$.
\end{itemize}
Then the induced subgraph $Z_{\mathcal V}$ of $Z$ with vertex set $Z^{(0)}-\mathcal V$ is 
$\delta'$--hyperbolic.

Moreover, for each $\lambda$, there exists $\mu=\mu(\lambda,\delta)$ such that the following holds.  

Let $Q\subseteq Z$ be a $(\lambda,\lambda)$--quasi-isometrically embedded induced subgraph of $Z$ such that 
$Q$ contains the star of $v$ whenever $v\in\mathcal V$ belongs to $Q$.  Then  the induced subgraph $Q_{\mathcal 
V}$ of $Z$ with 
vertex set 
$Q^{(0)}-\mathcal V\cap Q^{(0)}$ is $(\mu,\mu)$--quasi-isometrically 
embedded in $Z_{\mathcal V}$.
\end{lemma}

\begin{proof}[Proof of Lemma~\ref{lem:discrete_hyperbolic}]
We will equip $Z_{\mathcal V}$ with a hierarchically hyperbolic space structure $(Z_{\mathcal V},\mathfrak F)$ in which $\mathfrak F$ has 
empty orthogonality relation; from results in~\cite{HHS:quasiflats} we will then get that 
$Z_{\mathcal V}$ is $\delta'$--hyperbolic 
for some $\delta'$ depending on the HHS constants, which we shall see in this case depend only on $\delta$.  

Let $\mathfrak F=\{S\}\cup\mathcal V$.  The associated hyperbolic space $\mathcal C S$ is $Z$, and for each 
$v\in\mathcal V$, the associated space $\mathcal C v$ is 
$\link(v)$ (which is $\delta$--hyperbolic by hypothesis).  We declare $v\propnest S$ for all $v\in \mathcal V$; there is no other nontrivial nesting 
relation, and the orthogonality relation is empty. Therefore, two elements of $\mathfrak F$ are transverse if and only if 
they are distinct elements of $\mathcal V$.  The relation $\nest$ is easily seen to be a partial 
order with a unique maximal element and bounded chains (the complexity is $2$).

Define the projection $\pi_S:Z_{\mathcal V}\to \mathcal C S=Z$ to be the inclusion.  The projection $\pi_v:Z_{\mathcal V}\to \mathcal C v=\link(v)$ is 
defined as follows: given $x\in Z_{\mathcal V}$, consider all geodesics $\alpha$ in $Z$ from $x$ to $v$, and let $x'$ be the entry 
point of $\alpha$ in $\link(v)$.  Then $\pi_v(x)$ is the set of all such $x'$.  By Claim~\ref{claim:K} below, $\pi_v(x)$ has 
diameter bounded in terms of $\delta$ only.

(Here, we have used discreteness of $\mathcal V$ to ensure that $x'$ exists.  Indeed, discreteness ensures that $\link(v)\subset Z_{v}$, so that any geodesic in $Z$ from $x$ to $v$ must pass through $\link(v)$.)

\begin{claim}\label{claim:K}
 There exists $K=K(\delta)$ so that the following holds. Let $x,y\in Z_{\mathcal V}$ and let $\alpha,\beta$ be 
geodesics in $Z$ starting, respectively, at $x,y$ and intersecting $\link(v)$ only at their 
other respective endpoints $x',y'$. If $d_{\link(v)}(x',y')\geq K$, then any geodesic $\gamma$ in 
$Z$ from $x$ to $y$ contains $v$. 
 \end{claim}
\renewcommand{\qedsymbol}{$\blacksquare$}
\begin{proof}[Proof of Claim~\ref{claim:K}]
Let $\gamma$ be a geodesic in $Z$, joining $x$ to $y$, and not containing $v$. Then each vertex of $\gamma$ belongs to $Z_v$, i.e., $\gamma$ is a path in $Z_v$.  Moreover, since 
$\gamma$ is a geodesic of $Z$, it is again a geodesic of $Z_v$.  Fix a geodesic $q$ of $Z$ from $x'$ to $y'$; since 
$x',y'\in\link(v)$, we have $|q|\le 2$. 

Consider the quadrilateral in $Z$ with sides $\alpha,\beta,\gamma,q$.  By $\delta$--hyperbolicity of $Z$, this quadrilateral is $2\delta$--thin.  Let $\alpha'$ be the geodesic in $Z$ obtained by traveling backward along $\alpha$, starting from $x'$, for distance $\min\{10\delta+1,|\alpha|\}$, and define $\beta'$ analogously, using $y'$ and $\beta$.  By construction, $\alpha',\beta'$ are paths in $Z_v$, since $\alpha,\beta$ cannot contain $v$.  The endpoints $a,b$ of $\alpha',\beta'$ are $2\delta$--close to points $a',b'\in\gamma$ with $d_Z(a',b')\leq 24\delta + 4$.  If $|\alpha|\leq 10\delta+1$ then $a=a'$ and if $|\beta|\leq 10\delta+1$, we have $b=b'$.  Otherwise, $a,a'$ are both at least $8\delta+1$--far from $v$, so in any case a geodesic $[a,a']$ in $Z$ must lie in $Z_v$.  

Thus, using $\alpha',\beta',[a,a'],[b,b']$ and the part of $\gamma$ between $a'$ and $b'$, we construct a path of length at most $100(\delta+1)$ that joins $x'$ to $y'$ and lies in $Z_v$.

Hence, $d_{Z_v}(x',y')\leq 100(\delta+1)$.  Since $\link(v)$ is $(\delta,\delta)$--quasi-isometrically embedded in $Z_v$, we obtain $d_{\link(v)}(x',y')\leq \delta(100(\delta+1)+\delta)$, as required.   
\end{proof}
 \renewcommand{\qedsymbol}{$\Box$}

Moreover, the coarse map $\pi_v$ is $(K,K)$--coarsely Lipschitz, by Claim~\ref{claim:K}.  Indeed, if $x,y\in Z_{\mathcal V}$ are 
adjacent vertices, then the claim implies that their images lie at distance at most $K$ in $\link(v)$.

Let $\rho^v_S\subset Z$ be $\{v\}$.  Define $\rho^S_v:Z\to\link(v)$ by letting 
$\rho^S_v(x)=\pi_v(x)$ if $x\in Z_v$ and defining $\rho^S_v(w)$ arbitrarily when $w$ is in the open star of $v$.  The
consistency inequality for nesting holds by definition.

For $v,w\in \mathcal V$ distinct vertices, let $\rho^v_w=\pi_w(v)$. Let us verify the consistency inequality for transversality. Consider distinct $v,w\in\mathcal V$ and $x\in Z_\mathcal V$. Suppose that $d_{\link(v)}(\rho^w_v,\pi_v(x))\geq K$, for $K$ as in Claim \ref{claim:K}. Then any geodesic $\gamma$ from $x$ to $w$ passes through $v$. In particular, the set of entrance points in $\link(w)$ of geodesics from $x$ to $w$ is contained in the set of entrance points in $\link(w)$ for geodesics from $v$ to $w$. This implies $\pi_w(x)\subseteq \pi_w(v)$, so the two coarsely coincide.

Large links holds again by Claim \ref{claim:K}, which implies that given $x,y\in Z_{\mathcal V}$, the set of 
all $v\in\mathcal V$ on whose links $x,y$ have large projections are vertices of a geodesic in $X$ from $x$ to 
$y$.  

It remains to verify bounded geodesic image, partial realization, and 
uniqueness.  

\textbf{Bounded geodesic image:}  Let $\gamma$ be a geodesic in $Z$.  If $\gamma$ is disjoint from $\link(v)$, 
then $\diam(\rho^S_v(\gamma))\leq 
K$, by Claim~\ref{claim:K}.  This verifies bounded geodesic image.  

\textbf{Partial realization:} Let $p\in Z=\fontact S$ be the coordinate to realize.  If $p\in Z_{\mathcal V}$, 
let 
$x=p$.  If $p$ lies in the open star of $v$, let $x\in\link(v)$ be chosen arbitrarily.  Then 
$\dist_S(\pi_S(x),p)\le 1$, and no element of $\mathfrak F$ is transverse to, or contains, $S$.  If $p\in \mathcal C 
v=\link(v)$ is the coordinate to realize, let $x=p$.  Then $\dist_v(\pi_v(x),p)=0$ and $\dist_S(\pi_S(x),\rho^v_S)=1$. Also, $\dist_w(\rho^v_w,x)$ is bounded since $\pi_v$ is coarsely Lipschitz.  This verifies partial realization.  (Here, again, we used discreteness to ensure that $\link(v)\subset Z_{\mathcal V}$.)
 
\textbf{Uniqueness:}  Let $\kappa\ge0$ and fix $x,y\in Z_{\mathcal V}$.  Suppose that $d_Z(x,y)\leq \kappa$, and that $$d_{\link(v)}(\pi_v(x),\pi_v(y))\leq \kappa$$ for all $v\in\mathcal V$.  We need to bound $d_{Z_{\mathcal V}}(x,y)$ in terms of $\kappa$ and $\delta$.

Let $\eta$ be a geodesic 
of $Z_{\mathcal V}$ from $x$ to $y$.  If $\eta$ is a geodesic of $Z$, then $d_{Z_{\mathcal V}}(x,y)=d_Z(x,y)\leq\kappa$, and we are done. 
 So, assume that each geodesic 
$\gamma$ in $Z$ from $x$ to $y$ passes through some $v\in\mathcal V$, and fix such a $\gamma$.  

Let $v_1,\ldots,v_n\in\mathcal V$ be the vertices in $\mathcal V$ contained in $\gamma$.  Let $x_i,y_i$ be the entry and exit points of $\gamma$ in $\link(v_i)$, for $1\leq i\leq n$.  By discreteness, all of the $x_i,y_i$ lie in $Z_{\mathcal V}$, and are in particular all distinct from all $v_j$.  So, we can write
\[\gamma=\gamma_0p_1\gamma_1p_2\cdots\gamma_{n-1}p_n\gamma_n,\]

where each $\gamma_n$ is a (possibly trivial) geodesic in $Z_{\mathcal V}$ and each $p_i$ is the path of length $2$ from $x_i$ to $y_i$ passing through $v_i$.  Hence $$d_{Z_{\mathcal V}}(x,y)\leq \sum_{i=1}^nd_{\link(v_i)}(x_i,y_i)+\sum_{i=1}^n|\gamma_i|.$$

Note that each term of the left sum is at most $d_{\link(v_i)}(\pi_{v_i}(x),\pi_{v_i}(y))+2K$, where $K=K(\delta)$ bounds the diameter of $\pi_{v_i}$--projections of points.  So, by our assumption that all projections of $x,y$ are $\kappa$--close, we have

$$d_{Z_{\mathcal V}}(x,y)\leq n\kappa + n(2K-2) +d_Z(x,y)\leq (n+1)\kappa + n(2K-2),$$

since the $Z$--geodesic $\gamma$ is obtained by concatenating the $\gamma_i$ along with $n$ paths of length $2$.  Finally, each $v_i$ contributes at least $1$ to the length of $\gamma$, so $n\leq \kappa$.  We conclude that $d_{Z_{\mathcal V}}(x,y)\leq \kappa^2+\kappa(2K-1)$.  Since this depends only on $\kappa$ and $\delta$, uniqueness is verified.

Hence $(Z_{\mathcal V},\mathfrak F)$ is an HHS with empty orthogonality relation.  The constants implicit in the definition 
of an HHS depend only on $\delta$.  Then the proof of~\cite[Corollary 2.15]{HHS:quasiflats} (relying on Theorem 2.1 of~\cite{HHS:quasiflats}) 
implies that $Z_{\mathcal V}$ is a \emph{coarse median space} of rank $1$, where the constants/function in the definition of 
a coarse median space (see~\cite{Bowditch:coarse_median}) depend only on the HHS constants and hence only depend on $\delta$. 
Hence \cite[Theorem 2.1]{Bowditch:coarse_median} implies that 
$Z_{\mathcal V}$ is $\delta'$--hyperbolic, where $\delta'$ depends only on the coarse median constants and hence only on 
$\delta$.

\textbf{The ``moreover'' clause:} We now prove the ``moreover clause'' essentially by proving that $Q_{\mathcal 
V}$ is an HHS whose structure is compatible with that of $Z_{\mathcal V}$. Let $Q$ be an induced subgraph of 
$Z$ as in the statement.  Fix $v\in\mathcal 
V$.  Let $Q_v$ be the induced subgraph of $Q$ spanned by $Q^{(0)}-\{v\}$.  Define $q:Q_v\to2^{\link(v)}$ as 
follows.  Let $x\in Q_v$ be a vertex.  For each geodesic $\gamma$ in $Q$ from $x$ to $v$, let $x_\gamma$ be the 
entrance point of $\gamma$ in $\link(v)$.  Let $q(x)$ be the set of points $x_\gamma$, as $\gamma$ varies over 
the $Q$--geodesics from $x$ to $v$.  

For notational purposes, set $z=\pi_S:Z_{\mathcal V}\to 2^{\link(v)}$ as above.  Namely, given $x\in Z_{\mathcal V}$, consider all $Z$--geodesics $\alpha$ from 
$x$ to $v$, and for each $\alpha$, let $x_\alpha$ be the entrance point of $\alpha$ in $\link(v)$.  Let $z(x)$ 
be the union of the $x_\alpha$.   By Claim~\ref{claim:K}, $z(x)$ has diameter bounded in terms of $\delta$, 
i.e., $z$ is a (uniformly) well-defined coarse map.  We now show that $q$ is a well-defined coarse map that 
uniformly coarse coincides with the restriction of $z$ to $Q_v$.

Fix $x\in Q$ and let $\gamma,\alpha$ be as above.  By hypothesis, $\gamma$ is a 
$(\lambda,\lambda)$--quasigeodesic of $Z$ and hence $\epsilon$--fellow-travels in $Z$ with $\alpha$, where 
$\epsilon=\epsilon(\lambda,\delta)$.

Let $\alpha'$ be the subpath of $\alpha$ from $x$ to $x_\alpha$, so that $\alpha'$ is a geodesic of $Z$ 
avoiding $v$ and hence a geodesic of $Z_v$.  Let $\gamma'$ be the subpath of $\gamma$ from $x$ to $x_\gamma$, 
so that $\gamma'$ is a $(\lambda,\lambda)$--quasigeodesic of $Z$ avoiding $v$.

If $|\gamma'|<10^3\lambda\epsilon+\lambda^2=C$, then $|\alpha'|<C$, so 
$(\alpha')^{-1}\gamma'$ is a path in $Z_v$ of length $2C$ from $x_\alpha$ to $x_\gamma$.  
Since $\link(v)$ is $(\delta,\delta)$--quasi-isometrically embedded in $Z_v$, this gives an upper bound on 
$\dist_{\link(v)}(x_\alpha,x_\gamma)$ in terms of $\lambda,\epsilon,\delta$.

Suppose $|\gamma'|>C$.  Then $|\alpha'|>10^3\epsilon$.  Hence $\alpha',\gamma'$ contain points $a,c$ such 
that $\dist_Z(a,c)\leq\epsilon$ and $\dist_Z(a,v),\dist_Z(c,v)\in[100\epsilon,10^4\epsilon]$.  So, any 
$Z$--geodesic $\eta$ from $c$ to $a$ is a $Z_v$--geodesic of length at most $\epsilon$.  Hence, by travelling 
along $\gamma'$ from $x_\gamma$ to $c$, then along $\eta$, and then along $\alpha'$ from $a$ to $x_\alpha$, we 
obtain a path in $Z_v$ of length at most $10^4\epsilon+\lambda 10^4\epsilon+\lambda$.  As above, this means 
that $\dist_{\link(v)}(x_\alpha,x_\beta)$ is bounded above in terms of $\lambda,\epsilon,\delta$. This shows that $q(x)$ is a well-defined coarse map that uniformly coarsely coincides with $z(x)$ for $x\in Q$.

Now consider the HHS structure $(Z_{\mathcal V},\mathfrak F)$ constructed in the first part of the proof.  
Recall that the index set is $\mathfrak F=\{S\}\cup\mathcal V$, the hyperbolic space associated to $S$ is $Z$, 
and the hyperbolic space associated to $v\in\mathcal V$ is $\link(v)$.  The projection $\pi_S:Z_{\mathcal 
V}\to Z$ is the inclusion, and the projection $\pi_v:Z_{\mathcal V}\to\link(v)$ is the map $z$ discussed above.

Now, let $Q_{\mathcal V}$ be the induced subgraph of $Q$ produced by removing the vertices in $\mathcal V$.  
By 
the same argument that we used for $Z_{\mathcal V}$, the pair $(Q_{\mathcal V},\mathfrak F)$ is an HHS, where 
$\mathfrak F=\{S\}\cup\mathcal V$, the hyperbolic space associated to $S$ is $Q$, and the hyperbolic space for 
each $v$ is $\link(v)$.  The projection $Q_{\mathcal V}\to Q$ is again inclusion, and $\pi_v:Q_{\mathcal V}\to 
\link(v)$ is the map $q:Q_v\to 2^{\link(v)}$ defined above.  (We can apply the same argument we used for $Z$ 
because each $\link(v)$ is uniformly quasi-isometrically embedded in $Q_v$, since $\link(v)$ is 
$(\delta,\delta)$--quasi-isometrically embedded in $Z_v$ and $\link(v)\subset Q_v\subset Z_v$.)  Note that the 
constants for this HHS structure are different from those for $(Z_{\mathcal V},\mathfrak F)$, but they still 
depend only on $\delta$ and $\lambda$.

To conclude, we need to prove that $Q_{\mathcal V}\hookrightarrow Z_{\mathcal V}$ is a quasi-isometric 
embedding with constants depending on $\delta$ and $\lambda$.  We saw above that for each $v\in\mathcal V$, 
the projections $q,z$ to $\link(v)$ used in the two HHS 
structures coarsely coincide on $Q_{\mathcal V}$ (constants depend on $\delta,\lambda$).  Moreover, the inclusions 
$Q_{\mathcal V}\to 
Q,Q_{\mathcal V}\to Z$ obviously coincide.  

There exists $\lambda'=\lambda'(\lambda,\delta)$ such that the following hold for $a,b\in Q_{\mathcal V}$:
$$\dist_Q(a,b)\asymp_{\lambda',\lambda'}\dist_Z(a,b)$$  and, for each $v\in\mathcal V$, 
$$\dist_v(q(a),q(b))\asymp_{\lambda',\lambda'}\dist_v(z(a),z(b)).$$

Hence, by Theorem~\ref{thm:DF_plus_epsilon}, there exists $T_0$, depending on $\lambda'$ and the HHS constants 
for $Z_{\mathcal V}$, such that the following holds.  Let $T_1\geq T_0$.  Then there 
exists $\mu_1$, depending on $T_1$ and the HHS constants, such that $$\dist_{Z_{\mathcal 
V}}(a,b)\asymp_{\mu_1,\mu_1}\ignore{\dist_Q(a,b)}{T_1}+\sum_{v\in\mathcal V}\ignore{\dist_v(q(a),q(b)}{T_1}.$$

By the distance formula, applied to the HHS structure on $Q_{\mathcal V}$, the following holds.  For all 
sufficiently large $T_1$, we have $\mu_2$, depending only on the HHS constants for $Q_{\mathcal V}$ 
(which depend only on $\delta,\lambda$) such that $$\dist_{Q_{\mathcal 
V}}(a,b)\asymp_{\mu_2,\mu_2}\ignore{\dist_Q(a,b)}{T_1}+\sum_{v\in\mathcal V}\ignore{\dist_v(q(a),q(b)}{T_1}.$$

Hence there exists $\mu=\mu(\delta,\lambda)$ such that $\dist_{Q_{\mathcal V}}(a,b)\asymp_{\mu,\mu}\dist_{Z_{\mathcal 
V}}(a,b)$, as required.
\end{proof}

\subsection{$Y_\Delta$ is hyperbolic for $\mathcal C(\Delta)$ large}\label{subsec:inductive_hyperbolic}
Let $\Delta$ be a non-maximal simplex of $X$, and assume that $\diam(\mathcal C(\Delta))\geq\delta$.  Our main goal is to prove that $Y_\Delta$ is uniformly hyperbolic.  We also 
have a second goal, which is to prove that if $[\Delta]\nest[\Sigma]$ then 
$Y_\Delta\cap\mathcal C(\Sigma)$ is (uniformly) quasi-isometrically embedded in $Y_\Delta$.

We will apply the former conclusion throughout the rest of the paper.  The latter conclusion is used in the proof of 
Proposition~\ref{prop:hieromorphism}.  Here is the formal statement:

\begin{prop}\label{prop:Y_delta_hyperbolic}
For every $\delta, n$ there exists $\delta'$ so that the following holds. Let $(X,W,\delta,n)$ be a combinatorial HHS, and let $\Delta$ be a non-maximal simplex of $X$ such that $\diam(\mathcal C(\Delta))\geq\delta$. Then $Y_\Delta$ is $\delta'$-hyperbolic.

Moreover, let $\Sigma$ be a non-maximal simplex of $X$ with $[\Delta]\nest[\Sigma]$.  Then the inclusion 
$Y_{\Delta}\cap\mathcal C(\Sigma)\hookrightarrow Y_{\Delta}$ is a $(\delta',\delta')$-quasi-isometric embedding.
\end{prop}

The hyperbolicity will follow from Lemma~\ref{lem:Y_delta_hyperbolic} below.  The QI-embedding will follow from Lemma~\ref{lem:moreover}. 

For clarity, whenever we say that a 
constant is uniform we mean that it depends on $\delta,n$ only.

We now fix $\Delta$ as in the statement of Proposition~\ref{prop:Y_delta_hyperbolic} and define a sequence of spaces interpolating between $X^{+W}$ and $Y_{\Delta}$ as follows. 

\begin{defn}[Co-level]\label{defn:colevel}
Let $\Sigma$ be a non-maximal simplex.  We define the \emph{co-level} $cl[\Sigma]$ inductively as follows.  
First, define $cl[\emptyset]=0$.  

For $k\geq 0$, we say that $[\Sigma]$ has \emph{colevel at least $k+1$} if there exists $[\Sigma']$ with colevel at least $k$ and $[\Sigma]\propnest[\Sigma']$.  The co-level of $[\Sigma]$, denoted $cl[\Sigma]$, is the maximal $k$ such that $[\Sigma]$ has co-level at least $k$.
\end{defn}

We will be interested in the co-level of classes $[\Sigma]$ with $[\Delta]\nest [\Sigma]$ for our fixed 
simplex $\Delta$. Observe that for any simplices $[\Pi],[\Sigma]$, if $[\Pi]\propnest[\Sigma]$ then $cl[\Pi]>cl[\Sigma]$.

\begin{defn}\label{defn:coned_y}
 Let $\mathcal U=\{[\Sigma]:[\Delta]\nest[\Sigma]\}$. For $0\leq k\leq cl[\Delta]$, define $Y^k_{\Delta}$ to be 
the graph 
obtained from $Y_\Delta$ by connecting all pairs of vertices of $\link(\Sigma)\cap Y_\Delta$ for all 
$[\Sigma]\in \mathcal U$ 
with $cl[\Sigma]> k$.
\end{defn}

\begin{rem}(Heuristic picture of $Y^k_\Delta$)
 It might be useful to keep in mind the following rough picture of $Y^k_{\Delta}$. First, $Y^k_\Delta$ contains a quasi-isometric copy of each $\mathcal C(\Sigma)$ for $[\Sigma]$ of co-level $k$, and any two of these have bounded coarse intersection. In fact, these coarse intersections are best thought of to consist of coned-off copies of certain $\mathcal C(\Lambda)$ for $\Lambda$ of co-level strictly larger than $k$. Finally, $Y^k_{\Delta}$ does \emph{not} contain a copy of $\mathcal C(\Lambda)$ for $\Lambda$ of co-level strictly smaller than $k$, but rather a blown-up version of it similar to a $Y_\Lambda$-space for a link (this might become clearer when, in Section \ref{sec:hieromorphisms}, we discuss induced combinatorial HHS structures on links). We do not formulate the proof of hyperbolicity in these terms, but essentially we will show inductively that $Y^k_{\Delta}$ is hyperbolic relative to the aforementioned copies of the $\mathcal C(\Sigma)$ for $[\Sigma]$ of co-level $k$, with coned-off graph quasi-isometric to $Y^{k-1}_{\Delta}$.
\end{rem}

We will need the following auxiliary lemma, which is how we use the hypothesis on the diameter of $\mathcal C(\Delta)$.

\begin{lemma}\label{lem:Y_delta_containment}
Let $\Sigma,\Delta$ be non-maximal simplices of $X$ such that $\diam(\mathcal C(\Delta))\geq\delta$ and $[\Delta]\nest[\Sigma]$.  Then $\sat(\Sigma)\subset\sat(\Delta)$ and hence $Y_\Delta\subset Y_\Sigma$.
\end{lemma}

\begin{proof}
By Lemma~\ref{lem:SC_nesting}, there exists a simplex $\Pi$ of $\link(\Sigma)$ such that $[\Delta]=[\Sigma\star\Pi]$.  If $\Sigma'$ is a simplex with $\link(\Sigma)=\link(\Sigma')$, then the simplex $\Sigma'\star\Delta$ exists, and $\link(\Sigma'\star\Pi)=\link(\Sigma')\cap\link(\Pi)=\link(\Sigma\star\Pi)=\link(\Delta)$.  By definition, any $v\in\sat(\Sigma)$ is contained in some such $\Sigma'$, and hence in $\Sigma'\star\Pi$, from which it follows from the above that $v\in\sat(\Delta)$, as required.   
\end{proof}

We also need a graph $Z^k_\Delta$, which is a quasi-isometry model of $Y_\Delta^k$.  By definition, $Y^k_\Delta$ is obtained from $Y_\Delta$ by electrifying $\link(\Sigma)\cap Y_\Delta$ when $[\Delta]\nest[\Sigma]$ and $[\Sigma]$ has co-level at least $k+1$.  In $Z^k_\Delta$, for the $[\Sigma]$ with co-level exactly $k+1$, we replace the electrification by coning.  The reason for working with $Z^k_\Delta$ is that it enables the use of Lemma~\ref{lem:discrete_hyperbolic}, once we check that the links of the cone-points are hyperbolic and are quasi-isometrically embedded in the graph obtained from $Z^k_\Delta$ by removing them (i.e., in $Y^{k+1}_\Delta$).  This graph will be needed in the proof of Lemma~\ref{lem:Y_delta_hyperbolic}, and also in the proof of Lemma~\ref{lem:moreover}.

\begin{defn}\label{defn:Z-k}
Let $Z^{k}_\Delta$ be the graph obtained from $Y_\Delta$ by
\begin{itemize}
\item connecting all pairs of vertices of $\link(\Sigma)\cap Y_\Delta$ for all $[\Sigma]\in \mathcal U$ with $cl[\Sigma]> 
k+1$,
\item adding a vertex $v_{[\Sigma]}$ for each $[\Sigma]\in \mathcal U$ with $cl[\Sigma]=k+1$, and connecting any such 
$v_{[\Sigma]}$ to all vertices of $\link(\Sigma)\cap Y_\Delta$.
\end{itemize}
\end{defn}

We now inductively prove hyperbolicity of the $Y^k_\Delta$, in the lemma below. 

\begin{lemma}\label{lem:Y_delta_hyperbolic}
The following hold provided $\diam(\mathcal C(\Delta))\geq\delta$.
\begin{enumerate}
\item \label{item:qi}$Y^0_\Delta$ is uniformly quasi-isometric to $X^{+W}$.
\item \label{item:y_delta}$Y^{cl[\Delta]}_{\Delta}=Y_\Delta$.
\item \label{item:hyperbolic}$Y^k_\Delta$ is hyperbolic for $0\leq k\leq cl[\Delta]$, with 
uniform hyperbolicity constant.
\end{enumerate}
In particular, $Y_\Delta$ is uniformly hyperbolic.
\end{lemma}

\begin{proof}
Assertions~\eqref{item:qi} and~\eqref{item:y_delta} hold by construction.

For~\eqref{item:hyperbolic}, we prove hyperbolicity by induction on $k$. The case $k=0$ holds by 
item~\eqref{item:qi} and Definition~\ref{defn:combinatorial_HHS}. Suppose that $Y^{k}_\Delta$ is hyperbolic 
for some $k\geq1$ and $k< cl[\Delta]$. 

Since $Z^{k}_\Delta$ is uniformly quasi-isometric to $Y^{k}_\Delta$, we have that $Z^{k}_\Delta$ is 
uniformly hyperbolic. Moreover, $Y^{k+1}_\Delta$ is obtained from $Z^{k}_\Delta$ by removing the open star of 
all vertices $v_{[\Sigma]}$ described above (by construction, these vertices form a discrete set). In view of 
Lemma \ref{lem:discrete_hyperbolic}, it suffices to 
show that the link $\link(v_{[\Sigma]})$ of any $v_{[\Sigma]}$ in $Z^{k}_\Delta$ is hyperbolic and 
quasi-isometrically embedded in $Y^{k+1}_\Delta$.

\textbf{Hyperbolicity:} The vertex set of $\link(v_{[\Sigma]})$ is by definition $\link(\Sigma)\cap Y_\Delta$. We first construct a Lipschitz map $\psi: \link(v_{[\Sigma]})\to \mathcal C(\Sigma)$.  Since $\link(v_{[\Sigma]})=\link(\Sigma)\cap Y_\Delta$, whose vertex set is contained in $\mathcal C(\Sigma)$, we first define $\psi$ to be the inclusion on vertices.

Next, let $e$ be an edge of $\link(v_{[\Sigma]})$.  We claim that the endpoints of $e$ are joined by a path of uniformly bounded length in $\mathcal C(\Sigma)$, which suffices to extend $\psi$ to a (uniformly) Lipschitz map.  We check the claim by considering the various possibilities for $e$.  By the definition of $Z^k_\Delta$, one of the following holds: 
\begin{compactitem}
    \item $e$ is an edge of $Y_\Delta$ with endpoints in $\link(\Sigma)$, and hence is already an edge of $\mathcal C(\Sigma)$.  In this case, the uniformly bounded path is $e$ itself.

    \item There exists $\Sigma'\in \mathcal U$ with 
$cl[\Sigma']>cl[\Sigma]$, and the endpoints of $e$ are in $\link(\Sigma')$.  We claim that $\link(\Sigma')\cap \link(\Sigma)$ has uniformly bounded diameter in $\mathcal C(\Sigma)$.  Indeed, first note that since 
$cl[\Sigma']>cl[\Sigma]$, we have $[\Sigma]\not\nest[\Sigma']$.  In particular, $\sat(\Sigma')$ contains a 
vertex $\alpha$ not contained in $\sat(\Sigma)$, i.e., $\alpha\in Y_\Sigma$.  Now, every vertex of 
$\link(\Sigma)\cap\link(\Sigma')$ is joined by an edge of $Y_\Sigma$ to $\alpha$, so 
$diam_{Y_{\Sigma}}(\link(\Sigma)\cap\link(\Sigma'))\leq 2$.  Since $\mathcal C(\Sigma)\hookrightarrow 
Y_{\Sigma}$ is a uniform quasi-isometric embedding, we see that $\link(\Sigma)\cap\link(\Sigma')$ is uniformly 
bounded in $\mathcal C(\Sigma)$, which provides the desired path.
\end{compactitem}

Having thus shown that $\psi$ is uniformly Lipschitz, we show it is a quasi-isometry.  For this, we will define a map $\phi:\mathcal C(\Sigma)\to \link(v_{[\Sigma]})$ and show it is a uniformly Lipschitz quasi-inverse for $\psi$.

First, let $\phi$ be the identity on $\link(\Sigma)\cap Y_\Delta \subset \mathcal C(\Sigma)$.  Fix an arbitrary vertex $w\in \link(\Delta)\subseteq \mathcal C(\Sigma)$, and complete the definition of $\phi$ on vertices by sending every $v\in\link(\Sigma)-Y_\Delta$ to $w$.

Second, we check that $\phi$ is a quasi-inverse for $\psi$.  On vertices of $\link(\Sigma)\cap Y_\Delta$, this is clear; the two maps are literally inverses.  On the other hand, if $v\in\link(\Sigma)-Y_\Delta$, then by definition $\phi(v)=w\in\link(\Delta)$.  Now, since $v\not\in Y_\Delta$, then $v\in \sat(\Delta)$, so $v$ and $w$ are joined by an edge of $X$, since $w\in\link(\Delta)$.  Now, since $\link(\Delta)\subset\link(\Sigma)$, we have $w\in\mathcal C(\Sigma)$, and since $v\in \mathcal C(\Sigma)$, this edge is in $\mathcal C(\Sigma)$, and hence $\psi(\phi(v))=w$ is distance $1$ in $\mathcal C(\Sigma)$ from $v$.  This establishes that $\psi$ and $\phi$ are quasi-inverses.

We are left to show that $\phi$ is coarsely Lipschitz.  This reduces to showing that if  $v\in \link(\Sigma)-Y_\Delta$ is connected to $v'\in \link(\Sigma)\cap Y_\Delta$, then 
$v'$ and $w$ lie within bounded distance.  

First note, as above, that  $v\in \link(\Sigma)-Y_\Delta$ is equivalent to  $v\in 
\link(\Sigma)\cap \sat(\Delta)$. Since $\link(\Delta)\subset\link(\Sigma)$, every vertex of $\link(\Delta)$ is connected in $X$ to every vertex of $\Sigma$.  Together with the fact that $v\in\sat(\Delta)$, this implies that $\link(\Delta)\subset\link(\Sigma\star v)$.  In other words, $[\Delta]\nest[\Sigma\star v]$, i.e., $\Sigma\star v\in\mathcal U$. Also, clearly $[\Sigma\star v]\propnest [\Sigma]$, so that 
$cl[\Sigma\star v]>cl([\Sigma])$. From the definition of $Z^k_\Delta$, it now follows that any pair of vertices in $\link(\Sigma\star v)\cap Y_\Delta$ are 
connected in $\link(v_{[\Sigma]})$. 

Now we have two cases, reflecting that the edge of $\mathcal C(\Sigma)$ joining $v,v'$ is of one of two types:
\begin{compactitem}
    \item If $v$ and $v'$ are connected in $X$, then $v'\in \link(\Sigma\star v)\cap Y_\Delta$.  On the other hand, $w\in\link(\Delta)\subset\link(\Sigma\star v)\cap Y_\Delta$, so by the preceding discussion, $v'$ and $w$ are adjacent in $\link(v_{\Sigma})$, as required.

    \item The other possibility is that $v,v'$ are not connected in $X$, and are therefore joined by a $W$--edge of $\mathcal C(\Sigma)$.  By Definition~\ref{defn:combinatorial_HHS}.\eqref{item:C_0=C}, there exist maximal simplices $x,x'$ of $X$, connected in $W$, and, moreover, we have simplices $\Pi,\Pi'$, maximal in $\link(\Sigma)$, such that $x=\Sigma\star\Pi,x'=\Sigma\star\Pi'$, and $v\in\Pi,v'\in\Pi'$.

    We claim that there exists $v''\in \Pi\cap Y_\Delta$.  Indeed, if not, then $\Pi\subset\sat(\Delta)$.  Hence, for any $t\in\link(\Delta)^{(0)}\subset\link(\Sigma)$, we have a simplex $\Pi\star t\subset\link(\Sigma)$, contradicting maximality of $\Pi$.  Hence there is a vertex $v''\in \Pi\cap Y_\Delta$.

    Since $v''\in\Pi$, we have $v''\in\link(\Sigma)$.  Since $v\in\Pi$, there is an edge of $X$ joining $v,v''$, so $v''\in\link(\Sigma)\cap\link(v)=\link(\Sigma\star v)$.  Finally, we have chosen $v''\in Y_\Delta$.  So $v''\in\link(\Sigma\star v)\cap Y_\Delta$.

    Now, recall that, since $w$ and $v''$ are both in $\link(\Sigma\star v)\cap Y_\Delta$, they are joined by an edge of $\link(v_{[\Sigma]})$.  On the other hand, $v'$ and $v''$ are connected by a $W$--edge, since $x$ and $x'$ are $W$--adjacent.  Since $v',v''\in Y_\Delta$, we see that $v'$ is connected to $v''$ in $\link(v_{[\Sigma]})$.  Hence $v'$ and $w$ are at distance at most $2$ in $\link(v_{[\Sigma]})$, as required. 
\end{compactitem}

This completes the proof that $\psi$ is a uniform quasi-isometry, so by uniform hyperbolicity of $\mathcal C(\Sigma)$, which comes from Definition~\ref{defn:combinatorial_HHS}, we get the required hyperbolicity of $\link(v_{[\Sigma]})$.

\textbf{Quasi-isometric embedding of $\link(v_{[\Sigma]})$ in $Y^{k+1}_\Delta$:} To show that 
$\link(v_{[\Sigma]})$ is quasi-isometrically embedded in 
$Y^{k+1}_\Delta$ it suffices to find a coarsely Lipschitz extension $\hat{\psi}: Y^{k+1}_\Delta\to Y_\Sigma$ of 
$\psi$.  Indeed, given such a $\hat\psi$, consider the uniformly coarsely commutative diagram:
\begin{center}
$ 
 \begin{diagram}
  \node{\link(v_{[\Sigma]})}\arrow{e,t}{\psi}\arrow{s,J}\node{\mathcal C([\Sigma])}\arrow{s,J}\\
\node{Y^{k+1}_\Delta}\arrow{e,b}{\hat\psi}\node{Y_\Sigma}  
 \end{diagram}
$
\end{center}

\noindent where the inclusion on the right is a uniform quasi-isometric embedding (by 
Definition~\ref{defn:combinatorial_HHS}), the map $\psi$ is a uniform quasi-isometric embedding,  and the left 
inclusion is necessarily Lipschitz.  If we show that $\hat{\psi}$ is uniformly coarsely Lipschitz, the 
map $\link(v_{[\Sigma]})\to Y^{k+1}_\Delta$ will then be forced to be a uniform quasi-isometric 
embedding.

It is now that we use the assumption that $\diam(\mathcal C(\Delta))\geq\delta$, which gives us access to Definition~\ref{defn:combinatorial_HHS}.\eqref{item:chhs_delta} and its consequences.  Specifically, by Lemma~\ref{lem:Y_delta_containment}, we have $\sat(\Sigma)\subset\sat(\Delta)$.  It follows that the vertex set of $Y^{k+1}_\Delta$ is contained in the vertex set of $Y_\Sigma$, and at the level of vertex sets we declare $\hat{\psi}$ to be the inclusion.

The argument 
for why $\hat{\psi}$ gives a coarsely Lipschitz map is now similar to the argument for $\psi$ (it is still true 
that if $\link(\Sigma')\cap Y_\Sigma$ is not bounded, then $[\Sigma]$ is nested into $[\Sigma']$).  This 
completes the proof that $Y_\Delta$ is uniformly hyperbolic.
\end{proof}

\subsection{QI-embedding of $\mathcal C(\Sigma)\cap Y_\Delta$ in $Y_\Delta$}\label{subsec:QI-moreover}
We now turn to the QI-embedding part of Proposition~\ref{prop:Y_delta_hyperbolic}, which we restate for convenience:

\begin{lemma}\label{lem:moreover}
The inclusion $\mathcal C(\Sigma)\cap Y_\Delta\to Y_\Delta$ is a uniform quasi-isometric embedding, where $\Sigma$ is any non-maximal simplex of $X$ with $[\Delta]\nest[\Sigma]$.
\end{lemma}

The proof is similar in spirit to that of hyperbolicity, but a bit more technical.  Namely, we will construct a sequence of spaces embedded in the various hyperbolic spaces $Y^k_\Delta$ and prove inductively that they are quasi-isometrically embedded, with the last space in the sequence being $\mathcal C(\Sigma)\cap Y_\Delta$. As in the proof of hyperbolicity, we will also need intermediate spaces with suitable ``cone points'' in order to apply Lemma \ref{lem:discrete_hyperbolic}. Before the proof, we need some preliminary discussion to construct all the relevant spaces.

\subsubsection{Preliminaries on electrifying and coning off links inside $\mathcal C(\Sigma)\cap Y_\Delta$}  Recall that for $0\leq k\leq cl[\Delta]$, we have defined graphs $Y^k_\Delta$ and $Z^k_\Delta$ in Definition~\ref{defn:coned_y} and Definition~\ref{defn:Z-k} respectively.

Fix $[\Sigma]$ as in Lemma~\ref{lem:moreover}.  

\begin{defn}\label{defn:Q-and-R}
For $0\leq k\leq cl[\Delta]$, define $Q_\Delta^k$ to be the induced subgraph of $Y^k_\Delta$ spanned by $\link(\Sigma)\cap Y_\Delta$.  Define $R^k_\Delta$ to be the induced subgraph of $Z^k_\Delta$ spanned by the vertex set of $\link(\Sigma)\cap Y_\Delta$ and the cone-points $v_{[\Pi]}$ in $Z^k_{\Delta}$ associated to $[\Pi]$ for which $[\Pi]\nest[\Sigma]$.  
\end{defn}

\begin{rem}\label{rem:understanding-Q-R}
By Definition~\ref{defn:coned_y} and Definition~\ref{defn:Q-and-R}, the graph $Q^k_\Delta$ is obtained from $\mathcal C(\Sigma)\cap Y_\Delta$ by electrifying each subgraph of the form $\link(\Pi)\cap\link(\Sigma)\cap Y_\Delta$, where $[\Pi]$ is such that $[\Delta]\nest[\Pi]$ and $cl[\Pi]>k$.  In particular, $Q^{k-1}_\Delta$ is obtained from $Q^{k}_\Delta$ by electrifying the subgraphs $\link(\Pi)\cap\link(\Sigma)\cap Y_\Delta$ with $[\Delta]\nest[\Pi]$ and $cl[\Pi]=k$.

Meanwhile, the graph $R^{k-1}_\Delta$ is obtained from $Q^{k}_\Delta$ by coning off $\link(\Pi)\cap\link(\Sigma)\cap Y_\Delta$ for those $[\Pi]$ with $[\Delta]\nest[\Pi]$, $cl[\Pi]=k$, \emph{and} the additional property that $\link(\Pi)\subset \link(\Sigma)$.  So, if $v_{[\Pi]}$ is a cone-point of $Z^{k-1}_{\Delta}$ that lies in $R^{k-1}_{\Delta}$, then the entire star of $v_{[\Pi]}$ in $Z^{k-1}_{\Delta}$ lies in $R^{k-1}_{\Delta}$, which is a property we will need to apply Lemma \ref{lem:discrete_hyperbolic}.

But, if $cl[\Pi]=k$ and $[\Delta]\nest[\Pi]$ but $[\Pi]\not\nest[\Sigma]$, then $\link(\Pi)\cap\link(\Sigma)\cap Y_\Delta$ is electrified in $Q^{k-1}_\Delta$ but not coned off in $R^{k-1}_\Delta$. Hence, it is not immediately clear that $Q^{k-1}_\Delta$ is quasi-isometric to $R^{k-1}_\Delta$, but we will show later that it is.
\end{rem}

We relate the above spaces in the following commutative diagram, which we call the \emph{main diagram} (where one can roughly think of the bottom row as obtained by taking intersections with $\mathcal C(\Sigma)$):

\begin{center}
    $
    \begin{diagram}
    \node{ Y^k_\Delta}\arrow{e,t}{\ell_k}\node{ Z^{k-1}_\Delta}\arrow{e,t}{\eta}\node{ Y^{k-1}_\Delta}\\
    \node{ Q^k_\Delta}\arrow{n,l}{i_k}\arrow{e,b}{\ell_k}\node{ R^{k-1}_\Delta}\arrow{n,l}{j_{k-1}}\arrow{e,b}{\chi}\node{ Q^{k-1}_\Delta}\arrow{n,l}{i_{k-1}}
    \end{diagram}
    $
\end{center}

The maps are as follows:
\begin{itemize}
    \item $i_k$ and $i_{k-1}$ are inclusions; by Definition~\ref{defn:Q-and-R}, these are graph homomorphisms and hence $1$--Lipschitz;
    \item $j_{k-1}$ is an inclusion, and again, by Definition~\ref{defn:Q-and-R}, it is $1$--Lipschitz;
    \item $\ell_k$ are inclusions of induced subgraphs, as explained in Remark~\ref{rem:understanding-Q-R} and again $1$--Lipschitz; 
    \item $\eta$ is the inclusion on the image of $\ell_k$, and sends each $v_{[\Pi]}$ to an arbitrary vertex of $\link(\Pi)\cap Y_\Delta$ --- this is again $1$--Lipschitz since $v_{[\Pi]}$ is adjacent to every vertex of $\link(\Pi)\cap Y_\Delta$, and this subgraph is electrified in $Q^{k-1}_\Delta$ --- in fact, as noted in Section~\ref{subsec:inductive_hyperbolic}, $\eta$ is a uniform quasi-isometry;
    \item $\chi=i_{k-1}^{-1}\circ\eta\circ j_{k-1}$.  Note that if $v_{[\Pi]}\in j_{k-1}(R^{k-1}_\Delta)$, then $\link(\Pi)\cap Y_\Delta\subset\link(\Sigma)\cap Y_\Delta$, so $\eta(v_{[\Pi]})\in\link(\Sigma)$, whence $\chi$ is well-defined, and it is $1$--Lipschitz for the same reason $\eta$ is.
\end{itemize}

The next lemma about the $Q^k_\Delta$ supports the base case and the final step of an induction in the proof of Lemma~\ref{lem:moreover}.

\begin{lemma}\label{lem:Q-base-cases}
The graphs $Q^{k}_\Delta$ have the following properties:
\begin{enumerate}
    \item $Q^k_\Delta$ has uniformly bounded diameter for $k<cl[\Sigma]$.
    \item $Q^{cl[\Delta]}_\Delta=\mathcal C(\Sigma)\cap Y_\Delta$.
\end{enumerate}
\end{lemma}

\begin{proof}
The second assertion holds by definition: when $k=cl[\Delta]$, we start with $\mathcal C(\Sigma)\cap Y_\Delta$ and add no further edges.

The first assertion follows from the definition since the vertex set of $Q^k_\Delta$ is the vertex set of $\link(\Sigma)\cap Y_\Delta$, and $[\Delta]\nest[\Sigma]$, so when $k<cl[\Sigma]$, the definition says that any two vertices are joined by an edge.
\end{proof}

\subsubsection{$Y^k_\Delta$--bounded intersections of links} The next lemma will help us to control subsets $\link(\Pi)\cap\link(\Sigma)\cap Y_\Delta$ that are electrified in $Q^{k-1}_\Delta$, and coned off in $Z^{k-1}_\Delta$, but not coned off in $R^{k-1}_\Delta$.

\begin{lemma}\label{lem:electric-intersection}
Suppose that $0\leq k\leq cl[\Delta]$.  Let $\Pi$ be a non-maximal simplex of $X$ such that $[\Delta]\nest[\Pi]$, and $cl[\Pi]=k$, but $[\Pi]\not\nest[\Sigma]$.  

Then for all $x,y\in\mathcal C(\Pi)\cap\link(\Sigma)\cap Y_\Delta$, there exists a path $\gamma$ in $Y^k_\Delta$, with all vertices in $\mathcal C(\Pi)$, such that
\begin{itemize}
    \item $\gamma$ joins $x$ to $y$;

    \item $\gamma$ has length at most $\delta(\delta+2)+4$.
\end{itemize}
In particular, $\link(\Pi)\cap\link(\Sigma)\cap Y_\Delta$ has bounded diameter in the graph metric on $Y^k_\Delta$.
\end{lemma}

\begin{proof}
Fix a vertex $w_0\in\link(\Delta)\subset\link(\Pi)\cap\link(\Sigma)\cap Y_\Delta$.  Fix $x'\in\link(\Pi)\cap\link(\Sigma)\cap Y_\Delta$.

Suppose that there exists $v\in\link(\Pi)\cap\sat(\Delta)$ and that $x'$ is adjacent in $\mathcal C(\Pi)$ to $v$.  There is a simplex $v\star\Pi$  of $X$.  Now, since $v\in\sat(\Delta)$, we have $\link(\Delta)\subset\link(v)\cap\link(\Pi)$, so $v\star\Pi$ is non-maximal and $[\Delta]\nest[v\star\Pi]\propnest[\Pi]$.  Hence $cl[v\star\Pi]>k$, and it follows that any two vertices of $\link(v\star\Pi)$ are joined by an edge in $Y^k_\Delta$.

There are two possibilities to consider, according to the type of edge $e$ joining $x'$ to $v$:
\begin{itemize}
    \item If $e$ is an edge of $X$, then $x'\in\link(v)\cap\link(\Pi)$, so $x'\in\link(v\star\Pi)$, and thus $Y^k_\Delta$ has an edge $\alpha$ joining $x'$ to $w_0$.

    \item If $e$ is a $W$--edge (and there is no $X$--edge), then since $x',v\in\link(\Pi)$, Definition~\ref{defn:combinatorial_HHS}.\eqref{item:C_0=C} provides maximal simplices $\sigma,\tau$ of $\link(\Pi)$ such that $x'\in\sigma\star\Pi$, and $v\in\tau\star\Pi$, and $\sigma\star\Pi,\tau\star\Pi$ are $W$--adjacent.  Since $\tau\star\Pi$ is maximal, there exists $v'\in\tau\star\Pi-\sat(\Delta)$ (since saturations cannot contain maximal simplices).  Moreover, since $\mathcal C(\Delta)$ has diameter at least $\delta$, Lemma~\ref{lem:Y_delta_containment} implies that $\Pi\subset\sat(\Delta)$, so $v'\in\tau$.  Thus $x'$ is joined to $v'$ by a $W$--edge, and hence by an edge $\alpha_0$ of $Y^k_\Delta$.  And since $v'\in\link(v\star\Pi)$, there is an edge $\alpha_1$ of $Y^k_\Delta$ from $v'$ to $w_0$.  So $x'$ is joined to $w_0$ by the length--$2$ path $\alpha=\alpha_0\alpha_1$ in $ Y^k_\Delta$ whose vertices are in $\link(\Pi)$.  
\end{itemize}

Now let $x,y$ be as in the statement.  Since $[\Pi]\not\nest[\Sigma]$, we have $\sat(\Sigma)\not\subset\sat(\Pi)$.  Hence there exists $u\in\sat(\Sigma)\cap Y_\Pi$, so since $x,y\in\link(\Pi)$, there is a path of length $2$ in $Y_\Pi$ (consisting of $X$--edges) from $x$ to $y$.  Since $\mathcal C(\Pi)\to Y_\Pi$ is a $(\delta,\delta)$--quasi-isometric embedding, we get a path $\gamma'\subset\mathcal C(\Pi)$ from $x$ to $y$ with $|\gamma'|\leq \delta(\delta+2)$.  If $\gamma'\subset Y_\Delta$, we are done, taking $\gamma=\gamma'$.

Hence suppose that $\gamma'$ passes through $\sat(\Delta)$.  Let $v_0,v_1$ be the first and last vertices of $\sat(\Delta)$ along $\gamma'$ (it is possible that $v_0=v_1$).  Hence we have vertices $x,x',v_0,v_1,y',y$ in $\gamma'$, in that order, such that $x,x'$ and $y',y$ are joined by subgeodesics of $\gamma'$ lying in $Y_\Delta$, and $v_0,v_1\in\sat(\Delta)$, and the subgeodesics of $\gamma'$ from $x'$ to $v_0$ and from $v_1$ to $y'$ are edges of $\mathcal C(\Pi)$.  Then the first part of the proof shows that $x'$ is joined to $w_0$ by a path $\alpha_0$ in $Y_\Delta^k$, and $w_0$ is joined to $y'$ by a path $\alpha_1$ in $Y_\Delta^k$, such that $|\alpha_0|,|\alpha_1|\leq 2$, and all vertices in both paths are in $\link(\Pi)$.

Hence $x$ is joined to $y$ by a path $\gamma$ obtained by concatenating the subgeodesic of $\gamma'$ from $x$ to $x'$, the path $\alpha_0\alpha_1$, and the subgeodesic of $\gamma'$ from $y'$ to $y$.  So, $|\gamma|\leq \delta(\delta+2)+4$, as claimed.
\end{proof}

\subsubsection{Proof of Lemma~\ref{lem:moreover}}
To prove that $\mathcal C(\Sigma)\cap Y_\Delta\hookrightarrow Y_\Delta$ is a quasi-isometric embedding, we actually prove the following lemma:

\begin{lemma}\label{lem:inductive-Q-qie}
For all $0\leq k\leq cl[\Delta]$, the map $i_k:Q^k_\Delta\to Y^k_\Delta$ is a uniform quasi-isometric embedding.
\end{lemma}

\begin{proof}
First observe that the lemma holds for $0\leq k<cl[\Sigma]$ because $Q^k_\Delta$ has uniformly bounded diameter in this case, by Lemma~\ref{lem:Q-base-cases}.  So we only have to prove the lemma for $cl[\Sigma]\leq k\leq cl[\Delta]$.

Let $\ell=cl[\Delta]-cl[\Sigma]$.  We will argue by induction on $\ell$, i.e., for all $[\Sigma]$ with $cl[\Delta]-cl[\Sigma]=\ell$ and $[\Delta]\nest[\Sigma]$, we will prove the claim in the lemma for all $k$.

\textbf{Base case $\ell=0$:}  When $\ell=0$, the fact that $[\Delta]\nest[\Sigma]$ implies that $[\Delta]=[\Sigma]$. In this case we have $Y^{cl[\Delta]}_\Delta=Y_\Delta$ by Lemma \ref{lem:Y_delta_hyperbolic} and $Q^{cl[\Delta]}_\Delta=\mathcal C(\Delta)\cap Y_\Delta=\mathcal C(\Delta)$ by Lemma \ref{lem:Q-base-cases}. Hence, by Definition~\ref{defn:combinatorial_HHS}.\eqref{item:chhs_delta}, $i_{cl[\Delta]}$ is a $(\delta,\delta)$--quasi-isometric embedding.  Since we only need consider $k=cl[\Delta]=cl[\Sigma]$ in this case, we are done.

\textbf{Inductive step from $\ell-1$ to $\ell$:}  Let $\ell\geq 1$ and suppose that the lemma holds for any $[\Sigma']$ such that $[\Delta]\nest[\Sigma']$ and $cl[\Delta]-cl[\Sigma']<\ell$.  Suppose that $cl[\Delta]-cl[\Sigma]=\ell$.  

We now argue by induction on $k$ that $i_k$ is a uniform quasi-isometric embedding.  We saw above that this holds in the base cases $k<cl[\Sigma]$.  So, suppose $cl[\Sigma]\leq k\leq cl[\Delta]$ and suppose, by induction, that $i_{k-1}$ is a $(\delta_0,\delta_0)$--quasi-isometric embedding, where $\delta_0$ is a uniform constant.  Recall the main diagram:
\begin{center}
    $
    \begin{diagram}
    \node{ Y^k_\Delta}\arrow{e,t}{\ell_k}\node{ Z^{k-1}_\Delta}\arrow{e,t}{\eta}\node{ Y^{k-1}_\Delta}\\
    \node{ Q^k_\Delta}\arrow{n,l}{i_k}\arrow{e,b}{\ell_k}\node{ R^{k-1}_\Delta}\arrow{n,l}{j_{k-1}}\arrow{e,b}{\chi}\node{ Q^{k-1}_\Delta,}\arrow{n,l}{i_{k-1}}
    \end{diagram}
    $
\end{center}
in which all maps are $1$--Lipschitz, $\eta$ is a uniform quasi-isometry, and $i_{k-1}$ is a uniform quasi-isometric embedding.  We show below, in Claim~\ref{claim:more-bounded-intersections} and the surrounding discussion, that $\chi$ is a uniform quasi-isometry.  Inspecting the main diagram then shows that $j_{k-1}$ is a uniform quasi-isometric embedding.  

From here, we conclude as follows.  By construction, $j_{k-1}(R^{k-1}_\Delta)$ --- which is isomorphic to $R^{k-1}_\Delta$ via $j_{k-1}$ --- is an induced subgraph of $Z^{k-1}_{\Delta}$ that contains the full star of any cone-vertex of $Z^{k-1}_{\Delta}$ whenever it contains the cone-vertex.  Indeed, if $v_{[\Pi]}$ is a cone-vertex in $j_{k-1}(R^{k-1}_\Delta)$, then by Definition~\ref{defn:Q-and-R}, $\link(\Pi)\cap Y_\Delta\subset \link(\Sigma)\cap Y_\Delta$.  It then follows immediately from Lemma~\ref{lem:discrete_hyperbolic} that $i_k$ is a uniform quasi-isometric embedding.  The fact that $j_{k-1}$ is a quasi-isometric embedding is needed in order to satisfy the hypothesis of Lemma~\ref{lem:discrete_hyperbolic}, which requires a quasi-isometric embedding at the level of coned-off graphs.  The hypotheses involving the ``ambient'' graphs $Y^k_\Delta$ and $Z^{k-1}_\Delta$ were already checked during the proof of Lemma~\ref{lem:Y_delta_hyperbolic}.

\textbf{Proving that $\chi$ is a quasi-isometry:}  It remains to show that $\chi$ is a uniform quasi-isometry.  We have already seen that $\chi$ is $1$--Lipschitz, so it remains to construct a uniformly coarsely Lipschitz quasi-inverse $\bar\chi$.

On the subgraph $\mathcal C(\Sigma)\cap Y_\Delta$, which is a subgraph of both $Q^{k-1}_\Delta$ and $R^{k-1}_\Delta$, and which contains all vertices of $Q^{k-1}_\Delta$, we define $\bar\chi$ to be the identity.  To conclude that $\chi$ is a uniform quasi-isometry, it therefore suffices to bound the distance in $R^{k-1}_\Delta$ between any two vertices $x,y\in Q^{k-1}_\Delta$ that are joined by an edge in $Q^{k-1}_\Delta$ that is not an edge in $\mathcal C(\Sigma)\cap Y_\Delta$.

First, if there exists $\Pi$ such that $[\Delta]\nest[\Pi]$ and $cl[\Pi]>k$ and $x,y\in \link(\Pi)$, then $x,y$ are joined by an edge in $Q^k_\Delta$ and hence in $R^{k-1}_\Delta$, as required.

Second, if there exists $\Pi$ such that $[\Delta]\nest[\Pi]\nest[\Sigma]$ and $cl[\Pi]=k$ and $x,y\in \link(\Pi)$, then $R^{k-1}_{\Delta}$ contains a cone-point $v_{[\Pi]}$ connected by edges to both $x$ and $y$, and we are again done.

The remaining case is where there exists $\Pi$ such that $[\Delta]\nest[\Pi]$, and $cl[\Pi]=k$, but $[\Pi]\not\nest[\Sigma]$, and $x,y\in\link(\Pi)\cap\link(\Sigma)\cap Y_\Delta$.  This case is dealt with in the following claim, which is the whole reason we are arguing by induction on $\ell$ and $k$, and not merely on $k$.

\begin{claim}\label{claim:more-bounded-intersections}
For $x,y$ as above, there exists a uniform constant $C$ such that $\dist_{R^{k-1}_\Delta}(x,y)\leq C$.
\end{claim}
\renewcommand{\qedsymbol}{$\blacksquare$}
\begin{proof}[Proof of Claim~\ref{claim:more-bounded-intersections}]
We will produce a path in $Q^k_\Delta$, of uniformly bounded length, joining $x$ to $y$.  Since $Q^k_\Delta$ is a subgraph of $R^{k-1}_\Delta$, this suffices.

First, observe that $k<cl[\Delta]$.  Indeed, if $k=cl[\Delta]$, then $[\Pi]=[\Delta]$ since $cl[\Pi]=k$ and $[\Delta]\nest[\Pi]$.  But since $[\Delta]\nest[\Sigma]$, this implies $[\Pi]\nest[\Sigma]$, contradicting our choice of $x,y$ (this case was handled above).

We next reduce to the case where $y\in \link(\Delta)$ (which is contained in $\link(\Pi)\cap\link(\Sigma)\cap Y_\Delta$ since $[\Delta]\nest[\Sigma],[\Delta]\nest[\Pi]$).  Suppose that $x$ can be joined to some $x'\in\link(\Delta)$ by a path in $Q^k_\Delta$ of length bounded by a uniform constant $C_0$, and that the same is true with $x$ replaced by $y$ and $x'$ replaced by some $y'\in\link(\Delta)$.  Now, since $k<cl[\Delta]$, any two vertices of $\link(\Delta)$ are connected by an edge of $Y^k_\Delta$, and hence $x',y'$ are connected by an edge of $Q^k_\Delta$.  Thus $x,y$ are connected by a path of length at most $2C_0+1$ in $Q^k_\Delta$.  Since this bound is uniform, we are done.  Hence we can and shall assume $y\in\link(\Delta)$.

Next, since $[\Pi]\not\nest[\Sigma]$ and $cl[\Sigma] \leq k=cl[\Pi]$, there is no nesting relation between $[\Pi]$ and $[\Sigma]$.  As in the proof of Lemma~\ref{lem:electric-intersection}, but with $\Pi$ and $\Sigma$ switching roles, this implies that there exists $v\in\sat(\Pi)-\sat(\Sigma)=Y_\Sigma\cap\sat(\Pi)$.  Since $x,y\in\link(\Pi)$, we obtain a path of length $2$ that joins $x$ to $y$ (passing through $v$) and lies in $Y_\Sigma$.  Hence, by Definition~\ref{defn:combinatorial_HHS}.\eqref{item:chhs_delta}, there exists a geodesic $\alpha$ of $\mathcal C(\Sigma)$ that joins $x$ to $y$ and satisfies $|\alpha|\leq \delta(\delta+2)$.

If $\alpha\subset Y_\Delta$, then since $\mathcal C(\Sigma)\cap Y_\Delta$ is contained in $Q^k_\Delta$, and $\alpha$ has uniformly bounded length, we are done.  Therefore, we can assume that $\alpha$ contains a vertex $u$ of $\sat(\Delta)$.  Since $y\in\link(\Delta)$, there is an edge $f$ of $X$ from $u$ to $y$, so since $\alpha$ is a geodesic, this edge must be the terminal edge of $\alpha$.  Hence we have $\alpha=\alpha_0\cdot e\cdot f$, where $\alpha_0$ is a path in $\mathcal C(\Sigma)\cap Y_\Delta\subset Q^k_\Delta$ of length at most $\delta(\delta+2)$ and $e$ is an edge of $\mathcal C(\Sigma)\subset Q^k_\Delta$.

Let $x'$ be the vertex of $\alpha_0$ immediately preceding $u$, so that $e$ joins $x'$ to $u$.  

Note that $X$ contains a simplex $u\star\Sigma$, since $u\in\mathcal C(\Sigma)$.  Since $u\in\sat(\Delta)$ and $[\Delta]\nest[\Sigma]$, we have $\link(\Delta)\subset\link(u)\cap\link(\Sigma)=\link(u\star\Sigma)$, i.e., $[\Delta]\nest[u\star\Sigma]\propnest[\Sigma]$.  So, $cl[\Delta]\geq cl[u\star\Sigma]>cl[\Sigma]$, and $y\in\link(\Delta)$, so $y\in\link(u\star\Sigma)$.  

Now we analyze two cases, according to whether or not $e$ is an $X$--edge.

If $e$ is an $X$--edge, then $x'\in\link(u\star\Sigma)$. For later purposes, note that the following argument gives in particular that whenever $x'\in \link(u\star\Sigma)$ there exists a path of uniformly bounded length in $Q^k_\Delta$ from $x'$ to $y$.

Notice that $cl[\Delta]-cl[u\star\Sigma]<cl[\Delta]-cl[\Sigma]=\ell$, so by our induction hypothesis (from the ``outer'' induction, on $\ell$), the inclusion $\mathcal Q\hookrightarrow Y^k_\Delta$ is a $(\delta_1,\delta_1)$--quasi-isometric embedding, where $\delta_1$ is uniform and $\mathcal Q$ is the induced subgraph of $Y^k_\Delta$ spanned by the vertex set of $\mathcal C(u\star\Sigma)\cap Y^k_\Delta$.  Note that $\mathcal Q\subset Q^k_\Delta$ since $[u\star\Sigma]\nest[\Sigma]$.

Now, since $x,y\in\link(\Sigma)\cap\link(\Pi)\cap Y_\Delta$ and $[\Pi]\not\nest[\Delta]$, Lemma~\ref{lem:electric-intersection} provides a path $\beta$ in $Y^k_\Delta$, of length at most $\delta(\delta+2)+4$, that joins $x$ to $y$.  By traveling backward along $\alpha_0$, and then along $\beta$, we get a uniformly bounded path in $Y^k_\Delta$ from $x'$ to $y$.  Thus $\mathcal Q$ contains a path of uniformly bounded length (in terms of $\delta_1$) from $x'$ to $y$.  So, since $\mathcal Q$ is a subgraph of $Q^k_\Delta$, we have a path in $Q^k_\Delta$ that has uniformly bounded length and joins $x'$ to $y$.  Prepending $\alpha_0$ to this path completes the proof in this case.

Finally, suppose that $e$ is a $W$--edge, and that $x'$ is not joined to $u$ by an $X$--edge.  By Definition~\ref{defn:combinatorial_HHS}.\eqref{item:C_0=C}, there exist maximal simplices $\sigma,\tau$ of $\link(\Sigma)$ such that $x'\in\sigma\star\Sigma$ and $u\in\tau\star\Sigma$, and every vertex of $\sigma$ is joined to every vertex of $\tau$ by an edge of $\mathcal C(\Sigma)$.  Since saturations cannot contain maximal simplices, there exists $u'\in(\tau\star\Sigma)^{(0)}\cap Y_\Delta$.  Note that $u'$ cannot lie in $\sat(\Sigma)$, because $\sat(\Sigma)\subset\sat(\Delta)$ by Lemma~\ref{lem:Y_delta_containment}.  Hence $u'\in\tau\subset\link(\Sigma)$.

It follows that the concatenation of $\alpha_0$ with the edge of $\mathcal C(\Sigma)$ from $x'$ to $u'$ is a uniformly bounded path in $\mathcal C(\Sigma)\cap Y_\Delta$, so it suffices to find a uniformly bounded path from $u'$ to $y$.  But $u'$ is joined to $u$ by an $X$--edge, so that, as observed above, the previous argument applies with $u'$ replacing $x'$ to bound the distance in $Q^k_\Delta$ from $u'$ to $y$, and hence from $x$ to $y$.  This proves the claim.
\end{proof}
\renewcommand{\qedsymbol}{$\Box$}

\textbf{Conclusion:}  We have shown that $i_k$ is a quasi-isometric embedding.  The constants increased by a uniform amount (arising from Lemma~\ref{lem:discrete_hyperbolic}) in the inductive step from $k-1$ to $k$, but there are boundedly many such steps.  The constants also depended in a uniform way on the constants obtained at stage $\ell-1$, but there are boundedly many possible values of $\ell$ (namely the complexity).  Hence the constants are uniform.
\end{proof}

Finally:

\begin{proof}[Proof of Lemma~\ref{lem:moreover}]
Apply Lemma~\ref{lem:inductive-Q-qie} to the given $[\Sigma]$ in the case where $k=cl[\Delta]$; this shows that $i_k$ is a uniform quasi-isometric embedding.  But in view of Lemma~\ref{lem:Q-base-cases}, $i_k$ is the inclusion $\mathcal C(\Sigma)\cap Y_\Delta\hookrightarrow Y_\Delta$, as required.
\end{proof}

The proof of Proposition \ref{prop:Y_delta_hyperbolic} is now complete: Lemma~\ref{lem:Y_delta_hyperbolic} established uniform hyperbolicity of $Y_\Delta$, and Lemma~\ref{lem:moreover} established that $\mathcal C(\Sigma)\cap Y_\Delta\to Y_\Delta$ is a uniform quasi-isometric embedding when $[\Delta]\nest[\Sigma]$.\qed

We remind the reader that both conclusions of the proposition require that $\diam(\mathcal C(\Delta))\geq \delta$, but as we shall see below, we only need the proposition in that case.  Indeed, morally, we need hyperbolicity in many places for there to be a uniformly coarsely Lipschitz projection $Y_\Delta\to\mathcal C(\Delta)$, but we can define this arbitrarily when $\mathcal C(\Delta)$ is bounded.  We need the QI-embedding statement once, in the proof of Proposition~\ref{prop:hieromorphism}, which obviously holds when the simplex $\Sigma\star\Delta$ in the statement of the proposition has uniformly bounded augmented link.

\section{Combinatorial HHS structures on links, and hieromorphisms}\label{sec:hieromorphisms}
We continue to let $(X,W,\delta,n)$ be a combinatorial HHS.  In this section, we discuss induced combinatorial 
HHS structures on links of simplices. The goal is to show, roughly speaking, that the inclusion of a link in $X$ is compatible with the combinatorial HHS structures. This will allow us to perform inductive arguments, since the complexity of the HHS structure on a link is strictly smaller than that of $X$.

The main results of this section are Propositions \ref{prop:compatibility_1}, describing the combinatorial HHS structure on links, and \ref{prop:hieromorphism}, stating exactly what we mean by the inclusion being compatible with the combinatorial HHS structures.

Fix a non-maximal simplex $\Delta$ of $X$. Let $\mathcal C_0(\Delta)$ be the subgraph of $Y_\Delta$ 
obtained from $\link(\Delta)$ by adding an edge from $v$ to $w$ if the following holds: there exist 
simplices $x,y$ of $X$ such that $\Delta\star x,\Delta\star y$ are maximal simplices of $X$ that respectively 
contain $v$ and $w$ and are connected in $W$.  (Note that when $[\Delta]=[\emptyset]$, we have $\link(\Delta)=X$, and $\mathcal C_0(\Delta)$ is obtained from $X$ by joining $v,w$ whenever they are respectively contained in maximal simplices $x,y$ of $X$ with $x,y$ adjacent in $W$.  To see this, note that $[\Delta]=[\emptyset]$ implies that $\Delta=\emptyset$, since $\link(\emptyset)=X$.  Hence the maximal simplices of the form $\Delta\star x$ of $X$ are exactly the maximal simplices $x$ of $X$.)

The advantage of this definition for us is that 
it has better inheritance properties than $\mathcal C(\Delta)$. However, and this is the only use of Definition 
\ref{defn:combinatorial_HHS}.\ref{item:C_0=C}, it defines the same object:

\begin{lemma}\label{lem:C_0=C}
 Given a combinatorial HHS $(X,W)$ and a non-maximal simplex $\Delta$ of $X$, we have $\mathcal C_0(\Delta)=\mathcal C(\Delta)$.
\end{lemma}

\begin{proof}
 Clearly, $\mathcal C_0(\Delta)$ is a subgraph of $\mathcal C(\Delta)$, and the vertex sets are the same. So, it suffices to show that any edge $e$ of $\mathcal C(\Delta)$ is also an edge of $\mathcal C_0(\Delta)$. If $e$ is an edge of $X$, then this is clear. Otherwise, $e$ is an edge coming from $W$, meaning that the endpoints of $e$ are contained in $W$-adjacent maximal simplices. Definition \ref{defn:combinatorial_HHS}.\ref{item:C_0=C} provides $W$-adjacent maximal simplices containing the endpoints of the form required to yield an edge of $\mathcal C_0(\Delta)$.
\end{proof}

In the rest of this section we will most often use the notation $\mathcal C_0(\Delta)$, except in the proof of 
Proposition \ref{prop:hieromorphism}, where we need to use results from Section \ref{sec:y_hyp}.

\begin{defn}[Induced $\link(\Delta)$--graph]\label{defn:induced_link_graph}
Let $\Delta$ be a non-maximal simplex of $X$.  The \emph{induced 
$\link(\Delta)$--graph} is the graph $W^{\Delta}$ defined as follows.

The vertex set of $W^\Delta$ is the set of maximal simplices $\Sigma$ of $\link(\Delta)$. Notice that any such simplex $\Sigma$ has the property that $\Delta\star\Sigma$ is a maximal simplex of $X$, and hence a vertex of $W$.

Two simplices $\Sigma,\Sigma'$ of $\link(\Delta)$ are connected in $W^\Delta$ if and only 
if $\Delta\star\Sigma$ and $\Delta\star\Sigma'$ are connected in $W$.  By definition, $W^\Delta$ is a $\link(\Delta)$--graph in the sense of 
Definition~\ref{defn:X_graph}. 
\end{defn}

\begin{rem}[The induced map $W^\Delta\to W$]\label{rem:link_graph_homomorphism}
There is a natural injective graph morphism $W^\Delta\hookrightarrow W$ whose image is an induced 
subgraph.  Indeed, the maximal simplex $\Sigma$ of $\link(\Delta)$ is sent to the maximal 
simplex $\Delta\star\Sigma$ of $X$.  Hence we can view $W^\Delta$ as a subgraph of $W$.  
\end{rem}

\begin{defn}[The induced map $\link(\Delta)^{+W^\Delta}\to X^{+W}$]\label{defn:induced_map_link_graph_x_graph}
There is an induced simplicial map 
$\iota_\Delta:\duaug{\link(\Delta)}{W^\Delta}\to\duaug{X}{W}$ defined as follows.

For each $v\in\link(\Delta)^{(0)}$, let $\iota_\Delta(v)=v$, where $v$ is viewed as a vertex of $X$.  If 
$v,w\in\link(\Delta)^{(0)}$ are joined by an edge $e$ of $\link(\Delta)$, then $e$ is an edge of $X$ (and 
hence $X^{+W}$) and we 
let $\iota_\Delta(e)=e$.  If $\sigma,\sigma'$ are maximal simplices of $\link(\Delta)$ corresponding to vertices of 
$W^\Delta$, then by definition, $\sigma\star\Delta$ and $\sigma'\star\Delta$ are maximal simplices of $X$ 
corresponding to vertices of $W$.  If $\sigma,\sigma'$ are adjacent in $W^\Delta$, then $\sigma\star\Delta$ and 
$\sigma'\star\Delta$ are adjacent in $W$, by Definition~\ref{defn:induced_link_graph}.  Thus $\duaug{\link(\Delta)}{W^\Delta}$ contains the 1--skeleton $D$ 
of the join $\sigma\star\sigma'$ and $\duaug{X}{W}$ contains the $1$--skeleton $D'$ of the join 
$(\sigma\star\sigma')\star\Delta$, and we define $\iota_\Delta$ on $D$ to be the inclusion $D\to D'$.  
\end{defn}

\begin{rem}\label{rem:iota_Delta_C}
 Chasing the definitions, one sees that $\iota_\Delta(\link(\Delta)^{+W^\Delta})=\mathcal C_0(\Delta)$. In 
fact, the two graphs have the same vertex set, and in either case two vertices are connected by an edge if and 
only if they are contained in maximal simplices of the form $\Delta\star\Pi$ that are connected in $W$.
\end{rem}

\subsubsection*{Notation} We now revisit the notions associated to a combinatorial HHS in the context of $\link(\Delta)$. For each non-maximal simplex $\Sigma$ of $\link(\Delta)$, let $\link^\Delta(\Sigma)$ be the link of 
$\Sigma$ in $\link(\Delta)$, let $\sim_\Delta$ be the resulting equivalence relation (as in 
Definition~\ref{defn:simplex_equivalence}) on the simplices of $\link(\Delta)$, and let $\sat^\Delta(\Sigma)$ 
be the union of the vertex sets of simplices of $\link(\Delta)$ in the $\sim_\Delta$--class of $\Sigma$.  Let 
$Y^\Delta_\Sigma$ be the induced subgraph of $\duaug{\link(\Delta)}{W^\Delta}$ induced by the vertices of 
$\link(\Delta)-\sat^{\Delta}(\Sigma)$.

Let $\mathcal C_0^\Delta([\Sigma])$ be the subgraph of  $\duaug{\link(\Delta)}{W^{\Delta}}$ obtained from  
$\link^\Delta(\Sigma)$ by adding edges connecting $v,w$ if $v,w$ are vertices of maximal simplices of $\link(\Delta)$ of the form $\Sigma\star x,\Sigma\star y$ that are connected in $W^\Delta$.

Let $\mathfrak S_\Delta$ be the set of $\sim_\Delta$--classes of non-maximal simplices in $\link(\Delta)$.  
The 
relation $\sim_\Delta$ and the accompanying graphs $\link^\Delta(\Sigma),\sat^\Delta(\Sigma)$ (where $\Sigma$ 
is a 
simplex of $\link(\Delta)$) allow us to define relations $\nest,\orth,\transverse$ on $\mathfrak S_\Delta$ as 
in 
Definition~\ref{defn:nest_orth}.

The $\orth$ relation on $\mathfrak S_\Delta$ warrants some extra comment.  Let $\Sigma,\Sigma'$ be non-maximal 
simplices of 
$\link(\Delta)$.  Following Definition~\ref{defn:nest_orth}, we declare $\Sigma,\Sigma'$ to be 
orthogonal in 
$\mathfrak S_\Delta$ if and only if $\link^\Delta(\Sigma')\subseteq\link^\Delta(\link^\Delta(\Sigma))$.  Note 
that the 
left side is $\link(\Sigma')\cap\link(\Delta)$, and the right side is 
$\link(\link(\Sigma)\cap\link(\Delta))\cap\link(\Delta)$.  

\begin{lemma}[Induced map on index sets]\label{lem:map_on_simplices}
Let $\Sigma$ be a non-maximal (possibly empty) simplex of $\link(\Delta)$.  The assignment $\Sigma\mapsto 
\Sigma\star\Delta$ induces an injective map $\iota^*:\mathfrak S_\Delta\to\mathfrak S$ that preserves the 
$\nest,\orth,\transverse$ relations.
\end{lemma}

\begin{proof}
We first show that $\Sigma\sim_\Delta\Sigma'$ if and only if $\Sigma\star\Delta\sim \Sigma'\star\Delta$, 
showing that $\iota^*$ is well-defined and injective.  

Recall that $\link(\Sigma\star\Delta)=\link(\Sigma)\cap\link(\Delta)=\link^\Delta(\Sigma)$, and similarly for 
$\Sigma'$, so that we have  $\link(\Sigma\star\Delta)= \link(\Sigma'\star\Delta)$ if and only if 
$\link^\Delta(\Sigma)=\link^\Delta(\Sigma')$, as required.

\textbf{Preservation of nesting and non-nesting:} For the same reason as above, we have $\link(\Sigma\star\Delta)\subseteq \link(\Sigma'\star\Delta)$ if and only if $\link^\Delta(\Sigma)\subseteq\link^\Delta(\Sigma')$, that is, $[\Sigma]\nest[\Sigma']$ in $\mathfrak S^\Delta$ if and only if $[\Sigma\star\Delta]\nest[\Sigma'\star\Delta]$ in $\mathfrak S$.

\textbf{Preservation of orthogonality and non-orthogonality:}  Suppose that $\Sigma,\Sigma'$ are simplices of 
$\link(\Delta)$ satisfying 
$[\Sigma]\orth[\Sigma']$.  Then by definition, we have 
$$\link(\Sigma')\cap\link(\Delta)\subseteq\link(\link(\Sigma)\cap\link(\Delta))\cap\link(\Delta).$$  In other 
words, 
$\link(\Sigma'\star\Delta)\subseteq \link(\link(\Sigma\star\Delta))\cap\link(\Delta)$.  Hence 
$[\Sigma'\star\Delta]\orth[\Sigma\star\Delta]$ in $\mathfrak S$, as required.

Conversely, if $[\Sigma'\star\Delta]\orth[\Sigma\star\Delta]$ in $\mathfrak S$, then by definition 
$\link(\Sigma')\cap\link(\Delta)\subset\link(\link(\Sigma)\cap\link(\Delta))$, but since the left-hand side is contained in $\link(\Delta)$ we also have $$\link(\Sigma')\cap\link(\Delta)\subset\link(\link(\Sigma)\cap\link(\Delta))\cap\link(\Delta),$$ that is
$[\Sigma]\orth[\Sigma']$ 
in $\mathfrak S_\Delta$.
\end{proof}

\begin{cor}\label{cor:induced_finite_complexity}
If $\Delta$ is non-empty, then the complexity of $\link(\Delta)$ is strictly less than the complexity $n$ of $(X,W)$.
\end{cor}

\begin{proof}
Any $\nest$--chain in $\mathfrak S_{\Delta}$ maps via $\iota^*$ to a $\nest$--chain in $\mathfrak S$ of the 
same 
length, by Lemma~\ref{lem:map_on_simplices}.  Let $\Sigma$ be the $\nest$--maximal element of such a chain in 
$\mathfrak S_{\Delta}$.  Then $[\Sigma\star\Delta]$ is properly nested in $[\emptyset]$, which lies in 
$\mathfrak 
S-\iota^*(\mathfrak S_\Delta)$, whence the complexity of $\mathfrak S_\Delta$ is strictly less than the complexity $n$ of 
$\mathfrak S$. 
\end{proof}

We are aiming to show that $(\link(\Delta),W^\Delta,\delta',m)$ inherits a combinatorial HHS structure from 
$(X,W,\delta,n)$, where $\delta'$ is uniform, and $m<n$.  The preceding statements established what we need at 
the level of index sets.  Now we study the spaces $\mathcal C_0^\Delta(\Sigma)$.

\begin{lemma}\label{lem:curve_graph_iota_v2}
Let $\Sigma$ be a simplex of $\link(\Delta)$.  Then:
\begin{enumerate}
     \item\label{item:same_C} $\iota_\Delta$ restricts to an isomorphism of graphs from $\mathcal C_0^{\Delta}([\Sigma])$ to  $\mathcal C_0([\Sigma\star\Delta])$.
\item \label{item:containment_Y} Suppose that $diam(\mathcal C_0(\Sigma\star\Delta))\geq\delta$. Then 
$\iota_\Delta(Y_\Sigma^\Delta)= Y_{\Sigma\star\Delta}\cap \mathcal C_0(\Delta)$, and $\iota_\Delta$ restricts 
to a graph isomorphism from $Y_\Sigma^\Delta$ to $Y_{\Sigma\star\Delta}\cap \mathcal C_0(\Delta)$.

\end{enumerate}
\end{lemma}

\begin{proof}

 \textbf{Proof of item \eqref{item:same_C}.} Let $v$ be a vertex of $\mathcal C^{\Delta}([\Sigma])$, that is, 
$v$ is a vertex of $\link(\Delta)\cap \link(\Sigma)$. But then $v$ is a vertex of $\link(\Delta\star\Sigma)$, 
so that $\iota_\Delta(v)=v$ is a vertex of $\mathcal C([\Sigma\star\Delta])$. Notice also that, for similar 
reasons, every vertex of $\link(\Delta\star\Sigma)$ is in the image of $\iota_\Delta$, and hence 
$\iota_\Delta$ restricts to a bijection between the vertex sets of $\mathcal C^{\Delta}([\Sigma])$ and 
$\mathcal C([\Sigma\star\Delta])$. We now have to show that $v,w$ are connected by an edge in $\mathcal 
C_0^{\Delta}([\Sigma]))$ if and only if they are connected by an edge in $\mathcal C_0([\Sigma\star\Delta])$. 
Clearly, such $v,w$ are connected by an edge of $\link(\Delta)$ if and only if they are connected by an edge of 
$X$. Hence, suppose that $v,w$ are connected by one of the additional edges coming from $W^\Delta$, meaning 
that $v,w$ are vertices of maximal simplices $\Sigma\star x,\Sigma\star y$ of $\link(\Delta)$ that are 
connected in $W^\Delta$. But then, by definition of $W^\Delta$, $\Delta\star\Sigma\star x$ and 
$\Delta\star\Sigma\star y$ are connected in $W$, showing that $v,w$ are connected in $\mathcal 
C_0([\Sigma\star\Delta])$. The other implication is similar.

\textbf{Proof of item \eqref{item:containment_Y}.} Let us first check the statement at the level of vertex 
sets.

$\subseteq$: Let $x\in Y_\Sigma^\Delta$ be a vertex, which is to say 
that 
$x\in\link(\Delta)^{(0)}-\sat^\Delta(\Sigma)$.  Since $x\in \link(\Delta)$, we have 
$x\in\duaug{\link(\Delta)}{W^\Delta}$. 
 Clearly $x\in\mathcal C_0(\Delta)$, so we claim that $x\in Y_{\Sigma\star\Delta}$, i.e., that $x\in X^{(0)}-\sat(\Sigma\star\Delta)$.  Suppose to the 
contrary that there exists a simplex $\Pi'$ so that $x\in\Pi'$ and $\link(\Pi')=\link(\Sigma)\cap\link(\Delta)$.  

Then we know that $[\Sigma\star\Delta]\nest [\Delta\star x]$, since
$$\link(\Sigma\star\Delta)=\link(\Delta\star x)=\link(\Delta)\cap \link(x)\supseteq \link(\Delta)\cap \link(\Pi').$$
Hence, by Lemma \ref{lem:SC_nesting} there exists a simplex $\Pi''$ in $\link(\Delta\star x)$ with 
$[\Delta\star x\star \Pi'']=[\Sigma\star\Delta]$. In particular, $\link^\Delta(\Pi''\star 
x)=\link^\Delta(\Sigma)$, showing $x\in \sat^\Delta(\Sigma)$, a contradiction.

$\supseteq$: This is similar, but easier. Let $x\in Y_{\Sigma\star\Delta}\cap \mathcal C_0(\Delta)$ be a vertex. We have to show that $x\in\link(\Delta)^{(0)}-\sat^\Delta(\Sigma)$. If we had $x\in \sat^\Delta(\Sigma)$ then there would exist a simplex $\Pi$ of $\link(\Delta)$ so that $x\in \Pi$ and $\link^\Delta(\Pi)=\link^\Delta(\Sigma)$. But then $\link(\Delta\star\Pi)=\link(\Delta\star\Sigma)$, and we would have $x\in\sat(\Sigma\star\Delta)$, a contradiction.

Finally, we are left with arguing that the edge sets are also the same. Observe that two vertices $x,y$ of 
$Y_\Sigma^\Delta$ are connected by an edge if and only if they are connected by an edge in 
$\link(\Delta)^{+W^\Delta}$. By Remark \ref{rem:iota_Delta_C}, we have $\link(\Delta)^{+W^\Delta}=\mathcal 
C_0(\Delta)$. Also, since $Y_{\Sigma\star\Delta}\cap\mathcal C_0(\Delta)$ is an induced subgraph of $\mathcal C_0(\Delta)$, we have that 
vertices in $Y_{\Sigma\star\Delta}\cap \mathcal C_0(\Delta)$ are connected if and only if they are connected in 
$\mathcal C_0(\Delta)$. To sum up, two vertices $x,y$ of $Y_\Sigma^\Delta$ are connected by an edge if and only 
if they are connected by an edge in $\mathcal C_0(\Delta)$, if and only if they are connected in 
$Y_{\Sigma\star\Delta}\cap \mathcal C_0(\Delta)$, as required.
\end{proof}

\begin{prop}[Combinatorial HHS structure on links]\label{prop:compatibility_1}
 For each $\delta$ there exists $\delta''$ so that the following holds. Given a 
combinatorial HHS $(X,W,\delta,n)$ 
and a non-maximal, nonempty simplex $\Delta$ of $X$, there is $m<n$ so that 
$(\link(\Delta),W^\Delta,\delta'',m)$ is a combinatorial HHS. 
\end{prop}

\begin{proof}
The required bound on complexity of $\mathfrak S_\Delta$ comes 
from Corollary~\ref{cor:induced_finite_complexity}.

We will first check Condition \ref{item:C_0=C}, which (as in the proof of Lemma \ref{lem:C_0=C}) guarantees 
that $\mathcal C_0^\Delta(\Sigma)=\mathcal C^\Delta(\Sigma)$ for every simplex $\Sigma$ of $\link(\Delta)$.

Let $v,w\in\link^\Delta(\Sigma)$, and suppose that they are contained in $W^\Delta$-adjacent simplices 
$\Sigma_1,\Sigma_2$ of $\link(\Delta)$. By definition of the edges of $W^\Delta$, we have that 
$\Delta\star\Sigma_1$ and $\Delta\star\Sigma_2$ are $W$-adjacent (and clearly still contain $v,w$ 
respectively). Moreover, $v,w\in\link(\Delta\star\Sigma)$, so that by Condition \ref{item:C_0=C} for $(X,W)$ we 
have that $v,w$ are contained in $W$-adjacent maximal simplices $\Delta\star\Sigma\star 
\Gamma_1$,$\Delta\star\Sigma\star\Gamma_2$. This means that $v,w$ are contained in the $W^\Delta$-adjacent 
maximal simplices $\Sigma\star \Gamma_1, \Sigma\star \Gamma_2$ of $\link(\Delta)$, as required.

Next, we produce $\delta'$ so that each 
$\mathcal C_0^\Delta([\Sigma])$ is $\delta'$--hyperbolic, and 
$\mathcal C_0^\Delta([\Sigma])\hookrightarrow Y_\Sigma^\Delta$ is a 
$(\delta',\delta')$--quasi-isometric embedding.

By Lemma \ref{lem:curve_graph_iota_v2}.\ref{item:same_C}, $\mathcal C^\Delta_0([\Sigma])$ is isometric to $\mathcal C_0([\Sigma\star\Delta])$.  Hence $\mathcal C_0^\Delta([\Sigma])$ is $\delta$--hyperbolic.  

We now verify the quasi-isometric embedding part. It suffices to consider the case where 
$diam(\mathcal 
C^\Delta_0([\Sigma]))\geq\delta$. Hence, in view of Lemma 
\ref{lem:curve_graph_iota_v2}.\ref{item:containment_Y}, we can consider the following commutative diagram, 
where the horizontal arrows are restriction of $\iota_\Delta$ and the vertical arrows are inclusions:
\begin{center}
$
\begin{diagram}
\node{Y_\Sigma^\Delta}\arrow{e}\node{Y_{\Sigma\star\Delta}}\\
\node{\mathcal C_0^\Delta([\Sigma])}\arrow{n,J}\arrow{e,t}{\cong} \node{\mathcal C([\Sigma\star\Delta])}
\arrow{n,J}
\end{diagram}
$
\end{center}
the right vertical arrow is a $(\delta,\delta)$--quasi-isometric embedding 
by 
Definition~\ref{defn:combinatorial_HHS}, and the top horizontal 
arrow is a $1$--Lipschitz map, and the left vertical arrow is also $1$--Lipschitz.  
Since the bottom horizontal arrow is an isometry, there exists $\delta''$, depending 
only on $\delta$, so that the left vertical arrow is a 
$(\delta'',\delta'')$--quasi-isometric embedding.  (For instance, $\delta''=\max\{\delta,1\}$ suffices.)

It remains to verify condition~\eqref{item:chhs_join} from Definition~\ref{defn:combinatorial_HHS}, which is going to easily follow from the same condition for $(X,W)$.  Let 
$\Sigma,\Omega$ be non-maximal simplices of $\link(\Delta)$, and suppose that there exists some non-maximal 
simplex $\Gamma$ of $\link(\Delta)$ whose $\sim_\Delta$--class is nested into those of $\Sigma$ and $\Omega$, and $diam(\mathcal C^\Delta_0([\Gamma]))\geq\delta$.  
Then $[\Gamma\star\Delta]\nest[\Sigma\star\Delta],[\Omega\star\Delta]$ by Lemma \ref{lem:map_on_simplices}, and $diam(\mathcal C_0([\Delta\star\Gamma]))\geq \delta$ by Lemma \ref{lem:curve_graph_iota_v2}.\ref{item:same_C}.   By 
Definition~\ref{defn:combinatorial_HHS}.\eqref{item:chhs_join}, applied to the combinatorial 
HHS $(X,W)$, there exists a non-maximal simplex $\Pi$ of 
$\link(\Sigma\star\Delta)=\link(\Sigma)\cap\link(\Delta)$ such that 
$[\Pi\star(\Sigma\star\Delta)]\nest[\Omega\star\Delta]$.  Hence $\Pi$ is a simplex of $\link(\Delta)$ (so 
$[\Pi]\in\mathfrak S_\Delta$), and $\Pi$ is a simplex of $\link^\Delta(\Sigma)$ (as required by 
Definition~\ref{defn:combinatorial_HHS}.\eqref{item:chhs_join}).  Moreover, by the preceding nesting 
relation, we have $\link^\Delta(\Pi\star\Sigma)\subseteq\link^\Delta(\Omega)$, i.e., $\Pi\star\Sigma$ is nested 
in $\Omega$, with respect to the nesting in $\mathfrak S_\Delta$.

To check the remaining clause of Definition~\ref{defn:combinatorial_HHS}.\eqref{item:chhs_join}, suppose that 
$\Gamma$ is a non-maximal simplex of $\link(\Delta)$ with $\mathcal C_0^\Delta(\Gamma)\cong 
\mathcal C_0(\Gamma\star\Delta)$ having diameter at least $\delta$.  Then 
$[\Gamma\star\Delta]\nest[(\Sigma\star\Delta)\star\Pi]$, by 
Definition~\ref{defn:combinatorial_HHS}.\eqref{item:chhs_join}, applied to the combinatorial HHS 
$(X,W)$.  Hence, as above, we see that $[\Gamma]\nest[\Sigma\star \Pi]$, where equivalence classes 
and nesting are taken in $\mathfrak S_\Delta$.  This completes the proof.
\end{proof}

From Proposition~\ref{prop:compatibility_1} and the first part of Proposition~\ref{prop:Y_delta_hyperbolic}, we 
immediately obtain:

\begin{cor}\label{cor:induced_Y_delta_hyperbolic}
Let $(X,W,\delta,n)$ be a combinatorial HHS.  Then there exists $\delta'$, depending only on $\delta$ and $n$, 
so that the following holds.  Let $\Delta$ be a non-maximal simplex of $X$ and let $\Sigma$ be a non-maximal 
simplex of $\link(\Delta)$.  Then $Y^\Delta_\Sigma$ is $\delta'$--hyperbolic provided $\diam(\mathcal C^\Delta([\Sigma]))\geq\delta'$.
\end{cor}

\subsubsection*{Compatibility of structures} For convenience, we now recall some more notation from Section \ref{sec:setup_statement}, and what it yields in the case of $\link(\Delta)$ and $W^\Delta$. For each non-maximal simplex $\Sigma$ of $\link(\Delta)$, the projection 
$\pi^\Delta_{\Sigma}:W^\Delta\to\mathcal C^\Delta([\Sigma])$ is defined as follows: each vertex $w\in 
W^\Delta$ is a maximal simplex of $\link(\Delta)$, and in particular $w$ is not contained in 
$\sat^\Delta(\Sigma)$, since $\link^\Delta(\Sigma)\neq\emptyset$.  Thus $w\cap Y^\Delta_\Sigma$ is a complete 
graph, and we 
take $\pi^\Delta_{[\Sigma]}(w)$ to be the closest-point projection of $w\cap Y^\Delta_\Sigma$ on 
$\mathcal C^\Delta([\Sigma])$.  When $\Sigma,\Sigma'$ are non-orthogonal simplices of $\link(\Delta)$, 
$\rho^{[\Sigma]}_{[\Sigma']}$ is defined in Definition~\ref{defn:projections}.

Now we are ready to prove the second main proposition of this section:

\begin{prop}[Inclusions induce hieromorphisms]\label{prop:hieromorphism}
Let $X,W,\delta,n,\Delta,\Sigma$ be as in Proposition~\ref{prop:compatibility_1}.  Then the map 
$\Sigma\mapsto\Sigma\star\Delta$ induces an 
isometry $\mathcal C^\Delta([\Sigma])\to\mathcal C([\Sigma\star\Delta])$ so that the 
following diagrams uniformly coarsely commute 
whenever $\Sigma,\Sigma'$ are simplices of $\link(\Delta)$ for which 
$[\Sigma]\transverse[\Sigma']$, $[\Sigma]\propnest[\Sigma']$ or $[\Sigma']\propnest[\Sigma]$ in $\mathfrak 
S_\Delta$:\\
\begin{tabular}{ccc}
     $
\begin{diagram}
 \node{W^\Delta}\arrow{e,J}\arrow{s,l}{\pi^\Delta_{[\Sigma]}}\node{W}\arrow{s,r}{\pi_{[\Sigma\star\Delta]}}\\
 \node{\mathcal C^\Delta([\Sigma])}\arrow{e,t}{\cong}\node{\mathcal C([\Sigma\star\Delta])}
\end{diagram}
$
&\,\,\,\,\,\,\,\,\,\,\,\,\,\,\,\,\,\,\,\,\,
&
$
\begin{diagram}
\node{\mathcal C^\Delta([\Sigma])}\arrow{e,t}{\cong}\arrow{s,l}{\rho^{[\Sigma]}_{[\Sigma']}}\node{\mathcal C([
\Sigma\star\Delta])
} \arrow{s,r}{\rho^{[\Sigma\star\Delta]}_{[\Sigma'\star\Delta]}}\\
\node{\mathcal C^\Delta([\Sigma'])}\arrow{e,b}{\cong}\node{\mathcal C([\Sigma'\star\Delta])}
\end{diagram}
$
\end{tabular}
\end{prop}

\begin{proof}
First of all, let us recall that in view of Lemma \ref{lem:C_0=C} there is no difference between $\mathcal 
C_0$-spaces and $\mathcal C$-spaces. This allows us to use results from this section and from Section 
\ref{sec:y_hyp}.  Second, we may assume that $\diam(\mathcal C(\Sigma'\star\Delta))\geq \delta$, for otherwise uniform coarse commutation of both diagrams is immediate, since the maps are all uniformly coarsely constant in this case.

Hence, by Proposition~\ref{prop:Y_delta_hyperbolic}, there exists $\delta'=\delta'(\delta)$ so that $Y_{\Sigma\star\Delta}$ 
is 
$\delta'$--hyperbolic.  By (uniformly) enlarging $\delta'$, we have that $Y^\Delta_\Sigma$ is 
$\delta'$--hyperbolic for each non-maximal simplex $\Sigma$ of $\link(\Delta)$, by 
Corollary~\ref{cor:induced_Y_delta_hyperbolic} and Proposition~\ref{prop:Y_delta_hyperbolic}.

In view of Lemma \ref{lem:curve_graph_iota_v2}.\ref{item:containment_Y}, whenever $diam(\mathcal C_0^\Delta(\Sigma))\geq \delta$, the ``moreover'' clause of Proposition~\ref{prop:Y_delta_hyperbolic} says that there exists a uniform $\epsilon$ such 
that $Y^\Delta_\Sigma\hookrightarrow 
Y_{\Sigma\star\Delta}$ is an $(\epsilon,\epsilon)$--quasi-isometric embedding.

Now, as in the proof of Proposition~\ref{prop:compatibility_1}, provided that $diam(\mathcal C_0^\Delta(\Sigma))\geq \delta$, we have a 
commutative diagram
\begin{center}
$
\begin{diagram}
\node{Y_\Sigma^\Delta}\arrow{e,J}\node{Y_{\Sigma\star\Delta}}\\
\node{C^\Delta([\Sigma])}\arrow{n,J}\arrow{e,t}{\cong}\node{C([\Sigma\star\Delta])}
\arrow{n,J}
\end{diagram}
$
\end{center}
where all arrows are $(\epsilon',\epsilon')$--quasi-isometric embeddings, where 
$\epsilon'$ is uniform.  Let $p:Y_{\Sigma\star\Delta}\to \mathcal 
C([\Sigma\star\Delta])$ be coarse closest-point projection, and let 
$p^\Delta:Y^\Delta_\Sigma\to\mathcal C^\Delta([\Sigma])$ be the coarse 
closest-point projection.  

For convenience, in the following argument we identify $Y_\Sigma^\Delta$ and $C^\Delta([\Sigma])$ with their images under $\iota_\Delta$. Let $a\in Y^\Delta_\Sigma$.  Note that $p(a)\in\mathcal C(\Sigma\star\Delta)=\mathcal C^\Delta_\Sigma\subset Y^\Delta_\Sigma$.

Let $\gamma$ be a $Y^\Delta_\Sigma$--geodesic from $a$ to $p^\Delta(a)$.  So, 
$\gamma$ is a uniform quasi-geodesic of $Y_{\Sigma\star\Delta}$, and hence contains a point $b$ at uniformly 
bounded distance $C$ (in $Y^\Delta_{\Sigma}$) from $p(a)$, since $\gamma$ fellow-travels with a 
geodesic of $Y_{\Sigma\star\Delta}$ that starts at $a$ and ends at $p^\Delta(a)\in\mathcal C(\Sigma\star\Delta)$. 
 Now, travelling along $\gamma$ from $a$ to $b$ and then travelling along a $Y^\Delta_\Sigma$--geodesic from $b$ to $p(a)$ gives a 
path from $a$ to $p(a)$ of length at most $\dist_{Y^\Delta_\Sigma}(a,b)+C\geq|\gamma|-1$.  Hence 
$\dist_{Y^\Delta_{\Sigma}}(b,p^\Delta(a))\leq C+1$, so $\dist_{Y^\Delta_\Sigma}(p^\Delta(a),p(a))\leq 2C+1$. 
Hence we have a uniform bound on $\dist_{Y_{\Sigma\star\Delta}}(p(a),p^\Delta(a))$, and thus on 
$\dist_{\mathcal C(\Sigma\star\Delta)}(p(a),p^\Delta(a))$.  Thus $p^\Delta$ uniformly coarsely 
coincides with the restriction of $p$ to $Y^\Delta_\Sigma$.

Since
$\pi_{[\Sigma\star\Delta]},\pi^\Delta_{[\Sigma]}$, and 
$\rho^{[\Sigma']}_{[\Sigma]}$ and $\rho^{[\Sigma'\star\Delta]}_{[\Sigma\star\Delta]}$ were all defined  in terms of $p$ and $p_\Delta$ (and since when the target of the projection under consideration has bounded diameter the statement holds automatically), checking that the diagrams in the statement coarsely commute is now straightforward; we give the details below.

Let $w$ be a maximal simplex of $\link(\Delta)$, i.e., a vertex of $W^\Delta$.  So, the image of $w$ under $W^\Delta\hookrightarrow W$ is the maximal simplex $w\star\Delta$ of $W$.  Let $\Sigma$ be a non-maximal simplex of $\link(\Delta)$.

Observe that $(w\star\Delta)\cap Y_{\Sigma\star\Delta}=w\cap Y^\Delta_\Sigma$.  Indeed, if $a$ is a vertex of $(w\star\Delta)\cap Y_{\Sigma\star\Delta}$, and each vertex of $\Delta$ belongs to $\Sigma\star\Delta$ and hence to $\sat(\Sigma\star\Delta)$, we must have $a\in w$.  Since $w\subset\link(\Delta)$, we have by Lemma~\ref{lem:curve_graph_iota_v2}.\eqref{item:containment_Y} that $a\in Y^\Delta_\Sigma$.  Conversely, if $a\in w\cap Y^\Delta_\Sigma$, then by the same lemma, $a\in \link(\Delta)\cap Y_{\Sigma\star\Delta}$, and obviously $a\in w\star\Delta$.

Now, by definition, $$\pi^\Delta_{[\Sigma]}(w) = p^\Delta(w\cap Y^\Delta_\Sigma)$$ and $$\pi_{[\Sigma\star\Delta]}(w\star\Delta)=p((w\star\Delta)\cap Y_{\Sigma\star\Delta})=p(w\cap Y^\Delta_\Sigma),$$
so by uniform closeness of $p^\Delta$ and $p$, established earlier in the proof, the sets $\pi^\Delta_{[\Sigma]}(w)$ and $\pi_{[\Sigma\star\Delta]}(w\star\Delta)$ uniformly coarsely coincide.  Hence the diagram in the statement uniformly coarsely commutes.

Suppose that $\Sigma,\Sigma'$ are non-maximal simplices of $\link(\Delta)$ such that $[\Sigma]\transverse[\Sigma']$ (resp. $[\Sigma]\propnest[\Sigma']$).  Then $[\Sigma\star\Delta]\transverse[\Sigma'\star\Delta]$ (resp. $[\Sigma\star\Delta]\propnest[\Sigma'\star\Delta]$).

By definition, $\rho^{[\Sigma]}_{[\Sigma']} = \pi^\Delta_{[\Sigma']}(Y^\Delta_{\Sigma'}\cap \sat^\Delta(\Sigma))$ and $\rho^{[\Sigma\star \Delta]}_{[\Sigma'\star\Delta]}=\pi_{[\Sigma'\star\Delta]}(Y_{\Sigma'\star\Delta}\cap \sat(\Sigma\star\Delta))$.

If $a$ is a vertex of $\sat^\Delta(\Sigma)$, then $a\in \Pi$, for some simplex $\Pi$ of $\link(\Delta)$ with $\link(\Pi)\cap\link(\Delta)=\link(\Sigma)\cap\link(\Delta)$.  So, $\Pi\star\Delta$, and hence $a$, is contained in $\sat(\Sigma\star\Delta)$.  Moreover, if $a\in Y^\Delta_{\Sigma'}$, then $a\in Y_{\Sigma'\star\Delta}\cap \link(\Delta)$.  Thus, $Y^\Delta_{\Sigma'}\cap \sat^\Delta(\Sigma)$ is contained in $Y_{\Sigma'\star\Delta}\cap \sat(\Sigma\star\Delta)$.  Combining this with the above discussion relating $p^\Delta$ to $p$, we have that $\rho^{[\Sigma]}_{[\Sigma']}$ (computed in the combinatorial HHS structure on $\link(\Delta)$) is coarsely contained in, and hence coarsely coincides with, $\rho^{[\Sigma\star \Delta]}_{[\Sigma'\star\Delta]}$, as required.

Now suppose that $[\Sigma]\propnest[\Sigma']$.  The map $\rho^{[\Sigma']}_{[\Sigma]}$ is defined to be the 
restriction of $p^\Delta:Y^\Delta_{\Sigma}\to \mathcal C(\Sigma\star\Delta)$ to $\mathcal 
C(\Sigma'\star\Delta)\cap Y^\Delta_{\Sigma}$, and $\emptyset$ otherwise.  Meanwhile, 
$\rho^{[\Sigma'\star\Delta]}_{[\Sigma\star\Delta]}$ is the restriction of $p$ to $\mathcal 
C(\Sigma'\star\Delta)\cap Y_{\Sigma\star\Delta}$ and $\emptyset$ otherwise.  So, these maps coarsely agree on 
$\mathcal C(\Sigma'\star\Delta)\cap Y^\Delta_{\Sigma}$.  Since $\mathcal C(\Sigma'\star\Delta)\cap 
\sat^\Delta(\Sigma)$ is contained in $\mathcal C(\Sigma'\star\Delta)\cap \sat(\Sigma\star\Delta)$, we have that 
$\rho^{[\Sigma'\star\Delta]}_{[\Sigma\star\Delta]}(a)=\emptyset$ whenever $\rho^{[\Sigma']}_{[\Sigma]}(a)$ is 
defined and equal to $\emptyset$.  This completes the proof that the required diagrams coarsely commute.
\end{proof}

\section{$W$ is hierarchically hyperbolic}\label{sec:check_axioms}
Now we prove Theorem~\ref{thm:hhs_links}.  Throughout this section, $(X,W,\delta,n)$ is a combinatorial HHS.

\subsection{Strong bounded geodesic image}\label{sec:strong_big}
The following is the key lemma.  

\begin{lemma}[Strong BGI]\label{lem:strong_BGI}
For each $\delta,\delta'$, there exists $C$ such that the following holds.  Let 
$Z$ be a 
$\delta$--hyperbolic  graph and 
let $V$ be a nonempty subgraph of $Z$.  Let $L_V$ be an induced subgraph of $Z$  disjoint from $V$
with the property that $z\in L_V$ implies 
that $z$ is adjacent to each vertex of $V$. Let $Z_V$ be the induced subgraph of 
$Z$ with vertex set $Z^{(0)}-V^{(0)}$, and suppose $Z_V$ is 
$\delta'$--hyperbolic.

Suppose that $L_V$ is $\delta$--hyperbolic and 
$(\delta,\delta)$--quasi-isometrically embedded (hence quasiconvex with 
constant depending only on $\delta,\delta'$) in $Z_V$, and let $\pi:Z_V\to L_V$ 
be the 
coarse closest-point projection.  

Let $x,y\in Z^{(0)}$ and let 
$\gamma$ be a geodesic in $Z$ from $x$ to $y$.  Suppose that $\gamma\cap 
V=\emptyset$.  Then $\dist_{L_V}(\pi(x),\pi(y))\leq 
C$.
\end{lemma}

\begin{proof}
Let $x,y,\gamma$ be as in the statement.  Since $\gamma\cap V=\emptyset$, we have $\gamma\subset Z_V$, and 
$\gamma$ is necessarily a geodesic of $Z_V$.  Consider 
a geodesic quadrilateral formed by geodesics 
$[x,p],[p,q],[q,y]$ and $\gamma$ in $Z_V$, where $p\in\pi(x),q\in\pi(y)$.  Since $L_V$ is quasiconvex in $Z_V$, there exists a constant $K=K(\delta,\delta')$ such that $[p,q]$ lies in the $K$--neighborhood of $L_V$, as measured in $Z_V$ (and hence in $Z$).  

The quadrilateral is $2\delta'$--thin in $Z_V$.  Let $x'\in[x,p]$ and suppose that $x'$ is $2\delta'$--close to some point of $[p,q]$.  Then $x'$ is $(2\delta'+K)$--close to $L_V$, which is possible only if $\dist_{Z_V}(x',p)\leq 2(\delta'+K)$, for otherwise we contradict that $p\in\pi(x)$.  Similarly, the $2\delta'$--neighborhood of $[y,q]$ only intersects $[p,q]$ at points within $2(\delta'+K)$ of $q$.  But $[p,q]$ lies in the union of the $2\delta'$--neighborhoods of the other three sides of the quadrilateral.  So, either $\dist_{Z_V}(p,q)\leq 5(\delta'+K)$, and hence $\dist_{L_V}(\pi(x),\pi(y))\leq C$ for $C$ uniformly bounded, and we are done, or: a non-empty subpath of $[p,q]$, consisting of all but initial and terminal length--$5(K+\delta')$ segments, lies $2\delta'$--close to $\gamma$.  Hence we can assume that there exist $a,b\in\gamma$ such that 
$\dist_{Z_V}(a,\pi(x))\leq 10K\delta'$ and 
$\dist_{Z_V}(b,\pi(y))\leq10K\delta'$.

 Hence, for any $v\in V$, we have $\dist_Z(a,v)\leq 10K\delta'+1$ and 
$\dist_Z(b,v)\leq 10K\delta'+1$.  So, since $\gamma$ is a geodesic of $Z$, the 
subpath of $\gamma$ from $a$ to $b$ has length at most $20K\delta'+2$, i.e., 
$\dist_{Z_V}(a,b)\leq 20K\delta'+2$.  Hence $\dist_{Z_V}(\pi(x),\pi(y))\leq 
40K\delta'+2$.  Since $L_V$ is $(\delta,\delta)$--quasi-isometrically embedded 
in $Z_V$, this gives $\dist_{L_V}(\pi(x),\pi(y))\leq 
40K\delta\delta'+\delta^2+2\delta$, and 
we take the latter value to be $C$.
\end{proof}

Now fix $[\Delta]\in\mathfrak S$.  By assumption, $Z=X^{+W}$ is 
$\delta$--hyperbolic.  The subgraph $V=\sat(\Delta)$ has the property that 
$Z_V=Y_\Delta$.  Moreover, each vertex of $L_V=\mathcal 
C([\Delta])$ is joined by an edge of $Z$ to each vertex in $V$, and $L_V$ is 
$\delta$--hyperbolic and $(\delta,\delta)$--quasi-isometrically embedded in 
$Z_V$, because of Definition~\ref{defn:combinatorial_HHS}.  Finally, Proposition~\ref{prop:Y_delta_hyperbolic} 
ensures that $Z_V$ is 
$\delta'$--hyperbolic.

Hence, by Lemma~\ref{lem:strong_BGI}, we have:

\begin{lemma}\label{lem:super_bgi}
 There exists $C=C(\delta,\delta')$ so that the following holds. Let 
$[\Delta]\in\mathfrak S$ and let $x,y\in Y_\Delta$.  If $d_{\mathcal C 
([\Delta])}(x,y)> C$, then any geodesic $\gamma$ in $X^{+W}$ from $x$ to $y$ 
intersects 
$\sat(\Delta)$.
\end{lemma}

The preceding lemma will be used repeatedly in the next subsection.

\subsection{Proof of Theorem~\ref{thm:hhs_links}: verification of the 
HHS axioms}\label{subsec:verify_axioms_assuming_hyperbolic}

Before proceeding to the proof of Theorem~\ref{thm:hhs_links}, we remind the reader of the roles of the different spaces from Definition~\ref{defn:combinatorial_HHS}:

\begin{rem}[Connectedness, and the roles of the different spaces]\label{rem:connected}
Recall that, given a combinatorial HHS $(X,W)$, Definition~\ref{defn:combinatorial_HHS} does not assert that either $X$ or $W$ is connected.  In the case of $X$, the reader should bear in mind that this is by design: the objects about which we are making geometric claims and/or assumptions are $W$ and the various $Y_\Delta$ and $\mathcal C(\Delta)$, including $X^{+W}$.  

The complex $X$ does not function as a metric object, but rather as a ``database of links'' recording the index set (and nesting, orthogonality, and transversality relations) for an HHS structure on $W$.  We never refer to metric properties of $X$, only combinatorial properties.  

The spaces $\mathcal C(\Delta)$ are connected by assumption: Definition~\ref{defn:combinatorial_HHS}.\eqref{item:chhs_delta} asserts that these are hyperbolic, and in particular they are geodesic metric spaces (not extended metric spaces).

What we will have to prove is connectedness of $W$.  Indeed, our goal is to prove that $W$ --- a graph with the usual graph metric --- is an HHS, and, in particular, by the requirement in Definition~\ref{defn:HHS} that HHS are quasigeodesic metric spaces, this means proving that $W$ is connected.  The proof of Theorem~\ref{thm:hhs_links} is by induction on the complexity of the combinatorial HHS $(X,W)$.  
In the base cases, we will show that $W$ is either a single point, or $W=X^{+W}=\mathcal C(\emptyset)$, which is connected by assumption. 

In the inductive step, we will assume that, for each non-empty, non-maximal simplex $(\link(\Delta),W^\Delta)$, we have an HHS structure $(W^\Delta,\mathfrak S_\Delta)$.  Built into this is the assumption that $W^\Delta$ is connected since, by Definition~\ref{defn:HHS}, a graph whose graph-metric admits an HHS structure is a geodesic metric space, i.e., it is connected.  We will verify, as part of the proof that $(X,\mathfrak S)$ satisfies the ``Uniqueness'' axiom (Definition~\ref{defn:HHS}.\eqref{item:dfs_uniqueness}), that any two vertices of $W$ are at finite distance in the graph metric, which is to say that $W$ is connected.  We will do this by building paths.
\end{rem}

We are now ready for:

\begin{proof}[Proof of Theorem~\ref{thm:hhs_links}]
Fix $(X,W,\delta,n)$ as in the statement. 

For ``Large Links'' and ``Uniqueness'' we will have to proceed by induction on complexity, and assume that 
Theorem \ref{thm:hhs_links} holds for all the links of nonempty non-maximal simplices of the simplicial complex $X$. 

\textbf{Base cases:} First we explain the base cases $n\leq 1$.  First note that $\emptyset$ is always a simplex of $X$.  So the complexity $n=0$ only if $\emptyset$ is a maximal simplex, which means that $X=\emptyset$.  By Definition~\ref{defn:X_graph}, this implies that $W$ is a single point, and $(W,\emptyset)$ is a hierarchically hyperbolic space in this case, as required.

When the complexity $n=1$, the only non-maximal simplex is the empty set, so $X$ is a discrete set of vertices (0-simplices have to be maximal).  So, by Definition~\ref{defn:X_graph}, $W=X^{+W}$.  On the other hand, by the definition of a combinatorial HHS (Definition~\ref{defn:combinatorial_HHS}.\eqref{item:chhs_delta}) $X^{+W}$ is hyperbolic, which is to say that $W$ is hyperbolic.  In particular, $(W,\{[\emptyset]\})$ is an HHS, where the projection from $W$ to $\mathcal C(\emptyset)=X^{+W}$ is the identity.
 
\textbf{Inductive step:}  We refer the reader to Definition~\ref{defn:HHS}.  The proof will consist of verifying that $(W,\mathfrak S)$ satisfies each requirement in that definition.  

\emph{Inductive hypothesis:}  By Proposition~\ref{prop:compatibility_1}, there is a constant $\delta''=\delta''(\delta,n)$ such that $(\link(\Delta),W^\Delta,\delta'',m)$ is a combinatorial HHS whenever $\Delta$ is a nonempty non-maximal simplex of $X$, where $m<n$ is the complexity of $\link(\Delta)$.  Our inductive hypothesis is therefore that for each such $\Delta$, the pair $(W^\Delta,\mathfrak S_\Delta)$ is an HHS, with the constants in Definition~\ref{defn:HHS} depending only on $\delta''$ and $m$ (and hence bounded in terms of $\delta$ and $n$).  Moreover, we assume that the projections $\pi_\bullet$ and $\rho^\bullet$ are as given in Definition~\ref{defn:projections}, with $\link(\Delta)$ playing the role of $X$, and $W^\Delta$ playing the role of $W$, in that definition.  We assume that the relations are as in Definition~\ref{defn:nest_orth}, with $\link(\Delta)$ playing the role of $X$.

\emph{Verifying the HHS axioms:}  We now proceed to verify the axioms from Definition~\ref{defn:HHS}.  

First, the underlying space, $W$, is a  geodesic extended metric space, since it is a graph.  Definition~\ref{defn:combinatorial_HHS} requires $W$ to be a quasigeodesic space, which in our situation means that we have to show that $W$ is connected.  This is done below, when we verify the uniqueness axiom (Definition~\ref{defn:HHS}.\eqref{item:dfs_uniqueness}).  We will not require connectedness of $W$ for the earlier parts of the proof.  (The simplest argument that we are aware of for connectedness requires less precision, but it otherwise follows the same outline.)

Second, let $\mathfrak S$ be the associated index set, with relations described in Definition~\ref{defn:nest_orth}; for each $[\Delta]\in\mathfrak S$, the associated $\delta$--hyperbolic space is $\mathcal C([\Delta])$, and the various projections are as described in Definition~\ref{defn:projections}. 

We now verify the enumerated axioms from Definition~\ref{defn:HHS}:

\textbf{Projections:}  The projections $\pi_{[\Delta]}$ are coarsely 
Lipschitz coarse maps because that is true of coarse closest point 
projections to quasiconvex subsets of hyperbolic spaces, if $\mathcal C([\Delta])$ has diameter at least $\delta$ (hyperbolicity of 
$Y_{\Delta}$ holds by Proposition~\ref{prop:Y_delta_hyperbolic} in this case, and $\mathcal C([\Delta])$ is 
quasi-isometrically embedded in $Y_{\Delta}$ by assumption).  If $\diam(\mathcal C([\Delta]))\leq\delta$, then the claim is obvious. 

Let $v$ be a vertex of $\mathcal C(\Delta)$.  Then $v$ is contained in some maximal simplex $w$ of $X$, so 
$\pi_{[\Delta]}(w)$ contains $v$.  So, $\mathcal C(\Delta)=\cup_w\pi_{[\Delta]}(w)$, and in 
particular $\pi_{[\Delta]}$ has quasiconvex image.

\textbf{Nesting:} The relation $\nest$ is defined in Definition 
\ref{defn:nest_orth}, and it is clearly a partial order.  The maximal element 
is the equivalence class of $\emptyset$. The bounded sets and coarse maps $\rho^\bullet_\bullet$ are 
defined in Definition \ref{defn:projections}. 
If $[\Delta']\propnest [\Delta]$, then by definition $\link(\Delta')\subsetneq 
\link (\Delta)\subseteq Y_{\Delta}$. Notice that 
some vertex $v$ of $\Delta'$ is contained in $Y_{\Delta}$. Indeed, if we had 
$\Delta'\subseteq \sat(\Delta)$, then we would 
have $\link(\Delta)=\link(\sat(\Delta))\subseteq \link(\Delta')$ (see Remark 
\ref{rmk:lk_sat}), contradicting $\link(\Delta')\subsetneq \link 
(\Delta)$. Since $v\in \sat(\Delta')$, $\rho^{[\Delta']}_{[\Delta]}$ is 
non-empty. Moreover, $\sat(\Delta')$ has diameter at 
most $2$ in $Y_{\Delta}$ because any vertex of $\sat(\Delta')$ is connected to 
any vertex of $\link(\Delta')\subseteq 
\link(\Delta)\subseteq Y_\Delta$ (notice that $\link(\Delta')$ is non-empty 
because, by definition, $\Delta'$ is 
non-maximal). Hence $\rho^{[\Delta']}_{[\Delta]}\subset \mathcal C([\Delta])$ is 
uniformly bounded.

\textbf{Orthogonality:}  We defined $\orth$ in Definition 
\ref{defn:projections}.  It is symmetric because if $\link(\Delta')\subseteq 
\link(\link(\Delta))$ then $\link(\link(\Delta'))\supseteq 
\link(\link(\link(\Delta)))$, and then:

\begin{claim}\label{claim:triple_link}
 $\link(\Delta)=\link(\link(\link(\Delta)))$.
\end{claim}
\renewcommand{\qedsymbol}{$\blacksquare$}
\begin{proof}[Proof of Claim~\ref{claim:triple_link}]
 For any subcomplex $\Sigma$, chasing the definitions we see that 
$\Sigma\subseteq \link(\link(\Sigma))$.  For 
$\Sigma=\link(\Delta)$, we get the inclusion ``$\subseteq$''. For 
$\Sigma=\Delta$, and applying $\link$ to both sides, we get 
the inclusion ``$\supseteq$''.
\end{proof}
\renewcommand{\qedsymbol}{$\Box$}
\noindent Anti-reflexivity of $\orth$ follows from the fact that 
$\link(\Delta)$ is always disjoint from $\link(\link(\Delta))$,  and 
$\link(\Delta)$ is non-empty for any non-maximal simplex (recall that we are 
excluding maximal simplices).

Next, suppose that $[\Delta]\nest [\Delta']$ and $[\Delta']\orth [\Delta'']$. 
Then by the definitions,  $\link(\Delta)\subseteq 
\link(\Delta')\subseteq \link(\link(\Delta''))$, so $[\Delta]\orth [\Delta'']$. 
Also, if $[\Delta]\orth [\Delta']$, then 
$\link(\Delta)$ (which is non-empty) is contained in $\link(\link(\Delta'))$, 
which is disjoint from $\link(\Delta')$, so 
that $[\Delta]\not\nest [\Delta']$.

Now fix $[\Delta]$ and $[\Delta_1]\nest [\Delta]$, and suppose that  
$A=\{[\Delta_2]: [\Delta_2]\nest[\Delta], 
[\Delta_2]\orth [\Delta_1]\}$ is non-empty. We need to find $[\Delta']$ with 
$[\Delta']\propnest [\Delta]$ and 
$[\Delta_2]\nest[\Delta']$ for all $[\Delta_2]\in A$.

By definition, $[\Delta_2]\in A$ has $\link(\Delta_2)\subseteq \link(\Delta)$ 
and  $\link(\Delta_2)\subseteq 
\link(\link(\Delta_1))$. Consider $B=\bigcup_{[\Delta_2]\in A}\link(\Delta_2)$; 
we are looking for a simplex strictly 
containing $\Delta$ whose link contains $B$. We have 
$\link(B)=\bigcap_{[\Delta_2]\in A}\link(\link(\Delta_2))$, and also 
$\link(\Delta_1)\subseteq \link(\link(\Delta_2))$ for each $[\Delta_2]\in A$, so 
that any vertex $v$ in $\link(\Delta_1)$ 
(which exists) has:
\begin{itemize}
 \item $v\in \link(\Delta)$, because $\link(\Delta_1)\subseteq \link(\Delta)$, 
since $[\Delta_1]\nest [\Delta]$. In particular, the simplex 
$\Delta'=v\star\Delta$ is well-defined.
 \item $v\in \link(B)$, so that $B\subseteq \link(v)$. So, for any 
$[\Delta_2]\in A$, $\link(\Delta_2)\subseteq \link(\Delta')=\link(v)\cap 
\link(\Delta)$ (since $\link(\Delta_2)\subseteq \link(\Delta)$), i.e., 
$[\Delta_2]\nest [\Delta']$.
\end{itemize}

Since $\Delta\subsetneq \Delta'$, we have $[\Delta']\propnest [\Delta]$  (by 
Remark 
\ref{rem:set_theory}). Notice that the nesting is indeed proper 
because $v\in \link(\Delta)-\link(\Delta')$.

\textbf{Transversality and consistency:} Recall 
$\rho^{[\Delta']}_{[\Delta]}$ from Definition \ref{defn:projections}.  
If $[\Delta]\transverse 
[\Delta']$ then:
\begin{itemize}
\item $\sat(\Delta')\cap Y_\Delta\neq \emptyset$,  for otherwise we would have 
$\sat(\Delta')\subseteq \sat(\Delta)$, and 
hence $\link(\Delta)=\link(\sat(\Delta))\subseteq\link(\sat(\Delta'))= 
\link(\Delta')$ and $[\Delta]\nest[\Delta']$ (we used 
Remark \ref{rmk:lk_sat}).
\item $\link(\Delta')\cap Y_\Delta\neq \emptyset$,  for otherwise we would have 
$\link(\Delta')^{(0)}\subseteq \sat(\Delta)$, 
and $[\Delta]\orth[\Delta']$ since $\sat(\Delta)\subseteq \link(\link(\Delta))$ 
by Remark \ref{rmk:lk_sat}. So, 
$\sat(\Delta')\cap Y_\Delta$ has diameter at most $2$ in $Y_\Delta$, and 
hence $\rho^{[\Delta']}_{[\Delta]}$ is uniformly bounded.
\end{itemize}

Now, suppose that $d_{[\Delta]}(w,\rho^{[\Delta']}_{[\Delta]})\geq C$,  for some 
vertex $w$ of $W$ (corresponding to a maximal simplex $\Sigma_w$ of $X$), where $C$ is as in Lemma 
\ref{lem:super_bgi}. Then there is a geodesic $\gamma$ in $X^{+W}$ from 
$\sat(\Delta')-\sat(\Delta)$ to $w'$ intersecting 
$\sat(\Delta)$, where $w'\in 
\Sigma_w-\sat(\Delta)$. The 
minimal subgeodesic of $\gamma$ from $w'$ to $\sat(\Delta)$ does not intersect 
$\sat(\Delta')$. By Lemma \ref{lem:super_bgi} 
we must have $d_{[\Delta']}(w,\rho^{[\Delta]}_{[\Delta']})< C$.

\textbf{Consistency for nesting:}  Now suppose that 
$[\Delta']\propnest[\Delta]$.  Let $w\in W$.

Choose a vertex $w'\in X$ as follows. Suppose 
$w\in W^{(0)}$ corresponds to a maximal simplex $\Sigma_w$.  If 
every vertex of $\Sigma_w$ is in $\sat(\Delta)\cup\sat(\Delta')$, then since 
$\link(\Delta')\subset\link(\Delta)$, $X$ would contain 
$\Sigma_w\star\link(\Delta')$, contradicting maximality.  So $w'$ can be chosen 
in $\Sigma_w-(\sat(\Delta)\cup\sat(\Delta'))$.  So, $w'\in Y_\Delta$ and $w'\in 
Y_{\Delta'}$.

Let $p_\Delta:Y_\Delta\to\mathcal C([\Delta])$ and 
$p_{\Delta'}:Y_{\Delta'}\to\mathcal C([\Delta'])$ be closest-point projections.  

If $\mathcal C([\Delta'])$ has diameter at most $\delta$, then we are done, so assume that $\mathcal C([\Delta'])$ has diameter 
more than $\delta$.  Then by Lemma \ref{lem:SC_nesting}, we have 
$[\Delta']=[\Pi\star\Delta]$ for some simplex $\Pi$ of $\link(\Delta)$.

Let $\alpha$ be a geodesic in the $\delta'$--hyperbolic space $Y_\Delta$ joining $w'$ to $p_\Delta(w')\in\link(\Delta)$.

First suppose that $\alpha\cap\sat(\Delta')=\emptyset$. In this case we use the following lemma:

\begin{lemma}\label{Y_sats}
Let $\Sigma$ be a simplex of $\link(\Delta)$.  Then 
$Y_{\Sigma\star\Delta}\subset Y_{\Sigma}\cap Y_{\Delta}$.  
\end{lemma}

\begin{proof}
Suppose that 
$x\in\sat(\Delta)$, i.e., there exists a non-maximal simplex $\Pi$ of $X$ such 
that $\link(\Pi)=\link(\Delta)$ and $x\in \Pi$.  Since $\Sigma$ is a 
simplex 
of $\link(\Delta)$, there is a simplex $\Sigma\star\Pi$ in $X$, and 
$\link(\Sigma\star\Pi)=
\link(\Sigma)\cap\link(\Pi)=
\link(\Sigma)\cap\link(\Delta)=\link(\Sigma\star\Delta)$.  So $x\in 
\Pi\star\Sigma\sim_X\Delta\star\Sigma$, i.e., $x\in\sat(\Sigma\star\Delta)$, as 
required.  Hence $Y_{\Sigma\star\Delta}\subset Y_\Delta$. Similarly, 
$Y_{\Sigma\star\Delta}\subset Y_{\Sigma}$.  
\end{proof}

By Lemma~\ref{Y_sats}, $Y_{\Delta'}\subset Y_\Delta$.  Since $\alpha$ is a geodesic of $Y_\Delta$ and is 
entirely contained in $Y_{\Delta'}$, we have that $\alpha$ is also a geodesic of $Y_{\Delta'}$.  Now, 
$Y_{\Delta'}$ is $\delta'$--hyperbolic, and is obtained from $Y_{\Delta}$ by 
deleting $\sat(\Delta')-\sat(\Delta)$ (and any edges with at least one endpoint in $\sat(\Delta')$).  Now, 
$\mathcal C([\Delta'])$ 
is a subgraph of $Y_{\Delta'}$, each of whose vertices is adjacent in $Y_\Delta$ to each vertex of 
$\sat(\Delta')-\sat(\Delta)$.  Moreover, $\mathcal C([\Delta'])$ is $\delta$--hyperbolic and 
$(\delta,\delta)$--quasi-isometrically embedded in $Y_{\Delta'}$.  Hence, by Lemma~\ref{lem:strong_BGI}, $w'$ and $p(w')$ 
have $C$--close $p_{\Delta'}$--images on $\mathcal C([\Delta'])$, which is to say that 
$\diam(\pi_{[\Delta']}(w)\cup\rho^{[\Delta]}_{[\Delta']}(\pi_{[\Delta]}(w)))$ is uniformly bounded.

Next, suppose that $\alpha$ passes through some $v\in\sat(\Delta')-\sat(\Delta)$.  Since $\alpha$ is a geodesic from $w'$ to $p_\Delta(w')$, we have from the definition of $p_\Delta$ that 
$$\dist_{Y_\Delta}(w',p_\Delta(w'))=\dist_{Y_\Delta}(w',v)+\dist_{Y_\Delta}(v,p_\Delta(w'))\leq \dist_{Y_\Delta}(w,\link(\Delta))+1.$$
Since $v\in\sat(\Delta')$, for any $\ell\in\link(\Delta')\subset\link(\Delta)$, we have that $v$ is adjacent to $\ell$.  Fixing such an $\ell$ (which exists because $\Delta'$ is non-maximal), we have $\dist_{Y_\Delta}(w',\ell)\leq \dist_{Y_\Delta}(w',v)+1$. Hence 
$$\dist_{Y_\Delta}(w',\link(\Delta))\leq \dist_{Y_\Delta}(w',\ell)\leq \dist_{Y_\Delta}(w',v)+1\leq \dist_{Y_\Delta}(w',\link(\Delta))-\dist_{Y_\Delta}(v,p_\Delta(w'))+2.$$ This implies that $\dist_{Y_\Delta}(v,p_\Delta(w'))\leq 2$, which provides a uniform bound, in terms of $\delta$, on $\dist_{Y_\Delta}(p_\Delta(w'),p_\Delta(v))$.  This in turn gives a uniform bound on $\dist_{[\Delta]}(p_\Delta(w'),p_\Delta(v))$, in view of Definition~\ref{defn:combinatorial_HHS}.\eqref{item:chhs_delta}.  Finally, $p_\Delta(w')\in\pi_{[\Delta]}(w)$ and $v\in\sat(\Delta')\cap Y_\Delta$, so by Definition~\ref{defn:projections}, we have shown that $\dist_{[\Delta]}(\pi_{[\Delta]}(w),\rho^{[\Delta']}_{[\Delta]})$ is uniformly bounded, as required.

Finally, we need to check the following.  Suppose that $[\Delta']\propnest[\Delta]$, and $[\Delta]\propnest[\Sigma]$ or 
$[\Delta]\transverse[\Sigma]$, and  $[\Sigma]\not\perp[\Delta']$.  We claim that 
$\dist_{[\Sigma]}(\rho^{[\Delta]}_{[\Sigma]},\rho^{[\Delta']}_{[\Sigma]})$ is uniformly bounded, as required 
by item~\ref{item:dfs_transversal} of the definition of an HHS (Definition~\ref{defn:HHS}).  By definition, $\rho^{[\Delta]}_{[\Sigma]}=p_\Sigma(\sat(\Delta)\cap 
Y_{\Sigma})$ and $\rho^{[\Delta']}_{[\Sigma]}=p_\Sigma(\sat(\Delta')\cap Y_{\Sigma})$.  Since 
$[\Delta']\transverse[\Sigma]$ 
or $[\Delta']\propnest[\Sigma]$, there exists $v\in\link(\Delta')\cap Y_\Sigma$.  Now, $v$ is adjacent in $Y_\Sigma$ to each 
vertex in $\sat(\Delta')\cap Y_\Sigma$, and to each vertex in $\sat(\Delta)\cap Y_\Sigma$.  So, 
$\dist_{Y_\Sigma}(\sat(\Delta),\sat(\Delta'))\leq 2$, so, since $p_\Sigma$ is uniformly coarsely Lipschitz, 
$\dist_{[\Sigma]}(\rho^{[\Delta]}_{[\Sigma]},\rho^{[\Delta']}_{[\Sigma]})$ is uniformly bounded.

\textbf{Finite complexity:} This follows from 
Definition~\ref{defn:combinatorial_HHS}.(\ref{item:chhs_flag}).

\textbf{Bounded geodesic image:} Let $[\Delta]\propnest [\Delta']$. By definition of nesting, $\mathcal C ([\Delta])$ is (properly) contained in 
$\mathcal C ([\Delta'])$.  Let $E$ be a constant to be determined, and suppose 
that $\gamma$ is a geodesic in $\mathcal C ([\Delta'])$ that is disjoint from the 
$E$--neighborhood of $\mathcal C ([\Delta])$.  If $\gamma$ contains a vertex 
$v$ of $\sat(\Delta)$, then $\gamma$ is joined by an edge of $X$, and hence of 
$X^{+W}$, to some vertex $w\in\link(\Delta)\subset\link(\Delta')$.  Hence 
$\gamma$ passes through the $1$--neighborhood in $\mathcal C([\Delta'])$ of 
$\link(\Delta)$; by choosing $E>1$, this is impossible.  Hence $\gamma$ can be 
regarded as a geodesic in $Y_\Delta$ which is far from the quasiconvex subset 
$\mathcal C([\Delta])$ of the $\delta'$--hyperbolic space $Y_\Delta$.  Hence, 
provided $E$ is chosen sufficiently large in terms of $\delta'$ and $\delta$, 
the projections of the endpoints of $\gamma$ 
to $\mathcal C([\Delta])$ are $E$--close, as 
required. 

\textbf{Large links:} Let $[\Delta]\in\mathfrak S$.  Let $x,y\in W$. We need to produce a constant $E$, depending on $\delta$ and $n$ but independent of $[\Delta]$, such that there exist $[\Sigma_1],\ldots,[\Sigma_N]\propnest[\Delta]$ with the property that any $[\Sigma']\propnest[\Delta]$ with $\dist_{[\Sigma']}(x,y)>E$ satisfies $[\Sigma']\nest[\Sigma_i]$ for some $i$.  Moreover, we need to bound $N$ by a uniform linear function of $\dist_{[\Delta]}(x,y)$, and also bound $\dist_{[\Delta]}(x,\rho^{[\Sigma_i]}_{\Delta})$ above by $N$.  

First, suppose that $\Delta\neq\emptyset$.  Then by our induction hypothesis, $(W^\Delta,\mathfrak S^\Delta)$ is an HHS, and the projections for $W^\Delta$ 
coarsely coincide with those for $(X,W)$, as stated in Proposition \ref{prop:hieromorphism}.

The coordinates 
$(\pi_Y(x))_{Y\in\mathfrak S^\Delta}$ are consistent (with uniform constants), as we checked above (``Transversality and consistency'' and ``Consistency for nesting'').  Hence, using our induction hypothesis, we can apply the realization theorem for hierarchically hyperbolic spaces (Theorem~\ref{thm:realization}) to obtain a vertex $x'\in W^\Delta$ whose projections 
coarsely coincide with $\pi_Y(x)$ for all $Y\in\mathfrak S^\Delta$.  We can similarly construct $y'\in W^\Delta$ 
starting from the coordinates of $y$.  

In view of Lemma \ref{lem:SC_nesting}, every $[\Sigma]$ so that 
$[\Sigma]\nest [\Delta]$ and $\mathcal C(\Sigma)$ has diameter at least $\delta$ is in the image of the map 
$\iota^*$ from Lemma \ref{lem:map_on_simplices}.   Hence, it is easily seen that Large Links for $[\Delta],x,y$ 
follows from Large Links for $[\emptyset],x',y'$ in $(W^\Delta,\mathfrak S^\Delta)$ (up to increasing the 
constants).

Now we handle the case where $\Delta=\emptyset$.  Fix a geodesic $\alpha$ in 
$X^{+W}$ from 
$x'$ to $y'$, where $x',y'$ are vertices of $x,y$. For every non-maximal simplex $[\Sigma]$, by Lemma~\ref{lem:super_bgi} we have that either $\dist_{[\Sigma]}(x,y)$ is uniformly bounded by some $E$, or 
$\alpha$ contains a vertex $v$ in $\sat(\Sigma)$. Notice that $v\in\sat(\Sigma)$ implies $[\Sigma]\nest [v]$. Hence it suffices to let $\{\Sigma_i\}$ be the collection of all vertices of $\alpha$.

\textbf{Partial realization:}  Let $[\Delta_1],\ldots,[\Delta_k]\in\mathfrak S$ 
be pairwise-orthogonal, i.e., $\link(\Delta_i)\subset\link(\link(\Delta_j))$ for 
all $i\neq j$.  For each $i$, let $p_i\in\link(\Delta_i)$.  Since 
$p_i\in\link(\link(\Delta_j))$, we have that $p_1,\ldots,p_k$ are 
pairwise-adjacent vertices of $X$.  Hence there exists a maximal simplex $\Pi$ 
of $X$ containing $p_1,\ldots,p_k$.  Since $\Pi$ is a maximal simplex, it 
corresponds to a vertex $w$ of $W$.  For each $i$, we have 
$\Pi\cap\mathcal C([\Delta_i])\ni p_i$, so $p_i\in\pi_{[\Delta_i]}(w)$.

Next, suppose that $[\Sigma]\in\mathfrak S$ satisfies 
$[\Delta_i]\transverse[\Sigma]$ or $[\Delta_i]\propnest[\Sigma]$ for some $i$.  
Let $p:Y_\Sigma\to\mathcal C([\Sigma])$ be the coarse closest-point projection, so 
that $\rho^{[\Delta_i]}_{[\Sigma]}$ is by definition $p(\sat(\Delta_i)\cap 
Y_{\Sigma})$.  Now, if $[\Delta_i]\propnest[\Sigma]$, then 
$\link(\Delta_i)\subsetneq\link(\Sigma)$, so $p_i\in\link(\Sigma)$, whence 
$p_i\in Y_{\Sigma}$.  Hence $p(p_i)\subset\pi_{[\Sigma]}(w)$.  On the other hand, 
there exists $v\in\sat(\Delta_i)$ such that $v\not\in\sat(\Sigma)$, so 
$p(v)\subset\rho^{[\Delta_i]}_{[\Sigma]}$.  Since $p_i,v$ are adjacent in 
$Y_\Sigma$, we have that $p(p_i)$, and hence $\pi_{[\Sigma]}(w)$, lies uniformly 
close to $\rho^{[\Delta_i]}_{[\Sigma]}$.  When $[\Delta_i]\transverse[\Sigma]$, 
we again have some $v\in\sat(\Delta_i)\cap Y_\Sigma$, so $v$ is necessarily 
adjacent to $p_i$.  Hence $v$ lies at distance at most $2$ in $Y_\Sigma$ from 
every vertex of $\Pi\cap Y_\Sigma$.  Since $\Pi\cap Y_\Sigma$ contains at least 
one vertex, $\dist_{[\Sigma]}(p(v),p(\Pi\cap Y_\Sigma))$ is uniformly bounded, 
i.e., $\dist_{[\Sigma]}(\rho^{[\Delta_i]}_{[\Sigma]},\pi_{[\Sigma]}(w))$ is 
uniformly bounded, as required.

\textbf{Uniqueness axiom and connectedness of $W$:} In this part of the proof, we again use our induction assumption: links of nonempty non-maximal simplices carry a natural combinatorial HHS structure of strictly lower complexity, as in Proposition 
\ref{prop:compatibility_1}, which is furthermore compatible with the structure for $X$ as explained in 
Proposition \ref{prop:hieromorphism}.  These combinatorial HHS structures have their associated constants (from Definition~\ref{defn:combinatorial_HHS}) bounded in terms of $\delta$ and $n$, and in particular, independently of the link in question. 

Moreover, we are assuming by induction on complexity that 
Theorem \ref{thm:hhs_links} holds for all links of nonempty non-maximal simplices of $X$, i.e., for each nonempty non-maximal $[\Delta]$, we have an HHS structure $(W^\Delta,\mathfrak S_\Delta)$ where all of the constants from Definition~\ref{defn:HHS} depend only on the above combinatorial HHS constants, and are hence bounded uniformly in terms of $\delta,n$.  In particular, the function $\theta_u$ from Definition~\ref{defn:HHS}.\eqref{item:dfs_uniqueness}, can be taken to be the same for all HHS structures $(W^\Delta,\mathfrak S_\Delta)$ as $\Delta$ varies over nonempty non-maximal simplices of $X$.

Let $x,y$ be vertices in $W$, i.e., maximal simplices in $X$.  To prove connectedness amounts to proving that $\dist_W(x,y)$ is finite.  To prove uniqueness requires a strict strengthening of this, namely that $\dist_W(x,y)$ is bounded above by a fixed function of $\sup_{[\Delta]\in\mathfrak S}\dist_{[\Delta]}(x,y)$.  To do this, we will construct a path in $W$ from $x$ to $y$, and moreover, bound the length of this path in terms of $\sup_{[\Delta]\in\mathfrak S}\dist_{[\Delta]}(x,y)$.

First:

\begin{claim}\label{claim:finite-relevant}
There exists a uniform constant $E$ such that the set of $[\Delta]\in\mathfrak S$ such that $\dist_{[\Delta]}(x,y)>E$ is finite.  Hence $$\kappa=\kappa(x,y)=\sup_{[\Delta]\in\mathfrak S}\dist_{[\Delta]}(x,y)<\infty.$$
\end{claim}
\renewcommand{\qedsymbol}{$\blacksquare$}
\begin{proof}[Proof of Claim~\ref{claim:finite-relevant}]
By the large link axiom, which we have just verified above, there exists a constant $E$ such that the following holds.  Let $\alpha$ be a geodesic in $X^{+W}$ from $x'$ to $y'$, where $x'\in x,y'\in y$.  Then any $[\Delta]\in\mathfrak S-\{[\emptyset]\}$ for which  $\dist_{[\Delta]}(x,y)>E$ satisfies $[\Delta]\nest[v]$, where $v$ is one of the $\dist_{[\emptyset]}(x',y')+1$ vertices of $\alpha$.

By our induction hypothesis, for each such $v$, the pair $(W^v,\mathfrak S_v)$ is an HHS (with uniform constants and uniqueness functions), so by applying Theorem~\ref{thm:realization} (realization theorem) exactly as in the verification of the large link axiom (``Large links''), we obtain $x_v,y_v\in W^v$ whose projections to each $\mathcal C([\Delta]),[\Delta]\nest[v]$, uniformly coarsely coincide with $\pi_{[\Delta]}(x)$ and $\pi_{[\Delta]}(y)$ respectively.  Now, up to uniformly enlarging $E$, it follows, by, for example, an application of the distance formula (Theorem \ref{thm:distance_formula}) in the HHS $(W^v,\mathfrak S_v)$, that there are finitely many $[\Delta]\nest[v]$ with $\dist_{[\Delta]}(x,y)>E$.  Since there are finitely many choices for $v$, this proves the first assertion, from which the second assertion follows immediately.
\end{proof}
\renewcommand{\qedsymbol}{$\Box$}

Hence, to simultaneously prove connectedness of $W$ and the uniqueness axiom, it suffices to prove that, for 
any $\kappa\geq 0$ and any vertices $x,y\in W$ with $d_{[\Delta]}(x,y)\leq \kappa$ for all $[\Delta]\in\mathfrak S$, there is a path in $W$ that joins $x$ to $y$ 
and has length bounded in terms of $\kappa,\delta,n$.

Therefore, fix $\kappa\geq 0$.  Let $x,y$ be vertices of $W$ (that is, maximal simplices of $X$) so that 
$d_{[\Delta]}(x,y)\leq \kappa$ for all $[\Delta]\in\mathfrak S$. We need $\theta$, depending on 
$\kappa,\delta,n$ only, so that $x$ and $y$ can be joined by a path in $W$ of length at most $\theta$, which is to say that $d_W(x,y)\leq \theta$ and $x,y$ lie in the same path-component of $W$.

To achieve this, we will prove by induction on $k$ that for every $\kappa$ there exists $\theta(\kappa,\delta,n,k)$ so that whenever $d_{X^{+W}}(x,y)\leq k$ and $d_{[\Delta]}(x,y)\leq 
\kappa$ for all $[\Delta]\in\mathfrak S-\{[\emptyset]\}$, the vertices $x$ and $y$ can be joined by a path in 
$W$ of length at most $\theta$, and in particular, $d_W(x,y)\leq \theta$.

Note that we are inducting on distance in $\mathcal C(\emptyset)=X^{+W}$, which is connected by Definition \ref{defn:combinatorial_HHS}.\eqref{item:chhs_delta}.

\emph{Base case $k=0$:}  First suppose that $k=d_{X^{+W}}(x,y)=0$. Then $x$ and $y$ share a vertex $v$, and they are of the form $x'\star v$, $y'\star v$, for some simplices $x',y'$ in the link of $v$, and necessarily $x',y'$ are maximal simplices in $\link(v)$. Moreover, 
$$d_{\mathcal C^v([\Delta])}(x',y')\leq \kappa'(\kappa,\delta,n)$$ for any non-maximal simplex $\Delta$ of 
$\link(v)$, by Proposition \ref{prop:hieromorphism}. 

By induction on complexity, in $W^v$ there is path 
$x'=x'_1,\dots,x'_l=y'$ of maximal simplices of 
$\link(v)$ with $l\leq \theta(\kappa,\delta,n)$. Hence, in $W$ we have a path $x=x'_1\star v,\dots,x'_l\star v=y$ (notice 
that these are indeed maximal simplices of $X$), showing $d_{W}(x,y)\leq l$.
 
\emph{Inductive step:}  Now suppose that $k=d_{X^{+W}}(x,y)>0$, and consider a geodesic $\gamma$ of the $\delta$--hyperbolic space $X^{+W}$ such that $\gamma$ is of minimal length among geodesics
connecting a vertex $u$ of $x$ to a vertex of $y$.  (Recall that hyperbolicity, and in particular the property of being a geodesic space, holds for $X^{+W}$ by 
Definition~\ref{defn:combinatorial_HHS}.\eqref{item:chhs_delta}.)  We consider two cases:\\

{\bf Case 1.} Suppose that one of the following holds:
\begin{enumerate}
\item there is a vertex $t$ on $\gamma-\{u\}$ that lies in $\sat(\sigma)$ for some non-maximal simplex $\sigma$ with $[\sigma]\nest [u]$ and $diam(\mathcal C(\sigma))\geq \delta$.
\item $u$ is connected in $X$ to the second vertex $t$ of $\gamma$.
 \end{enumerate}
 Consider the closest $t$ to $y$ satisfying either condition. 

Any $\sigma$ as above has the property that $[\sigma]\nest [t]$ since $t\in \sat(\sigma)$. Hence, in view of 
Definition \ref{defn:combinatorial_HHS}.(\ref{item:chhs_join}), there exists a simplex 
$\tau$ in the link of $u$ so that $[\tau\star u]\nest [t]$ and any $[\omega]$ with $diam(\mathcal C(\omega))\geq \delta$ 
nested into both $[u]$ and $[t]$ is nested into $[\tau\star u]$. For later purposes, we can and will pick $\tau=t$ if $t$ is 
the second vertex of $\gamma$ and $u,t$ are connected by an edge of $X$. Set $\tau'=\tau\star u$.

Before proceeding, observe that in either of the two itemized situations, $\tau'$ is a non-maximal simplex.  
Indeed, in the first case, where $\tau'$ comes from 
Definition~\ref{defn:combinatorial_HHS}.\eqref{item:chhs_join}, this is because of 
Remark~\ref{rem:join_nonmax}.  In the second case, where $\tau'=t\star u$, we can assume that $\tau'$ is 
non-maximal by the following argument.  If $\tau'$ is maximal, then $\dist_{X^{+W}}(\tau',y)\leq k-1$, so by 
the inductive hypothesis, $\dist_W(\tau',y)\leq\theta(\kappa,\delta,n,k-1)$, and in particular $y,\tau'$ lie in the same 
component of $W$.  

Also, 
$\dist_{X^{+W}}(\tau',x)=0$, so 
$\dist_W(x,\tau')\leq\theta(\kappa,\delta,n,0)$, and in particular $x,\tau'$ are in the same component of $W$.  

Transitivity of the binary relation ``are connected by a path in $W$'' on the vertex set of $W$ shows that $x,y$ are in the same component of $W$, so $\dist_W(x,y)$ is finite and the triangle inequality bounds $\dist_W(x,y)$, as required.  So, for the remainder of the argument, we can and shall assume that $\tau'$ is non-maximal.

\begin{claim}\label{claim:pleasant_z}
 There exists $C=C(\delta,n,\kappa)$ so that there is a maximal simplex $z''$ in $\link(\tau')$ with the property that, for 
$z=z''\star \tau'$, we have  $d_{\mathcal C([\Delta])}(z,y)\leq C$ for all non-maximal, non-empty simplices $\Delta$ of $X$.
\end{claim}

\renewcommand{\qedsymbol}{$\blacksquare$}
\begin{proof}
Let us consider $W^{\tau'}$, which is an HHS by Proposition \ref{prop:compatibility_1} and induction on $n$. 
Moreover, the 
projections to the various hyperbolic spaces for $W^{\tau'}$ can be computed from those for $W$, by 
Proposition~\ref{prop:hieromorphism}. 

We can now apply the 
realization theorem (Theorem~\ref{thm:realization}) to the coordinates of $y$ (which we verified to be 
consistent when we checked consistency) and find a 
maximal simplex $z''$ in $\link(\tau')$ with the property that $d_{\mathcal C([\Delta])}(z''\star \tau',y)\leq C$ for all 
$[\Delta]$ nested into $[\tau']$. Consider now some $[\Delta]$ not nested into $[\tau']$, and $\Delta\neq\emptyset$. 
Moreover, we can assume $diam(\mathcal C(\Delta))\geq \delta$. There are a few cases.

If $[\Delta]$ is not nested into $[u]$, then $u$ is not in $\sat(\Delta)$, and the projections of $y$ and $z$ to $\mathcal C([\Delta])$ coarsely coincide with the projection of $u$, as required.

Suppose that $[\Delta]$ is nested into $[u]$. Let $\gamma'$ be the final subgeodesic of $\gamma$ starting at $t$. If $\gamma'$ does not intersect $\sat(\Delta)$, then the projections of $y$ and $z$ to $\mathcal C([\Delta])$ coarsely coincide and we are done. Hence, suppose that $\gamma'$ intersects $\sat(\Delta)$. Since we chose $t$ as close as possible to $y$, the intersection must consist of $t$ only. Hence, $t$ lies in $\sat(\Delta)$, which implies that $[\Delta]$ is nested into $[t]$. Since $[\Delta]$ is nested into both $[u]$ and $[t]$, it is nested into $\tau'$, and in this case the bound $C$ exists by construction as stated above.
\end{proof}

Let $z'$ be a (necessarily) maximal simplex in $\link(u)$ so that $z''\star \tau'=z'\star u$ (so that $z=z'\star u$), and 
write $x=x'\star u$.

Now consider the HHS $(W^u,\mathfrak S_u)$, which exists by induction.  By assumption, for all non-maximal simplices $\Delta$, we have $\dist_{[\Delta]}(x,y)\leq\kappa$, and by Claim~\ref{claim:pleasant_z}, $\dist_{[\Delta]}(y,z)\leq C$, from which the triangle inequality gives $\dist_{[\Delta]}(x,z)\leq \kappa + C$.  Hence $\dist_{\mathcal C^u(\Delta)}(x',z')\leq \kappa+C$ whenever $[\Delta]\nest[u]$.  So, by the uniqueness axiom in the HHS $(W^u,\mathfrak S_u)$, the distance between $x'$ and $z'$ in $W^u$ can be bounded in terms of $\delta,n,\kappa$.  

In other words, in $W^u$ there is a path $x'=x'_1,\dots,x'_l=z'$ of 
maximal simplices of $\link(u)$ with $l\leq \theta(\kappa,\delta,n)$. Hence, there is a path $x=x'_1\star u,\dots,x'_l\star 
u=z$ in $W$.

\begin{claim}\label{claim:2-subcases-uniqueness} 
Either $d_{X^{+W}}(z,y)<d_{X^{+W}}(x,y)$ or $d_{X^{+W}}(z,y)=d_{X^{+W}}(x,y)$ and the first edge 
of $\gamma$ is not in $X$, while the first edge of some geodesic of minimal length from $z$ to $y$ is in $X$.
\end{claim}

\begin{proof}
Since $[\tau']$ is nested into $[t]$, we have in particular that any vertex of $z$ not in $\tau'$ is connected 
in $X$ to $t$. That is, $z$ is within distance $1$ of $t$ in $X^{+W}$, and more precisely either $t$ is a 
vertex of $z$ or $t$ is connected in $X$ to a vertex of $z$. Hence we have $d_{X^{+W}}(z,y)<d_{X^{+W}}(x,y)$ 
unless $t$ is the second vertex of $\gamma$ and \emph{not a vertex of $z$}. But then in this case the first edge of 
$\gamma$ is not an edge of $X$: if that edge $e$ was in $X$, since $t$ is an endpoint of $e$ we would have 
$t=\tau$. This is a contradiction: recall that $\tau\star u=\tau'$, so $t=\tau$ implies $t\in\tau'\subset z''\star\tau'=z$, contradicting that $t$ is not a vertex of $z$.  

On the other hand, there is a path of length $d_{X^{+W}}(x,y)$ from $z$ to $y$ starting with an edge of $X$. 
If this path is not a geodesic, then $d_{X^{+W}}(z,y)<d_{X^{+W}}(x,y)$, and otherwise the other possibility 
holds.
\end{proof}

If $\dist_{X^{+W}}(z,y)<\dist_{X^{+W}}(x,y)=k$, then by induction on $k$, there is a path in $W$, of length bounded by a constant depending only on $\kappa,k,\delta,n$, that joins $z$ to $y$.  Otherwise, by Claim~\ref{claim:2-subcases-uniqueness}, $\dist_{X^{+W}}(z,y)=k$, and some minimal-length $X^{+W}$--geodesic from $z$ to $y$ starts with an edge $e$ of $X$.  Now, $e$ is contained in some maximal simplex $z_1$ of $X$, and $z,z_1$ share a vertex, namely the initial vertex of $e$.  Noting that $\dist_{X^{+W}}(z_1,y)<\dist_{X^{+W}}(x,y)$, we argue as above to uniformly bound the distance in $W$ from $z_1$ to $y$, and apply the base case to produce a path in $W$ from $z$ to $z_1$ and a bound on $\dist_W(z,z_1)$.  Hence, by concatenating paths, we again see that $z,y$ are in the same component of $W$ and at distance bounded in terms of $\kappa,\delta,n,k$.  

Finally, we earlier exhibited a path of bounded length in $W$ from $x$ to $z$, so by concatenating again, we see that $x,y$ are in the same component of $W$, and obtain a bound on $\dist_W(x,y)$, as required.\\

{\bf Case 2.}  For each non-maximal simplex $\sigma$ with $[\sigma]\nest [u]$ and $diam(\mathcal C(\sigma))\geq \delta$, $\gamma-\{u\}$ does not intersect $\sat(\sigma)$. Moreover,  the first edge of $\gamma$ is not an edge of $X$.

Let $v$ be the second vertex of $\gamma$. In this case, there exist maximal simplices $p$ and $q$ that are joined by an edge 
in $W$, $p$ contains $u$ and $q$ contains $v$.

\begin{claim}
 There is $C=C(\delta,n,\kappa)$ so that for every $[\Delta]$ nested into $[u]$ we have $d_{\mathcal C([\Delta])}(p,x)\leq C$.
\end{claim}

\begin{proof}
Consider any $[\Delta]$ nested into $[u]$, and we can assume $diam(\mathcal C(\Delta))\geq \delta$.  The projection of $p$ to 
$\mathcal C(\Delta)$ coarsely coincides with the projection of $q$, since we checked that projections are coarsely Lipschitz 
and $p,q$ are adjacent in $W$, and in turn the projection of $q$ coarsely coincides with that of $v$, since the latter is 
well-defined by the hypothesis about $\gamma$ avoiding saturations. For the same reason plus Lemma \ref{lem:super_bgi}, the 
projection of $v$ coarsely coincides with that of $y$. Finally, by hypothesis the projection of $y$ coarsely coincides with 
that of $x$, and hence we are done.
\end{proof}
\renewcommand{\qedsymbol}{$\Box$}

Write $x=x'\star u$ and $p=p'\star u$.  Since $W^u$ is an HHS, in $W^u$ there is a path $x'=x'_1,\dots,x'_l=p'$ 
of length $l$ of maximal simplices of $\link(u)$ with $l\leq \theta(\kappa,\delta,n)$. Hence, there is a path 
$x=x'_1\star u,\dots,x'_l\star 
u=p$ in $W$. In particular, there is a path of length $l+1$ from $x$ to $q$. Clearly, 
$d_{X^{+W}}(q,y)<d_{X^{+W}}(u,y)$, and we 
are done.\\

\textbf{Conclusion:}  We have shown that the combinatorial HHS $(X,W)$ gives rise to a hierarchically hyperbolic space 
$(W,\mathfrak S)$.  It remains to prove the statement about group actions.  Let $G$ act on $X$.  By hypothesis, there are 
finitely many $G$--orbits of links of non-maximal simplices.  
Since $\mathfrak S$ corresponds $G$--equivariantly and bijectively to this set of links, the action of $G$ on $\mathfrak S$ 
is cofinite, as required by Definition~\ref{defn:HHG}.  

The action on $X$ also induces an action on 
the set of maximal simplices of $X$, and hence on the vertex set of $W$.  By hypothesis, this action extends to a 
metrically proper, cobounded action on $W$, as required by Definition~\ref{defn:HHG}.  

Let $\Delta$ be a non-maximal simplex, i.e., $[\Delta]\in\mathfrak S$.  For any 
$g\in G$, the automorphism $g:X\to X$ restricts to an isomorphism $g:\link(\Delta)\to\link(g\Delta)$.  In 
particular, 
the $G$--action on $\mathfrak S$ induced by the action on $X$ preserves $\nest,\orth,\transverse$.

Moreover, if $v,w\in\link(\Delta)$ are contained in $W$--adjacent maximal simplices $\sigma,\tau$, then $gv,gw$ are contained 
in $W$--adjacent maximal simplices $g\sigma,g\tau$, so the isomorphism $g:\link(\Delta)\to\link(g\Delta)$ extends to an 
isometry $g:\mathcal C(\Delta)\to\mathcal C(g\Delta)$.   

Similarly, since $g Sat(\Delta)=Sat(g\Delta)$, we get an isometry 
$g:Y_\Delta\to Y_{g\Delta}$.  Because $\pi_{[\Delta]}$ was defined in terms of coarse closest-point projection 
$Y_\Delta\to\mathcal C(\Delta)$, we have $\pi_{g[\Delta]}(gw)=g(\pi_{[\Delta]}(w))$ for all $w\in W$.  Similarly, if 
$[\Sigma]\propnest[\Delta]$ or $[\Sigma]\transverse[\Delta]$, then 
$\rho^{g[\Sigma]}_{g[\Delta]}=g(\rho^{[\Sigma]}_{[\Delta]})$.  Thus the action of $G$ on the HHS $(X,W)$ has all the 
properties listed in Definition~\ref{defn:HHG}, so $(G,\mathfrak S)$ is an HHG.
\end{proof}

\section{Hierarchical hyperbolicity from actions on hyperbolic complexes}\label{sec:propaganda}

A motivating application of Theorem \ref{thm:hhs_links} is to groups acting on hyperbolic simplicial complexes with 
finite stabilizers of maximal simplices. In this setup, at the cost of reducing the generality, one can even 
formulate conditions that imply those in Theorem \ref{thm:hhs_links} and do not need to refer to an 
$X$--graph; Theorem \ref{thm:propaganda} achieves just that. We start with two preliminary definitions:

\begin{defn}\label{defn:chhs_action}
 We say that the group $G$ \emph{acts on} the combinatorial HHS $(X,W)$ if it acts on $X$ in such a way that 
the action on maximal simplices extends to $W$. We say that the action on $(X,W)$ is proper (resp. cocompact) 
if the action on $W$ is proper (resp. cocompact).
\end{defn}

\begin{defn}[Hyperbolic $H$--space]\label{defn:hyp_H_space}
Let $H$ be a group acting on a simplicial complex $Y$.  We say that $Y$ is a \emph{hyperbolic $H$--space} if there is a 
graph $Y'$, consisting of $Y^{(1)}$ together with a set of additional edges $\mathcal E$ called the \emph{additional edge set}, such that $Y'$ is hyperbolic, the 
action of $H$ on $Y^{(1)}$ extends to an action of $H$ on $Y'$, and $\mathcal E$ contains finitely many $H$--orbits of edges.
\end{defn}

\begin{rem}\label{rem:qi_type}
 Varying the additional edge set in the definition above does not change the quasi-isometry type. In particular, if $Y$ is already connected, then it is a hyperbolic $H$--space if and only if it is already hyperbolic.
\end{rem}

Now we can state the theorem:

\begin{thm}\label{thm:propaganda}
  Let $G$ be a group acting cocompactly on the flag simplicial complex $X$, and suppose that the stabilizer of each 
maximal simplex is finite.  Suppose, in addition, that for all non-maximal simplices $\Delta,\Sigma$ of $X$:
  
\begin{enumerate}[(A)]
\item\label{item:stab_hyp} $\link(\Delta)$ is a hyperbolic $Stab_G(\link(\Delta))$-space quasi-isometrically 
embedded in $$\mathcal E\cup \left(X-\bigcup_{\link(\Sigma)=\link(\Delta)}\Sigma\right),$$ where $\mathcal E$ 
is the additional edge set,

\item\label{item:propaganda_join} $\link(\Delta)\cap\link(\Sigma)=\link(\Delta\star\Pi)\star\Pi'$ for some simplices $\Pi,\Pi'$ of $\link(\Delta)$,

\item \label{item:link_conn} $\link(\Delta)$ is connected unless $\Delta$ is a codimension--$1$ face of a maximal 
simplex.

\end{enumerate}

Then $G$ acts properly and cocompactly on a combinatorial HHS $(X,W)$. In particular, $G$ is a hierarchically hyperbolic group.  Moreover, $W$ can be chosen with the properties that
\begin{itemize}
    \item any two maximal simplices of $X$ corresponding to $W$--adjacent vertices share a codimension--$1$ face;
    \item $\mathcal C(\Delta)$ contains finitely many $\stabilizer_G(\link(\Delta))$--orbits of edges for each simplex $\Delta$ which is a codimension--$1$ face of a maximal simplex.
\end{itemize}
\end{thm}

\begin{rem}[Metric in Condition~\eqref{item:stab_hyp}]\label{rem:condition_c}
For clarity, in item \ref{item:stab_hyp}, the quasi-isometric embedding statement refers to the restriction to 
$\link(\Delta)$ of the hyperbolic metric built into  the assumption that it is a hyperbolic 
$Stab_G(\link(\Delta))$-space; recall that all 
such metrics are naturally quasi-isometric to each other. Since $\link(\Delta)$ is contained 
in $X-\bigcup_{\link(\Sigma)=\link(\Delta)}\Sigma$, it makes sense to add the additional edges to 
the latter graph. Finally, notice that, in view of condition \ref{item:link_conn}, for $\Delta$ not 
codimension--$1$ in a maximal simplex we could have more simply stated that $\link(\Delta)$ is hyperbolic and 
quasi-isometrically embedded into $X-\bigcup_{\link(\Sigma)=\link(\Delta)}\Sigma$.

Condition~\eqref{item:stab_hyp} is playing the role of Definition~\ref{defn:combinatorial_HHS}.\eqref{item:chhs_delta}: adding $\mathcal E$ amounts to augmenting the links, and $X-\cup_{\link(\Sigma)=\link(\Delta)}\Sigma$ is playing the role of $Y_\Delta$.
\end{rem}

\begin{rem}\label{rem:condition-B}
 Condition \eqref{item:propaganda_join} holds for curve graphs. Since we do not need this fact explicitly, we will only mention that it can be proven using arguments similar to those used in Claim \ref{claim:sat_structure}, and we leave the details to the interested reader. An interesting case to keep in mind, which shows the reason behind the  ``$\star\Pi'$'', is when $\Delta$ and $\Sigma$ are each a curve, and these curves fill the complement of some other curve. The intersection of the links consists of the latter curve only.
\end{rem}

\begin{rem}[Improper actions]\label{rem:improper}
 We will only use the condition on stabilizers of maximal simplices being finite to get properness of the action of $G$ on $W$. Dropping that condition, we can still conclude the following:
 \begin{itemize}
  \item $G$ acts coboundedly by HHS automorphisms on $W$, i.e., the HHS 
structure $(W,\mathfrak S)$, and the $G$--action on $W$, satisfy everything from Definition~\ref{defn:HHG} 
except for the requirement that the action of $G$ on $W$ is proper;
  \item  $Cay(G,S\cup\{H_i\})$ is equivariantly quasi-isometric to $W$, where the $H_i$ are stabilizers of representatives of orbits of maximal simplices and $S$ is a finite subset of $G$ so that $S\cup\{H_i\}$ is a generating set.
 \end{itemize}
\end{rem}

\begin{proof}[Proof of Theorem~\ref{thm:propaganda}]
We will construct $W$ as in the statement, and it will then follow from 
Theorem~\ref{thm:hhs_links} that $G$ is a hierarchically hyperbolic group.

\textbf{Construction of $W$:}  A simplex $\Delta$ of $X$ is \emph{almost-maximal} if $\Delta$ is a codimension--$1$ face of 
a maximal simplex of $X$.

Let $\{\link(\Delta_1),\dots,\link(\Delta_k)\}$ contain exactly one element of each $G$--orbit of links of almost-maximal 
simplices.  (Note that this is not quite the same as taking a list of $G$--orbit representatives of almost-maximal simplices 
and then taking links, since multiple $G$--distinct almost-maximal simplices can have the same link.)

For each $i\leq k$, there is an additional edge set $\mathcal E^i=\{e^i_1,\ldots,e^i_{\ell(i)}\}$ of edges such that, by adding edges $ge^i_j$ to 
$\link(\Delta_i)$ (for $j\leq \ell(i)$ and $g\in\stabilizer_G(\link(\Delta_i))$), we obtain a hyperbolic graph; this set of 
additional edges exists by condition~\eqref{item:stab_hyp}.  For each $i,j$, let $v^i_j,w^i_j$ denote the endpoints of 
$e^i_j$.

We are free to replace the hyperbolic graph obtained from $\link(\Delta_i)$ in this way by any other 
$\stabilizer_G(\link(\Delta_i))$--equivariantly quasi-isometric graph with finitely many 
$\stabilizer_G(\link(\Delta_i))$--orbits of ``additional'' edges.  For later use we make the following specific choice for the additional edges:

\begin{rem}\label{rem:how-we-added-edges}
By hypothesis, $\link(\Delta_i)$ (with the metric obtained by adding the edges $e_j^i$ and their translates) is 
quasi-isometrically embedded in $\mathcal E^i\cup (X-Sat(\Delta_i))$.  Suppose that $v,w\in\link(\Delta_i)$ 
lie at distance at most $2$ in 
$\mathcal E^i\cup X-Sat(\Delta_i)$.  Then $v,w$ lie at uniformly bounded distance (denoted $B$) in 
$\link(\Delta_i)$ (with the extra edges).  

Notice that $\stabilizer(\link(\Delta_i))$, and in fact even $\stabilizer(\Delta_i)$, acts with finitely many orbits of vertices on 
$\link(\Delta_i)$, since vertices $v',w'$ in $\link(\Delta_i)$ are in the same $\stabilizer(\Delta_i)$-orbit 
if $\Delta_i\star v',\Delta_i\star w'$ are in the same $G$-orbit of simplices with a marked vertex.

Hence $\stabilizer(\link(\Delta_i))$ acts cocompactly on 
$\link(\Delta_i)\cup \stabilizer_G(\link(\Delta_i))\cdot\{e^i_1,\ldots,e^i_{\ell(i)}\}$.
Moreover, $\link(\Delta_i)\cup \stabilizer_G(\link(\Delta_i))\cdot\{e^i_1,\ldots,e^i_{\ell(i)}\}$ is locally finite since $\stabilizer_G(\link(\Delta_i))$ acts on $\link(\Delta_i)$ with finite stabilizers of vertices (since $\Delta_i$ is almost-maximal, and hence $\Delta_i\star v$ is maximal for any vertex $v$ in $\link(\Delta_i)$). In particular, there are finitely many orbits of 
paths of length at most $B$.  Hence, by adding finitely many more $\stabilizer_G(\link(\Delta_i))$--orbits of edges to 
$\link(\Delta_i)$, we can and shall assume that any such $v,w$ are joined by a $\stabilizer_G(\link(\Delta_i))$--translate of 
an edge in 
$\{e^i_j\}$.  This assumption will simplify matters in the proof of Claim~\ref{claim:C_0=C} below.    
\end{rem}

Now, define an $X$--graph $W$ as follows:
\begin{itemize}
     \item The vertex set of $W$ is the set of maximal simplices of $X$.
     \item If $x,y$ are maximal simplices of $X$, we join $x,y$ by an edge if there exists $i\leq k,j\leq\ell(i)$ and $g\in 
G$ such that $x=g(\Delta_i'\star v^i_j)$ and $y=g(\Delta_i'\star w^i_j)$, where $\Delta_i'$ is an almost-maximal simplex 
with $\link(\Delta'_i)=\link(\Delta_i)$.  By construction, $W$--adjacent maximal simplices must intersect in a 
common 
almost-maximal simplex, which verifies the first item in the ``moreover'' clause of Theorem~\ref{thm:propaganda}.
\end{itemize}

\textbf{The $G$--action on $W$:}  The $G$--action on the set of maximal simplices of $X$ induces a $G$--action 
on $W$ by graph 
automorphisms, by the construction of $W$.  Moreover, since $G$ acts cocompactly on $X$, there are finitely 
many $G$--orbits of 
vertices in $W$.  Every edge has the form $\{g(\Delta_i\star v^i_j),g(\Delta_i\star w^i_j)\}$, where $\Delta_i$ is one of 
finitely many almost-maximal simplices and $(v^i_j,w^i_j)$ is one of finitely many pairs of vertices, so $W$ has finitely 
many $G$--orbits of edges.  Hence $G$ acts cocompactly on $W$.  

Since we have added finitely many orbits of edges to each $\link(\Delta_i)$, and stabilizers of maximal simplices are 
finite, it follows that $W$ is locally finite.  Since $G$ acts on $W$ with finite vertex stabilizers (because maximal 
simplices in $X$ have finite stabilizers by hypothesis), it then follows that $G$ acts properly on $W$.  

Hence, to complete the proof, it suffices to verify that $(X,W)$ is a combinatorial HHS.

\textbf{$(X,W)$ is a combinatorial HHS:} First, $W$ is an $X$--graph by construction.

Also, observe that condition \ref{item:propaganda_join} is stronger than Definition 
\ref{defn:combinatorial_HHS}.(\ref{item:chhs_join}). In fact, any $\Gamma$ whose link is not a non-trivial join (e.g., if 
$\diam(\mathcal C_0(\Gamma))\geq 3$)  and so that $\link(\Gamma)\subseteq\link(\Delta)\cap\link(\Sigma)$ cannot intersect 
$\Pi'$, and so we have $\link(\Gamma)\subseteq \link(\Delta\star\Pi)$. (Moreover, $\link(\Delta\star 
\Pi)\subseteq\link(\Sigma)$ by definition, so $[\Delta\star\Pi]\nest[\Sigma]$, as required.)

By definition, $X$ is a flag complex. The remaining parts of Definition~\ref{defn:combinatorial_HHS} are checked in the following claims.

\begin{claim}
 $X$ has finite complexity.
\end{claim}

\renewcommand{\qedsymbol}{$\blacksquare$}
\begin{proof}
If we had $\Pi'=\emptyset$ in condition \ref{item:propaganda_join}, then inclusion of links would yield 
reverse inclusion of representative simplices, and the proof would be straightforward from finite dimension. 
Since $\Pi'$ could be non-empty (which happens even in curve graphs), we need some understanding of join 
structures on links.

 For a simplex $\Delta$, let $\Theta_\Delta$ be any simplex in $\link(\Delta)$ so that $\link(\Delta)$ is a 
join of $\Theta_\Delta$ and some sub-complex, and $\Theta_\Delta$ is maximal with this property. For 
convenience, for a simplex $\Delta$ we let $\#\Delta=|\Delta^{(0)}|$. Define $c(\Delta)=(\# \Gamma_\Delta, -\# 
\Theta_\Delta)$, where $\Gamma_\Delta$ is any simplex with $\link(\Gamma_\Delta)=\link(\Delta\star 
\Theta_\Delta)$ with the maximal number of vertices among all choices.
 
 We claim that if $\link(\Delta)\subsetneq \link(\Sigma)$, then $c(\Sigma)<c(\Delta)$ in the lexicographic order. Since $G$ acts on $X$ cocompactly, 
$\dimension X<\infty$, so this readily implies finite complexity.

First, let us show $\link(\Delta\star \Theta_\Delta)\subseteq \link(\Sigma\star\Theta_\Sigma)$. Notice that $\Theta_\Sigma\cap \link(\Delta)\subseteq \Theta_\Delta$, since if we had a vertex $v$ in $\Theta_\Sigma\cap \link(\Delta)- \Theta_\Delta$, then we could add it to $\Theta_\Delta$ to form a larger simplex with the property that $\link(\Delta)$ is a join of the simplex and some sub-complex, contradicting maximality of $\Theta_\Delta$. Consider now a vertex $v\in \link(\Delta\star \Theta_\Delta)$, that is $v\in\link(\Delta)-\Theta_\Delta$. Then $v\in \link(\Sigma)$, and it cannot lie in $\Theta_\Sigma$, so $v\in \link(\Sigma\star\Theta_\Sigma)$, as required.

If $\link(\Delta\star \Theta_\Delta)\subsetneq \link(\Sigma\star\Theta_\Sigma)$, then condition \ref{item:propaganda_join} implies that we can write $\link(\Gamma_\Delta)=\link(\Gamma_\Sigma\star\Pi)\star\Pi'$. By maximality of $\Theta_\Delta$, we have $\Pi'=\emptyset$, so we must have $\Pi\neq\emptyset$. Hence, $\# \Gamma_\Sigma<\#\Gamma_\Sigma\star\Pi\leq \#\Gamma_\Delta$, and hence $c(\Sigma)<c(\Delta)$ (regardless of the number of vertices of $\Theta_\Delta,\Theta_\Sigma$).

If $\link(\Delta\star \Theta_\Delta)= \link(\Sigma\star\Theta_\Sigma)$, then $\#\Gamma_\Delta=\#\Gamma_\Sigma$, so we have to show $\# \Theta_\Delta<\# \Theta_\Sigma$.

Let $v\in\Theta_\Delta$ be a vertex. Then $v\in \link(\Delta)$, and in particular it lies either in $\link(\Sigma\star\Theta_\Sigma)$ or in $\Theta_\Sigma$. However, in the present situation, the former cannot occur since $v\notin \link(\Delta\star \Theta_\Delta)$. Hence, $\Theta_\Delta\subseteq \Theta_\Sigma$. If we had equality, then we would have
$$\link(\Delta)=\link(\Delta\star\Theta_\Delta)\star\Theta_\Delta=\link(\Sigma\star\Theta_\Sigma)\star \Theta_\Sigma=\link(\Sigma),$$
but we are assuming $\link(\Delta)\subsetneq \link(\Sigma)$. Hence, $\# \Theta_\Delta<\# \Theta_\Sigma$, as required.
\end{proof}

\begin{claim}\label{claim:join_max} Let $\Delta$ be a simplex, let $v,w\in\link(\Delta)$ be vertices, and let $\Sigma\star 
v,\Sigma\star w$ be maximal simplices of $X$. Then there exists a simplex $\Pi$ in $\link(\Delta)$ so  $\Delta\star\Pi\star 
v, \Delta\star\Pi\star w$ are maximal simplices of $X$.
\end{claim}

\begin{proof} 
If $v=w$, then we can just take any maximal simplex $\Pi\star v$ in $\link(\Delta)$, so suppose that this is 
not the case. 

Write $\link(\Delta)\cap\link(\Sigma)=\link(\Sigma\star\Pi_0)\star\Pi'_0$, as in condition 
\ref{item:propaganda_join}. Note that $\Sigma$ is almost maximal since, for instance, $\Sigma\star v$ is maximal.  Hence, $\link(\Sigma)$ is discrete, and thus $\link(\Delta)\cap\link(\Sigma)$ is discrete.  

We claim that $\Pi_0'=\emptyset$.  If not, then since $\link(\Delta)\cap\link(\Sigma)$ is discrete, $\Pi_0'$ is a single vertex and $\link(\Sigma\star\Pi_0)=\emptyset$.  This contradicts that $\link(\Sigma)\cap\link(\Delta)$ contains two distinct vertices, namely $v$ and $w$.  Hence $\Pi_0'=\emptyset$.  

Next, observe that $\Pi_0=\emptyset$.  Indeed, if not, then $\Sigma\star\Pi_0$ is a simplex whose link contains $v$ and $w$, so $\Pi_0\star v\star \Sigma$ is a simplex properly containing the $\Sigma\star v$, which was assumed to be maximal, a contradiction.

Hence $\link(\Delta)\cap\link(\Sigma)=\link(\Sigma),$ which is to say that $\link(\Sigma)\subseteq 
\link(\Delta)$.

Now, write $\link(\Sigma)=\link(\Delta)\cap\link(\Sigma)=\link(\Delta\star\Pi)\star\Pi'$, again using condition \ref{item:propaganda_join}. Since $\link(\Sigma)$ is discrete, we have $\Pi'=\emptyset$. But then $\Pi$ is the simplex we wanted (notice that $\link(\Delta\star\Pi)$ is discrete, so $\Delta\star\Pi\star v$, $\Delta\star\Pi\star w$ are maximal simplices).
\end{proof}

\begin{claim}\label{claim:C_0=C}
     Let $\Delta$ be a simplex of $X$ and let $\sigma_v,\sigma_w$ be $W$--adjacent maximal simplices of $X$ respectively containing 
vertices $v,w\in\link(\Delta)$ that are distinct and non-adjacent in $X$.  Then there exist simplices $\Sigma_v,\Sigma_w$ of $\link(\Delta)$ such that 
$\Sigma_v\star\Delta$ and $\Sigma_w\star\Delta$ are maximal simplices respectively containing $v,w$, and 
$\Sigma_v\star\Delta,\Sigma_w\star \Delta$ are adjacent in $W$.
\end{claim}

\begin{proof}[Proof of Claim~\ref{claim:C_0=C}]
Let $v,w$ be as in the statement.  Then there exists $\Delta_i$ such that $\sigma_v=g(\Delta_i\star v^i_j)$ and 
$\sigma_w=g(\Delta_i\star w^i_j)$ for some $g\in G,j\leq \ell(i)$.

Notice that since $v,w$ are distinct and not adjacent, neither $v$ nor $w$ is contained in $g\Delta_i$.  (Indeed, if $v\in g\Delta_i$, then $v\in g\Delta_i\star gw^j_i=g(\Delta_i\star w^i_j)$, and hence $v$ is adjacent in $X$ to $w$, contradicting our hypothesis.)

     Hence, $gv^i_j=v,gw^i_j=w$.  Since $\sigma_v,\sigma_w$ share the 
almost-maximal face $g\Delta_i$, Claim~\eqref{claim:join_max} provides a simplex $\Pi$ of $\link(\Delta)$ such that 
$\Delta\star\Pi\star v$ and $\Delta\star \Pi\star w$ are maximal simplices of $X$.  

Suppose that $g\Delta_i$ is not contained in $Sat(\Delta\star \Pi)$.  Then there is a path of length $2$ in $X$ from $v$ to 
$w$ that avoids $Sat(\Delta\star\Pi)$.  So, $\Delta\star\Pi\star v$ and $\Delta\star\Pi\star w$ are $W$--adjacent, because 
of how we added extra edges to the link of the almost-maximal simplex $\Delta\star\Pi$ --- see Remark~\ref{rem:how-we-added-edges}.  Hence we are done, with 
$\Sigma_v=\Pi\star v$ and $\Sigma_w=\Pi\star w$.

Otherwise, suppose that $g\Delta_i\subseteq Sat(\Delta\star\Pi)$.  So, $[\Delta\star\Pi]\nest [g\Delta_i]$.  
Hence, there exist simplices $\Pi'$ and $\Pi''$ such that 
$\link(\Delta\star\Pi)=\link(g\Delta_i\star\Pi')\star\Pi''$.  But $\Delta\star\Pi$, being almost-maximal, has 
discrete link, so $\Pi''=\emptyset$.  Thus there 
exists $\Pi'$ (a simplex of $\link(g\Delta_i)$) with $[\Delta\star\Pi]=[g\Delta_i\star\Pi']$.  Since $g\Delta_i$ is 
almost-maximal and $g\Delta_i\star\Pi'$ is necessarily non-maximal, we have $\Pi'=\emptyset$.  So, 
$[\Delta\star\Pi]=[g\Delta_i]$.  By definition, this means that $\link(\Delta\star\Pi)=\link(g\Delta_i)$.  So, the extra 
edges added to $\link(\Delta\star\Pi)$ --- which are determined by the link of $\Delta\star\Pi$, not the simplex itself --- 
are exactly the edges determined by $g\Delta_i$, so $\Delta\star\Pi\star v, \Delta\star\Pi\star w$ are $W$-adjacent.   
Again, we are done, with $\Sigma_v=\Pi\star v$ and $\Sigma_w=\Pi\star w$. 
\end{proof}

Claim~\ref{claim:C_0=C} says that $(X,W)$ satisfies condition~\eqref{item:C_0=C} from 
Definition~\ref{defn:combinatorial_HHS}.  We now verify condition~\eqref{item:chhs_delta}, which has two parts, the second of which verifies the second item in the ``moreover'' clause of the theorem:

\begin{claim}\label{claim:hyperbolic_C}
There exists $\delta$ such that $\mathcal C(\Delta)$ is $\delta$--hyperbolic for each non-maximal simplex $\Delta$ of $X$, and it is moreover obtained from $\link(\Delta)$ by adding finitely many $Stab_G(\link(\Delta))$-orbits of edges.
  
\end{claim}
\begin{proof}[Proof of Claim~\ref{claim:hyperbolic_C}]
Fix $\Delta$.  We first prove the statement about finitely many $\stabilizer_G(\link(\Delta))$-orbits of edges.  When $\link(\Delta)$ is in the same $G$--orbit as $\link(\Delta_i)$ for some $i$, we have by construction that vertices of $\link(\Delta)$ are joined by a $W$--edge only if they are joined by a $\stabilizer_G(\link(\Delta))$--translate of one of the finitely many additional edges in $\mathcal E^i$.  And, by discreteness of $\link(\Delta)$ in this case, no two vertices are joined by an $X$--edge.  Hence there are finitely many $\stabilizer_G(\link(\Delta))$--orbits of edges in $\mathcal C(\Delta)$ in this case.

In general, if $\Delta$ is not (up to equivalence of links) almost maximal, we argue as follows.  By Claim~\ref{claim:C_0=C}, vertices $v,w\in\link(\Delta)$ are $\mathcal C(\Delta)$ adjacent only if they are either adjacent in $X$ (i.e. in $\link(\Delta)$), or belong to $W$--adjacent maximal simplices of the form $\Delta\star\Sigma_v,\Delta\star\Sigma_w$.  In the latter case, the proof of Claim~\ref{claim:C_0=C} shows that there is a simplex $\Pi$ of $\link(\Delta)$ such that $\Delta\star\Sigma_v=\Delta\star\Pi\star v$ and $\Delta\star\Sigma_w=\Delta\star\Pi\star w$, and $[\Delta\star \Pi]=[g\Delta_i]$ for some $i\leq k$, and some $g\in G$, and moreover there exists $j\leq \ell(i)$ such that $gv^j_i=v,gw^j_i=w$.  

Recall that $\stabilizer_G(\Delta)\leq \stabilizer_G(\link(\Delta))$ acts cocompactly on $\link(\Delta)$.  In particular, $gv^j_i$ belongs to one of finitely many $\stabilizer_G(\link(\Delta))$--orbits of maximal simplices $\Pi\star v$ of $\link(\Delta)$.  Since $W$ is locally finite, it follows that the edge of $W$ from $\Delta\star\Pi\star v$ to $\Delta\star\Pi\star w$ belongs to one of finitely many $\stabilizer_G(\link(\Delta))$--orbits, and hence the edge from $v$ to $w$ in $\mathcal C(\Delta)$ also belongs to one of finitely many $\stabilizer_G(\link(\Delta))$--orbits.

Hence it remains to show that 
$\mathcal C(\Delta)$ is hyperbolic.  

By the proof of Lemma~\ref{lem:C_0=C} and the fact that we have already verified 
condition~\eqref{item:C_0=C} from Definition~\ref{defn:combinatorial_HHS}, we know that $\mathcal C_0(\Delta)=\mathcal C(\Delta)$ for any simplex $\Delta$, so it suffices to show that $\mathcal C_0(\Delta)$ 
is hyperbolic.  If $\Delta$ is almost-maximal, then the choice of edges implies that 
$\mathcal C_0(\Delta)=\mathcal C(\Delta)$ is hyperbolic.

Otherwise, $\link(\Delta)$ is connected, by condition~\eqref{item:link_conn}.  Since $\link(\Delta)$ can be made hyperbolic 
by adding finitely many $\stabilizer_G(\link(\Delta))$--orbits of edges, by condition~\eqref{item:stab_hyp}, $\link(\Delta)$ 
is quasi-isometric to a hyperbolic graph and is therefore hyperbolic.  So, it suffices to show that $\mathcal 
C(\Delta)$ is quasi-isometric to $\link(\Delta)$.  

Now, the inclusion $\link(\Delta)\to \mathcal C(\Delta)$ is Lipschitz and bijective on vertex sets, so we need to show that 
the inverse map on vertex sets is coarsely Lipschitz.  

Suppose that $v,w\in \mathcal C(\Delta)$ are adjacent.  If $v,w$ are adjacent in $\link(\Delta)$, then we are done.  So, 
suppose that $v,w$ belong to $W$--adjacent maximal simplices $x,y$.  Then $x\cap y=\Sigma$ is an almost-maximal simplex, 
because of how the edges in $W$ were defined. Now, if $\Sigma\subseteq Sat(\Delta)$, then 
$\link(\Delta)\subseteq\link(\Sigma)$.  But since $\Sigma$ is almost-maximal, $\link(\Sigma)$ contains no 
$X$--edges, whence 
$\link(\Delta)$ also contains no $X$--edges.  But this means that $\Delta$ is almost-maximal, a contradiction. 
 So $\Sigma$ contains a 
vertex $u$ of $X-Sat(\Delta)$, so $v,u,w$ is a path of length $2$ in $X-Sat(\Delta)$ from $v$ to $w$.  Since $\link(\Delta)$ 
is quasi-isometrically embedded in $X-Sat(\Delta)$, by condition~\eqref{item:stab_hyp} and the fact that $\link(\Delta)$ is connected (see Remark \ref{rem:qi_type}), this means that $v,w$ lie at 
uniformly bounded distance in $\link(\Delta)$, as required.
\end{proof}

The final claim is:

\begin{claim}\label{claim:qi_embedding}
There exists $\delta$ such that $\mathcal C(\Delta)$ is $(\delta,\delta)$--quasi-isometrically embedded in $Y_\Delta$, for 
all non-maximal simplices $\Delta$ of $X$.     
\end{claim}
\begin{proof}[Proof of Claim~\ref{claim:qi_embedding}]

It suffices to show the claim for a fixed $\Delta$.

Let $Z_\Delta$ be obtained from 
$X-Sat(\Delta)$ by connecting $\mathcal C(\Delta)$-adjacent vertices of $\link(\Delta)$. By the second part of Claim \ref{claim:hyperbolic_C} and condition \eqref{item:stab_hyp}, $\mathcal C(\Delta)$ is quasi-isometrically embedded into $Z_\Delta$.

We now show that $Z_\Delta\hookrightarrow Y_{\Delta}$ is a uniform quasi-isometric embedding.

Since $Z_\Delta$ is a subgraph of $Y_\Delta$ and the inclusion $Z_\Delta\hookrightarrow Y_\Delta$ is bijective on vertices, 
it suffices to show that if $e$ is an edge of $Y_\Delta$ that is not an edge of $Z_\Delta$, then the endpoints $v,w$ of $e$ 
are uniformly close in $Z_\Delta$.  Any such $v,w$ are contained in maximal simplices $x,y$ of $W$.  Now, by construction of 
$W$, we have that $x\cap y$ is an almost-maximal simplex.  If $(x\cap y)^{(0)}\not\subseteq Sat(\Delta)$, then $v,w$ are 
joined by a path of length $2$ in $X-Sat(\Delta)\subset Z_\Delta$, and we are done.  Otherwise, $(x\cap y)^{(0)}\subset 
Sat(\Delta)$.  Hence $\link(\Delta)\subseteq\link(x\cap y)$ so, since $x\cap y$ is almost-maximal, we have that 
$[\Delta]=[x\cap y]$.  Hence $v,w$ are joined by an edge of $\mathcal C(\Delta)$ (coming from the hyperbolic $G$--graph 
structure on $\link(\Delta)=\link(x\cap y)$).  Thus $v,w$ are adjacent in $Z_\Delta$, as required.
\end{proof}
\renewcommand{\qedsymbol}{$\Box$}
This completes the proof that Definition~\ref{defn:combinatorial_HHS}.\eqref{item:chhs_delta} holds for $(X,W)$, and hence 
completes the proof. 
\end{proof}

\section{Mapping class group quotients}\label{sec:MCG}

We now state Theorem~\ref{thm:MCG_quotient}, describing hierarchically hyperbolic quotients of mapping class groups.  In 
this section we discuss the various consequences of the theorem, while the proof is given in Section 
\ref{sec:proof_quotients}.

In this section, and the next section, we abuse language slightly: we
often use the same notation for a simplicial graph and for the flag
complex determined by the graph.  In particular, we do not distinguish
between the curve graph $\mathcal C(S)$ of a surface $S$, and the
curve complex $\mathcal C(S)$.  For example, when talking about the
metric on $\mathcal C(S)$, we mean the graph metric on the
$1$--skeleton; when talking about a combinatorial HHS $(X,W)$ with
$X=\mathcal C(S)$, we mean the full curve complex.

Recall that the complexity of a connected orientable surface $S$ of
finite type is $\xi(S)=3\genus(S)+p(S)-3$, where $\genus(S)$ is the genus and
$p(S)$ is the number of punctures.

\begin{thm}\label{thm:MCG_quotient}
 Let $S$ be a connected orientable surface of finite type of complexity at least 2.
Let $F\subseteq 
MCG(S)$ be any finite set, and let $Q<MCG(S)$ be a convex-cocompact subgroup. If all 
hyperbolic groups are residually finite, then the following holds.

For all $-1 \leq i < \xi(S)$ there exists a normal subgroup $N_i\triangleleft MCG(S)$ such that the quotient 
$\phi \co MCG(S)\to MCG(S)/N_i=\bar{G}_i$ has the 
following properties:
 \begin{enumerate}[(I)]
  \item\label{item:larg_inj_radius} (\textbf{Large injectivity radius.}) $\phi|_{F}$ is injective.
  \item \label{item:HH_structure} (\textbf{Explicit HHS structure.}) The action of $\bar G_i$ on $\mathcal C(S)/N_i$ satisfies 
the hypotheses of Theorem \ref{thm:propaganda}, so that $\bar{G}_i$ is a hierarchically hyperbolic group.  More precisely, 
$\bar G_i$ acts properly and cocompactly on a 
combinatorial HHS $(\mathcal C(S)/N_i,W)$, and the corresponding HHS structure $(\bar G_i,\mathfrak S_{N_i})$ satisfies:
\begin{itemize}
 \item the map $b:\mathfrak S_{N_i}\to\mathfrak S^{\geq 1}/N_i$ from Definition \ref{defn:b} below is well-defined and a bijection,  where $\mathfrak S^{\geq 
1}$ is the set of isotopy classes subsurfaces of $S$ without annular or thrice-punctured sphere components;
 \item two elements $U,V\in\mathfrak S_{N_i}$ are nested 
(resp. orthogonal) if and only if $b(U),b(V)$ have representatives in $\mathfrak S^{\geq 1}$ that are nested (resp. 
disjoint);
 \item there exists $B$ so that for any element $U$ of $\mathfrak S_{N_i}$  such that $b(U)$ has a representative of 
complexity at most $i$, we have that $\mathcal C(U)^{(0)}$ is finite and $\diam(\mathcal C(U))\leq B$.
\end{itemize}
\item\label{item:poc} (\textbf{Convex-cocompact injects.}) $\phi|_Q$ is injective and the orbit maps of $Q$ to $\mathcal 
C(S)/N_i$ are quasi-isometric embeddings; in particular $\bar{G}_i$ is infinite.
 \end{enumerate}
\end{thm}

The existence of a bijection $b$ with the stated properties is sufficient for many applications; in practice we will use the following construction (the existence of lifts of simplices of $\mathcal C(S)/N_i$ to $\mathcal C(S)$ is part of the claim in the theorem that $b$ is well-defined):

\begin{defn}\label{defn:b}
 Given any non-maximal simplex $\Delta$ of $C(S)/N_i$, consider a lift $\lift\Delta$ to $\mathcal C(S)$. The vertex set of the 
link of $\lift\Delta$ in $\mathcal C(S)$ consists of all curves (regarded as vertices of $\mathcal C(S)$) contained in a subsurface that we denote $S_{\lift\Delta}$. Define $b([\Delta])=[S_{\lift\Delta}]_{N_{i}}$, where $[\cdot]_{N_{i}}$ denotes the $N_{i}$--orbit.
\end{defn}

The groups $\bar{G}_i$ will be constructed inductively, with $\bar{G}_{i+1}$ being a quotient of $\bar{G}_i$, and $\bar{G}_{-1}$ being a quotient by suitable powers of Dehn twists. In particular, $N_{-1}$ is a normal subgroup generated by normal powers of Dehn twists, and $N_i<N_{i+1}$ for all $i$. Roughly speaking, for each $i$, we have to ensure that subgroups of the mapping class group coming from subsurfaces of complexity $i$ become finite, and to do so we mod out finite-index subgroups of those subgroups. More information on all this is provided at the beginning of Section \ref{sec:proof_quotients}.

\begin{rem*}\label{rem:resid_finite_hypothesis}
The residual finiteness hypothesis will be applied to particular 
hyperbolic groups which arise in the proof.
\end{rem*}

Theorem~\ref{thm:MCG_quotient} will be proven in Section~\ref{subsec:main_theorem}. 
We now establish the corollaries stated in 
the introduction. Additionally, in the case of the 
closed genus 2 surface, we prove  
Theorem~\ref{thm:MCG_quotient} without needing to assume 
residual finiteness for hyperbolic groups.

For Theorem~\ref{thm:MCG_quotient}, the case of $i=-1$ warrants extra 
focus, since there the mapping classes being quotiented are those 
supported in annuli, namely Dehn twists. Accordingly, in this case we 
provide a more explicit description of the kernel $N_{-1}$, which we 
state as Theorem~\ref{thm:MCG_mod_Dehn}. (In the general case, we 
prefer to keep the statement of Theorem \ref{thm:MCG_quotient} 
more concise rather than adding a more detailed 
description of $N_i$.) Also, note that in this case there is no 
residual finiteness assumption required.

\begin{thm}\label{thm:MCG_mod_Dehn}
 Let $S$ be a connected orientable surface of finite type of complexity at least 2.  Given $K\geq 1$, denote by $DT_K$ the normal subgroup 
generated by all $K$-th powers of Dehn twists. There exists $K_0\ge1$ so that, for any multiple $K$ of $K_0$, 
$MCG(S)/DT_K$ is an infinite hierarchically hyperbolic group. More precisely, given $F,Q$ as in 
Theorem~\ref{thm:MCG_quotient}, all 
conclusions of Theorem \ref{thm:MCG_quotient} hold with $i=-1$ and $N_{-1}=DT_K$ for any sufficiently large multiple of 
$K_0$.
\end{thm}

\begin{proof}
The proof follows verbatim the proof of 
Theorem~\ref{thm:MCG_quotient}, applied in the case $i=-1$, with the following modifications:
\begin{itemize}
     \item We do not need the choice of $H$ in Lemma \ref{lem:fi}, only the coloring of the subsurfaces 
described in the lemma, and the following fact: there is $K'_0>0$ so that any 
$K'_0$--th power of an element of $MCG(S)$ preserves the coloring.
    \item The choice of $N=N_{-1}$ 
made in Notation \ref{not:correct_powers} below can be replaced by choosing $\Gamma^\theta_Y<\langle \tau_Y\rangle$ to be 
$\langle \tau_Y^K\rangle$ for a suitably large multiple $K$ of $K'_0$.  
\end{itemize}
Here, and in Corollary \ref{cor:genus_2_non_intro}, we take $K_0$ to be such a suitably large multiple of 
$K'_0$. 
\end{proof}

In the special case of Theorem~\ref{thm:MCG_mod_Dehn} for the 
the genus-two closed surface, an even stronger conclusion holds which we 
obtain in the following.

\begin{cor}\label{cor:genus_2_non_intro}
  There exists $K_0\ge1$ so that for all non-zero multiples $K$ of $K_0$, the following holds. The 
quotient $MCG(\Sigma_2)/DT_K$ is hyperbolic relative to an infinite index subgroup commensurable to the product of two 
$C'(1/6)$-groups, where $DT_K$ denotes the normal subgroup generated by all $K$-th powers of Dehn twists.
\end{cor}

\begin{proof}
 The peripheral subgroup of the relatively hyperbolic structure will
 be the image $H$ in $G_K$ of the stabilizer $\lift H$ of a fixed
 curve $\gamma$ that cuts $\Sigma_2$ into two $\Sigma_{1,1}$
 subsurfaces.  Note that $\lift H$ is virtually a central extension by
 a Dehn twist of a product of virtually free groups, which are
 isomorphic to the mapping class group of $\Sigma_{1,1}$.  By
 \cite[Proposition 5.8]{DFDT}, provided $K$ is sufficiently large, we
 have the following.  The subgroup $H$ arises from $\lift H$ as the
 quotient by the subgroup generated by $K$-th powers of Dehn twists
 around $\gamma$ and curves contained in one of the $\Sigma_{1,1}$.
 In particular, $H$ is commensurable to the product of two groups, each of which is the quotient of a free 
group by $K$-th powers of certain elements, and
 this finite collection of elements is independent of $K$.  In
 particular, up to increasing $K$, $H$ is commensurable to the product
 of two finitely presented $C'(1/6)$-groups.
 
 Note that $H$ has infinite index, for example because $G_K$ is 
 acylindrically hyperbolic \cite[Theorem 3.1]{DFDT}, and 
thus cannot be commensurable to a product of infinite groups.
 
 We are left to check relative hyperbolicity, for which we use \cite{Russell:rel_HHS}.  Recall that the index set of the HHS 
structure is $\mathfrak S^{\geq 1}/DT_K$, with elements being nested if and only if they have nested representatives, and 
similarly for orthogonality. The only orthogonal pairs in $\mathfrak S^{\geq 1}$ are pairs $Y_i,W_i$ of surfaces homeomorphic 
to $\Sigma_{1,1}$. The set of surfaces $\{Y_i\cup W_i\}$ satisfies the following two properties:
 \begin{itemize}
  \item Whenever $U,V\in\mathfrak S^{\geq 1}$ satisfy $U\orth V$, there exists $i$ so that $U,V\subseteq Y_i\cup W_i$,
  \item For $i\neq j$, there is no $U\in \mathfrak S^{\geq 1}$ so 
  that $U\subseteq Y_i\cup W_i$ and $U\subseteq Y_j\cup W_j$.
 \end{itemize}
Therefore, the analogous properties hold for $\{Y_i\cup W_i\}/DT_K\subseteq \mathfrak S^{\geq 1}/DT_K$, that is,  
$(G_K,\mathfrak S^{\geq 1}/DT_K)$ has \emph{isolated orthogonality} in the sense of \cite[Definition 
4.1]{Russell:rel_HHS}.  Hence, by \cite[Theorem 4.3]{Russell:rel_HHS}, $G_K$ is hyperbolic relative to $H$.
\end{proof}

\begin{rem}
 We believe that that quotients of $MCG(\Sigma_2)$ by suitable large powers of Dehn twists around 
\emph{separating} curves are hyperbolic relative to subgroups which are virtually a direct product of free 
groups. However, we cannot use quotients of curve graphs to witness this, since, roughly, those quotients are 
expected to have an HHS structure that still has annular curve graphs corresponding to non-separating curves, 
and these annular curve graphs are not ``visible'' in the curve graph of $\Sigma_2$.
\end{rem}

\begin{rem}\label{rem:Francesco_trick}
Corollary \ref{cor:genus_2_non_intro} implies that $MCG(\Sigma_2)$ is
fully residually non-elementary hyperbolic; we now explain why, and
then provide a different argument for this fact relying on results of
a very different nature dating back to \cite{Picard1,Picard2}.

 Let $K_0$ be as in Theorem \ref{thm:MCG_mod_Dehn}.  Up to replacing 
 $K_0$ with a multiple, we can assume this constant is large enough 
 to satisfy the hypothesis of   
\cite[Theorem 3.1, Proposition 5.8]{DFDT}. Let $F\subseteq MCG(\Sigma_2)$ be finite. Using Theorem 
\ref{thm:MCG_mod_Dehn}, choose a non-zero multiple $K$ of $K_0$ so that $\phi|_F$ is injective, where $\phi:MCG(\Sigma_2)\to 
MCG(\Sigma_2)/DT_K=G_K$ is the quotient map. 
 
 The peripheral subgroup of the relatively hyperbolic structure on $G_K$ from Corollary \ref{cor:genus_2_non_intro} is residually finite by residual finiteness of 
$C'(1/6)$-groups, which follows from applying 
\cite[Theorem 1.2]{Wise:cubulating_small_can}, \cite[Theorem 
1.1]{vhak}, and \cite[Theorem 
4.4]{HaglundWise:special}. 

We can then apply the relatively hyperbolic Dehn filling theorem \cite{Os-perfill,GrMa-perfill} to 
construct a non-elementary hyperbolic quotient of $G_K$ where $F$ embeds, as we wanted.

The simpler argument is based on the following observation that was pointed out to us by Francesco Fournier Facio. Suppose that the residually finite group $G$ has a non-elementary hyperbolic quotient $H$. Then $G$ is fully residually non-elementary hyperbolic. To see this, fix a finite subset $F$ of $G$, and consider a finite quotient $Q$ of $G$ where $F$ embeds. Then $G$ also maps to $H\times Q$, the image is non-elementary hyperbolic, and $F$ embeds.

This observation can be applied to $MCG(\Sigma_2)$. In fact, mapping class groups are residually finite \cite{Grossman:res_fin_MCG}, and a non-elementary hyperbolic quotient of $MCG(\Sigma_2)$ can be constructed as follows, as pointed out to us by Ian Agol. First, $MCG(\Sigma_2)$ maps onto $MCG(\Sigma_{0,6})$ by modding out by the hyperelliptic involution (see e.g. \cite[Proposition 3.3]{genus_2_linear}). In turn, $MCG(\Sigma_{0,6})$ maps onto the fundamental group of a finite-volume complex hyperbolic orbifold, see \cite[Theorem 0.2]{Thurston:quotients} or references therein. Such fundamental group has a non-elementary hyperbolic quotient, say by the relatively hyperbolic Dehn filling theorem.

As a side remark, we note that \emph{pure} mapping class groups of
punctured spheres with at least 4 punctures map onto the
non-elementary hyperbolic group $PMCG(\Sigma_{0,4})$ via repeated use
of the Birman exact sequence.  It would be interesting to know if 
this could be promoted to the 
the full mapping class groups of punctured spheres.
\end{rem}

In the proof of Corollary \ref{cor:fully_residually_hyperbolic} we 
explain how to construct hyperbolic quotients of the mapping class 
group; we then employ this construction in the remaining 
corollaries.

\begin{cor}\label{cor:fully_residually_hyperbolic}
     Let $S$ be a connected orientable surface of finite type of complexity at least 2.  If all 
hyperbolic groups are residually finite, then $MCG(S)$ is fully residually non-elementary hyperbolic.
\end{cor}

\begin{proof}
 Let $F\subseteq MCG(S)$ be finite, and let $Q$ be any
 convex-cocompact subgroup of $MCG(S)$ isomorphic to the free group on
 2 generators (for instance, by \cite{Fujiwara:Freesubgroups2}, any 
 sufficiently high powers of a pair of independent pseudo Anosovs 
 will yield such a $Q$).
 
 By Theorem \ref{thmi:MCG}, there is a hierarchically hyperbolic
 quotient $\bar G$ of $MCG(S)$ such that: $MCG(S)\to\bar G$ is
 injective on $F$ and $Q$; $Q$ quasi-isometrically embeds into
 $\bar G$; and, all hyperbolic spaces in the HHS structure are bounded, except 
 for the one space associated to the $\nest$--maximal domain 
 (the quotient of $\mathcal C(S)$). Thus, by
 \cite[Corollary 2.15]{HHS:quasiflats}, $\bar G$ is hyperbolic, and
 since $Q$ embeds, $\bar G$ is non-elementary.
\end{proof}

\begin{cor}\label{cor:separable}
     Let $S$ be a connected orientable surface of finite type of complexity at least 2.  If all hyperbolic groups are residually finite, then 
every convex-cocompact subgroup of $MCG(S)$ is separable.
\end{cor}

\begin{proof}
We will show below that all \emph{torsion-free} convex-cocompact
subgroups are separable.  This is sufficient, by the following
argument.  Let $Q$ be convex-cocompact and let $\widehat Q$ be a
finite index torsion-free subgroup of $Q$, which exists by, say,
intersecting $Q$ with a torsion-free finite-index subgroup of
$MCG(S)$.  Then $\widehat Q$ is closed in the profinite topology,
which implies that its cosets are also closed.  So, $Q$ is closed
since it is a finite union of closed sets.  This reduces the claim to
the case where $Q$ is torsion-free, which we now address.

 Let $Q$ be a convex-cocompact torsion-free subgroup, and let $g\in MCG(S)-Q$.  We consider two cases:
 
\emph{Non-pseudo-Anosov case.} First suppose that $g$ is reducible or
periodic (that is, it acts with bounded orbits on $\mathcal
C(S)$).  Construct a hyperbolic quotient $\bar G$ of $MCG(S)$ as in
the proof of Corollary \ref{cor:fully_residually_hyperbolic}, for
$F=\{1,g\}$ and our given $Q$.  The image $\bar Q$ of $Q$ is
quasi-convex in $\bar G$, and the image $\bar g$ of $g$ is
non-trivial.  Since $g$ has bounded orbits in $\mathcal C(S)$, we have
that $\bar g$ has bounded orbits in the quotient of $\mathcal C(S)$.
Hence, $\bar g$ has finite order, and in particular it is not in $\bar
Q$, since $\bar Q$ is also torsion-free.  In view of the fact that we
are assuming that all hyperbolic groups are residually finite, by
\cite[Theorem 0.1]{AGM} we can find a finite quotient of $\bar
G$ in which the image of $\bar g$ is not in the image of $Q$.
 
 \emph{Pseudo-Anosov case.} For sufficiently large $n>0$, the subgroup $\langle g^n, Q\rangle$ is convex-cocompact and is 
naturally isomorphic to $\langle g^n\rangle*Q$ (see e.g.~\cite[Theorem M]{RST:morse}).  Construct a hyperbolic 
quotient $\bar G$ of $MCG(S)$, as in the proof of 
Corollary \ref{cor:fully_residually_hyperbolic}, except using the convex-cocompact subgroup $\langle g^n, Q\rangle$.  Since 
$\langle g^n, Q\rangle$ is quasi-isometrically embedded in $\bar G$, so is $Q$.  Moreover, $\bar g\not\in\bar Q$ since 
$\langle g^n\rangle*Q\to\bar G$ is a quasi-isometric embedding.  We conclude as above using \cite[Theorem 
0.1]{AGM}.
\end{proof}

\begin{cor}\label{cor:basically_omnipotence}
Let $S$ be a connected orientable surface of finite type of complexity at least 2.  If all hyperbolic groups are residually finite, then the 
following holds. Let $g,h\in MCG(S)$ be pseudo-Anosovs with no common proper power, and let $q\in \mathbb Q_{>0}$. Then there 
exists a finite group $G$ and a homomorphism $\psi:MCG(S)\to G$ so that $ord(\psi(g))/ord(\psi(h))=q$, where $ord$ denotes 
the order.
\end{cor}

\begin{proof}
For sufficiently large $n$, the elements $g^n,h^n$ freely generate a convex-cocompact free subgroup $Q$.  Construct $\bar G$ 
as in the proof of Corollary \ref{cor:fully_residually_hyperbolic}.  Hence we have a hyperbolic quotient $\bar G$ of 
$MCG(S)$ where the images $\bar g,\bar h$ of $g$ and $h$ have infinite order and have no common proper power.
 
We can now quotient $\bar G$ by suitable (large) powers of $\bar g$ and $\bar h$ to find a further hyperbolic quotient 
$\check{G}$ where the images of $g$ and $h$ satisfy the condition on the orders as in the statement. Using residual 
finiteness of $\check{G}$ we finally find the finite quotient of $MCG(S)$ that we were looking for.
\end{proof}

\section{Proof of Theorem \ref{thm:MCG_quotient}}\label{sec:proof_quotients}

\subsection{Outline}\label{subsec:MCG_proof_outline}

We start with a rough outline of the proof of Theorem
\ref{thm:MCG_quotient}, in which we will take successive 
quotients of $MCG(S)$.  

\subsubsection{First quotient: Dehn twists} We start by describing the
first quotient, which is the quotient of $MCG(S)$ by the normal
subgroup generated by suitable powers of Dehn twists.  In this
outline, we denote this normal subgroup by $N$.  We will check
hierarchical hyperbolicity of $MCG(S)/N$ by considering its action
on $\mathcal C(S)/N$ and applying Theorem \ref{thm:propaganda}.  It 
was already proven in \cite{DFDT} that $\mathcal C(S)/N$ is 
hyperbolic; here we will further develop the technology from 
\cite{DFDT} to gain additional information about $\mathcal C(S)/N$.

The key tool will be lifts.  In particular, the way that hyperbolicity
of $\mathcal C(S)/N$ is proven in \cite{DFDT} is by showing that
geodesic triangles in $\mathcal C(S)/N$ can be lifted to geodesic
triangles in $\mathcal C(S)$ (compare with Proposition
\ref{prop:omni_lift} below).  
An important tool for doing such 
lifting is a version of ``Greendlinger's Lemma.'' 
Roughly, this provides us with a normal form in which every term 
contributes a large projection to some domain and forces the lift to 
travel near some specific vertex in the curve graph. 
The tool through which we obtain our Greendlinger's Lemmas, 
Lemma \ref{lem:short_less_complex_base},  is 
that of a \emph{composite rotating family} in the sense of 
\cite{D_PlateSpinning}.

The key generalization we provide here is that, rather than lifting 
triangles, here we lift more
general objects, namely \emph{generalized $m$--gons}.  A generalized
$m$--gon is roughly a concatenation of simplices and geodesics in
links; we formalize this idea in Definition~\ref{defn:generalized_m_gon}.

We will show that generalized $m$--gons can be lifted, provided that
$m$ is not too large. This will be the main tool to reduce various
statements about links in $\mathcal C(S)/N$ to statements about links
in $\mathcal C(S)$, which can be verified by curves-on-surfaces
considerations.

To prove hierarchical hyperbolicity of $MCG(S)/N$ we use this type of argument repeatedly; see Subsections 
\ref{subsec:MCG_quotient_lemmas} and \ref{subsec:main_theorem}. 

The condition of Theorem \ref{thm:propaganda} which requires the 
most new work to check is the quasi-isometric embedding requirement in 
Theorem~\ref{thm:propaganda}.\eqref{item:stab_hyp}. In order to check
that condition, we use certain concatenations of geodesics in links 
which we call
\emph{approach paths} below.  These will give rise to the most
complicated generalized $m$--gons that we will consider.

Finally, choosing sufficiently large powers of Dehn twists allows one
to make sure that a given finite set embeds in the quotient. 
Moreover, this  
allows one to preserve the contracting directions, which are 
characterized by having bounded projections to all proper 
subsurfaces; these are the ``convex-cocompact'' directions. 
This is essentially
because the version of Greendlinger's Lemma mentioned above says that
nontrivial elements of $N$ create large subsurface projections.

\subsubsection{Further quotients} So far, we found that the quotient of
$MCG(S)$ by suitable powers of Dehn twists is hierarchically hyperbolic.  To pass to further hierarchically 
hyperbolic quotients, and eventually to a hyperbolic quotient, we use a similar method, approximately speaking. 

At any given stage, we have a hierarchically hyperbolic group, and at the ``bottom level
of the hierarchy'' we have hyperbolic groups (in the case of the first
quotients, we had $\mathbb Z$ subgroups). In each of the hyperbolic groups, we take
sufficiently deep finite-index subgroups, and quotient by those.  In
these later stages, in order to find the appropriate finite-index
subgroups, we use the hypothesis that hyperbolic groups are residually finite.

Once again, in this setting we establish a composite rotating family,
a Greendlinger's Lemma, and the ability to lift.  However, this 
approximate description hides several technical difficulties, some quite
serious, as we now explain.

As mentioned above, we use the combinatorial/geometric structure of $\mathcal C(S)$ to prove various properties 
of the quotient via lifting. This means that either one has to make sure that all those properties pass to the 
quotients of the curve graphs, or to take lifts to $\mathcal C(S)$ for the further quotients as well. We choose the second 
option.

The ability to first lift to the
previous quotient of $\mathcal C(S)$, and from there to $\mathcal
C(S)$ is an essential aspect of our induction. 
Accordingly, we establish this as
Proposition
\ref{prop:technical_MCG_quotient}.\eqref{item:lift_to_CS}-\eqref{item:lift_splx_induction}
below.

More generally, we collect all the properties that are required for the 
inductive hypothesis in Proposition \ref{prop:technical_MCG_quotient}; a
couple of them follow from the others, but we found it helpful to have
all of them collected in a single place.

Even given the ability to lift, just checking that the aforementioned sufficiently deep finite-index subgroups define a composite rotating family requires a
significant amount of work, since relatively straightforward lifting
arguments are not sufficient.  This is one of a few places where we
found it efficient to insist that the kernels of our quotients are
contained in a carefully chosen finite-index subgroup of the mapping
class group (as required by Proposition
\ref{prop:technical_MCG_quotient}.\eqref{item:deep_BBF}), see Lemma
\ref{lem:fi}. The subgroup we use is contained in the one constructed by 
Bestvina--Bromberg--Fujiwara \cite[Section~5]{BBF}. 
We do not think that choosing this subgroup so specifically is
strictly needed, but the strategies we are aware of to get around
using it are significantly more complicated than using the subgroup.

\subsubsection{Structure of the section}
We now explain how the rest of the section is organized.  In
Subsection \ref{subsec:coloring}, we construct the ambient
finite-index subgroup of $MCG(S)$.  Then, after recalling the
definitions of composite projection system and composite rotating
family in Subsection \ref{subsec:CPS_CRF}, we set up the case of the
first quotient in Subsection \ref{subsec:inductive_setup}. In Subsection
\ref{subsubsec:higher_complexity}, we set up the induction in 
Proposition \ref{prop:technical_MCG_quotient}.  
After that, we can describe the
composite projection system and composite rotating family for the
further quotients, in Subsection \ref{subsec:composite_induction}.
From that point on (and \textbf{only} from that point on), the
proofs in the case of Dehn twists and in the case of the further
quotients are largely the same, and are done together, with the
occasional digression where the two cases are treated differently.

\subsubsection{Warning to the reader} Some of the combinatorial properties of quotients of curve graphs that we 
will encounter correspond to topological properties that can be stated in terms of subsurfaces, curves, etc. 
While this might help with intuition and to motivate why they are relevant, we emphasize that very often it 
will not be straightforward at all to relate properties of quotients of curve graphs to topological properties, 
and we will have to rely on the combinatorial HHS viewpoint, taking advantage of topological arguments only 
after lifting to the curve graph. In particular, in all of our statements we can only use combinatorial, rather 
than topological, language.

\subsection{Coloring subsurfaces}\label{subsec:coloring}
Throughout the proof of Theorem~\ref{thm:MCG_quotient}, we will use a
strengthened form of the coloring of the subsurfaces of $S$
constructed in~\cite{BBF}.

\begin{lemma}[Enhanced BBF subgroups]\label{lem:fi}
 There exists a finite coloring $\mathfrak S=\mathfrak S_1\sqcup...\sqcup \mathfrak S_t$ of the collection $\mathfrak S$ of 
all subsurfaces of $S$, so that distinct elements with the same color overlap.
 
 Moreover, for every integer $q>0$ there exists a finite-index torsion-free normal subgroup $H_q$ of $MCG(S)$ so that
 \begin{enumerate}
  \item\label{item:invariant_color} the coloring is $H_q$-invariant, meaning $H_q\mathfrak S_j=\mathfrak S_j$ for all $j$;
  \item\label{item:multi-orbit} for every subsurface $Y$ with at least one component which is not an annulus or a pair of 
pants, and any curve $\gamma$ on $Y$, there exists a curve $\alpha$ on 
$Y$ such that no $H_{q}$--translate of $\alpha$ is disjoint 
from, or isotopic to, $\gamma$;
  \item\label{item:Dehn_powers} for every Dehn twist $\tau\in MCG(S)$ we have $H_q\cap \langle\tau\rangle \leq 
\langle\tau^q\rangle$.
 \end{enumerate}
\end{lemma}

\begin{proof}
Bestvina--Bromberg--Fujiwara constructed a coloring of $\mathfrak S$ 
with the property that distinct elements with the same color overlap 
\cite[Proposition 5.8]{BBF}.  Additionally, the coloring they 
construct has the property that there is a normal, torsion-free finite-index subgroup $H^0$ of $MCG(S)$ so that 
the colors are exactly the $H^0$--orbits of the induced action on 
$\mathfrak S$. Thus, this coloring satisfies item~\eqref{item:invariant_color} as well as the requirements 
before it.

(Strictly speaking, Bestvina-Bromberg-Fujiwara produced a coloring of the set of \emph{connected} subsurfaces 
by $H^0$--orbits with no two subsurfaces of the same color being disjoint.  This extends to all of $\mathfrak 
S$: just color each disconnected subsurface by its $H^0$--orbit.  If $Y$ is disconnected, then $Y$ and $gY$ 
cannot be disjoint, for $g\in H^0$, because then each component of $Y$ would be disjoint from its 
$g$--translate.)

Starting from the BBF coloring, 
we will pass to increasingly deep finite-index subgroups of $MCG(S)$. 
We begin by noting that each of the 
enumerated properties in the statement is stable under passing to further finite-index subgroups.

We now arrange for the coloring to satisfy 
item~\eqref{item:multi-orbit}.  Consider a subsurface $Y$ as in the 
statement, and any curve $\gamma$ on $Y$. Let $g$ be a mapping class that is supported on $Y$ and does not 
stabilize $\gamma$.  By taking $g$ to be, for example, a partial pseudo-Anosov supported on $Y$, we can pass 
to a positive power and assume 
that our $g$ with the preceding property also satisfies $g\in H^0$.

By \cite[Theorem 1.4]{LeiningerMcReynolds}, $\stabilizer(\gamma)$ is separable in $MCG(S)$. Hence $H^0$ has a 
finite-index subgroup $H^{Y,\gamma}$ such that $\stabilizer(\gamma)\cap H^0\leq H^{Y,\gamma}$ but $g\not \in 
H^{Y,\gamma}$.  It follows that $H^{Y,\gamma}g\cap \stabilizer(\gamma)=\emptyset$.

Let $\alpha=g\gamma$.  Then for any $h\in H^{Y,\gamma}$, we have that $h\alpha=hg\gamma$ is in the same 
$H^0$--orbit as $\gamma$.  On the other hand, $h\alpha\neq \gamma$, for otherwise we would have $hg\in 
H^{Y,\gamma}\cap\stabilizer(\gamma)$, which is impossible.  Thus $h\alpha$ and $\gamma$ must intersect, by the 
defining property of $H^0$, applied to annular subsurfaces.

The above paragraph holds for any fixed $\gamma$ and $Y$. Since there 
are only finitely many $H^0$--orbits of pairs $(Y,\gamma)$, we can 
conclude by taking a finite intersection of the $H^{Y,\gamma}$, where $(Y,\gamma)$ varies over 
orbit-representatives. Denote this intersection by $H^1$ and we now 
have a subgroup which satisfies item~\eqref{item:multi-orbit}, as 
desired.

Let us now fix $q$ and arrange item \eqref{item:Dehn_powers}.  
We will construct a finite-index normal subgroup $H^2<MCG(S)$ so that for every Dehn 
twist $\tau\in MCG(S)$, we have $H^2\cap \langle\tau\rangle < 
\langle\tau^q\rangle$. Then, setting $H_q=H^1\cap H^2$, we will have that $H_q$ 
is a finite-index 
subgroup of $H^1$, which ensures that properties \ref{item:invariant_color} 
and \ref{item:multi-orbit} hold as well.

Since there are finitely many conjugacy classes of Dehn twists, it suffices to show that for any given Dehn twist 
$\tau$ there is a normal finite-index subgroup $H^\tau$ so that $H^\tau\cap \langle\tau\rangle < \langle\tau^q\rangle$ (so 
that we can take a finite intersection of the $H^\tau$ for all $\tau$ in a complete list of conjugacy representatives).

Fix a Dehn twist $\tau$.  There are two cases:
\begin{itemize}
     \item Suppose that $\tau$ is a Dehn twist around a non-separating curve.  Let $S=S_{g,p}$, where $g$ is the genus and 
$p$ is the number of punctures.  

If $p=0$, then the usual action of $MCG(S)$ on $H_1(S,\mathbb Z/q\mathbb Z)$ gives a homomorphism $\Psi:MCG(S)\to 
Sp(2g,\mathbb Z/q\mathbb Z)$.  By \cite[Proposition 6.3]{FarbMargalit}, $\Psi(\tau)$ has order $q$.  We 
let $H^\tau=\mathrm{ker}(\Psi)$.

Suppose $p\geq 1$.  Let $S'=S_{g,p-1}$, so that $S$ is obtained from
$S'$ by removing a point $x$.  Let $F:PMCG(S)\to PMCG(S')$ denote the
surjection in the Birman exact sequence, where $PMCG(S)\leq MCG(S)$ is
the finite index subgroup fixing each puncture.  Then $F(\tau)$ is a
Dehn twist $\tau'$ in $S'$ around a non-separating curve.  By
induction, $MCG(S')$ has a finite-index subgroup $H^{\tau'}$ such that
$H^{\tau'}\cap\langle\tau'\rangle<\langle (\tau')^q\rangle$.  The
subgroup $H^{\tau}=F^{-1}(H^{\tau'})$, which has finite index in
$PMCG(S)$ and hence in $MCG(S)$, has the desired
property.

\item Now suppose $\tau$ is a Dehn twist around a separating curve.    It suffices to consider the case that $q$ is a power 
of a prime $p$, since in general we can take the intersection of the 
finite-index subgroups dictated by the prime factorization of $q$. By \cite[Theorem 1.2]{Paris:residual}, $\tau$ lies in a 
finite-index normal subgroup $H^3$ of $MCG(S)$ which is residually $p$. In particular, 
there is a finite quotient $H^3/N$ of $H^3$ so that $\tau$ maps to an element of order $q'$ with $q|q'$, and $[H^3:N]$ is a 
power of $p$. Since the intersection of finitely many subgroups of index a power of $p$ also has index a power of $p$, we can 
take $N$ to be characteristic in $H^3$, whence normal in $MCG(S)$, and set $H^\tau=N$.    
\end{itemize}
This completes the proof of the lemma.
\end{proof}

\subsection{Composite projection systems and composite rotating families}\label{subsec:CPS_CRF}
Now we recall two definitions from~\cite{D_PlateSpinning} that we will need below.  Specifically, the following combines 
\cite[Definitions 1.1-1.2]{D_PlateSpinning}. The reader can keep in mind as $\bbY_*$ the family of curves on a surface, which can be split into finitely many families with the property that curves in each family pairwise intersect by a result of \cite{BBF}. Also, $\Act(Y)$ is the set of curves that intersect the curve $Y$, and $d_Y$ denotes the distance between subsurface projections. (This was the motivating example for this notion in \cite{D_PlateSpinning}.)

\begin{defn}[Composite projection system]\label{defn:CPS}
        Let $\bbY_*$ be a countable set equipped with a finite coloring $\bbY_* =\sqcup_{j=1}^m \bbY_j$.  For each 
$Y\in\bbY_*$, let $j(Y)$ denote the value $j$ for which $Y\in\bbY_j$.

		 A \emph{composite projection system} on a countable set
		 $\bbY_*$ is the data consisting of: a constant $\theta>0$; a 
		 family of subsets, one 
		 for each $Y\in\bbY_*$ denoted $\Act(Y)\subset \bbY_*$ (called the
		 \emph{active set} for $Y$) such that $\bbY_{j(Y)}
		 \subset\Act(Y)$; and a family of functions $d_Y \co (
		 \Act(Y)\setminus \{Y\} \times \Act(Y) \setminus \{Y\}) \to
		 \bbR_+$, satisfying the following whenever all quantities are
		 defined:
         
        \begin{enumerate}
          \item[\textbf{(CPS1)}]  \label{item:cps_symmetry}     (symmetry) $d_Y (X,Z) = d_Y (Z,X)$ for all 
$X,Y,Z$;
          \item[\textbf{(CPS2)}] (triangle inequality)
       $d_Y(X,Z) + d_Y(Z,W) \geq d_Y(X,W) $ for all $X,Y,Z, W$; \label{item:cps_triangle} 
         \item[\textbf{(CPS3)}] (Behrstock
       inequality) $\min\{ d_Y(X,Z), d_Z (X,Y) \} \leq \theta
       $ for all $X,Y,Z$;\label{item:cps_behrstock} 
          \item[\textbf{(CPS4)}] (properness) $\{Y\in \bbY_j, d_Y (X,Z) >\theta\}$ is
            finite for all $X,Z$; \label{item:cps_properness} 
          \item[\textbf{(CPS5)}] (separation) $d_Y (Z,Z) \leq \theta$ for all 
$Z,Y$.\label{item:cps_separation} 
       \end{enumerate}
         The map $\Act$ is required to satisfy three further properties:  
         \begin{enumerate}
           \item[\textbf{(CPS6)}] \label{item:cps_symmact}(symmetry in action) $X\in \Act(Y)$ if and
         only if $Y\in \Act(X)$;
           \item[\textbf{(CPS7)}] \label{item:cps_closeact}(closeness in inaction) if   $X\notin \Act(Z)$, for
         all $Y \in \Act(X) \cap \Act(Z)$, $ d_Y( X,Z )\leq \theta
         $; 
         \item[\textbf{(CPS8)}] \label{item:cps_finitefill}(finite filling)  for all $\calZ\subset \bbY_*$, 
there
           is a finite collection 
            of elements $X_j$ in $\calZ$ such that $\cup_j
            \Act(X_j)$ covers  $\cup_{ X\in \calZ} \Act(X)$.
         \end{enumerate}
An \emph{automorphism} of a composite projection system is a 
bijection $g\co\bbY_*\to\bbY_*$ such that:
\begin{itemize}
     \item $g$ preserves each $\bbY_j$;
     \item for all $Y\in\bbY_*$, we have $\Act(gY)=g(\Act(Y))$;
     \item for all $Y$ and all $X,Z\in\Act(Y)$, we have $d_{g(Y)}(g(X),g(Z))=d_Y(X,Z)$.
\end{itemize}
\end{defn}

The following is a variant of the notion introduced in \cite[Definition 2.1]{D_PlateSpinning}.

\begin{defn}[$(\Theta_{P},\Theta_{Rot})$--Composite rotating 
	family]\label{defn:ThetaCRF}
         Let $\Theta_{P}$ and $\Theta_{Rot}$ be constants. 
		 A \emph{$(\Theta_{P},\Theta_{Rot})$--composite rotating family} on a composite projection system
    endowed with an action of a group $G$  by automorphisms    
    is a  family of subgroups $\Gamma_Y,
    Y\in \bbY_*$ such that   
    \begin{enumerate}
      \item[\textbf{(CRF1)}] \label{item:crf_inf}for all $X\in \bbY_*, \Gamma_X < G_X= {\rm
     Stab}_G(X)$, is an infinite group of rotations around $X$, with proper isotropy, meaning that for all $R>0$ and $Y\in \Act(X)$ the set $F_Y^X(R)=\{\gamma\in \Gamma_X: d_X(Y,\gamma Y)<R\}$ is finite.
      \item[\textbf{(CRF2)}] \label{item:crf_equi}for all $g\in G$, and all $X\in \bbY_*$ , one
    has  $\Gamma_{gX}= g\Gamma_X g^{-1}$,
    \item[\textbf{(CRF3)}] \label{item:crf_commute}if $X\notin \Act(Z)$ then $\Gamma_X$ and $\Gamma_Z$ 
commute, 
      \item[\textbf{(CRF4)}] \label{item:crf_rot}for all $j$, for all $X,Y, Z \in \bbY_j$, if $d_Y(X,Z)
        \leq \Theta_P $
        then for all
        $g\in \Gamma_Y\setminus\{1\}$,  $d_Y(X,gZ) \geq \Theta_{Rot}$. 
    \end{enumerate}
\end{defn}

	\begin{rem}\label{rem:comrotfamilyconstants} Dahmani's original
	definition of a composite rotating family, \cite[Definition
	2.1]{D_PlateSpinning}, doesn't use the metric $d_{Y}$ from
	the composite projection system, but rather a perturbation 
	which differs from $d_{Y}$ by a bounded amount.
	
	Our definition above relies on two constants 
	$\Theta_{P}$ and $\Theta_{Rot}$. In \cite[Definition 
	2.1]{D_PlateSpinning}, instead the constants  $\Theta_{P}$ and 
	$\Theta_{Rot}$ are fixed depending only on the constant $\theta$ 
	from the composite projection system. For suitably chosen values 
	of $\Theta_{P}$ and $\Theta_{Rot}$ a 
	$(\Theta_{P},\Theta_{Rot})$--composite rotating family is a 
	composite rotating family in the sense of Dahmani's 
	\cite[Definition 2.1]{D_PlateSpinning}. 
	
	Specifically, in the notation of 
	\cite[\S~1.2.1]{D_PlateSpinning} we can choose 
	$\Theta_{P}=c_{*}+21m\kappa+\kappa$ and 
	$\Theta_{Rot}>2c_{*}+2\Theta_{P}+20(\kappa+\Theta)+\kappa$, where 
	all these constants are functions of 
	the constant $\theta$ 
	from the composite projection system. These two 
	relations differ from the analogous ones in 
	\cite[\S~1.2.1]{D_PlateSpinning} by an additional term of 
	$\kappa$ to take into account the perturbation of the metric by 
	at most $\kappa$, and $\kappa\geq\theta$.
\end{rem}

Accordingly, with a slight abuse of terminology we define:
\begin{defn}[Composite rotating family]\label{defn:CRF} 
	A \emph{composite rotating family} is a 
	$(\Theta_{P},\Theta_{Rot})$--composite rotating family with 
	constants as in Remark~\ref{rem:comrotfamilyconstants}.
\end{defn}

\begin{rem*}
    The definition of proper isotropy given above is taken from \cite[Definition 2.3]{DFDT}, and it is weaker than the corresponding definition from \cite[\S 1.2.2]{D_PlateSpinning}, which requires that the finite set from Definition~\ref{defn:ThetaCRF}.(CRF1) is independent of $Y\in \Act(X)$, that is, $F_Y^X(R)=F^X(R)$. We now explain why, despite this, all results from \cite{D_PlateSpinning} still apply to a $(\Theta_{P},\Theta_{Rot})$--composite rotating family defined as above (that is, using proper isotropy from \cite{DFDT}).

    \begin{itemize}
        \item First of all, proper isotropy is only needed in \cite{D_PlateSpinning} to the extent that it is needed to prove \cite[Lemma 1.4]{D_PlateSpinning}; there are no other uses of $F^X(R)$ (simply denoted $F(R)$ in \cite{D_PlateSpinning}). In said lemma, the only relevant value of $R$ is $R=10\kappa$ for $\kappa$ as in Remark \ref{rem:comrotfamilyconstants}. Hence, it suffices to show that for a $(\Theta_{P},\Theta_{Rot})$--composite rotating family, we can take $F^X(10\kappa)=\{1\}$.

        \item Spelling out the above, we have to show that for all $Y\in \Act(X)$ and $\gamma\in \Gamma_X-\{1\}$ we have $d_X(Y,\gamma Y)>10\kappa$; this would follow directly from Definition~\ref{defn:ThetaCRF}.(CRF4) except that we have to consider $Y\in \Act(X)$ but not necessarily in $\bbY_{j(X)}$. However, given such $Y$, we use that $\Gamma_Y$ is infinite, as stipulated in Definition~\ref{defn:ThetaCRF}.(CRF1), to produce $\alpha\in\Gamma_Y-F^Y_X(10\kappa)$.  By definition, $d_Y(X,\alpha X)\geq10\kappa>\theta$.  Hence, by Definition~\ref{defn:CPS}.(CPS3), we have $d_X(Y,\alpha X)\leq\theta\leq\Theta_P$.  Since $G$ preserves each $\bbY_i$, we have $X,\alpha X\in \bbY_{j(X)}$.

        \item Now let $\gamma\in\Gamma_X-\{1\}$.  Then by Definition~\ref{defn:ThetaCRF}.(CRF4), we have $d_X(\alpha X,\gamma\alpha X)\geq\Theta_{rot}$.  By equivariance, we have $d_X(\gamma\alpha X,\gamma Y)\leq \theta$, so by Definition~\ref{defn:ThetaCRF}.(CPS2), we get $d_X(Y,\gamma Y)\geq\Theta_{rot}-2\theta>10\kappa$, in view of our choice of $\Theta_{rot}$.
    \end{itemize}
This shows that a $(\Theta_{P},\Theta_{Rot})$--composite rotating family defined as above satisfies all of the statements about composite rotating families established in \cite{D_PlateSpinning} and in \cite{DFDT}, and hence we can freely apply both sets of results below.  It also explains why the results from \cite{D_PlateSpinning} used in \cite{DFDT} apply in the setting of the latter paper, where $\Theta_{rot}$ is always assumed to be sufficiently large.
\end{rem*}

\subsection{Setting up the induction}\label{subsec:inductive_setup}
Theorem~\ref{thm:MCG_quotient} is proven by induction on $i$.  

\subsubsection{The annular case $i=-1$}\label{subsubsec:annular_setup}  
The base case, where $i=-1$, will be verified almost 
identically to the inductive step, but the notation used in the proof has a slightly different meaning in the base case.  

\begin{notation}\label{not:annuli}
When $i=-1$, we will use the following notation:
  \begin{itemize}
  \item[]
  \item $N_{-2}=\{1\}$,
  \item $X=\mathcal C(S)$,
  \item $\bbY_*=(\mathcal C(S))^{(0)}$
  
  \item $d_Y$ denotes the distance in the annular curve graph $\mathcal C(Y)$ of $Y$.
 \end{itemize}
\end{notation}

\begin{rem} The collection $\bbY_*$ above provides a composite 
	projection system as shown in \cite[\S 3]{D_PlateSpinning}, where 
	for all $W,Y,Z\in\bbY_{*}$ the distance between $Y$ and $Z$ as 
	measured in $W$ is $d_{W}(\pi_W(Y), \pi_W(Z))$, where $\pi_W$ is the (annular) subsurface projection.
\end{rem}

In the case $i=-1$, we have the following:

\begin{lemma}[Composite rotating family, annular case]\label{lem:i=-2_rotating_family}
For every $\theta>0$ there exists $\theta_0>0$ so that the following holds.  

Let $\tau_1,\dots,\tau_k$ be a complete list of 
conjugacy representatives of Dehn twists in $MCG(S)$. For each $j\leq k$, let $\theta_j$ be a positive multiple of 
$\theta_0$. For 
$Y\in\bbY_*$, let $\Gamma^\theta_Y=\langle\tau_Y^{\theta_{j(Y)}}\rangle$, where $\tau_Y$ is the Dehn twist around $Y$, and 
$j(Y)$ has the property that $\tau_Y$ is conjugate to $\tau_{j(Y)}$.

Then the subgroups $\Gamma^\theta_Y$ form a composite rotating family on the composite projection system $\bbY_*$. Moreover,
$$\min\{d_{\mathcal C(Y)}(x,\gamma x): Y\in\bbY_*,\gamma\in\Gamma_Y^{\theta}-\{1\},x\in\mathcal C(Y)\}>\theta.$$

Finally, set $N= \langle\langle \Gamma^\theta_Y\rangle\rangle$. If $\Sigma$ is a simplex of $\mathcal C(S)$ so that the complement of the multicurve $\Sigma^{(0)}$ has one 
complexity-1 component $Y$ while all others are pairs of pants, then $\stabilizer(\Sigma)/(N\cap \stabilizer(\Sigma))$ is an 
infinite hyperbolic group acting with finite point-stabilizers on $\link(\Sigma)/ (N\cap 
\stabilizer(\Sigma))$.
\end{lemma}

\begin{rem}\label{rem:suppress}
 We remind the reader of the convention for (combinatorial) hierarchically hyperbolic spaces: when writing 
distances in a hyperbolic space/link, we suppress the $\pi_\bullet$ notation for projection maps.  So, e.g. 
$\dist_{\mathcal C(Y)}(x,y)$ means $\dist_{\mathcal C(Y)}(\pi_Y(x),\pi_Y(y))$.  See 
Notation~\ref{notation:suppress_pi}.
\end{rem}

\begin{proof}
Except the last conclusion, the proof is identical to the one give by 
Dahmani in \cite[Section 3]{D_PlateSpinning}, with two minor changes: 
first, while in Dahmani's case all the Dehn twists are raised to the 
same power, ours are allowed to vary; and, second our 
coloring was chosen in a more specific way.

To verify the last conclusion, we may thus apply the versions of 
results from \cite{DFDT} which allow for variable powers of Dehn 
twists and our particular coloring, as noted above.

For the last conclusion, with a suitable choice of powers of Dehn
twists we can apply \cite[Proposition 5.8]{DFDT}.  This says that
$N\cap \stabilizer(\Sigma)$ is generated by the powers of Dehn twists
supported in $Y$ and those supported around the curves of
$\Sigma^{(0)}$.  More precisely, \cite[Proposition 5.8]{DFDT} applies
to stabilizers of vertices in the curve graph, but we can apply it
inductively passing to subsurfaces.  But then
$\stabilizer(\Sigma)/(N\cap \stabilizer(\Sigma))$ is virtually a
quotient $\hat{G}$ of $MCG(Y)$ by powers of Dehn twists, which is a
hyperbolic group by \cite[Theorem
6.8.(1)]{DFDT}.

We now argue that the action of $\stabilizer(\Sigma)/(N\cap
\stabilizer(\Sigma))$ on $\link(\Sigma)/ (N\cap \stabilizer(\Sigma))$
has finite stabilizers. These stabilizers are quotients of
stabilizers of vertices of $\link(\Sigma)$ since every element of the stabilizer is the image in the quotient by an element in the stabilizer of an orbit, but this representative can be multiplied by an element of $N\cap \stabilizer(\Sigma)$ to ensure that it lies in a vertex stabilizer. Quotients of stabilizers of vertices are virtually
generated by commuting Dehn twists, and therefore become finite after
quotienting by $N\cap \stabilizer(\Sigma)$.
\end{proof}

The following lemma is crucial to show that we can lift various objects from quotients of curve graphs to curve graphs. Roughly, it says that for $N$ the normal subgroup generated by large powers of Dehn twists, two vertices being in the same $N$-orbit can be witnessed by a large annular projection that can be ``shortened'' with a suitable power of a Dehn twist which lies in $N$. This kind of ``shortening'' lemmas are often referred to as ``Greendlinger lemmas'' in analogy with the the Greendlinger lemma for classical small cancellation groups which allows one to shorten a word representing the identity using a relator.

\begin{lemma}[``Greendlinger lemma'', annular case]\label{lem:short_less_complex_base} 
There exists a diverging function $\frak T$ so that the following holds for $\theta>0$.  Let 
$N=N_{-1}=\langle\langle 
\{\Gamma^\theta_Y\}\rangle\rangle$.  Then there is a well-ordered set $\mathfrak C$, and an assignment 
$\gamma\in N\mapsto c(\gamma)\in \mathfrak C$, with $c(1)$ 
the minimal element of $\mathfrak C$. 

Moreover, for all $\gamma\in 
N-\{1\}$ and all simplices $\Delta$ of $X$, there is $Y\in\bbY_*$ and $\gamma_Y\in\Gamma_Y^\theta$ so that 
$c(\gamma_Y \gamma)<c(\gamma)$ and either
  \begin{itemize}
   \item $\Delta\subseteq \fix(\Gamma_Y^\theta)$, or $\gamma \Delta\subseteq \fix(\Gamma_Y^\theta)$, or
   \item $d_Y(\Delta,\gamma \Delta)>\frak T(\theta)$ (and the quantity is defined).
  \end{itemize}
 \end{lemma}

\begin{proof}
	We note that the proof of this result 
	is a minor variation of the proof of \cite[Corollary 4.6]{DFDT}. 
	The key difference is that here we must consider simplices 
	rather than just single vertices.
	
	The proof in the present annular case is identical, verbatim, to 
	the proof in the general case  
	(Lemma~\ref{lem:short_less_complex}), so we postpone the proof 
	until then. 
\end{proof}

Lemma \ref{lem:i=-2_rotating_family} will be used in combination with Lemma~\ref{lem:fi} to choose a finite-index subgroup $H$ of $MCG(S)$, 
and an $H$--invariant coloring of the subsurfaces, as follows.

\begin{notation}\label{not:correct_powers}
Below $MCG(S)$ is given its usual HHS structure as in \cite[Section 
11]{HHS_II}. Let $C$ be the constant from
Definition~\ref{defn:HHS}.(\ref{item:dfs:bounded_geodesic_image})(bounded
geodesic image) for $MCG(S)$, see \cite[Section 11]{HHS_II} or 
\cite[Theorem~3.1]{MM_II}. Let $\kappa$ be given by Lemma
\ref{lem:conv_cocpt=cobounded} applied to the subgroup $Q<MCG(S)$ from
Theorem \ref{thm:MCG_quotient}.  Fix some $\theta>0$ such that:

 \begin{itemize}
   \item $\frak T(\theta)>( 6\genus(S)+2p(S)-3)C$,
    \item $\frak T(\theta)> \max\{d_{\mathcal C(Y)}(fx,g x):Y\in\bbY_*, f,g\in F, x\in X\}$,
    \item $\frak T(\theta)>\kappa+C$.
 \end{itemize}
Let $\theta_0$ be as in Lemma 
\ref{lem:i=-2_rotating_family}, for the given $\theta$.  Fix from now on $H=H_{\theta_0}$ as in Lemma \ref{lem:fi}. For $Y\in\bbY_*$, 
let $\Gamma^\theta_Y=H\cap\langle\tau_Y\rangle<\langle\tau_Y^{\theta_0}\rangle$; notice that Lemma 
\ref{lem:i=-2_rotating_family} 
applies to $\Gamma^\theta_Y$ by the construction of $H$. For each $Y$, let $\Gamma_Y=\Gamma_Y^\theta$.  Set 
$N=N_{-1}=\langle\langle \{\Gamma_Y\}\rangle\rangle$.

At this point, we also fix an $H$--invariant coloring $\mathfrak S=\bigsqcup_j\mathfrak S_j$ as in Lemma~\ref{lem:fi}, 
which will remain the same at all stages of the induction.
\end{notation}

We emphasize that $H$ will remain the same at all subsequent stages of
the induction, even though we have thus far only defined $\Gamma^\theta_Y$
in cases where $Y$ is an annulus.

\subsection{Inductive conditions}\label{subsubsec:higher_complexity}
We now set the notation for $i\geq -1$.  Suppose we have constructed $\bar G_i=MCG(S)/N_i$ for a given $F$ and 
 $Q$, at all 
complexities up to $i$.  We will make further assumptions for the inductive step. To state these 
 we need the following definition.
 
 \begin{defn}[Generalized $m$--gon]\label{defn:generalized_m_gon}
     A \emph{generalized $m$--gon} in a simplicial graph is a sequence $\tau_0,\dots,\tau_{m-1}$ so that:
     \begin{itemize}
 \item Each $\tau_j$ is either a simplex (\textbf{type S}), together with non-empty sub-simplices $\tau_j^\pm$, or a 
geodesic in $\link(\Delta_j)$ for some (possibly empty) simplex $\Delta_j$ (\textbf{type G}) with 
endpoints $\tau^\pm_j$.
 \item $\tau^+_j=\tau^-_{j+1}$ (indices are taken modulo $m$). 
\end{itemize}
(The second bullet implies that $\tau_j\cap\tau_{j+1}$ is non-empty.)
\end{defn}

The main inductive statement is the following proposition, which is a more precise version of Theorem \ref{thm:MCG_quotient}; many of the additional points are required to inductively obtain composite projection systems. We denote pointwise stabilizers by $\pstab$. Also, given a simplicial map $q:Y\to Z$ of simplicial complexes, a \emph{lift} of an ordered simplex $\Delta$ of $Z$ (that is, a simplex with an ordering on its vertices) is an ordered simplex $\Sigma$ of $Y$ so that $q(\Sigma)=\Delta$, and the map is order-preserving at the level of vertices. We will denote ordering on vertices on the subscripts, e.g., $(v_0,\dots,v_k)$ if the vertex set of the simplex consists of the $v_j$. A lift of a simplex is a lift of the simplex with any order on its vertices. A \emph{lift} of a generalized $m$-gon $\tau=\tau_0,\dots,\tau_{m-1}$ in $Z$ is an $m$-gon $\tau'_0,\dots,\tau'_{m-1}$ in $Y$ with $q(\tau'_j)=\tau_j$, $\tau'_j$ is \textbf{type S}/\textbf{type G} if and only if $\tau_j$ is, 
and if $\tau_j$ is a 
geodesic in $\link(\Delta_j)$ then $\tau'_j$ is a geodesic in $\link(\Sigma_j)$ for some lift $\Sigma_j$ of $\Delta_j$.

 \begin{prop}\label{prop:technical_MCG_quotient}
Let $S$, $F$, and $Q$ be as in Theorem \ref{thm:MCG_quotient}, with
$S$ having genus $\genus(S)$ with $p(S)$ punctures.  For $-1\leq i\leq
3\genus(S)+p(S)-4$, there exists a quotient $\phi:MCG(S)\to
MCG(S)/N_i=\bar{G}_i$ such that properties Theorem
\ref{thm:MCG_quotient}.\eqref{item:larg_inj_radius}-\eqref{item:poc}
hold.  Moreover, the following additional properties also hold, where $q:\mathcal C(S)\to \mathcal C(S)/N_i$ is the quotient map:
 
  \begin{enumerate}[(I)]
  \setcounter{enumi}{3}
  \item\label{item:deep_BBF} $N_i<H$, where $H<MCG(S)$ is as in Notation \ref{not:correct_powers}, and $N_{-1}<N_i$.
 
  \item\label{item:hyperelliptic} For all distinct $f,g\in F$ either $f^{-1}g$ has finite order or there exists a vertex $x$ 
of $\mathcal C(S)/N_i$ so that $f(x)\neq g(x)$.
 \item\label{item:lift_to_CS} For $m\leq 
 \max\{4,6\genus(S)+2p(S)-3\}$, any generalized $m$-gon in $\mathcal 
C(S)/N_i$ 
can be lifted to $\mathcal C(S)$.
  \item\label{item:lift_splx_induction} For every ordered simplex $\Delta$ of $C(S)/N_i$ there is a unique $N_i$-orbit of 
lifts $\lift\Delta$ in $\mathcal C(S)$,  and for 
any such lift we have $q(\sat(\lift\Delta))=\sat(\Delta)$. 
\item \label{item:BGI} There exists $C_i$ with the following property.  Let $\Delta$ be a simplex of $\mathcal C(S)/N_i$, and let $v_0,\dots,v_k$ be a geodesic of 
$\link(\Delta)$.  Suppose that for some simplex $\Sigma$ of $\mathcal C(S)/N_i$ we have that $\dist_{\mathcal C(\Sigma)}(v_0,v_k)$ is defined and at least $C_i$. Then 
there exists $i$ so that $\Delta^{(0)}\cup\{v_i\}$ is contained in $Sat(\Sigma)$.
  \item\label{item:bottom_hyp} If an element $b([\Delta])\in \mathfrak S^{\geq 1}/N_i$, for $[\Delta]\in\mathfrak S_{N_i}$, has a representative of 
complexity $i+1$,  
then $\pstab(\sat(\Delta))$ is hyperbolic and acts properly and cocompactly on $\mathcal 
C(\Delta)$.
 \end{enumerate}
 \end{prop}

 \begin{rem}\label{rem:we-used-propaganda}
 The combinatorial HHS $(\mathcal C(S)/N_i,W)$ from Theorem~\ref{thm:MCG_quotient}.\eqref{item:HH_structure} is obtained by applying Theorem~\ref{thm:propaganda} to the action of $\bar G_i$ on $\mathcal C(S)/N_i$; see Section~\ref{subsec:main_theorem}.  
 \end{rem}
 
\subsubsection{Guide to the proof of Proposition \ref{prop:technical_MCG_quotient} (hence Theorem \ref{thm:MCG_quotient})}\label{subsec:guide_to_proof}
 \begin{conv}\label{conv:induction}
  From now and until the end of the Section we assume that either:
  \begin{itemize}
   \item $i=-1$, with the notation from Notation \ref{not:annuli},
   \item $i>-1$, and Proposition holds with $i$ replaced by $i-1$.
  \end{itemize}
 \end{conv}
 
 \begin{rem}\label{rem:0step} For $i=0$, the composite projection system is empty and 
	 thus the quotient we are taking is trivial in the sense that we 
	 are just quotienting by the trivial subgroup. The reason this is 
	 empty boils down to the fact that the complexity $0$ subsurfaces 
	 are thrice punctured spheres and thus for each the curve graph 
	 is empty.
	 Thus, since the composite projection system is empty, statements 
	 involving the set $\bbY_{*}$ are all vacuously true.  	
 \end{rem}

 Here is a list of where properties \eqref{item:larg_inj_radius}--\eqref{item:lift_splx_induction} are verified:

 \begin{itemize}
	 \item Items \eqref{item:larg_inj_radius} and \eqref{item:hyperelliptic} are proven together in Subsection \ref{subsec:preserve}.

	   \item Item \eqref{item:HH_structure} is the content of Subsection \ref{subsec:main_theorem}.

	     \item Item \eqref{item:poc} is also proven in Subsection \ref{subsec:preserve}.

  \item Item \eqref{item:deep_BBF} is simply a restatement of the 
  assumptions we made in Notation \ref{not:annuli}, for $i=-1$, and Notation
  \ref{not:not_annuli}, for $i>-1$.
  
  \item Item \eqref{item:lift_to_CS} and the first part of Item \eqref{item:lift_splx_induction} hold by Proposition 
\ref{prop:omni_lift}, with the second part of Item \eqref{item:lift_splx_induction} being Lemma \ref{lem:quotient_sats}.
 \end{itemize}

 The remaining properties follow from properties \eqref{item:larg_inj_radius}--\eqref{item:lift_splx_induction}, as shown in 
the following lemmas.
 
 \begin{lemma}\label{lem:BGI_quotient}
  Assume properties
  \eqref{item:larg_inj_radius}--\eqref{item:lift_splx_induction} hold for
  our current $i$.  Then property~\eqref{item:BGI} holds for our current
  $i$.
 \end{lemma}

\begin{proof}
 By property \eqref{item:HH_structure}, we are working in a combinatorial HHS with underlying simplicial complex $\mathcal 
C(S)/N_i$.
 
If some vertex of $\Delta^{(0)}$ was not in $Sat(\Sigma)$, then there would a path of length $2$ in $Y_\Sigma$ from $v_0$ to 
$v_k$ and thus $\dist_{\mathcal C(\Sigma)}(v_0,v_k)$ would be  
uniformly bounded. This is a contradiction when $C$ is 
sufficiently large. Thus all vertices in $\Delta^{(0)}$ are contained 
in $Sat(\Sigma)$, which yields that $[\Sigma]\nest[\Delta]$.

Using the 
combinatorial HHS structure of $\link(\Delta)$ (Propositions \ref{prop:compatibility_1} and \ref{prop:hieromorphism}) we have 
that that the geodesic must intersect $Sat(\Sigma)$ by Lemma \ref{lem:super_bgi}, as required.
\end{proof}

\begin{lemma}\label{lem:fix_sat_hyperbolic}
Assume properties 
\eqref{item:larg_inj_radius}--\eqref{item:lift_splx_induction} for 
our current $i$.  Then property \eqref{item:bottom_hyp} holds for our current $i$.
\end{lemma}

\begin{proof}
Recall that $\pstab$ denotes pointwise stabilizers.

 We use the combinatorial HHS $(\mathcal C(S)/N_i,W)$, which exists by item \eqref{item:HH_structure}. We will check that 
$\pstab(\sat(\Delta))$ is hierarchically hyperbolic using the technology of Section \ref{sec:hieromorphisms}.
 
Recall the graph $W^\Delta$ defined in Definition \ref{defn:induced_link_graph}, which has vertex set the set of maximal simplices of $\link(\Delta)$, and has the property that $(\link(\Delta),W^\Delta)$ is a combinatorial HHS. It is readily seen that 
$\pstab(\sat(\Delta))$ acts on $W^\Delta$. We now 
check that this action is proper and cocompact.\\
 
\textbf{Properness:}  Let $\Sigma$ be a vertex of $W^\Delta$, that is, a maximal simplex of $\link(\Delta)$.  
Then 
$\stabilizer_{\pstab(\sat(\Delta))}(\Sigma)\leq\stabilizer_{\bar G_i}(\Delta\star\Sigma)$, and the latter group is finite, since 
$\bar G_i$ acts properly on $W$ again by item \eqref{item:HH_structure}. Hence $\pstab(\sat(\Delta))$ acts properly on $W^\Delta$.\\

\textbf{Cocompactness for $i=-1$:} We check that there are finitely many $\pstab(\sat(\Delta))$--orbits of edges in $W^\Delta$. The proof is different depending on $i$.

For $i=-1$ (that is, $\bar G_i$ is a quotient of $MCG(S)$ by powers of Dehn twists), no simplex of $\mathcal C(S)/N_i$ corresponds (via the map $b$ from Definition~\ref{defn:b}) to a subsurface of complexity $i+1$, since $\mathfrak S^{\geq 1}$ does not contain any complexity--$0$ surfaces.  So the lemma holds vacuously for $i=-1$.\\

\textbf{Cocompactness for $i=0$:} The argument for cocompactness when $i=0$ is more complicated. In this case, Remark~\ref{rem:0step} says that $\bar G_0=\bar G_{-1}$ and $\mathcal C(S)/N_0=\mathcal C(S)/N_{-1}$.  Our simplex $\Delta$ has the property that $b([\Delta])$ is an orbit of complexity--$1$ subsurfaces.

We will use the following claim about $\bar G_{-1}=\bar G_0$ and $N_{-1}=N_0$ later as well, so we state it separately, together with the assumptions we need to prove it.

\begin{claim}\label{claim:finite-stabs-stab-link}
    Assume that for $i=-1$ the following items from Proposition~\ref{prop:technical_MCG_quotient} hold:
    \begin{itemize}
        \item item \eqref{item:deep_BBF};
        \item item \eqref{item:lift_to_CS};
        \item the part of \eqref{item:lift_splx_induction} about uniqueness of $N_{-1}$-orbits.
    \end{itemize}
    
    Let $\Delta$ be an almost-maximal simplex of $\mathcal C(S)/N_{-1}$, and let $\Gamma_\Delta$ be the quotient of $\stabilizer(\link(\Delta))$ by the kernel of the action on $\link(\Delta)$. Then $\Gamma_\Delta$ has finite vertex-stabilizers in $\link(\Delta)$.
\end{claim}
\renewcommand{\qedsymbol}{$\blacksquare$}
\begin{proof}
Fix a vertex $v$ of $\link(\Delta)$.  The simplex $\Delta\star v$ is maximal because of the assumption that $\Delta$ is almost-maximal.  By the third bullet point, we can fix, once and for all, lifts $\lift\Delta$ and $\lift v$ of $\Delta$ and $v$ such that $\mathcal C(S)$ contains the simplex $\lift\Delta\star\lift v$.

\textbf{The subsurface $Y$:}  The $0$--skeleton of our fixed lift $\lift\Delta$ is a multicurve in $S$ whose complement is a subsurface we denote by $Y'$.  By assumption, $\Delta\star v$ is a maximal simplex of $\mathcal C(S)/N_{-1}$, so the simplex $\lift\Delta\star\lift v$ of $\mathcal C(S)$ is maximal.  Indeed, if it was not maximal, we could extend it to a bigger simplex, and distinct vertices in a common simplex of $\mathcal C(S)$ are not in the same $N_{-1}$--orbit, by the first bullet point, resulting in a simplex in $\mathcal C(S)/N_{-1}$ strictly containing $\Delta\star v$. Hence the vertex set of $\lift \Delta\star\lift v$ is a pants decomposition of $S$.  Therefore, $\{\lift v\}$ is a pants decomposition of a complexity--$1$ component $Y$ of $Y'$.  The curves in $\partial Y$ belong to $\lift\Delta$, and the remaining curves of $\lift\Delta$ form a pants decomposition of the complement of $Y$.

\textbf{Keeping track of groups:}  Recall that $\phi: MCG(S)\to \bar G_{-1}$ denotes the quotient map.  Let $\omega:\stabilizer(\link(\Delta))\to Sym(\link(\Delta))$ be the action on $\link(\Delta)$, and we identify $\Gamma_\Delta$ with $\image(\omega)$.   

We now argue that we have
$$q(\link(\lift\Delta))=\link(\Delta).$$

Indeed, given a vertex $w$ of $\link(\Delta)$, we can lift $\Delta\star w$ to a 
simplex of the form $\lift\Delta\star\lift w$, and hence $q(\lift w)=w$. This shows the containment ``$\supseteq$'', the other containment follows from a use of the first bullet point similar to the argument above.

Since $q(\link(\lift\Delta))=\link(\Delta)$, the homomorphism $\phi$ restricts on $\stabilizer_{MCG(S)}(\link(\lift\Delta))$ to a homomorphism to $\stabilizer(\link(\Delta))$. 

Next, observe that $\link(\lift\Delta)$ corresponds to the set of curves in $Y$, and hence $\stabilizer_{MCG(S)}(\link(\lift\Delta))$ is contained in $\stabilizer_{MCG(S)}(\partial Y)$.  Let $\Lambda_0$ be the kernel of the action of $\stabilizer_{MCG(S)}(\link(\lift\Delta))$ on the set of boundary curves of $Y$, so that $\Lambda_0$ has finite index in $\stabilizer_{MCG(S)}(\partial Y)$.  We summarize this with the following diagram,
$$\Lambda_0\hookrightarrow\stabilizer_{MCG(S)}(\link(\lift\Delta))\stackrel{\phi}{\longrightarrow}\stabilizer(\link(\Delta))\stackrel{\omega}{\longrightarrow}\Gamma_\Delta,$$
where the rightmost arrow is surjective and the leftmost arrow has finite-index normal image.

\textbf{Reducing to a claim about lifting elements:}  We need to show that $\omega(\stabilizer(\link(\Delta))\cap\stabilizer(v))$ is finite.  To do this, we will show that for each $g\in\stabilizer(\link(\Delta))\cap\stabilizer(v)$, there exists $\lift g\in\stabilizer_{MCG(S)}(\link(\lift\Delta))\cap\stabilizer_{MCG(S)}(\lift v)$ such that $\phi(\lift g)=g$.  We first explain why the latter statement suffices to prove the former.

First, $\Lambda_0\cap\stabilizer_{MCG(S)}(\lift v)$ has finite index in $\stabilizer_{MCG(S)}(\link(\lift\Delta))\cap\stabilizer_{MCG(S)}(\lift v)$, so we can fix (independently of $g$) finitely many elements $\lift g_1,\ldots,\lift g_k\in\stabilizer_{MCG(S)}(\link(\lift\Delta))\cap\stabilizer_{MCG(S)}(\lift v)$ representing all cosets of $\Lambda_0\cap\stabilizer_{MCG(S)}(\lift v)$.  

So, for some $j$, we have $\lift g=\lift g_j\lift h$, where $\lift h\in\Lambda_0\cap\stabilizer_{MCG(S)}(\lift v)$.  Now, $\lift h$ stabilizes each component of $\partial Y$, and stabilizes the curve $\lift v$, so $\lift h$ has the same action on $\mathcal C(Y)$ as does its restriction to $Y$ (which is defined since $\lift h$ stabilizes all boundary components).  Hence, $\phi(\lift h)$ acts in the same way as  $\phi(\tau^\ell)$ on $\link(\Delta)$, i.e. $\omega\circ\phi(\hat h)=\omega\circ\phi(\tau^\ell)$ (by a simple argument using $q(\link(\lift\Delta))=\link(\Delta)$; here $\tau$ is the $\hat v$--Dehn twist).

Since $\tau$ has a positive power contained in $N_{-1}$ --- see Notation \ref{not:correct_powers} --- we conclude that there are finitely many possibilities for $\omega\circ\phi(\hat h)$.  Since $\omega(g)=\omega\circ\phi(\lift g_i)\cdot\omega\circ\phi(\lift h)$, there are therefore only finitely many possibilities for $g$, i.e. the stabilizer of $v$ in $\Gamma_\Delta$ is finite, as required.

So, it remains to produce the lift $\lift g$ of $g$.

\textbf{Auxiliary vertices in $\link(\Delta)$:}  There exists a vertex $w\in\link(\Delta)-\{v\}$.  Indeed, by the first bullet point, $N_{-1}$ is contained in the subgroup $H<MCG(S)$ from Notation~\ref{not:correct_powers} (where $H$ is in particular a subgroup provided by Lemma~\ref{lem:fi}).  By Lemma~\ref{lem:fi}.\eqref{item:multi-orbit}, we have the following.  Recall that $\lift v$ is a curve on the complexity--$1$ subsurface $Y$, so by the lemma, there is a curve $\lift w$ on $Y$ such that $\lift v$ intersects (and differs from) every $H$--translate of $\lift v$.  In particular, $\lift v$ and $\lift w$ are not in the same $N_{-1}$--orbit, so the image $w$ of $\lift w$ in $\mathcal C(S)/N_{-1}$ is different from $v$.  On the other hand, $\lift w\in\link(\lift\Delta)$, so $w=q(\lift w)\in\link(\Delta)$.

\textbf{Lifting elements by lifting $4$--gons:} Fix $w\in\link(\Delta)-\{v\}$ and recall that $g\in\stabilizer(\link(\Delta))\cap\stabilizer(v)$.  We have a generalized $4$--gon $\Delta\star v=\Delta\star gv,gv\star g\Delta,g\Delta\star w,w\star\Delta$.  (The first equality $gv=v$ is because $g\in\stabilizer(v)$ and the simplex $g\Delta\star w$ exists because $w\in\link(\Delta)=g\link(\Delta)=\link(g\Delta)$ where the second equality holds since $g\in\stabilizer(\link(\Delta))$.)  Using the second bullet point, and the the third bullet point, we can lift this generalized $4$--gon to $\mathcal C(S)$ in such a way that our lift contains $\lift\Delta$ and $\lift v$ as the lifts of $\Delta$ and $v$.  

We claim that there exists $\lift g\in MCG(S)$ such that $\phi(\lift g)=g$ and such that $\lift g(\lift \Delta\star\lift v)$ is the lift of $g\Delta\star gv=g\Delta\star v$ appearing in the lift of our generalized $4$--gon.  Indeed, for any $\hat g_0\in\phi^{-1}(g)$, the simplex $\hat g_0(\lift\Delta\star\lift v)$ is a lift of $g\Delta\star v$.  Moreover, by the third bullet point, there exists $n_0\in N_{-1}$ such that $n_0\hat g_0(\lift\Delta\star\lift v)$ is the lift of $g\Delta\star v$ appearing in the lift of the generalized $4$--gon, and we take $\lift g=n_0\lift g_0$.

The lifted $4$--gon tells us that $\lift g(\lift \Delta\star\lift v)$ and $\lift\Delta\star \lift v$ have a common vertex which is a lift of $v=gv$.  The unique lift of $v$ in $\lift\Delta\star\lift v$ is $\lift v$,  so $\lift v\in \lift g(\lift\Delta\star \lift v)$.  We cannot have $\lift v\in\lift g\lift\Delta$, since that would imply $v\in g\Delta$, contradicting that $v=gv\in\link(g\Delta)$.  Hence $\lift v=\lift g\lift v$, i.e. $\lift g\in\stabilizer_{MCG(S)}(\lift v)$.  

On the other hand, $\lift\Delta^{(0)}$ and $\lift g\lift\Delta^{(0)}$ are both multicurves on $S$, both disjoint from the curves $\lift v$ and $\lift w$, as can be seen by examining the lifted generalized $4$--gon.  By definition, $Y$ is the complexity--$1$ component of the complement of $\lift\Delta^{(0)}$, so $\lift gY$ is the complexity--$1$ component of the complement of $\lift g\lift\Delta^{(0)}$.

Now, since $v\neq w$, the lifts $\lift v$ and $\lift w$ of $v$ and $w$ have to be distinct.  Since $\lift v$ and $\lift w$ are both in $\link(\lift\Delta)$, both of these vertices, regarded as curves, are in $Y$ (and they are not boundary curves of $Y$ because the boundary curves of $Y$ all belong to $\lift \Delta$).  Since $Y$ has complexity $1$, and $\lift v,\lift w$ are distinct, these curves fill $Y$.  

On the other hand, $\lift v$ and $\lift g\lift w$ are in the link of $\lift g\lift \Delta$ (by considering the lifted $4$--gon), so, as curves, they belong to $\lift gY$.  Indeed, they are intersecting curves in the complement of $\lift g\lift\Delta^{(0)}$, so they must belong to a common component of that complement, and since they fill a complexity--$1$ subsurface, the component they belong to must be $\lift gY$.  Hence $Y$ is a subsurface of $\lift gY$, i.e. $Y=\lift gY$.  Since $\lift g$ stabilizes $Y$, it must stabilize $\partial Y$. 

 We have produced $\lift g$ with $\phi(\lift g)=g$ and $\lift g\in\stabilizer_{MCG(S)}(\partial Y)\cap\stabilizer_{MCG(S)}(\lift v)$.  As explained above, this completes the proof of the claim.
\end{proof}
\renewcommand{\qedsymbol}{$\Box$}

Now we complete the proof of cocompactness for $i=0$.

Denote by $\Lambda_\Delta$ the quotient of $\pstab(\sat(\Delta))$ by the kernel of its action on $\link(\Delta)$. Both $\Lambda_\Delta$ and the group $\Gamma_\Delta$ from Claim \ref{claim:finite-stabs-stab-link} can be thought of as subgroups of $Sym(\link(\Delta))$, where we have $\Lambda_\Delta<\Gamma_\Delta$. We now claim that $\Lambda_\Delta$ is in fact a finite-index subgroup of $\Gamma_\Delta$. This suffices to prove cocompactness since, by item \eqref{item:HH_structure} and the ``moreover'' part of Theorem \ref{thm:propaganda}, there are finitely many $\stabilizer(\link(\Delta))$--orbits of edges in $\mathcal C(\Delta)$. Hence, the same is true for $\pstab(\sat(\Delta))$.

To show the claim, we use the fact if two groups $\Lambda<\Gamma$ act faithfully on a non-empty set, with $\Gamma$ having finite stabilizers and $\Lambda$ having finitely many orbits, then $\Lambda$ has finite index in $\Gamma$. 

Fix a lift $\lift\Delta$ of $\Delta$ to $\mathcal C(S)$, and let $Y$ be the complexity--$1$ component of the complement in $S$ of the multicurve $\lift\Delta^{(0)}$, which exists since $b([\Delta])$ is an orbit of complexity--$1$ subsurfaces (recall Definition \ref{defn:b}). 

Since $Y$ has complexity $1$, the simplex $\lift \Delta$ is almost-maximal in $\mathcal C(S)$, so $\Delta$ is almost-maximal in $\mathcal C(S)/N_{-1}$.  Thus Claim~\ref{claim:finite-stabs-stab-link} implies that $\Gamma_\Delta$ has finite stabilizers.

On the other hand, the kernel of the action of $\stabilizer_{MCG}(Y)$ on the complement of $Y$ is contained in $\pstab(\sat(\lift\Delta))$.  This kernel acts with finitely many orbits of curves on $Y$. Since, by item \eqref{item:lift_splx_induction}, we have  $q(\sat(\lift\Delta))=\sat(\Delta)$ (which implies that $\pstab(\sat(\lift\Delta))$ maps to $\pstab(\sat(\Delta))$), we conclude that the same holds for $\pstab(\sat(\Delta))$.  This completes the proof of cocompactness for $i=0$.\\

\textbf{Cocompactness for $i>0$:}  Now suppose $i>0$.  Now the simplex $\Delta$ of $\mathcal C(S)/N_i$ corresponds to a subsurface of complexity $i+1>1$.

Suppose that $\sigma,\tau$ are $W^\Delta$--adjacent maximal simplices of $W^\Delta$.  By Lemma~\ref{lem:C_0=C}, this occurs if and only if $\sigma\star\Delta$ and $\tau\star\Delta$ are adjacent in $W$.  By the first item in the ``moreover'' clause of Theorem~\ref{thm:propaganda}, there exists a common codimension--$1$ face $\eta\star\Delta$ of $\tau\star\Delta$ and $\sigma\star\Delta$.  Since $i+1>1$, the simplex $\Delta$ is not almost-maximal, so $\eta\neq\emptyset$, and therefore $[\eta\star\Delta]\propnest[\Delta]$.  From Theorem~\ref{thm:MCG_quotient}.\eqref{item:HH_structure}, first and second bullet points, $b([\eta\star\Delta])$ is an orbit of subsurfaces of complexity at most $i$, so by the third bullet point of the same statement, $\mathcal C([\eta\star\Delta])$ is finite.  Hence, fixing $\eta$, there are only finitely many possibilities for $\tau\star\Delta$ and $\sigma\star\Delta$ and thus finitely many possibilities for the $W$--edge joining them.

Thus, to show that the action of $\pstab(\sat(\Delta))$ on $W^\Delta$ has finitely many orbits of edges, it suffices to show that the simplex $\eta$ of $\link(\Delta)$ belongs to one of finitely many $\pstab(\sat(\Delta))$ orbits.  

We argue as in the previous case that $\pstab(\sat(\Delta))$ acts with finitely many orbits of simplices on $\link(\Delta)$.  By assumption, we can fix a lift $\lift\Delta$ of $\Delta$ to $\mathcal C(S)$.  By \eqref{item:lift_splx_induction}, $q(\sat(\lift\Delta))=\sat(\Delta)$, so it suffices to show that $\pstab(\sat(\lift\Delta))$ acts with finitely many orbits of simplices on $\link(\lift\Delta)$.  Now, let $Y$ be the union of the non-pants components of the complement in $S$ of the multicurve $\lift\Delta^{(0)}$.  Then $\link(\lift\Delta)^{(0)}$ is the vertex set of $\mathcal C(Y)$, and $\lift\Delta^{(0)}$ consists of the boundary curves of $Y$ plus a pants decomposition of $S-Y$.  Now, $\stabilizer_{MCG(S)}(Y)$ acts on $Y$ with finitely many orbits of multicurves, i.e. it acts with finitely many orbits of simplices on $\link(\lift\Delta)$.  Hence the same is true for the finite index subgroup $\Lambda$ of $\stabilizer_{MCG(S)}(Y)$ that stabilizes each boundary curve of $Y$.  But since any $\lift g\in\Lambda$ can be restricted to $Y$, there exists $\lift h\in MCG(S)$ that acts as the identity on $S-Y$ and as $\lift g$ on $\mathcal C(Y)$.  Now, $\lift h\in\pstab(\sat(\lift\Delta))$, and it follows that the latter group also has finitely many orbits of simplices in $\link(\lift\Delta)$, as required.\\

\textbf{Conclusion:}  Hence, $\pstab(\sat(\Delta))$ acts properly and cocompactly on the combinatorial HHS $(\link(\Delta), 
W^\Delta)$, and it is therefore a hierarchically hyperbolic group, by Theorem~\ref{thm:hhs_links}. 

To show that $W^\Delta$ and hence $\pstab(\sat(\Delta))$ are hyperbolic, it suffices to show that this HHS structure on $W^\Delta$ has at most one unbounded hyperbolic space, namely $\mathcal C(\Delta)=\mathcal C_0^\Delta([\emptyset])$ (see Lemma~\ref{lem:curve_graph_iota_v2}.\eqref{item:same_C} for the latter equality).  Indeed, let $\Sigma$ be any non-empty, non-maximal simplex of $\link(\Delta)$.  By Lemma~\ref{lem:curve_graph_iota_v2}.\eqref{item:same_C} and Lemma~\ref{lem:C_0=C} (saying roughly that links are the same for $X$ as for $\link(\Delta)$), we have to bound the diameter of $\mathcal C([\Delta\star\Sigma])$.  But since $\Sigma$ is nonempty, we have $[\Delta\star\Sigma]\propnest[\Delta]$. This implies that $b([\Delta\star\Sigma])$ has a representative $U$ of complexity at most $i$, since  $b([\Delta])$ has a representatives $Y$ with $U$ properly nested in $Y$, so that the complexity of $U$ is strictly lower than that of $Y$ (recall that $b$ preserves nesting, see item \eqref{item:HH_structure}).

We can now apply the third bullet point of item \eqref{item:HH_structure} (recall that we are hypothesizing this property for the current $i$) to conclude that $\mathcal C([\Delta\star\Sigma])$ is uniformly bounded.  Hence $W^\Delta$ and $\pstab(\sat(\Delta))$ are hyperbolic.
\end{proof}

\subsection{Verifying composite properties, and Greendlinger's Lemma}\label{subsec:composite_induction}
In this subsection only, we are working in the case $i>-1$, i.e., we are doing the inductive step laid out in 
Convention~\ref{conv:induction}.  The base case, $i=-1$, was already handled when we verified that the 
subgroups $\Gamma_Y$ used to form $N_{-1}$ form a composite rotating 
family, etc.  Many of the additional 
points in Proposition~\ref{prop:technical_MCG_quotient} were introduced to enable arguments in this subsection.

Using our inductive hypothesis, we consider the combinatorial 
 hierarchically hyperbolic space $(\mathcal C(S)/N_{i-1},W)$. We let $X=\mathcal C(S)/N_{i-1}$.
 
 We let $\bbY_*$ denote the collection of equivalence classes 
 $[\Delta]$ of simplices $\Delta$ of $X$ so that: 
 \begin{itemize}
  \item $b([\Delta])=N_{i-1}Y$ for some subsurface $Y$ of $S$ of complexity $i$ (recall that $b$ is the 
bijection between equivalence classes of simplices and $N_{i-1}$--orbits of subsurfaces),
  \item $\fontact ([\Delta])$ is unbounded.
 \end{itemize}

 In view of \eqref{item:deep_BBF}, the coloring of $\mathfrak S$ descends to a coloring of $\mathfrak 
S/N_{i-1}$, so 
that $\bbY_*$ is partitioned into finitely many countable 
 families.

Before showing that $\bbY_*$ defines a composite projection system, we need a preliminary lemma.

\begin{lemma}\label{lem:sat=link_link}
 Let $[\Delta]\in \bbY_*$. Then $\link(\link(\Delta))^{(0)}=\sat(\Delta)$. In particular, $[\Sigma]\in \bbY_*$ satisfies $[\Delta]\orth[\Sigma]$ if and only if $\link(\Sigma)^{(0)}\subseteq \sat(\Delta)$.
\end{lemma}

\begin{proof}
The containment $\sat(\Delta)\subseteq \link(\link(\Delta))$ holds 
for any simplex in any complex (by definition of the saturation), so it suffices to prove 
the other containment.

 Since $\mathcal C(\Delta)$ is unbounded, we can find vertices $v,w\in\link(\Delta)$ that are not connected by 
any path in $\link(\Delta)$ of length less than 3. Let $x$ be a vertex of $\link(\link(\Delta))$. We have a 
generalized $4$-gon $\Delta\star v, v\star x,x\star w, w\star \Delta$, which we can lift to $\mathcal C(S)$. 
 We denote by $\lift \Delta, \lift v$, etc. the various lifts.   In terms of curves, the multicurve 
$\Delta^{(0)}$ is disjoint 
from the curves $\lift v$ and $\lift w$. Moreover, from the hypothesis that $v$ and $w$ are sufficiently far 
in $\link(\Delta)$, we see that $\lift v\cup\lift w$ fills the (necessarily unique, since $\link(\Delta)$ is unbounded) 
component of the complement of
$\Delta^{(0)}$ which is not a pair of pants. Since $\lift x$ is disjoint from $\lift v$ and $\lift w$, we then 
see that $\lift x\in \sat(\lift\Delta)$ (that is, it is part of a pants decomposition of the complement of the 
surface filled by $\lift v\cup\lift w$, as is $\Delta^{(0)}$). Since $q(\sat(\lift\Delta))=\sat(\Delta)$, we 
have $x\in\sat(\Delta)$, as required.
\end{proof}

\begin{lemma}\label{lem:compositeprojectionsystem} 
$\bbY_*$ defines a composite projection system.
\end{lemma}

\begin{proof} 
	We will use the HHS structure associated to $(X,W)$ via Theorem~\ref{thm:hhs_links} and 
the inductive hypothesis that $(X,W)$ is a combinatorial HHS.
	
	First, since all the elements of $\bbY_*$ are associated to 
	equivalence classes of 
	domains in $\frakS$ of the same complexity, no nesting can occur 
	and thus any pair of these associated domains is either 
	transverse or orthogonal.
	
	Accordingly, to each $Y\in\bbY_*$, we define the set $Act(Y)$ to be the set of all elements 
of $\bbY_*-\{Y\}$ which are associated to 
	domains transverse to $Y$.  We will almost always work with $\actminus Y$.  The symmetry in action axiom 
(Definition~\ref{defn:CPS}.(CPS6)) is immediate.
		
	Recall that projections in a combinatorial HHS are defined in
	Definition \ref{defn:projections}.  Consider $Y\in\bbY_*$ and
	$W,Z\in \actminus{Y}$, and define
	$d_{Y}(W,Z)=d_{Y}(\rho^W_Y,\rho^Z_Y)$.

	Symmetry (Definition~\ref{defn:CPS}.(CPS1)) and the triangle inequality 
(Definition~\ref{defn:CPS}.(CPS2)) follow immediately from the
	fact that $d_{\fontact Y}$ is a distance function.  The
	Behrstock inequality (Definition~\ref{defn:CPS}.(CPS3)) follows from 
Definition~\ref{defn:HHS}.\eqref{item:dfs_partial_realization},\eqref{item:dfs_transversal}. 
 Properness (Definition~\ref{defn:CPS}.(CPS4)) follows from the distance formula 
(Theorem~\ref{thm:distance_formula}) and Definition~\ref{defn:HHS}.\eqref{item:dfs_partial_realization}.

	The separation axiom (Definition~\ref{defn:CPS}.(CPS5)) holds trivially since $d_{Y}$ is a distance 
	function. 
	
	The closeness in inaction axiom (Definition~\ref{defn:CPS}.(CPS7)) holds, since if $W\notin\act(Z)$ 
	and $Y\in\act(Z)\cap\act(W)$, then $W$ and $Z$ are orthogonal and 
	$Y$ is transverse to both of them. This implies that 
	$d_{Y}(W,Z)$ is uniformly bounded, by~\cite[Lemma 1.5]{HHS:boundary}.
	
	The finite filling axiom (Definition~\ref{defn:CPS}.(CPS8)) turns 
	out to be the hardest. We first record the 
following which allows us to verify the active sets in this axiom via links 
and saturations:
	\begin{claim}
	 Let $\mathcal W\subseteq \bbY_*$. Then, for $[\Sigma]\in\bbY_*$, we have $[\Sigma]\in 
\bigcup_{w\in\mathcal W}\act(w)$ if and only if $\link(\Sigma)\nsubseteq \bigcap_{[\Delta]\in \mathcal 
W}\sat(\Delta)$.
	\end{claim}

\renewcommand{\qedsymbol}{$\blacksquare$}
	\begin{proof}
	Passing to the complements, we show that $[\Sigma]\in\bigcap_{w\in\mathcal W}\act(w)^c$ if and only 
if $\link(\Sigma)\subseteq \bigcap_{[\Delta]\in \mathcal W}\sat(\Delta)$.
	
	This is just because  we have $[\Sigma]\notin \act([\Delta])$ $\iff$ $[\Sigma]\orth[\Delta]$, and in view of Lemma \ref{lem:sat=link_link} this is equivalent to $\link(\Sigma)\subseteq \sat(\Delta)$.
	\end{proof}

	In view of the claim it suffices to show that for any $\mathcal W\subseteq \bbY_*$ there exists a finite 
subcollection $\{[\Delta_1],\dots,[\Delta_n]\}\subseteq \mathcal W$ so that 
$\bigcap_j\sat(\Delta_j)=\bigcap_{[\Delta]\in \mathcal W}\sat(\Delta)$. This readily follows from the 
following 
claim.

\begin{claim}
\label{claim:sat_seq}
There does not exist an infinite sequence $[\Delta_1],[\Delta_2], \dots$ in $\bbY_*$ so that $$\bigcap_{j\leq 
n+1}\sat(\Delta_j)\subsetneq \bigcap_{j\leq n}\sat(\Delta_j)$$ for all $n$.
\end{claim}

\begin{proof}
First, we show that for every $[\Delta]\in\bbY_*$ there are simplices $\Pi(\Delta)$ and $\Pi'(\Delta)$ so that we have 
$\sat(\Delta)=\link(\Pi(\Delta))\star \Pi'(\Delta)$.

We record the following observation about curve graphs:

\begin{claim}\label{claim:sat_structure}
Let $\Sigma$ be a simplex in $\mathcal C(S)$.  Then there exist simplices $\Sigma',\Sigma''$ of $\mathcal 
C(S)$ such that $\sat(\Sigma)=\link(\Sigma')\star\Sigma''$.
\end{claim}

\begin{proof}
First of all, the vertex set of $\link(\Sigma)$ consists of all essential curves in the complement of $\Sigma$, that is, it coincides with the vertex set of the curve graph of the (possibly disconnected) subsurface $Y$ of $S$ consisting of all components of the complement of $\Sigma$ which are not pairs of pants. Moreover, $\Sigma$ is the simplex whose vertex set consists of
\begin{itemize}
 \item all curves of $S$ isotopic to a boundary component of $Y$ in $S$ (note that two boundary components of $Y$ can be isotopic in $S$), and
 \item a pants decomposition of the union $Z$ of the non-annular components of the complement of $Y$. 
\end{itemize}
We let $\Sigma''$ be the simplex with vertex set the curves as in the first item.

 Any simplex with the same link as $\Sigma$ admits the same description, and moreover any essential curve in $Z$ can be completed to a pants decomposition. This implies that $\sat(\Sigma)$ consists of the vertex set of $\Sigma''$ together with all essential curves of $Z$. We can then pick $\Sigma'$ to be any simplex whose link has vertex set the set of essential curves of $Z$; namely, $\Sigma'$ has vertex set consisting of all curves of $S$ isotopic to a boundary component of $Z$ in $S$, together with a pants decomposition of $Y$.
\end{proof}

Consider any lift $\lift\Delta$ of $\Delta$ to $\mathcal 
C(S)$. Claim~\ref{claim:sat_structure} implies that $\sat(\lift\Delta)=\link(\Delta')\star\Delta''$ for some 
simplices $\Delta',\Delta''$ of $\mathcal 
C(S)$.

Also notice that $q(\link(\Delta'))=\link(q(\Delta'))$.  Indeed, the containment 
``$\subseteq$'' is clear, and the other 
containment follows from the fact that for any vertex $v$ in $\link(q(\Delta'))$ we can lift $v\star q(\Delta')$ to a 
simplex containing $\Delta'$ (by existence and uniqueness of the orbit of lifts of simplices, 
item~\eqref{item:lift_splx_induction}). In particular, we get that 
$\sat(\Delta)=q(\sat(\lift\Delta))=\link(q(\Delta'))\star q(\Delta'')$, and we are done.

Now suppose, for a contradiction, that we have a sequence $[\Delta_1],[\Delta_2], \dots$ as in the 
statement of Claim \ref{claim:sat_seq}.  By induction, suppose we proved that 
$I_n=\bigcap_{j\leq n}\sat(\Delta_j)=\link(\Pi_n)\star\Pi'_n$ for some simplices $\Pi_n,\Pi'_n$ with $\Pi_n\subseteq 
\Pi_{n+1}$. The base case $n=1$ is given by the argument above. 

Then $I_{n+1}=I_n\cap \link(\Pi(\Delta_n))\star 
\Pi'(\Delta_n)=\link(\Pi_n)\star\Pi'_n\cap  \link(\Pi(\Delta_n))\star \Pi'(\Delta_n)$ is readily seen to be of the required 
form in view of condition \ref{item:propaganda_join} of Theorem \ref{thm:propaganda}, which holds for $X$ by inductive 
hypothesis (Theorem \ref{thm:MCG_quotient}.\ref{item:HH_structure}). Since $I_{n+1}\subsetneq I_n$, for each $n$ we have either $\Pi'_{n+1}\subsetneq 
\Pi'_n$ or, if not, $\Pi_n\subsetneq \Pi_{n+1}$. This is easily seen to imply that the $\Pi_n$ have arbitrarily many 
vertices, contradicting the finite dimensionality of $X$.  This proves the claim.
\end{proof}
 \renewcommand{\qedsymbol}{$\Box$}
This completes the proof of Lemma \ref{lem:compositeprojectionsystem}.
\end{proof}

\begin{lemma}\label{lem:compositerotatingfamily}
For every $\theta>0$ the following holds. For any $Y=[\Delta]\in\bbY_*$ we can choose a finite-index subgroup 
$\Gamma_Y^\theta<\pstab(\sat(\Delta))$ contained in $H/N_{i-1}<\bar G_{i-1}=MCG(S)/N_{i-1}$ so that
the subgroups $\Gamma_Y^\theta$ form a composite rotating family on the composite projection system $\bbY_*$. 
Moreover,
$$\min\{d_{Y}(x,\gamma x): Y\in\bbY_*,\gamma\in\Gamma_Y-\{1\},x\in\mathcal C(Y)\}>\theta.$$
\end{lemma}
 
\begin{proof} In Lemma~\ref{lem:compositeprojectionsystem} we 
	established that  $\bbY_*$ is a composite 
	projection system.
	
	Definition~\ref{defn:ThetaCRF}.(CRF1): Roughly, this says that $\Gamma_Y^\theta$ needs to be an infinite group acting on $\bbY_*$ fixing $Y$ and, for any $R$, having only finitely many elements that move some point of $\mathcal C(Y)$ at most distance $R$.

Any sufficiently deep finite-index normal subgroup $\Gamma_Y^\theta$ satisfies this by 
inductive hypothesis 
\eqref{item:bottom_hyp} (which guarantees that $\Gamma_Y^\theta$ acts properly on an unbounded graph, so that it 
must be infinite).
	
	Definition~\ref{defn:ThetaCRF}.(CRF2): It suffices to choose the $\Gamma^\theta_Y$ as follows. For each $Y=[\Delta]$ 
in a given set of representatives of $\bar G_{i-1}$-orbits in $\bbY_*$, choose a normal subgroup $\Gamma^\theta_Y$ of 
$\pstab(\sat(\Delta))$, and extend the choice to all $Y$.
	
	Definition~\ref{defn:ThetaCRF}.(CRF3): For this, we use that $\Gamma_Y^\theta$ is contained in $H/N_{i-1}$. We 
first make a preliminary 
claim.
	
	\begin{claim}
	 Let $[\Delta]\orth[\Sigma]$, and let $g\in \pstab(\sat(\Delta))$, $h\in \pstab(\sat(\Sigma))$. Then $gh(x)=hg(x)$ for all $x\in \sat(\Delta)\cup\link(\Delta)$, and $gh(x)=hg(x)=g(x)$ for all $x\in\link(\Delta)$.
	\end{claim}
 \renewcommand{\qedsymbol}{$\blacksquare$}

	\begin{proof}
	First, we prove that $h(\sat(\Delta))\subseteq \sat(\Delta)$.  The same argument also shows 
$g(\sat(\Sigma))\subseteq \sat(\Sigma)$.
	
	Let $v\in \sat(\Delta)$, that is, $v$ is a vertex of a simplex $\Delta'$ with 
$\link(\Delta')=\link(\Delta)$. We have to show that $hv\in \sat(\Delta)$. We have that $hv$ is a vertex of 
$h\Delta'$, so we have $hv\in \sat(\Delta)$ provided that $\link(h\Delta')=\link(\Delta')$. But 
$\link(\Delta')=\link(\Delta)\subseteq\sat(\Sigma)$ by Lemma \ref{lem:sat=link_link}, and since $h\in 
\pstab(\sat(\Sigma))$, we have that $h$ fixes $\link(\Delta')$. Hence, $\link(h\Delta')=\link(\Delta')$, as 
required.
	
	Now let $x\in \sat(\Delta)$. Then $h(x)\in \sat(\Delta)$, so that $g(h(x))=h(x)$. On the other hand, 
$g(x)=x$, so that $h(g(x))=h(x)$. This shows that $gh(x)=hg(x)$ for all $x\in \sat(\Delta)$. Similarly, 
$gh(x)=hg(x)=g(x)$ for all $x\in \sat(\Sigma)$. To conclude, just notice that 
$\sat(\Delta)\cup\link(\Delta)\subseteq \sat(\Delta)\cup \sat(\Sigma)$ since $\link(\Delta)\subseteq 
\sat(\Sigma)$ by Lemma \ref{lem:sat=link_link}.
	\end{proof}
	
	Suppose that $X\notin \act(Z)$, where $X=[\Sigma]$, $Z=[\Delta]$; we have to show that $\Gamma^\theta_X$ commutes 
with $\Gamma^\theta_Z$. Then $X\orth Z$, since there is no nesting relation between distinct elements of $\bbY_*$. For $g\in 
\pstab(\sat(\Delta)), h \in \pstab(\sat(\Sigma))$, by the claim above we have that $gh$ and $hg$ act in the same way on 
$\sat(\Delta)\cup \link(\Delta)$, so we are done once we prove the following:
	
	\begin{claim}
	 Let $a\in H/N_{i-1}\leq\bar G_{i-1}$ act trivially on $\sat(\Delta)\cup \link(\Delta)$. Then $a$ is the 
identity.
	\end{claim}

	\begin{proof}
	Let $\phi:MCG(S)\to \bar G_{i-1}$ be the quotient map. Let $P$ be a maximal simplex of $X$ containing 
$\Delta$.
The vertex set of $P$ is contained in $\sat(\Delta)\cup \link(\Delta)$, so $a$ fixes $P$ 
pointwise. We can lift $P$ to a simplex $\lift P$ of $\mathcal C(S)$, and since there is a unique $N_{i-1}$--orbit 
of such 
lifts, there is $\lift a\in H$ so that $\phi(\lift a)=a$ and $\lift a$ fixes $\lift 
P$ 
pointwise.
	
	We claim that $\lift P$ is a maximal simplex of $\mathcal C(S)$.  Indeed, suppose to the contrary that $\lift P$ is 
properly contained in a simplex $\lift P'$.  Then, since $X$ is simplicial, $\lift P'$ projects to a simplex $P'$ of $X$ 
strictly containing $P$, a contradiction.  So, $\lift P$ is maximal.

	In terms of curves, $\lift a$ fixes the pants decomposition $\lift P^{(0)}$, so that $\lift a$ is a product of powers 
of Dehn twists around the curves of $\lift P^{(0)}$. Since $N_{-1}\cap \langle \tau\rangle=H\cap \langle 
\tau\rangle$ and 
$N_{-1}<N_{i-1}<H$ (by induction hypothesis \eqref{item:deep_BBF}), we see that $\lift a\in N_{i-1}$, so that $a$ is 
the identity, as 
required.	
\end{proof}
 \renewcommand{\qedsymbol}{$\Box$}

	Definition~\ref{defn:ThetaCRF}.(CRF4): Since $\Gamma^\theta_Y$ acts properly on $\mathcal C(Y)$, and there are 
finitely many $\bar 
G_{i-1}$-orbits in $\bbY_*$, choosing sufficiently deep subgroups for orbit representatives as in the proof of (CRF2) 
ensures the required property. The ``moreover'' part follows similarly from $\pstab(\sat(\Delta))$ acting properly and cocompactly on $\mathcal C(\Delta)$ (inductive hypothesis \eqref{item:bottom_hyp}), so that it suffices to choose deep subgroups for orbit representatives.
\end{proof}

 \begin{lemma}[``Greendlinger'']\label{lem:short_less_complex} 
There exists a diverging function $\frak T$ so that the following holds. Let $\theta>0$ and let 
$\Gamma^\theta_Y=\Gamma_Y$ be as 
in Lemma \ref{lem:compositerotatingfamily}.

 Then, for $N=\langle\langle \{\Gamma^\theta_Y\}\rangle\rangle$, the following holds.  There is a well-ordered set $\mathfrak 
C$, and an assignment $N\ni\gamma\mapsto c(\gamma)\in \mathfrak C$, with $c(1)$ minimal in $\mathfrak C$ and the following additional properties.

For all $\gamma\in 
N-\{1\}$ and all simplices $\Delta$ of $X$, there is $Y\in\bbY_*$ and $\gamma_Y\in\Gamma_Y^\theta$ so that $c(\gamma_Y 
\gamma)<c(\gamma)$ and either
  \begin{itemize}
   \item $\Delta\subseteq \fix(\Gamma_Y)$, or $\gamma \Delta\subseteq \fix(\Gamma_Y)$, or
   \item $d_Y(\Delta,\gamma \Delta)>\frak T(\theta)$ (and the quantity is defined).
  \end{itemize}
 \end{lemma}
 
 \begin{rem}
 \label{rem:fix_sat}
  For $Y\in\bbY_*$, $\sat(Y)$ is fixed pointwise by $\Gamma_Y$, so that $\pi_Y(\Delta)$ is defined for any 
simplex $\Delta$ of $X$ not contained in $\fix(\Gamma_Y)$.
 \end{rem}

 \begin{proof}[Proof of Lemma~\ref{lem:short_less_complex}]
Independently of $\theta$, we now choose, for each simplex $\Delta$ of $X$ and each $\bbY_k$, some 
$Y^\Delta_k\in\bbY_k$. We 
make 
the choice as follows.

First, we 
choose 
$H/N_{i-1}$-orbit-representatives $\Delta_j$ of simplices, and make the choice for each 
of them. Then, for any other simplex $\Delta$, we choose a group element $g\in H/N_{i-1}$ mapping the suitable orbit 
representative, say $\Delta_j$, to the given simplex $\Delta$, and set $Y^\Delta_k=gY^{\Delta_j}_k$ (since $H/N_{i-1}$ preserves 
the colors, this is a well-defined element of $\bbY_k$). Set
 
 $$D=\max_{j,k} \sup\{d_Y(\Delta_j,Y^{\Delta_j}_k): Y\in \bbY_k,\Delta_j\nsubseteq \sat(\Delta_j)\},$$
 which is finite since there are finitely many colors and finitely many orbits (since $H/N_{i-1}$ has finite index in $\bar G_{i-1}$, 
and $\bar G_{i-1}$ acts cocompactly on $X$), and the distance formula takes finite values.

 We refer the reader to ~\cite[Theorem 4.1]{DFDT} (which follows by combining Dahmani's construction in~\cite[Section 2.4.2]{D_PlateSpinning} with~\cite[Prop. 
2.13, Lemma 2.16, Lemma 2.17]{D_PlateSpinning}), where the subgroup $N$ in Lemma~\ref{lem:short_less_complex} is described as an increasing union of subgroups $N=\bigcup_\alpha N_\alpha$, where $\alpha$ varies over countable ordinals. 

For $\gamma\in N$, we denote by $\alpha(\gamma)$ the smallest ordinal such that $\gamma$ is conjugate into $N_{\alpha(\gamma)}$. Two properties of $\alpha(\gamma)$ observed in \cite{DFDT} are that $\alpha(\gamma)$ is never a limit ordinal and $\alpha(\gamma)=0$ if and only if $\gamma=1$. Moreover, $N_\alpha$ for $\alpha$ a successor ordinal has a certain amalgamated product decomposition used in \cite[Definition 4.3]{DFDT} to define $n(\gamma)$, which is the length of the cyclic normal form in $N_{\alpha(\gamma)}$ for the conjugacy class of $\gamma$. Following \cite{DFDT}, we consider the complexity $c(\gamma)$ given by $(\alpha(\gamma),n(\gamma))$, ordered lexicographically.

After conjugation by a suitable element $h$ and replacing $\Delta$ with $h^{-1}\Delta$, 
one can assume that $\gamma \in N_{\alpha(\gamma)}$.

Let $Y=v$ and $\gamma_Y=\gamma_v$ be as in \cite[Proposition 4.5, fourth and second bullet]{DFDT} with 
$\nu=Y^\Delta_{j(\gamma)}$, so that $\gamma_Y\gamma$ has shorter cyclic normal form. We argue that this implies 
that 
$c(\gamma_Y\gamma)<c(\gamma)$. Indeed, it the length of the new normal form is still 
greater than $2$, we have $\alpha(\gamma_Y\gamma)=\alpha(\gamma)$ but the second coordinates satisfy $n(\gamma_Y\gamma)<n(\gamma)$. If not, the length is reduced to $1$ (or $0$)  and 
$\gamma_Y 
\gamma$ is actually conjugate into   $N_{\alpha(\gamma)-1}$, so that $\alpha(\gamma_Y\gamma)<\alpha(\gamma)$.
In either case we have $c(\gamma_Y\gamma)<c(\gamma)$, as required.

We are left to show that one of the alternatives applies. If $\Delta,\gamma \Delta\nsubseteq \fix(\Gamma_Y)$ (otherwise we 
are done), then we have that $$d_Y(\Delta,\gamma \Delta)\geq d_Y(Y^\Delta_{j(\gamma)},\gamma 
Y^\Delta_{j(\gamma)})-2D,$$ as 
required; we can take $\frak T(\theta)=\theta-2D$.
 \end{proof}

\begin{notation}\label{not:not_annuli}
   We now choose $\theta$ sufficiently large to ensure that
   \begin{itemize}
   \item $\frak T(\theta)>( 6\genus(S)+2p(S)-3)C_{i-1}$, for $C_{i-1}$ as in Proposition 
\ref{prop:technical_MCG_quotient}.\eqref{item:BGI},
    \item $\frak T(\theta)> \max\{d_{\mathcal C(Y)}(fx,g x):Y\in\bbY^*, f,g\in F, x\in X\}$,
    \item $\frak T(\theta)>\kappa+C_{i-1}$, for $\kappa$ as in Lemma \ref{lem:conv_cocpt=cobounded} applied to 
the image of $Q$ in $\bar G_{i-1}$.
   \end{itemize}
We set $\Gamma_Y=\Gamma^\theta_Y$ as in Lemma \ref{lem:short_less_complex}, and set $N=\langle\langle 
\{\Gamma_Y\}\rangle\rangle$.  We emphasize that $N$ is the kernel of the quotient $\bar G_{i-1}\to 
\bar G_{i}$, and is not the same as $N_i$.
\end{notation}

\subsection{Lifting}\label{subsec:lifting_gons}
In this subsection, we are back in the setting of Convention~\ref{conv:induction}.  In other words, either 
$i=-1$ (base case) or $i>-1$ (inductive step).  The goal in this section is to lift generalized $m$--gons from 
$X/N$ to $X$.  This will typically be used in conjunction with our inductive hypothesis (see Proposition~\ref{prop:technical_MCG_quotient}), to lift all the way 
back up to $\mathcal C(S)$.

(For clarity, in the proposition we recall parts of the definitions of lifts.)

 \begin{prop}[Lifting generalized $m$--gons through $q$]\label{prop:omni_lift}
 For every $m\leq 6\genus(S)+2p(S)-3$ the following hold.
  \begin{enumerate}
   \item\label{item:lift_splx} For every simplex $\Delta$ of $X/N$,  together with an order $(v_0,\dots,v_k)$ on the vertices 
there exists a unique $N$-orbit of simplices $\Sigma$ in $X$, each with an order $(w_0,\dots,w_k)$ on its vertices, such 
that $\Sigma$ is a lift of the ordered simplex $\Delta$ (meaning that $q(w_j)=v_j$ for $0\leq j\leq k$).
   \item \label{item:lift_geod_link}
 Given a simplex $\Delta$ in $X/N$, a  lift $\Sigma$ of $\Delta$ in $X$, and a geodesic $\gamma$ in the link of $\Delta$, we 
have that $\gamma$ can be lifted to a geodesic in the link of $\Sigma$.
 
\item \label{item:lift_gon} Any generalized $m$-gon $\tau=\tau_0,\dots,\tau_{m-1}$ in $X/N$ can be lifted to a generalized 
$m$-gon $\tau'_0,\dots,\tau'_{m-1}$ in $X$ so that $\tau'_j$ is \textbf{type S}/\textbf{type G} if and only if $\tau_j$ is, 
and if $\tau_j$ is a 
geodesic in $\link(\Delta_j)$ then $\tau'_j$ is a geodesic in $\link(\Sigma_j)$ for some lift $\Sigma_j$ of $\Delta_j$.
  \end{enumerate}
 \end{prop}
 
 \begin{proof}
 Consider a generalized $m$-gon $\tau=\tau_0,\dots,\tau_{m-1}$. If $\tau_j$ is of \textbf{type G}, and $\Delta_j$  is the 
corresponding simplex, let $d_j$ be the number of vertices of $\Delta_j$. If $\tau_j$ is \textbf{type S}, let $d_j$ be the 
number of vertices of $\tau_j$ minus $2$. Finally, set $d(\tau)=\max_j\{d_j\}$. 

\textbf{Good lifts:}  We say that $\tau_j$ of \textbf{type S} has a \emph{good lift} if it can be lifted, and that $\tau_j$ 
of \textbf{type G} has a \emph{good lift} if it can be lifted to a geodesic in the link of a lift of the corresponding 
simplex.
 
 \textbf{Structure of the proof:}  The proofs of the 3 items are interlaced, and more specifically we will prove the 
following claims, for any $d\geq 2$:
 \begin{enumerate}[(a)]
 \item If Item \ref{item:lift_splx} holds whenever $\Delta$ has at most $d$ vertices, then Item \ref{item:lift_geod_link}  
holds whenever $\Delta$ has at most $d-2$ vertices.
 \item If Item \ref{item:lift_splx} holds whenever $\Delta$ has at most $d$ vertices, then Item \ref{item:lift_gon}  holds 
whenever $d(\tau)\leq d-2$.
 \item If Item \ref{item:lift_gon} holds whenever $d(\tau)\leq d-2$ and Item \ref{item:lift_splx}  holds whenever $\Delta$ 
has at most $d$ vertices, then Item \ref{item:lift_splx} holds whenever $\Delta$ has at most $d+1$ vertices.
 \end{enumerate}

$(a)-(c)$ feed into an induction that proves all 3 items, with the base case specified below and corresponding to Item \ref{item:lift_splx} for $\Delta$ with at most two vertices (from which by $(a)$ and $(b)$ we can deduce Items \ref{item:lift_geod_link} and \ref{item:lift_gon} for $\Delta$ with at most 2 vertices, which in turn by $(c)$ yield Item \ref{item:lift_splx} for $\Delta$ with at most 3 vertices, etc.).

\textbf{Base case:}  The base cases are the following:
\begin{itemize}
     \item For each vertex $\bar v$ of $X/N$, there is a unique $N$--orbit of vertices $v\in X$ such that $q(v)=\bar v$, 
just because of how $X/N$ is defined.

    \item For each edge $\bar e$ of $X/N$, with endpoints $\bar v,\bar w$, the fact that $N$ acts simplicially on $X$ (in 
such a way that stabilizers of simplices fix them pointwise) implies that there is a unique orbit of edges $e$ of $X$, with 
endpoints $v,w$, such that $q(v)=\bar v,q(w)=\bar w$.
\end{itemize}

We now prove $(a)-(c)$, keeping $(b)$ last since it is the hardest.
 
 \textbf{Proof of $(a)$:} Let $v_0,v_1,\dots$ be the sequence of vertices along $\gamma$.  Then $\Delta\star v_0$ is a 
simplex of $X/N$, which can be lifted to a simplex $\Sigma'\star v'_0$ of $X$. By the uniqueness clause of Item 
\ref{item:lift_splx}, up to applying an element of $N$, we can assume $\Sigma'=\Sigma$, so that we lifted the first 
vertex of the geodesic in the appropriate link. Suppose that we lifted the sequence of vertices $v_0,\dots,v_k$ to 
$v'_0,\dots,v'_k$. Since $\Delta\star v_k\star v_{k+1}$ is a simplex of $X/N$, it can be lifted to $X$, and similarly to the 
argument for $v_0$, the lift can be chosen to be of the form $\Sigma\star v'_k\star v'_{k+1}$, so that we can lift $v_k$ to 
$v'_k$. Inductively, we can the lift the whole geodesic $\gamma$.
 
 \textbf{Proof of $(c)$:} Consider a simplex $\Delta$ with at most $d+1$ vertices.  By the base case of the induction, we 
can assume that $\Delta$ has at least 3 vertices, so $\Delta=\Delta'\star v\star w$, for some non-empty simplex $\Delta'$. We 
can think of $\Delta$ as a generalized $3$-gon $\Delta'\star v,v\star w,\Delta'\star w$, with all sides of 
\textbf{type S}, and each 
$d_j\leq d-2$. Hence, the 3-gon can be lifted, which also provides a lift of $\Delta$. 

For the uniqueness 
part, suppose that 
$\Delta$, endowed with an order on the vertices, has two lifts not in the same $N$-orbit, and consider an 
arbitrary 
codimension--$1$ face $\Delta''$, with opposite vertex $u$. We know that $\Delta''$ has a unique orbit of lifts, so that we see 
that there exist two lifts of $\Delta$ of the form $\Sigma\star u_1,\Sigma\star u_2$ that are not in the same 
$N$-orbit, but $u_1,u_2$ are (as two lifts of $\Delta$ can be made to coincide on the lifts of $\Delta''$ that they contain). But this means that there exists a vertex $t$ in $\Sigma$ so that $t\star u_1,t\star u_2$ 
are not in the same $N$-orbit. However, their endpoints are, and hence either the projections of the edges to $X/N$ yield 
loops in $X/N$, or the two projections together form a bigon. In both cases we contradict Item \ref{item:lift_gon} (for $m=1$ 
or $m=2$), since the loops/bigon cannot be lifted to $X$, which is a simplicial graph. For later purposes, we note that we 
also just proved:
 
 \begin{lemma}\label{lem:X_K_simplicial}
 $X/N$ is a simplicial graph.
\end{lemma}
 
 \textbf{Proof of $(b)$:} By $(a)$ we know that Item \ref{item:lift_geod_link} holds whenever $\Delta$ has at 
most $d-2$ vertices,  and from this we 
see that if we have a generalized $m$-gon $\tau=\tau_0,\dots,\tau_{m-1}$ with $d(\tau)\leq d-2$ then we can lift it to an 
``open'' generalized $m$-gon $\eta=\eta_0,\dots,\eta_{m-1}$, which is defined in the same way as a generalized $m$-gon except 
that we do not require $\eta^+_{m-1}=\eta^-_{0}$. We call $\eta^-_{0}$ and $\eta^+_{m-1}$ the initial and terminal marking. 
However, we still have that $\eta^+_{m-1}$ and $\eta^-_{0}$ are in the same $N$-orbit, since they are both lifts of 
$\tau^-_0$, say $\eta^+_{m-1}=g\eta^-_{0}$, for some $g\in N$. Assume that among all possible choices of lifts 
and elements $g\in N$ with $g\eta^-_{0}=\eta^+_{m-1}$, we picked one 
that minimizes the complexity $c(g)$ 
from Lemma \ref{lem:short_less_complex} or Lemma \ref{lem:short_less_complex_base}. If 
$g\eta^-_{0}=\eta^-_{0}$, we are done.

Otherwise, we use Lemma \ref{lem:short_less_complex} or Lemma \ref{lem:short_less_complex_base} (depending on whether $i>-1$ or not) to change the lift; let $Y\in\bbY_*$ and $\gamma_Y\in \Gamma_Y$ be as in the lemma for $\gamma=g$ and $\Delta=\eta^-_0$. In particular, the ``complexities'' satisfy $c(\gamma_Yg)<c(g)$.

First, suppose $\eta^-_0\subseteq \fix(\Gamma_Y)$. Then we can apply $\gamma_Y$ to all the lifts, contradicting minimality of $c(g)$ since $(\gamma_Y g) \gamma_Y\eta^-_0= \gamma_Y g\eta^-_0= \gamma_Y \eta^+_{m-1}$ ($\gamma_Y\eta^-_0$ and $\gamma_Y \eta^+_{m-1}$ being the new initial and terminal markings).

Second, suppose $g\eta^-_0=\eta^+_{m-1}\subseteq \fix(\Gamma_Y)$. Then $\eta^+_{m-1}=\gamma_Yg\eta^-_{0}$, again contradicting minimality of $c(g)$ (without even changing the lifts).

Lastly, suppose $d_Y(\eta^-_0,\eta^+_{m-1})> \frak T(\theta)$. Assume that some $\eta^+_k$ is contained in 
$\fix(\Gamma_Y)$.  In this case, we can replace the 
lifts $\eta_j$ for $j>k$ with the lifts $\gamma_Y\eta_j$, and get a 
new open generalized $m$-gon, $\gamma_Y g$ maps the initial marking to the terminal marking, contradicting minimality of $c(g)$.

Suppose instead that for each $k$ there is $x_k\in\eta^+_k$ and $x_{-1}\in\eta^-_{0}$ so that $x_k\notin \fix(\Gamma_Y)$. 
Then we have $\frak T(\theta)<d_Y(x_{-1},gx_{-1})\leq \sum d_Y(x_j,x_{j+1})$. By the choice of $\theta$ in 
Notation \ref{not:not_annuli} large compared to the BGI constant, we see that there must be a type G $\eta_k$, 
with corresponding simplex $\Delta_k$, and 
$v_k\in\eta_k$ so that $\Delta_k\cup \{v_k\}$ is contained $\fix(\Gamma_Y)$ (here we are using Proposition \ref{thm:MCG_quotient}.\eqref{item:BGI} together with Remark \ref{rem:fix_sat}). We can replace the lifts 
$\eta_j$ for $j>k$ with the lifts $\gamma_Y\eta_j$, as well as replacing the terminal path of $\eta_k$ starting at 
$v_k$ also with $\gamma_Y\eta_k$. Then, we conclude as before.
 \end{proof}

\subsection{Supporting lemmas for Theorem~\ref{thm:MCG_quotient}}\label{subsec:MCG_quotient_lemmas}
Recall that $X=\mathcal C(S)/N_{i-1}$, and $N$ is the kernel of the quotient $\bar G_{i-1}\to \bar G_i$, so that the map $\mathcal C(S)\to \mathcal C(S)/N_i=X/N$ factors as $\mathcal C(S)\to X\to X/N$.

\begin{conv}\label{conv:lift_all_the_way}
 Since generalized $m$-gons can be lifted from $X/N$ to $X$, and from $X$ to $\mathcal C(S)$, they can be 
lifted from $X/N$ to $\fontact(S)$. All the lifts in this 
subsection are of the latter type. When we cite Proposition \ref{prop:omni_lift} in this subsection, we 
will always use it together with lifting from $X$ to $\mathcal C(S)$.
\end{conv}

\begin{lemma}\label{lem:quotient_links}
 For every simplex $\Delta$ of $X/N$, and any simplex $\Sigma$ of $\mathcal C(S)$ that is a lift of $\Delta$  we have that 
$q(\link(\Sigma))=\link(\Delta)$. Moreover, if $i=-1$, then $\link(\Delta)=\link(\Sigma)/(N\cap Stab(\Sigma))$.
\end{lemma}

\begin{proof}
 Fix a lift $\Sigma$ of $\Delta$. Given a vertex $v$ of $\link(\Delta)$, we can lift $\Delta\star v$ to a 
simplex in $\mathcal C(S)$, and since there is a unique orbit of lifts of $\Delta$, we can choose the lift to 
contain $\Sigma$, and therefore be of the form $\Sigma\star\lift v$. Then $q(\lift v)=v$, and similarly we can 
show that edges of $\link(\Delta)$ arise from edges in the link of $\Sigma$.

Let us now prove the moreover part. We have to show that if two vertices $v,w$ of $\link(\Sigma)$ are 
$N$-translates, then they are $(N\cap \stabilizer(\Sigma))$-translates. This is because $v,w$ being 
$N$-translates implies that $\Sigma\star v$, $\Sigma\star w$ are lifts of the same simplex, and in particular 
the simplices are in the same $N$-orbit. This implies that there 
exists $h\in N$ that stabilizes $\Sigma$ and 
maps $v$ to $w$, as required.
\end{proof}

\begin{rem}[Connected links]\label{rem:quotient_links}
A consequence of Lemma  \ref{lem:quotient_links} is that all simplices of $X/N$ have connected links except 
co-dimension 1 faces in maximal simplices (since this holds in $\mathcal C(S)$).
\end{rem}

\begin{defn}[Approach path]\label{defn:simple_hierarchy}
An \emph{approach path} in $X/N$ is a sequence of paths $\gamma_1,\dots,\gamma_m$ and simplices 
$\Delta_1,\dots,\Delta_{m+1}$, so that
\begin{itemize}
\item the endpoint of $\gamma_j$ is the starting point of $\gamma_{j+1}$,
 \item $\gamma_j$ is a geodesic in the link of  a (possibly empty) simplex 
$\Delta_j$,
 \item the endpoint of $\gamma_m$ is in the link of $\Delta_{m+1}$,
 \item $\Delta_j$ is a proper sub-simplex of $\Delta_{j+1}$.
\end{itemize}
We say that the approach path starts (resp. ends) at $x$ if $x$ is the starting point of $\gamma_0$ (resp. endpoint of 
$\gamma_m$). We call $\Delta_{m+1}$ the \emph{terminal simplex}, and resulting path the concatenation of the $\gamma_j$.
\end{defn}

\begin{rem}\label{rem:max_depth}
 $m$ as above is bounded by $3\genus(S)+p(S)-3$ (that is, the complexity of $S$), which equals the maximal 
number of 
vertices of a simplex in $X/N$ by Proposition \ref{prop:omni_lift}.\ref{item:lift_splx}.
\end{rem}

\begin{lemma}
 Given a vertex $x\in X/N$ and a simplex $\Delta$ of $X/N$ so that $x\notin Sat(\Delta)$, there exists an approach path that starts at $x$ and has terminal simplex $\Delta'$ so that $[\Delta']=[\Delta]$.
\end{lemma}

\begin{proof}
 Consider $x$ and $\Delta$ as in the statement, and pick any geodesic $\gamma'_1$ in $X/N=\link(\emptyset)$ that intersects $Sat(\Delta)$ only at its endpoint $v_1$; notice that $[\Delta]\nest [v_1]$. Consider the subgeodesic $\gamma_1$ of $\gamma'_1$ obtained removing the last edge, and set $\Delta_2=v_1$. If $[\Delta_2]=[\Delta]$, we are done.
 
 Otherwise, inductively, suppose that we have an approach path $\gamma_1,\dots,\gamma_j$ starting at $x$ and terminating at $x_j\in \link(\Delta_{j+1})$ and $[\Delta]\propnest [\Delta_{j+1}]$. In particular, $\Delta_{j+1}$ has connected link (see Remark \ref{rem:quotient_links}), so that we can consider a shortest geodesic $\gamma'_{j+1}$ in $\link(\Delta_{j+1})$ to $\link(\Delta)$. If $\gamma'_{j+1}$ does not intersect $\sat(\Delta)$, we conclude by setting $\gamma_{j+1}=\gamma'_{j+1}$, and otherwise we can find an initial subgeodesic $\gamma_{j+1}$ of $\gamma'_{j+1}$ that intersects $Sat(\Delta)$ only at its endpoint $v_{j+1}$. We set $\Delta_{j+2}=\Delta_j\star v_{j+1}$, notice that we are done if $[\Delta_{j+2}]=[\Delta]$, and otherwise reapply the inductive procedure. This terminates by Remark \ref{rem:max_depth}.
\end{proof}

\begin{lemma}\label{lem:entrance_retraction}
 Let $\Delta$ be a simplex of $X/N$, and endow $\link(\Delta)$ with any metric induced by adding finitely many 
$\stabilizer(\link(\Delta))$--orbits of edges, as in Definition 
\ref{defn:hyp_H_space}. Let $\rho:X/N-Sat(\Delta)\to\link(\Delta)$ map the vertex $x$ to the endpoint 
of an arbitrary approach path that starts at $x$ and has terminal simplex equivalent to $\Delta$. Then $\rho$ is coarsely 
Lipschitz, and $\rho$ restricts to the identity on $\link(\Delta)$.
\end{lemma}

\begin{proof}
The fact that $\rho$ restricts to the identity on $\link(\Delta)$ is immediate from the definition of $\rho$.

 It is enough to show that adjacent vertices $v,w$ map uniformly close under $\rho$. To this end, form a $(2m+3)$-gon,  where 
$m\leq 3\genus(S)+p(S)-3$, using approach paths starting at $v$ and $w$, a single edge from $v$ to $w$, and 
 simplices 
$\Delta\star\rho(v),\Delta\star\rho(w)$.
 
 Consider a lift of this $(2m+3)$-gon (notice that $2m+3\leq 6\genus(S)+2p(S)-3$), which contains a lift 
$\Sigma$ of 
$\Delta$.  Notice that all sides of \textbf{type G} of the 
lifted $(2m+3)$-gon are geodesics in links of simplices of $\mathcal C(S)$ with strictly fewer vertices than $\Sigma$, and in particular 
they are geodesics in links of simplices not equivalent to $\Sigma$. Moreover, no vertex on a side of 
\textbf{type G} is in 
$Sat(\Sigma)$, since the image of a vertex of $Sat(\Sigma)$ is in $Sat(\Delta)$ (this is because such vertex is contained in 
a simplex $\Sigma'$ in $Sat(\Sigma)$ equivalent to $\Sigma$, and using Lemma \ref{lem:quotient_links} we see that 
$q(\Sigma')$ is equivalent to $\Delta$). In particular, by the Bounded Geodesic Image Theorem \cite[Theorem 3.1]{MM_II} (recall that we are in $\mathcal C(S)$) all sides have bounded subsurface 
projection to $\link(\Sigma)$, providing a bound on the distance between the lifts of $\rho(v)$ and $\rho(w)$ and hence on the 
distance between $\rho(v)$ and $\rho(w)$.  This concludes the proof.
\end{proof}

\subsection{Checking hierarchical hyperbolicity}\label{subsec:main_theorem}
We now have all the tools to prove the main conclusion of Theorem~\ref{thm:MCG_quotient}, namely hierarchical 
hyperbolicity of the quotient group.  We apply Convention \ref{conv:lift_all_the_way} about lifting to 
$\mathcal C(S)$ rather than to $X$ in this subsection as well, with exceptions clearly marked.

\begin{proof}[Proof of Theorem \ref{thm:MCG_quotient}.\eqref{item:HH_structure}]
We check that the action of $\bar{G}_i=\bar G_{i-1}/N$ on 
$X/N$ satisfies the hypotheses of Theorem \ref{thm:propaganda}. That is, we take $N_i$ to be the kernel of the map $MCG(S)\to \bar G_i$.

First, $X/N$ 
is simplicial by Lemma \ref{lem:X_K_simplicial}, and $\bar G_{i-1}/N$ acts on $X/N$ by simplicial automorphisms, and the action 
is cocompact since the action of $MCG(S)$ on $\mathcal C(S)$ is cocompact.
 
In view of Proposition \ref{prop:omni_lift}.\ref{item:lift_splx}, any maximal simplex $\Delta$ of $X/N$ is the  projection 
of some maximal simplex $\Sigma$ of $\mathcal C(S)$, which represents a unique $N_{i}$--orbit.  Since $\Sigma$ is maximal, 
$\stabilizer_{MCG(S)}(\Sigma)$ contains a finite-index abelian subgroup $A$ generated by powers by Dehn twists around the curves 
corresponding to the vertices of $\Sigma$.  Now, if $\bar g\in\stabilizer_{MCG(S)/N_{i}}(\Delta)$, let $g$ represent the left 
coset $\bar g$ of $N_{i}$, so $g\Sigma=h\Sigma$ for some $h\in N_{i}$.  In other words, $\bar g=\bar g'$ for some 
$g'\in\stabilizer_{MCG(S)}(\Sigma)$, i.e., $\stabilizer_{MCG(S)/N_{i}}(\Delta)$ is contained in the image of 
$\stabilizer_{MCG(S)}(\Sigma)$. Now, $A$ has finite image in $MCG(S)/N_{i}$, so since $g'$ represents one of finitely many cosets of 
$A$ in $\stabilizer_{MCG(S)}(\Sigma)$, we see that  $\stabilizer_{MCG(S)/N_{i}}(\Delta)$ is finite, as required.                                                                        It remains to check conditions~\eqref{item:stab_hyp}--\eqref{item:link_conn} from Theorem~\ref{thm:propaganda}.\\
 
\textbf{Proof of~\eqref{item:stab_hyp}:} Let $\Delta$ be a non-maximal simplex of $X/N$.  If $\Delta$ is not almost-maximal, 
then lift $\Delta$ to a simplex $\Sigma$ of $\mathcal C(S)$, note that $\link(\Sigma)$ is 
hyperbolic (since it is either a non-trivial join or the curve graph of a surface of complexity at least 2), 
and deduce 
that $\link(\Delta)$ is hyperbolic since we can lift triangles in $\link(\Delta)$ by Proposition 
\ref{prop:omni_lift}.\ref{item:lift_gon}.

\textbf{Almost-maximal $\Delta$ case:}  If $\Delta$ is almost-maximal, we divide into cases according to the value of $i$.  For $i=-1$, the group $\bar G_{-1}$ is obtained as a quotient by powers of Dehn twists.  The quotient $\bar G_0$ is equal to $\bar G_{-1}$ (see Remark~\ref{rem:0step}).  For $i>-1$, we are taking further proper quotients.\\ 

\emph{The case $i=-1$:}  Let $i=-1$.  Let $\Gamma_\Delta$ denote the image of $\stabilizer(\link(\Delta))$ in the group of permutations of $\link(\Delta)$.  We will show that $\Gamma_\Delta$ acts with finitely many orbits and with finite stabilizers on $\link(\Delta)$, and that $\Gamma_\Delta$ is a hyperbolic group.  Together, these facts imply that $\link(\Delta)$ is a hyperbolic $\stabilizer(\link(\Delta))$--space, as required.

We first show that $\Gamma_\Delta$ acts with finite stabilizers.  Since we are in the case $i=-1$ and $\Delta$ is almost-maximal, this will follow from Claim~\ref{claim:finite-stabs-stab-link} once we check that the three itemized assumptions in that claim hold in the present situation.  The first assumption is that Proposition~\ref{prop:technical_MCG_quotient}.\eqref{item:deep_BBF} holds for $i=-1$, but this is immediate from the choice of $N_{-1}$ in Notation~\ref{not:correct_powers}.  The second assumption is that Proposition~\ref{prop:technical_MCG_quotient}.\eqref{item:lift_to_CS} holds; this is because of Proposition~\ref{prop:omni_lift}, which also implies that the first part of Proposition~\ref{prop:technical_MCG_quotient}.\eqref{item:lift_splx_induction} holds.  So, we can apply Claim~\ref{claim:finite-stabs-stab-link} and conclude that $\Gamma_\Delta$ acts on $\link(\Delta)$ with finite stabilizers.

Next, we verify that $\Gamma_\Delta$ acts with finitely many orbits.  Using Proposition~\ref{prop:omni_lift}, we have a lift $\lift\Delta$ of $\Delta$ to $\mathcal C(S)$.   By Lemma~\ref{lem:quotient_links}, we have that $q(\link(\lift\Delta))=\link(\Delta)$, and therefore, since $\stabilizer_{MCG(S)}(\lift\Delta)\leq \stabilizer_{MCG(S)}(\link(\lift\Delta))$, we have $\phi(\stabilizer_{MCG(S)}(\lift\Delta))\leq \stabilizer(\link(\Delta))$.  Hence, letting $L$ denote the image of $\phi(\stabilizer_{MCG(S)}(\lift\Delta))$ in the permutation group of $\link(\Delta)$, we have $L\leq \Gamma_\Delta$.

Since $\Delta$ is almost-maximal in $\mathcal C(S)/N_{-1}$, the simplex $\lift\Delta$ is almost-maximal in $\mathcal C(S)$, so its $0$--skeleton is a multicurve whose complement has a single non-pants component, a complexity--$1$ subsurface of $S$ denoted $Y$.  The vertices in $\link(\lift\Delta)$ are the curves on $Y$, and the subgroup $\stabilizer_{MCG(S)}(\lift\Delta)$ acts on $Y$, with finitely many orbits of curves.  Hence $\phi(\stabilizer_{MCG(S)}(\lift\Delta))$ acts on $\link(\Delta)$ with finitely many orbits.  In other words, $L$, and hence $\Gamma_\Delta$, has finitely many orbits in $\link(\Delta)$.  

It remains to show that $\Gamma_\Delta$ is hyperbolic.  Now, we have shown that $L<\Gamma_\Delta$ acts with finitely many orbits on $\link(\Delta)$, and that $\Gamma_\Delta$, and hence also $L$, acts on $\link(\Delta)$ with finite stabilizers.  By construction, the actions of $L$ and $\Gamma_\Delta$ are faithful.  So $L$ has finite index in $\Gamma_\Delta$.  Hence, to show hyperbolicity of $\Gamma_\Delta$, we just have to show hyperbolicity of $L$.

Now, Lemma~\ref{lem:i=-2_rotating_family} implies that $\phi(\stabilizer_{MCG(S)}(\lift\Delta))$ is hyperbolic, so to get hyperbolicity of $L$, it is enough to show that the kernel of the action of $\phi(\stabilizer_{MCG(S)}(\lift\Delta))$ on $\link(\Delta)$ is finite.  But, by Lemma~\ref{lem:quotient_links}, $\link(\Delta)=\link(\lift\Delta)/(N\cap \stabilizer_{MCG(S)}(\lift\Delta))$, and so by Lemma~\ref{lem:i=-2_rotating_family}, $\stabilizer_{MCG(S)}(\lift\Delta)/(N\cap \stabilizer_{MCG(S)}(\lift\Delta))$ acts on $\link(\Delta)$ with finite point-stabilizers, as needed.\\

\emph{The case $i=0$:}  For $i=0$, we have that $N_0=N_{-1}$, and $G_0=G_{-1}$.  As before, $\Delta$ is an almost-maximal simplex; all that has changed is the relationship between $i$ and the complexity of the subsurface $Y$ obtained by lifting $\Delta$ to $\mathcal C(S)$: that complexity is still $1$, but now $1=i+1$.  However, the exact same argument as in the case $i=-1$ shows that $\Gamma_\Delta$ is a hyperbolic group acting with finitely many orbits and finite stabilizers on $\link(\Delta)$, so again $\link(\Delta)$ is a hyperbolic $\stabilizer(\link(\Delta))$--space, as required.\\ 

\emph{The case $i=1$:} Suppose that $i=1$.  Since $\Delta$ is almost-maximal,  the associated subsurface $Y$ obtained above as the non-pants component of the complement of the multicurve corresponding to a lift $\lift\Delta$ of $\Delta$ still has complexity $1$.  

Let $\Delta'$ be the image of $\lift\Delta$ under the quotient $\mathcal C(S)\to\mathcal C(S)/N_0=X$.  Then $\Delta'$ is a lift of $\Delta$ to $X$; this is one of the exceptions to the convention on lifting all the way to $\mathcal C(S)$.

Now, by induction, we can assume that the action of $\bar G_0$ on $\mathcal C(S)/N_0$ satisfies all the conclusions of Proposition~\ref{prop:technical_MCG_quotient}, and in particular item \eqref{item:bottom_hyp}, which says in particular that $\pstab(\sat(\Delta'))$ acts with finitely many orbits of vertices on $\link(\Delta')$, since $Y$ has complexity $1$ and $\lift\Delta$ is also a lift of $\Delta'$ to $\mathcal C(S)$.

Now, the kernel $N$ of $\bar G_0\to\bar G_1$ contains a finite-index subgroup of $\pstab(\sat(\Delta'))$.  Indeed, $[\Delta']\in\bbY$ since $Y$ has complexity $1$ (see the beginning of Section~\ref{subsec:composite_induction} for the definition of $\bbY$, which we are applying in the case where $X=\mathcal C(S)/N_0$).  By the definition of $N$ in Notation~\ref{not:not_annuli}, $N$ therefore contains a finite-index subgroup $\Gamma^\theta_{[\Delta']}$ of $\pstab(\sat(\Delta'))$.  By the previous part of the argument, it follows that $N\cap \stabilizer(\link(\Delta'))$ acts with finitely many orbits of vertices on $\link(\Delta')$, and therefore the image of  $\link(\Delta')$ under the quotient map $X\to X/N$ is finite.  Since the conclusion of Lemma~\ref{lem:quotient_links} about links mapping to links also applies to $X\to X/N$, we conclude that $\link(\Delta)$ is finite.  In particular, it is a hyperbolic $\stabilizer(\link(\Delta))$--space.

\emph{The case $i>1$:}  Now suppose that $i>1$.  Since $\Delta$ is almost-maximal, the subsurface $Y$ has complexity $1\leq i-1$, so $\link(\Delta)$ has finite vertex set by the fact that the third bullet of 
Theorem \ref{thm:MCG_quotient}.\eqref{item:HH_structure} holds inductively.\\

 This completes the proof of the $\stabilizer(\link(\Delta))$--hyperbolicity clause of Theorem~\ref{thm:propaganda}.\eqref{item:stab_hyp} about $\stabilizer(\link(\Delta))$--hyperbolicity.  The quasi-isometric embedding clause follows from the existence of the coarse retraction $\rho$ provided by Lemma \ref{lem:entrance_retraction}.\\
 
\textbf{Proof of~\eqref{item:link_conn}:} This is Remark~\ref{rem:quotient_links}.\\
 
\textbf{Proof of~\eqref{item:propaganda_join}:} Let $\Sigma,\Delta$ be simplices of $X/N$.  Recall 
that we need to show 
that there exist simplices $\Pi,\Pi'$ of $\link(\Delta)$ such that
\begin{enumerate}[(i)]
     \item $\link(\Delta)\cap\link(\Sigma)=\link(\Delta\star\Pi)\star\Pi'.$\label{item:B_restate}
\end{enumerate}

We will in fact prove the following:

\begin{claim}\label{claim:2_bullets}
 For each $\Delta,\Sigma$ one of the following holds:
\begin{itemize}
 \item There exists a vertex $v$ in $\link(\Delta)$ so that $\link(\Delta)\cap\link(\Sigma)\subseteq Star(v)$, or
 \item $\link(\Delta)\subseteq \link(\Sigma)$.
\end{itemize}
\end{claim}

{\bf Induction on co-level.}
We now show how to conclude the proof given the claim.

First, notice that \eqref{item:B_restate} holds whenever $\Delta$ is maximal.  Consider some pair $\Delta,\Sigma$, and suppose that \eqref{item:B_restate} holds for any pair $\Delta',\Sigma'$ for which 
$\Delta'$ has strictly more vertices than $\Delta$. If the second bullet holds, then we can set $\Pi=\Pi'=\emptyset$, and 
\eqref{item:B_restate} holds for $\Delta,\Sigma$. If the first bullet holds, then consider the simplices $\Pi_0,\Pi'_0$ 
obtained from \eqref{item:B_restate} applied to $\Delta\star v, \Sigma$. If $v\in\link(\Sigma)$, then we can set 
$\Pi=\Pi_0\star v$, $\Pi'=\Pi'_0\star v$. If not, we set $\Pi=\Pi_0\star v$, $\Pi'=\Pi'_0$. In either case, \eqref{item:B_restate} holds, we spell out the first case, the other being very similar:
$$\link(\Delta)\cap\link(\Sigma)=(\link(\Delta)\cap Star(v))\cap\link(\Sigma)=(\link(\Delta\star v)\star v)\cap \link(\Sigma)=$$
$$(\link(\Delta\star v)\cap\link(\Sigma))\star v= \link(\Delta\star v\star \Pi_0)\star\Pi'_0\star v.$$

It remains to prove the claim:

\renewcommand{\qedsymbol}{$\blacksquare$}
\begin{proof}[Proof of Claim \ref{claim:2_bullets}]
If $\Delta$ is maximal, the second bullet holds, so we assume that this is not the case. Also, if $\link(\Delta)\cap\link(\Sigma)=\emptyset$, then we can take as $v$ any vertex in $\link(\Delta)$.  So, from now on we assume $\link(\Delta)\cap\link(\Sigma)\neq\emptyset$.

Let $\Lambda$ be a maximal simplex of $\link(\Delta)\cap\link(\Sigma)$.  Let $\lift\Delta\star\lift\Lambda$ be a lift of the 
simplex $\Delta\star\Lambda$, provided by Proposition~\ref{prop:omni_lift}. In terms of curves, $(\lift\Lambda)^{(0)}\cup(\lift\Delta)^{(0)}$ is a multicurve, and as such it can be  completed to a 
pants decomposition $M_{\Lambda,\Pi}=(\lift\Lambda)^{(0)}\cup(\lift\Delta)^{(0)}\cup(\lift\Pi)^{(0)}$.

{\bf Reducing to $\Pi=\emptyset$.}
Let $v\in\link(\Delta)\cap\link(\Sigma)$. We claim that $v\in\link(\Delta\star\Pi)$. We can  lift the generalized $4$-gon 
$\Delta\star\Lambda,\Lambda\star \Sigma,\Sigma\star v, v\star\Delta$ to $\mathcal C(S)$, obtaining the lifts $\lift\Delta,$ etc., and in 
fact we can assume $\lift\Delta,\lift\Lambda$ coincide with the previously chosen lifts. If we had $v\notin 
\link(\Delta\star\Pi)$, then we would also have $\lift v\notin \link(\lift\Delta\star\lift\Pi)$, and also $\lift 
v\notin\lift\Lambda$. But $\lift v\in\link(\lift\Delta)$, so $\lift v\notin\link(\lift\Pi)$. 
Moreover, by maximality of $\Lambda$, we also have $\lift v\notin \link(\lift\Lambda)$. Hence, as a curve, $\lift v$ 
intersects both $(\lift\Lambda)^{(0)}$ and $(\lift\Pi)^{(0)}$ (but not $(\lift\Delta)^{(0)}$ nor $(\lift\Sigma)^{(0)}$). We now 
make a multicurve $\sigma$ from $\lift v$, which we assume to be in minimal position with respect to $M_{\Lambda,\Pi}$, by 
considering the boundary of a regular neighborhood of $\lift v\cup \lift\Lambda^{(0)}$. Notice that $\sigma$ is disjoint 
from $\lift\Sigma^{(0)}$, $\lift\Lambda^{(0)}$, and $\lift\Delta^{(0)}$. Also, at least one component $\sigma_0$ of $\sigma$ 
is not parallel to $\lift\Lambda^{(0)}$ (since it intersects $\lift\Pi^{(0)}$ non-trivially), and this contradicts the 
maximality of $\Lambda$. This proves that $\link(\Delta)\cap\link(\Sigma)\subseteq \link(\Delta\star\Pi)$.
If $\Pi$ is non-empty, the first bullet holds with $v$ any vertex of $\Pi$. Hence, we now assume $\Pi=\emptyset$.

{\bf The case $\Pi=\emptyset$.}
We now know that $(\lift\Lambda)^{(0)}\cup(\lift\Delta)^{(0)}$ is a pants decomposition,  and that in fact this holds for any 
maximal simplex $\Lambda$ of $\link(\Delta)\cap\link(\Sigma)$. If all such maximal simplices share a vertex $v$, then the 
first bullet holds for this $v$.

Otherwise, for each 
vertex $w$ of a fixed maximal simplex $\Lambda$ in $\link(\Delta)\cap\link(\Sigma)$, we can find another such simplex 
$\Theta$ not sharing $w$ with $\Lambda$. We claim that for each curve $\delta$ in $\lift\Delta^{(0)}$ which is the boundary curve of a 
pair of pants not all of whose boundary curves are in $\lift\Delta^{(0)}$, $\Sigma$ has a lift avoiding that curve. In fact, 
we can lift a generalized $4$-gon $\Delta\star\Lambda,\Lambda\star \Sigma,\Sigma\star \Theta, \Theta\star\Delta$, where 
$\Theta$ and $\Lambda$ do not share a vertex corresponding to one of the boundary curves of a pair of pants as 
above. Then, 
the lift of $\Sigma$ will not intersect $\delta$, for otherwise it will have to intersect some curve either in the lift of 
$\Lambda$ or in the lift of $\Theta$.

Let now $\Delta_0$ be the sub-simplex of $\Delta$ consisting of all  vertices $q(\delta)$ for $\delta$ as above. Then each 
vertex of $\Sigma$ is in $Star(\Delta_0)$. Hence, $\Delta_0,\Lambda$, and $\Sigma$ are contained in a common simplex, which 
can be lifted to a simplex in $\mathcal C(S)$. We can also arrange the corresponding lifts of $\Delta_0$ and $\Lambda$ so that the lift 
of $\Lambda$ is $\lift\Lambda$ and the lift of $\Delta_0$ has vertex set consisting of all curves $\delta$ as 
above.

This gives a lift $\lift\Sigma'$ of 
$\Sigma$ which is contained in the union of the pairs of pants in the complement of 
$(\lift\Lambda)^{(0)}\cup(\lift\Delta)^{(0)}$, all of whose boundary curves are in $\lift\Delta^{(0)}$. This implies that 
$\link(\lift\Delta)\subseteq \link(\lift\Sigma')$, and hence $\link(\Delta)\subseteq \link(\Sigma)$, that is, the second 
bullet.

This concludes the proof of Claim \ref{claim:2_bullets}.
\end{proof} 
\renewcommand{\qedsymbol}{$\Box$}

We now prove the statement about the index set $\mathfrak S_{N_i}$ of the HHS structure on $\bar G_i$. Recall that $\mathfrak S_{N_i}$  is 
the set of equivalence classes of non-maximal simplices.
We first recall the bijection $b$.  Given any non-maximal simplex $\Delta$ of $X/N$, consider a lift $\lift\Delta$ to $\mathcal C(S)$. The vertex set of the 
link of $\lift\Delta$ in $\mathcal C(S)$ consists of all curves (regarded as vertices of $\mathcal C(S)$) contained in a subsurface that we denote $S_{\lift\Delta}$. Define $b([\Delta])=[S_{\lift\Delta}]_{N_{i}}$, where $[\cdot]_{N_{i}}$ denotes the $N_{i}$--orbit. Notice that choosing a different lift yields a subsurface in the same $N_{i}$--orbit, since all lifts of $\Delta$ are in the same $N_{i}$--orbit.

We now complete the proof that $b$ is well-defined.  What is left to prove is that equivalent simplices yield the same orbit.

\begin{claim}\label{claim:no_link_join}
 No link is a join of a non-empty simplex and some subcomplex.  In particular, no link consists of a single non-empty simplex.
\end{claim}

\renewcommand{\qedsymbol}{$\blacksquare$}
\begin{proof}[Proof of Claim~\ref{claim:no_link_join}]
Consider 
the link of $\Delta$ 
and some vertex $v\in\link(\Delta)$, and let us show that there is a vertex $w$ of $\link(\Delta)$ not 
connected to, or equal to, $v$. Consider a lift $\lift\Delta\star\lift v$ of $\Delta\star v$. By Lemma \ref{lem:fi}.\eqref{item:multi-orbit} there 
exists a vertex $\lift w$ of $\link(\lift\Delta)$ so that $\lift v$ and $\lift w$ are not in the same 
$H$--orbit, and no $H$--translate of $\lift w$ is adjacent to $\lift v$. Since $N_{i}<H$, this means that the 
image $w$ of $\lift w$, which is in $\link(\Delta)$, is distinct from $v$ and not connected to $v$, as 
required.
\end{proof}
\renewcommand{\qedsymbol}{$\Box$}

In view of condition \ref{item:propaganda_join}, we now see that $\link(\Sigma)\subseteq \link(\Delta)$ if 
and only if there exists a simplex $\Pi$ in $\link(\Delta)$ so that $\link(\Sigma)=\link(\Delta\star \Pi)$. This also holds 
in $\mathcal C(S)$, and hence we get that $[\Delta]\nest[\Sigma]$ if and only if there are $S_{\lift\Delta}, S_{\lift\Sigma}$ as above 
with $S_{\lift\Delta}$ nested into $S_{\lift\Sigma}$. This implies that if $\link(\Delta)=\link(\Sigma)$ then 
$S_{\lift\Delta}\subseteq S_{\lift\Sigma}\subseteq gS_{\lift\Delta}$ for some $g\in N_{i}$, and since $S_{\lift\Delta}, 
gS_{\lift\Delta}$ have the same complexity they need to coincide, showing $S_{\lift\Delta}= S_{\lift\Sigma}$, as we wanted.

Notice that we also showed that $[\Delta]\nest[\Sigma]$ if and only if the corresponding orbits contain nested 
representatives. We are only left to show the analogous statement for orthogonality/disjointness, which we will reduce to 
the nesting statement.

Consider a non-maximal simplex $\Delta$ of $X/N$, and a lift $\lift\Delta$. Then $\lift\Delta$ contains the boundary 
multicurve $(\lift\Delta_0)^{(0)}$ of $S_{\lift\Delta}$. Denote $\Delta_0=q(\lift\Delta_0)$, and pick any maximal simplex 
$\Lambda$ in $\link(\Delta)$. We claim that a simplex $\Sigma$ satisfies $[\Delta]\orth[\Sigma]$ if and only if 
$[\Sigma]\nest[\Delta_0\star\Lambda]$. This implies the orthogonality/disjointness statement, in view of the nesting 
statement.

Fix a lift $\lift\Lambda$ of $\Lambda$ contained in $\link(\hat \Delta)$, and notice that the vertex set forms a pants 
decomposition of $S_{\lift\Delta}$. If $[\Sigma]\nest[\Delta_0\star\Lambda]$, then there is a lift $\lift\Sigma$ so that 
$\link(\lift\Sigma)\subseteq \link(\lift\Delta_0\star\lift\Lambda)$. Any curve disjoint from 
$(\lift\Delta_0\star\lift\Lambda)^{(0)}$ is disjoint from curves in $S_{\lift\Lambda}$, so that $\link(\lift\Sigma)$ and 
$\link(\lift\Delta)$ form a join. The same then holds for $\Delta,\Sigma$, showing $[\Delta]\orth[\Sigma]$, as required.

Suppose now that $[\Delta]\orth[\Sigma]$. Our goal is to show that any vertex $v\in\link(\Sigma)$ lies in 
$\link(\Delta_0\star\Lambda)$. For later use, note that there is a well-defined simplex $\Lambda\star v$ in view of the definition of orthogonality, since $\Lambda\subset\link(\Delta)$.

We claim that we can find another maximal simplex $\Theta$ in $\link(\Delta)$ so that any lift of $\Theta$ contained in $\link(\lift\Delta)$ has vertex set which, together with $(\lift\Lambda)^{(0)}$ fills $S_{\lift\Delta}$. We will do so by showing that we can find $\Theta$ in $\link(\Delta)$ such that for every vertex $v$ of $\Lambda$ there exists a vertex $w$ of $\Theta$ such that $v$ and $w$ are not connected. In particular, as curves, any lifts of $v$ and $w$ intersect, from which we see that the vertices of any lift of $\Theta$ contained in $\link(\lift\Delta)$, as curves, intersect all curves of $(\lift\Lambda)^{(0)}$. Since $(\lift\Lambda)^{(0)}$ is a pants decomposition of $S_{\lift\Delta}$, this yields the required filling statement. Now, by Claim \ref{claim:no_link_join}, $\link(\Delta)$ does not consist of the maximal simplex $\Lambda$ only, so that there exists a vertex $w_1$ in $\link(\Delta)$ which is not connected to some vertex of $\Lambda$ (recall that $\Lambda$ is a maximal simplex of $\link(\Delta)$). If $w_1$ is not connected to any vertex of $\Lambda$, we can simply complete $w_1$ to a maximal simplex $\Theta$ of $\link(\Delta)$. If not, let $\Lambda_2=\Lambda\cap \link(w_1)$. Again by Claim \ref{claim:no_link_join}, used as above, there exists $w_2\in\link(\Delta\star w_1)$ which is not connected to some vertex of $\Lambda_2$. If it is not connected to any vertex, we can take $\Theta$ with vertex set containing $w_1,w_2$, and otherwise we continue for finitely many steps.

 Consider now any vertex $v\in\link(\Sigma)$, which we need to show is in $\link(\Delta_0\star\Lambda)$. We 
can lift a generalized 4-gon $\Delta_0\star\Lambda,\Lambda\star v,v\star\Theta,\Theta\star\Delta_0$, with the lift 
of $\Delta_0\star\Lambda$ being $\lift\Delta_0\star\lift\Lambda$. We then see that, as a curve, the lift 
$\lift v$ is disjoint from $S_{\lift\Delta}$ (since $\lift\Theta$ and $\lift\Lambda$ fill $S_{\lift\Delta}$), so that $v$ lies either in $\link(\Delta_0\star\Lambda)$ or in $\Delta_0$. 

We are left to argue that $v$ does not lie in $\Delta_0$. Suppose that this was the case. Again by \ref{claim:no_link_join}, there is a vertex $w$ in $\link(\Sigma)$ which is not connected to $v$), and in particular any lift $\lift w$ intersects $\lift v$. Applying the above argument with $w$ replacing $v$, we would find a lift $\lift w$ that cuts $S_{\lift\Delta}$ since it intersects its boundary, a contradiction.

Finally, we complete the proof of Theorem \ref{thm:MCG_quotient}.\eqref{item:HH_structure} by proving that if 
$b([\Delta])=N_iY$ for $Y$ of complexity at most $i$, then $\fontact(\Delta)$ has finite vertex set. This 
suffices to get a uniform diameter bound since there are finitely many orbits of simplices. 

\emph{This is the only place in this subsection where lifts are taken to be in $X$ rather than in $\mathcal 
C(S)$.} First, if $b([\Delta])=N_iY$ for $Y$ of complexity less than $i$, then the link of any lift 
$\lift\Delta$ of $\Delta$ in $X$ has finitely many vertices, by the inductive hypothesis. Since the conclusion 
of Lemma \ref{lem:quotient_links} about links mapping to links applies to the map $X\to X/N$ as well, with the 
same proof, we conclude that the link of $\Delta$ has finite vertex set, as required.

Suppose now $b([\Delta])=N_iY$ for $Y$ of complexity $i$. Consider a lift $\lift\Delta$ of $\Delta$ in $X$. By 
inductive hypothesis (specifically, Proposition \ref{prop:technical_MCG_quotient}.\eqref{item:bottom_hyp}), 
$\pstab(\sat(\lift\Delta))$ acts cocompactly on $\link(\lift\Delta)$. This readily implies that 
$\pstab(\sat(\Delta))$ acts cocompactly on $\link(\Delta)$, by Lemma \ref{lem:quotient_links} (used as above). 
Now, $\pstab(\sat(\Delta))$ is a quotient of $\pstab(\sat(\lift\Delta))$ by a finite-index subgroup by construction 
of the composite rotating family, and we are done. 
\end{proof}

We conclude this subsection with the proof Proposition \ref{prop:technical_MCG_quotient}.\eqref{item:lift_splx_induction}, which we now have the tools to prove:

\begin{lemma}\label{lem:quotient_sats}
 Let $\Delta$ be a non-maximal simplex of $X/N$, and consider a lift $\lift\Delta$ to $\mathcal C(S)$. Then $q(\sat(\lift\Delta))=\sat(\Delta)$.
\end{lemma}

\begin{proof}
 We use that the bijection $b$ defined in Definition \ref{defn:b} is well-defined. Let $v\in\sat(\Delta)$. Then 
$v$ is a vertex of some $\Delta'$ with the same link as $\Delta$. 

Since $b$ is well-defined, $\Delta'$ must have a lift $\lift\Delta'$ such that the vertex set of 
$\link(\lift\Delta')$ consists of all curves contained $S_{\lift\Delta}$. That is, $\lift\Delta'$ has the same 
link as $\lift \Delta$, so that $\lift\Delta'\subseteq \sat(\lift\Delta)$ and 
$v\in\Delta'=q(\lift\Delta')\subseteq q(\sat(\lift\Delta))$. We just proved $\sat(\Delta)\subseteq 
q(\sat(\lift\Delta))$.
 
 Let $v\in q(\sat(\lift\Delta))$. Then $v=q(\lift v)$ for some vertex $\lift v$ of a simplex $\lift\Delta'$ 
with the same link as $\lift\Delta$. This implies that $\Delta'=q(\lift\Delta')$ has the same link as 
$\Delta$, by Lemma \ref{lem:quotient_links}. Hence $v\in\sat(\Delta)$, so $
q(\sat(\lift\Delta))\subseteq \sat(\Delta)$.
\end{proof}

\subsection{Preservation properties}\label{subsec:preserve}

\begin{lemma}\label{lem:some_rotation}
 Let $\gamma\in N$ and $x\in X^{(0)}$. Then either $\gamma x=x $ or there exists $Y\in \bbY_*$ so that 
$\dist_Y(x,\gamma x)>\frak T(\theta)$.
\end{lemma}

\begin{proof}
Since $\mathfrak C$ is well-ordered, we can argue by (transfinite) induction on $c(\gamma)$.  The statement 
holds for $\gamma=1$, that is, for the minimal element of $\mathfrak C$.

Suppose $\gamma\neq 1$ and suppose $\gamma x\neq x$. Let 
$\gamma_Y$ be as in Lemma 
\ref{lem:short_less_complex}/Lemma \ref{lem:short_less_complex_base}.  If the second conclusion of Lemma 
\ref{lem:short_less_complex}/Lemma \ref{lem:short_less_complex_base} applies, then we are done.  

Otherwise, suppose $x\in Fix(\Gamma_Y)$. Then, $\gamma_Y\gamma x\neq x$ (since $\gamma_Y^{-1}x=x$), and 
$c(\gamma_Y\gamma)<c(\gamma)$, so by induction there exists $W\in\bbY_*$ so that $\dist_W(x,\gamma_Y\gamma 
x)>\frak T(\theta)$.  Hence, we get $\dist_{\gamma_Y^{-1}W}(x,\gamma x)>\frak T(\theta)$ (we used 
$\gamma_Y^{-1}x=x$ 
again).

The case $\gamma x\in Fix(\Gamma_Y)$ is similar: $\gamma_Y\gamma x\neq x$ since $\gamma_Y\gamma x=\gamma x$, 
and $c(\gamma_Y\gamma)<c(\gamma)$, so by induction there exists $W\in\bbY_*$ so that $\dist_W(x,\gamma_Y\gamma 
x)>\frak T(\theta)$.  Hence, $\dist_W(x,\gamma x)>\frak T(\theta)$, again because $\gamma_Y\gamma x=\gamma x$.
\end{proof}

\begin{proof}[Proof of Theorem \ref{thm:MCG_quotient}.\eqref{item:larg_inj_radius} and Proposition \ref{prop:technical_MCG_quotient}.\eqref{item:hyperelliptic}]
 Let $f,g\in F$ be distinct. We have to prove that $f^{-1}g\notin N_{i}$. If $f^{-1}g$ has finite order, then 
$f^{-1}g$ cannot be in $H$, since $H$ is torsion-free, and hence in particular not in $N_{i}$. If not, let 
$x\in X$ so that $f(x)\neq g(x)$, which exists by induction for $i>-1$, while for $i=-1$ it exists because 
infinite order elements of $MCG(S)$ act non-trivially on $\mathcal C(S)$. Then by Lemma 
\ref{lem:some_rotation}, if we had $f^{-1}g\in N_{i}$ we would have some $Y\in\bbY_*$ so that 
$\dist_Y(f(x),g(x))>\frak T(\theta)$.  This contradicts the choice of $\theta$ in Notation \ref{not:not_annuli}.
\end{proof}

\begin{lemma}\label{lem:cobounded_geod}
Let $\kappa$ be as in Lemma~\ref{lem:conv_cocpt=cobounded}. Let $x,y\in X$ be so that 
$\dist_{[\Delta]}(x,y)\leq \kappa$ for all 
non-empty non-maximal simplices $\Delta$ of $X$. Let $[x,y]$ be a geodesic from $x$ to $y$.  Then $q|_{[x,y]}$ 
is an isometric embedding.
\end{lemma}

\begin{proof}
 Consider a lift $[x,y']$ to $\fontact(S)$ of a geodesic $[\bar x,\bar y]$ in $X/N$ connecting the 
images $\bar x,\bar y$ of $x,y$. Then $y'=\gamma y$ for some $\gamma\in N$. Since $\mathfrak C$ is 
well-ordered, we can choose $\gamma$ to have minimal $c(\gamma)$ among: 
\begin{itemize}
 \item the $N_i$-orbits of the pair $(x,y)$ 
and of the geodesic $[x,y]$,
\item all lifts $[x,y']$ of $[\bar x,\bar y]$,
\item   all choices of $\gamma$ with $y'=\gamma 
y$.  
\end{itemize}

We claim that $\gamma=1$, which will show that $\dist_X(x,y)=\dist_{X/N}(\bar x,\bar y)$, which readily 
implies the desired conclusion.
 
 Suppose $\gamma\neq 1$. Consider $Y\in\bbY_*$ and $\gamma_Y\in \Gamma_Y$ with $c(\gamma_Y\gamma)<c(\gamma)$ as 
in Lemma \ref{lem:short_less_complex}/Lemma \ref{lem:short_less_complex_base}.  There are three cases.

\begin{itemize}
 \item If $y\in \fix(\Gamma_Y)$, then we can replace $[x,y]$ with $\gamma_Y[x,y]$, and $[x,y']$ with 
$\gamma_Y[x,y']$. Then $\gamma_Y\gamma$ maps the second 
endpoint of the original geodesic (i.e., $\gamma_Y y=y$) to the second endpoint of the lift (i.e., $\gamma_Yy'$), meaning $\gamma_Y y'=\gamma_Y\gamma y$.  This 
contradicts minimality of $c(\gamma)$.

\item If $y'\in Fix(\Gamma_Y)$, then $y'=\gamma_Yy'= \gamma_Y\gamma y$, contradicting minimality of $\gamma$.

\item Otherwise, $\dist_Y(y',y)>\frak T(\theta)$. If $x\in 
\fix(\Gamma_Y)$, we can replace $[x,y]$ with $\gamma_Y[x,y]$, 
and $[x,y']$ with $\gamma_Y[x,y']$, and we contradict minimality of $c(\gamma)$ since $\gamma_Y 
y'=\gamma_Y\gamma y$. Otherwise, $\dist_Y(x,y)$ is well-defined and $\leq \kappa$. By the choice of $\theta$ 
(much larger than $\kappa$) in Notation \ref{not:annuli} or Notation \ref{not:not_annuli}, we get that 
$\dist_Y(x,y')$ is well-defined and large enough that $[x,y']$ intersects $\fix(\Gamma_Y)$ (by the defining 
property of $C_{i}$, see Proposition \ref{prop:technical_MCG_quotient}.\eqref{item:BGI}).
 
We can now replace a terminal subgeodesic $[v,y']$ of the lift $[x,y']$ starting at some $v\in\fix(\Gamma_Y)$ 
by its translate $\gamma_Y[v,y']$, thereby obtaining a new lift. This contradicts minimality of $c(\gamma)$ 
since $\gamma_Yy'=\gamma_Y\gamma y$.
\end{itemize}
This completes the proof of the lemma.
\end{proof}

We conclude by proving the remaining statement in Theorem~\ref{thm:MCG_quotient}.

\begin{proof}[Proof of Theorem \ref{thm:MCG_quotient}.\eqref{item:poc}]
Recall that we need to prove that the map $\phi$ restricts to an injective map on $Q$, and orbit maps from $Q$ 
to $\mathcal C(S)/N_i$ are quasi-isometric embeddings.

We first prove the quasi-isometric embedding statement. By the choice of $\kappa$ (coming from Lemma \ref{lem:conv_cocpt=cobounded}) and Lemma 
\ref{lem:cobounded_geod}, for any $x_0\in X$ and $g\in Q$, we have that any geodesic $[x_0,gx_0]$ in $X$ 
projects to a geodesic $[\bar x_0,g\bar x_0]$ in $X/N$. Since the length of $d_X(x_0,gx_0)$ is comparable up to 
multiplicative and additive constants with the word length of $g$, by induction, the same holds for 
$\dist_{X/N}(\bar x_0,g\bar x_0)$. This suffices to show that $Q$-orbit maps are quasi-isometric embeddings.

To show the injectivity statement, note that $Q$-orbit maps being quasi-isometric embeddings implies that the 
kernel of $\phi|_Q$ is finite. But this implies that the kernel of $\phi|_Q$ must be trivial since it 
is contained in $N_i$, whence in $H$, which is torsion-free by construction (see Lemma \ref{lem:fi}). Hence, 
$\phi|_Q$ is injective.
\end{proof}

\bibliographystyle{alpha}
\bibliography{links}

\newcommand{\etalchar}[1]{$^{#1}$}
\begin{thebibliography}{BFG{\etalchar{+}}07}

\bibitem[AAS07]{AAS:MCG_not_rel_hyp}
James~W. Anderson, Javier Aramayona, and Kenneth~J. Shackleton.
\newblock An obstruction to the strong relative hyperbolicity of a group.
\newblock {\em J. Group Theory}, 10(6):749--756, 2007.

\bibitem[AB23]{AbbottBehrstock:conjugator}
Carolyn Abbott and Jason Behrstock.
\newblock Conjugator lengths in hierarchically hyperbolic groups.
\newblock {\em Groups Geom. Dyn.}, 17(3):805--838, 2023.

\bibitem[ABD21]{ABD}
Carolyn Abbott, Jason Behrstock, and Matthew~Gentry Durham.
\newblock Largest acylindrical actions and stability in hierarchically
  hyperbolic groups.
\newblock {\em Trans. Amer. Math. Soc. Ser. B}, 8:66--104, 2021.
\newblock With an appendix by Daniel Berlyne and Jacob Russell.

\bibitem[AF19]{AramayonaFunar}
Javier Aramayona and Louis Funar.
\newblock Quotients of the mapping class group by power subgroups.
\newblock {\em Bull. Lond. Math. Soc.}, 51(3):385--398, 2019.

\bibitem[AGM09]{AGM}
Ian Agol, Daniel Groves, and Jason~Fox Manning.
\newblock Residual finiteness, {QCERF} and fillings of hyperbolic groups.
\newblock {\em Geom. Topol.}, 13(2):1043--1073, 2009.

\bibitem[Ago13]{vhak}
I.~Agol.
\newblock The virtual {H}aken conjecture.
\newblock {\em Doc. Math.}, 18:1045--1087, 2013.
\newblock With an appendix by Agol, Daniel Groves, and Jason Manning.

\bibitem[ANS{\etalchar{+}}24]{AbbottNgSpriano:uniform}
Carolyn~R. Abbott, Thomas Ng, Davide Spriano, Radhika Gupta, and Harry Petyt.
\newblock Hierarchically hyperbolic groups and uniform exponential growth.
\newblock {\em Math. Z.}, 306(1):Paper No. 18, 33, 2024.

\bibitem[BB01]{genus_2_linear}
Stephen~J. Bigelow and Ryan~D. Budney.
\newblock The mapping class group of a genus two surface is linear.
\newblock {\em Algebr. Geom. Topol.}, 1:699--708, 2001.

\bibitem[BBF15]{BBF}
Mladen Bestvina, Ken Bromberg, and Koji Fujiwara.
\newblock Constructing group actions on quasi-trees and applications to mapping
  class groups.
\newblock {\em Publ. Math. Inst. Hautes \'{E}tudes Sci.}, 122:1--64, 2015.

\bibitem[BDM09]{BDM:thick}
Jason Behrstock, Cornelia Dru\c{t}u, and Lee Mosher.
\newblock Thick metric spaces, relative hyperbolicity, and quasi-isometric
  rigidity.
\newblock {\em Math. Ann.}, 344(3):543--595, 2009.

\bibitem[BF92]{BF:combination}
M.~Bestvina and M.~Feighn.
\newblock A combination theorem for negatively curved groups.
\newblock {\em J. Differential Geom.}, 35(1):85--101, 1992.

\bibitem[BFG{\etalchar{+}}07]{BFGM}
Uri Bader, Alex Furman, Tsachik Gelander, Nicolas Monod, et~al.
\newblock Property ({T}) and rigidity for actions on banach spaces.
\newblock {\em Acta mathematica}, 198(1):57--105, 2007.

\bibitem[BHS17a]{HHS_III}
Jason Behrstock, Mark~F. Hagen, and Alessandro Sisto.
\newblock Asymptotic dimension and small-cancellation for hierarchically
  hyperbolic spaces and groups.
\newblock {\em Proc. Lond. Math. Soc. (3)}, 114(5):890--926, 2017.

\bibitem[BHS17b]{HHS_I}
Jason Behrstock, Mark~F. Hagen, and Alessandro Sisto.
\newblock Hierarchically hyperbolic spaces, {I}: {C}urve complexes for cubical
  groups.
\newblock {\em Geom. Topol.}, 21(3):1731--1804, 2017.

\bibitem[BHS19]{HHS_II}
Jason Behrstock, Mark Hagen, and Alessandro Sisto.
\newblock Hierarchically hyperbolic spaces {II}: {C}ombination theorems and the
  distance formula.
\newblock {\em Pacific J. Math.}, 299(2):257--338, 2019.

\bibitem[BHS21]{HHS:quasiflats}
Jason Behrstock, Mark~F. Hagen, and Alessandro Sisto.
\newblock Quasiflats in hierarchically hyperbolic spaces.
\newblock {\em Duke Math. J.}, 170(5):909--996, 2021.

\bibitem[Bir74]{Birman:Braids}
J.~Birman.
\newblock {\em Braids, links, and mapping class groups}, volume~82 of {\em
  Annals of Math. Studies}.
\newblock Princeton University Press, 1974.

\bibitem[Bow12]{Bow:rel_hyp}
B.~H. Bowditch.
\newblock Relatively hyperbolic groups.
\newblock {\em Internat. J. Algebra Comput.}, 22(3):1250016, 66, 2012.

\bibitem[Bow13]{Bowditch:coarse_median}
Brian~H. Bowditch.
\newblock Coarse median spaces and groups.
\newblock {\em Pacific J. Math.}, 261(1):53--93, 2013.

\bibitem[Bow20]{Bowditch:WP_rigid}
Brian~H. Bowditch.
\newblock Large-scale rank and rigidity of the {W}eil-{P}etersson metric.
\newblock {\em Groups Geom. Dyn.}, 14(2):607--652, 2020.

\bibitem[BR20]{BerlaiRobbio}
Federico Berlai and Bruno Robbio.
\newblock A refined combination theorem for hierarchically hyperbolic groups.
\newblock {\em Groups Geom. Dyn.}, 14(4):1127--1203, 2020.

\bibitem[BR22]{BerlyneRussell}
Daniel Berlyne and Jacob Russell.
\newblock Hierarchical hyperbolicity of graph products.
\newblock {\em Groups Geom. Dyn.}, 16(2):523--580, 2022.

\bibitem[Bri10]{Bridson:MCG_not_CAT(0)}
Martin~R. Bridson.
\newblock Semisimple actions of mapping class groups on {${\rm CAT}(0)$}
  spaces.
\newblock In {\em Geometry of {R}iemann surfaces}, volume 368 of {\em London
  Math. Soc. Lecture Note Ser.}, pages 1--14. Cambridge Univ. Press, Cambridge,
  2010.

\bibitem[BRW17]{BridsonReidWilton}
Martin~R. Bridson, Alan~W. Reid, and Henry Wilton.
\newblock Profinite rigidity and surface bundles over the circle.
\newblock {\em Bull. Lond. Math. Soc.}, 49(5):831--841, 2017.

\bibitem[CMM21]{CMM:normal}
Matt Clay, Johanna Mangahas, and Dan Margalit.
\newblock Right-angled {A}rtin groups as normal subgroups of mapping class
  groups.
\newblock {\em Compos. Math.}, 157(8):1807--1852, 2021.

\bibitem[Dah18]{D_PlateSpinning}
Fran{\c{c}}ois Dahmani.
\newblock The normal closure of big {D}ehn twists, and plate spinning with
  rotating families.
\newblock {\em Geom. Topol.}, 22:4113--4144, 2018.

\bibitem[DDLS20]{Veech:hhs}
Spencer Dowdall, Matthew~G Durham, Christopher~J Leininger, and Alessandro
  Sisto.
\newblock Extensions of veech groups are hierarchically hyperbolic.
\newblock {\em arXiv preprint arXiv:2006.16425}, 2020.

\bibitem[DG18]{DG:recognize_DF}
Fran\c{c}ois Dahmani and Vincent Guirardel.
\newblock Recognizing a relatively hyperbolic group by its {D}ehn fillings.
\newblock {\em Duke Math. J.}, 167(12):2189--2241, 2018.

\bibitem[DHS17]{HHS:boundary}
Matthew~Gentry Durham, Mark~F. Hagen, and Alessandro Sisto.
\newblock Boundaries and automorphisms of hierarchically hyperbolic spaces.
\newblock {\em Geom. Topol.}, 21(6):3659--3758, 2017.

\bibitem[DHS20]{HHS:corrigendum}
Matthew~Gentry Durham, Mark~F. Hagen, and Alessandro Sisto.
\newblock Correction to the article {B}oundaries and automorphisms of
  hierarchically hyperbolic spaces.
\newblock {\em Geom. Topol.}, 24(2):1051--1073, 2020.

\bibitem[DHS21]{DFDT}
Fran\c{c}ois Dahmani, Mark Hagen, and Alessandro Sisto.
\newblock Dehn filling {D}ehn twists.
\newblock {\em Proc. Roy. Soc. Edinburgh Sect. A}, 151(1):28--51, 2021.

\bibitem[DT15]{DurhamTaylor:stability}
Matthew Durham and Samuel~J Taylor.
\newblock Convex cocompactness and stability in mapping class groups.
\newblock {\em Algebraic \& Geometric Topology}, 15(5):2837--2857, 2015.

\bibitem[DT19]{DT-DecidingIsom}
Fran\c{c}ois Dahmani and Nicholas Touikan.
\newblock Deciding isomorphy using {D}ehn fillings, the splitting case.
\newblock {\em Invent. Math.}, 215(1):81--169, 2019.

\bibitem[(Ed95]{Kirby_list}
Rob~Kirby (Ed.).
\newblock Problems in low-dimensional topology.
\newblock In {\em Proceedings of Georgia Topology Conference, Part 2}, pages
  35--473. Press, 1995.

\bibitem[FM02]{FarbMosher}
Benson Farb and Lee Mosher.
\newblock Convex cocompact subgroups of mapping class groups.
\newblock {\em Geom. Topol.}, 6:91--152, 2002.

\bibitem[FM12]{FarbMargalit}
Benson Farb and Dan Margalit.
\newblock {\em A primer on mapping class groups}, volume~49 of {\em Princeton
  Mathematical Series}.
\newblock Princeton University Press, Princeton, NJ, 2012.

\bibitem[Fuj15]{Fujiwara:Freesubgroups2}
Koji Fujiwara.
\newblock Subgroups generated by two pseudo-{A}nosov elements in a mapping
  class group. {II}. {U}niform bound on exponents.
\newblock {\em Trans. Amer. Math. Soc.}, 367(6):4377--4405, 2015.

\bibitem[Fun99]{Funar}
Louis Funar.
\newblock On the {TQFT} representations of the mapping class groups.
\newblock {\em Pacific J. Math.}, 188(2):251--274, 1999.

\bibitem[GM08]{GrMa-perfill}
D.~Groves and J.~F. Manning.
\newblock Dehn filling in relatively hyperbolic groups.
\newblock {\em Israel J. Math.}, 168:317--429, 2008.

\bibitem[Gro75]{Grossman:res_fin_MCG}
Edna~K. Grossman.
\newblock On the residual finiteness of certain mapping class groups.
\newblock {\em J. London Math. Soc. (2)}, 9:160--164, 1974/75.

\bibitem[Hae16]{Haettel}
Thomas Haettel.
\newblock Higher rank lattices are not coarse median.
\newblock {\em Algebraic \& Geometric Topology}, 16(5):2895--2910, 2016.

\bibitem[Ham05]{Hamenstadt}
Ursula Hamenst{\"a}dt.
\newblock Word hyperbolic extensions of surface groups.
\newblock {\em arXiv preprint math/0505244}, 2005.

\bibitem[HMS24]{HHS:artin}
Mark Hagen, Alexandre Martin, and Alessandro Sisto.
\newblock Extra-large type {A}rtin groups are hierarchically hyperbolic.
\newblock {\em Math. Ann.}, 388(1):867--938, 2024.

\bibitem[HQR22]{HQR:big}
Camille Horbez, Yulan Qing, and Kasra Rafi.
\newblock Big mapping class groups with hyperbolic actions: classification and
  applications.
\newblock {\em J. Inst. Math. Jussieu}, 21(6):2173--2204, 2022.

\bibitem[HRSS22]{HHS:graph-man}
Mark Hagen, Jacob Russell, Alessandro Sisto, and Davide Spriano.
\newblock Equivariant hierarchically hyperbolic structures for 3-manifold
  groups via quasimorphisms.
\newblock {\em arXiv preprint arXiv:2206.12244}, 2022.

\bibitem[Hum92]{Humphries:powersdehn}
Stephen~P. Humphries.
\newblock Normal closures of powers of {D}ehn twists in mapping class groups.
\newblock {\em Glasgow Math. J.}, 34(3):313--317, 1992.

\bibitem[HW08]{HaglundWise:special}
Fr\'{e}d\'{e}ric Haglund and Daniel~T. Wise.
\newblock Special cube complexes.
\newblock {\em Geom. Funct. Anal.}, 17(5):1551--1620, 2008.

\bibitem[Iva06]{Ivanov_list}
Nikolai~V. Ivanov.
\newblock Fifteen problems about the mapping class groups.
\newblock In {\em Problems on mapping class groups and related topics},
  volume~74 of {\em Proc. Sympos. Pure Math.}, pages 71--80. Amer. Math. Soc.,
  Providence, RI, 2006.

\bibitem[KK13]{KK1}
Sang-hyun Kim and Thomas Koberda.
\newblock Embedability between right-angled {A}rtin groups.
\newblock {\em Geom. Topol.}, 17(1):493--530, 2013.

\bibitem[KK14]{KK2}
Sang-Hyun Kim and Thomas Koberda.
\newblock The geometry of the curve graph of a right-angled {A}rtin group.
\newblock {\em Internat. J. Algebra Comput.}, 24(2):121--169, 2014.

\bibitem[KL96]{KapovichLeeb:actions}
Michael Kapovich and Bernhard Leeb.
\newblock Actions of discrete groups on nonpositively curved spaces.
\newblock {\em Math. Ann.}, 306(2):341--352, 1996.

\bibitem[KL08]{KentLeininger:convex_cocompact}
Autumn~E Kent and Christopher~J Leininger.
\newblock Shadows of mapping class groups: capturing convex cocompactness.
\newblock {\em Geometric and Functional Analysis}, 18(4):1270--1325, 2008.

\bibitem[KM12]{KahnMarkovic}
Jeremy Kahn and Vladimir Markovic.
\newblock Immersing almost geodesic surfaces in a closed hyperbolic three
  manifold.
\newblock {\em Ann. of Math. (2)}, 175(3):1127--1190, 2012.

\bibitem[LM07]{LeiningerMcReynolds}
Christopher~J. Leininger and D.~B. McReynolds.
\newblock Separable subgroups of mapping class groups.
\newblock {\em Topology Appl.}, 154(1):1--10, 2007.

\bibitem[MM99]{MM_I}
Howard~A. Masur and Yair~N. Minsky.
\newblock Geometry of the complex of curves. {I}. {H}yperbolicity.
\newblock {\em Invent. Math.}, 138(1):103--149, 1999.

\bibitem[MM00]{MM_II}
H.~A. Masur and Y.~N. Minsky.
\newblock Geometry of the complex of curves. {II}. {H}ierarchical structure.
\newblock {\em Geom. Funct. Anal.}, 10(4):902--974, 2000.

\bibitem[MS22]{MangioniSisto}
Giorgio Mangioni and Alessandro Sisto.
\newblock Rigidity of mapping class groups mod powers of twists.
\newblock {\em arXiv preprint arXiv:2212.11014}, 2022.

\bibitem[Nic13]{Nica}
Bogdan Nica.
\newblock Proper isometric actions of hyperbolic groups on {$L^p$}-spaces.
\newblock {\em Compos. Math.}, 149(5):773--792, 2013.

\bibitem[Osi07]{Os-perfill}
D.~V. Osin.
\newblock Peripheral fillings of relatively hyperbolic groups.
\newblock {\em Invent. Math.}, 167(2):295--326, 2007.

\bibitem[Par09]{Paris:residual}
Luis Paris.
\newblock Residual {$p$} properties of mapping class groups and surface groups.
\newblock {\em Trans. Amer. Math. Soc.}, 361(5):2487--2507, 2009.

\bibitem[Pic81]{Picard1}
{\'E}mile Picard.
\newblock Sur une extension aux fonctions de deux variables du probl{\`e}me de
  riemann relatif aux fonctions hyperg{\'e}om{\'e}triques.
\newblock In {\em Annales scientifiques de l'{\'E}cole Normale Sup{\'e}rieure},
  volume~10, pages 305--322, 1881.

\bibitem[Pic85]{Picard2}
{\'E}mile Picard.
\newblock Sur les fonctions hyperfuchsiennes provenant des s{\'e}ries
  hyperg{\'e}om{\'e}triques de deux variables.
\newblock In {\em Annales scientifiques de l'{\'E}cole Normale Sup{\'e}rieure},
  volume~2, pages 357--384, 1885.

\bibitem[PS23]{PetytSpriano}
Harry Petyt and Davide Spriano.
\newblock Unbounded domains in hierarchically hyperbolic groups.
\newblock {\em Groups Geom. Dyn.}, 17(2):479--500, 2023.

\bibitem[Rei06]{Reid:separable}
Alan~W. Reid.
\newblock Surface subgroups of mapping class groups.
\newblock In {\em Problems on mapping class groups and related topics},
  volume~74 of {\em Proc. Sympos. Pure Math.}, pages 257--268. Amer. Math.
  Soc., Providence, RI, 2006.

\bibitem[RS20]{RobbioSpriano}
Bruno Robbio and Davide Spriano.
\newblock Hierarchical hyperbolicity of hyperbolic-2-decomposable groups.
\newblock {\em arXiv preprint arXiv:2007.13383}, 2020.

\bibitem[RST22]{RST:morse}
Jacob Russell, Davide Spriano, and Hung~Cong Tran.
\newblock The local-to-global property for {M}orse quasi-geodesics.
\newblock {\em Math. Z.}, 300(2):1557--1602, 2022.

\bibitem[RST23]{RST}
Jacob Russell, Davide Spriano, and Hung~Cong Tran.
\newblock Convexity in hierarchically hyperbolic spaces.
\newblock {\em Algebr. Geom. Topol.}, 23(3):1167--1248, 2023.

\bibitem[Rus22]{Russell:rel_HHS}
Jacob Russell.
\newblock From hierarchical to relative hyperbolicity.
\newblock {\em Int. Math. Res. Not. IMRN}, (1):575--624, 2022.

\bibitem[Sco78]{scott1978subgroups}
Peter Scott.
\newblock Subgroups of surface groups are almost geometric.
\newblock {\em J. London Math. Soc. (2)}, 17(3):555--565, 1978.

\bibitem[Thu98]{Thurston:quotients}
William~P. Thurston.
\newblock Shapes of polyhedra and triangulations of the sphere.
\newblock In {\em The {E}pstein birthday schrift}, volume~1 of {\em Geom.
  Topol. Monogr.}, pages 511--549. Geom. Topol. Publ., Coventry, 1998.

\bibitem[Vok22]{Vokes}
Kate~M. Vokes.
\newblock Hierarchical hyperbolicity of graphs of multicurves.
\newblock {\em Algebr. Geom. Topol.}, 22(1):113--151, 2022.

\bibitem[Wis04]{Wise:cubulating_small_can}
D.~T. Wise.
\newblock Cubulating small cancellation groups.
\newblock {\em Geom. Funct. Anal.}, 14(1):150--214, 2004.

\bibitem[Yu05]{Yu}
Guoliang Yu.
\newblock Hyperbolic groups admit proper affine isometric actions on
  {$l^p$}-spaces.
\newblock {\em Geom. Funct. Anal.}, 15(5):1144--1151, 2005.

\end{thebibliography}
\end{document}